\documentclass[a4paper,10pt,twoside]{article}
\input art.sty
\usepackage{a4wide}
\SetMathAlphabet\mathcal{normal}{U}{rsfs}{m}{n} 

\fancypagestyle{plain}{ 
\fancyhead{} 
}

\title{Foncteur de Picard d'un champ algébrique}
\date{\empty}

\def\Pic{\hbox{\rm Pic}}
\def\piczar{\Pic_{\X/S \, \rm{(Zar)}}}
\def\picet{\Pic_{\X/S \, \text{(\'Et)}}}
\def\piceth{\Pic_{\X/S \, \emph{(\'Et)}}}
\def\picfppf{\Pic_{\X/S \, \rm{(fppf)}}}
\def\pic{\Pic_{\X/S}}

\def\piczero{\Pic_{\X/S}^0}
\def\champic{\mathcal{P}ic}
\def\liset{\text{{\rm Lis-\'et}}}
\def\lisets{\text{\underline{Lis-\'et}}}

\def\gll{\text{{\rm Llc}}}
\def\et{\text{\rm {\'Et}}}
\def\fppf{\text{\rm fppf}}
\def\fppfc{\text{\rm fppfc}}
\def\pl{\text{\rm pl}}

\def\A{\mathbb{A}}
\def\X{\mathcal{X}}
\def\Y{\mathcal{Y}}
\def\U{\mathcal{U}}

\def\inv{\hbox{\rm Inv}\,}
\def\aff{\hbox{{\rm (Aff/}S)}}
\def\Xt{\widetilde{\X}}
\def\Yt{\widetilde{\Y}}

\def\Lt{\widetilde{\Lc}}
\def\Mt{\widetilde{\Mc}}
\def\defm{\hbox{\rm Defm}}

\renewcommand{\fleche}{%
\xymatrix@C=1pc{\ar[r] &}}
\def\flechelongue{%
\xymatrix{\ar[r] &}}

\entrymodifiers={!!<0pt,0.7ex>+}
\begin{document}

\maketitle
\vskip -2cm
\begin{abstract}
In this article we study the Picard functor and the Picard stack of an algebraic stack. We give a new and direct proof of the representability of the Picard stack. We prove that it is quasi-separated, and that the connected component of the identity is proper when the fibers of $\X$ are geometrically normal. We study some examples of Picard functors of classical stacks. In an appendix, we review the lisse-etale cohomology of abelian sheaves on an algebraic stack.
\end{abstract}

\tableofcontents \vfill
\section{Introduction}

L'objet principal de cet article est l'étude du foncteur de Picard et du champ de Picard d'un champ algébrique. Notre but est de généraliser au cas des champs algébriques les résultats classiques valables dans le cadre des schémas, que l'on peut trouver dans \cite{poly_Kleiman} ou \cite{BRL}.

Si $\X$ est un champ algébrique sur un schéma $S$, on appelle champ de Picard de $\X$, et l'on note $\champic(\X/S)$, le champ classifiant des faisceaux inversibles sur $\X$, \emph{i.e.} le champ dont la catégorie fibre en $U$ est la catégorie des faisceaux inversibles sur $\X\times_S U$. On définit aussi le foncteur de Picard relatif de la manière suivante. On commence par définir un foncteur $P_{\X/S}$ qui à tout $U$ associe
$$P_{\X/S}(U)=\frac{\Pic(\X\times_S U)}{\Pic(U)}.$$
On note alors $\piczar$ (resp. $\picet, \picfppf$) le faisceau associ\'e \`a
$P_{\X/S}$ pour la topologie de Zariski (resp. \'etale, fppf)\footnote{Lorsque nous parlerons \emph{du} foncteur de Picard de $\X/S$, il s'agira de $\picfppf$. Nous le noterons $\pic$ s'il n'y a pas d'ambiguïté.}. Nous démêlons dans la première partie les liens étroits qui unissent les différents foncteurs de Picard et le champ de Picard.

La première question qui se pose est celle de l'algébricité du champ de Picard (donc par~\ref{comparaison_champ_foncteur} celle de la représentabilité du foncteur de Picard). Elle a été résolue par Aoki.

\begin{thm}[Aoki, \cite{Aoki_Hom}~5.1]
\label{algebricite_champ_Picard}
Si $\X$ est propre et plat sur $S$, alors le champ $\champic(\X/S)$ est un champ algébrique au sens d'Artin.
\end{thm}

Cependant ce résultat apparaît comme un cas particulier d'un théorème plus général et plus difficile. Aoki montre en effet dans \cite{Aoki_Hom} que si $\X$ et $\Y$ sont deux champs algébriques, alors sous de bonnes hypothèses le champ $\fHom(\X,\Y)$ est un champ algébrique au sens d'Artin. La d\'emonstration fait appel \`a un certain nombre de r\'esultats non triviaux concernant les d\'eformations de morphismes d'espaces alg\'ebriques (\cite{Illusie_CCD}), les d\'eformations de morphismes repr\'esentables de
champs alg\'ebriques (\cite{Olsson_defm}), les d\'eformations de champs
alg\'ebriques (\cite{Aoki_defm}) et la th\'eorie du complexe
cotangent (\cite{Illusie_CCD}, \cite{LMB} et
\cite{Olsson_Sheaves_on_Artin_stacks}). Le cas du champ de Picard s'en déduit en prenant pour $\Y$ le champ $\bgm$.
Nous en donnons au paragraphe~\ref{algebricite} une preuve directe, plus rapide et plus concrète.

Nous traitons au paragraphe~\ref{quasi-separation} la question de la quasi-séparation. Le résultat principal est le suivant :

\begin{thm}
\label{quasi_separation_du_champ_de_Picard}
Soient $S$ un sch\'ema noeth\'erien, et $\X$ un $S$-champ alg\'ebrique propre,
plat, et cohomologiquement plat en dimension z\'ero. Alors le champ de Picard
$\champic(\X/S)$ est quasi-s\'epar\'e (donc est algébrique au sens de \cite{LMB}).
\end{thm}

Nous étudions ensuite la composante neutre du foncteur de Picard. Le résultat des EGA sur la composante connexe des fibres le long d'une section (\cite{EGA4_3}~(15.6.5)) que l'on utilise usuellement pour la définir ne s'applique pas tel quel lorsque $\pic$ est un espace algébrique. Notre premier travail a donc été de généraliser ce résultat au cas des espaces algébriques afin de s'assurer que la construction décrite précédemment fournirait bien un sous-espace algébrique ouvert $\piczero$ du foncteur de Picard.

Une fois la composante neutre définie, nous nous posons la question de sa propreté sur $S$. Nous aboutissons au résultat suivant, qui généralise les résultats analogues pour les schémas.

\begin{thm}
\label{composante_neutre_propre}
Soit $\X$ un $S$-champ algébrique propre, plat et cohomologiquement plat en dimension zéro. On suppose que $\X$ est à fibres géométriquement normales, que $\piczero$ est un ouvert de $\Pic_{\X/S}$ et qu'il est de plus séparé et de type fini sur $S$. Alors $\piczero$ est propre sur $S$.
\end{thm}

Nous montrons d'abord le théorème sur un corps (\emph{cf.}~\ref{composante_neutre_propre_sur_k}). La démonstration repose sur une étude des $\Z$-torseurs sur un champ normal : nous commençons par montrer que si $\X$ est un champ algébrique localement noethérien et normal alors le groupe $H^1(\X,\Z)$ est trivial.

Enfin, nous illustrons le texte en étudiant le foncteur de Picard de quelques champs algébriques classiques : l'espace de modules des courbes elliptiques, le champ des racines $n^{\text{ièmes}}$ d'un faisceau inversible, et les courbes tordues d'Abramovich et Vistoli.

\medskip
\noindent
{\sc Revue de la cohomologie lisse-étale}

Le défaut de fonctorialité du topos lisse-étale des champs algébriques est maintenant bien connu : si $f : \X \fleche \Y$ est un morphisme de champs algébriques, le foncteur $f^{-1}$ n'est pas toujours exact. Ceci a pour conséquence fâcheuse que la \og machine \fg\ SGA4 ne s'applique pas toujours et qu'il faut par conséquent travailler un peu plus finement pour obtenir un certain nombre de propriétés d'apparence pourtant élémentaire sur la cohomologie lisse-étale des champs algébriques. Ce travail a été fait en grande partie par Olsson et Laszlo pour le cas des faisceaux quasi-cohérents (voir \cite{Olsson_Sheaves_on_Artin_stacks}) ou des coefficients finis (voir \cite{Laszlo_Olsson_Six_operationsI}). Mais les faisceaux abéliens qui ne jouissent pas d'une telle structure ont pour l'instant été laissés de côté.

Or le faisceau abélien que l'on rencontre le plus souvent lorsque l'on s'intéresse au foncteur de Picard n'est pas quasi-cohérent : il s'agit de $\gm$. Nous avons donc été amené à démontrer au fil de nos travaux diverses propriétés. Leur nombre croissant nous a finalement conduit à les regrouper dans un chapitre à part. Cette annexe est donc une sorte de \og fourre-tout \fg\ cohomologique qui, bien loin de prétendre à l'exhaustivité, se contente au contraire de répondre aux strictes exigences des autres chapitres.  Nous espérons tout de même que nous aurons ainsi contribué à combler une lacune de la littérature existante.

Je remercie chaleureusement Laurent Moret-Bailly pour la patience et la disponibilité dont il a su faire preuve dans son travail de direction, et pour les relectures minutieuses de ma thèse, dont cet article est issu.

\subsection*{Conventions}
\vskip -1.5mm
Suivant \cite{LMB}, sauf mention expresse du contraire, tous les champs algébriques (\emph{a fortiori} tous les schémas et tous les espaces algébriques) seront quasi-séparés. Un champ \og algébrique \fg\ non-quasi-séparé, \emph{i.e.} un champ dont la diagonale est représentable et localement de type fini, et qui admet une présentation lisse, sera appelé \og champ algébrique au sens d'Artin\fg. Nous renvoyons aux annexes pour le formalisme cohomologique.

\section{Foncteurs de Picard et champ de Picard}

\subsection{Description cohomologique}
Comme dans le cas des schémas, on peut calculer le groupe de Picard de $\X$ à l'aide de la cohomologie de $\gm$. Rappelons tout d'abord le résultat de descente ci-dessous.

\begin{sousprop}[\cite{LMB}~(13.5)]
\label{descente_fidelement_plate}
Soit $f : \Y\fleche \X$ un morphisme fidèlement plat de champs algébriques, que l'on
suppose de plus quasi-compact ou localement de présentation finie. Alors la
catégorie des faisceaux inversibles sur $\X$ est équivalente à
la catégorie ${\rm Desc}(\Y,\X)$ décrite de la manière suivante. Un objet est un couple $(\Lc,\alpha)$ où $\Lc$ est un
faisceau inversible sur $\Y$ et
$\alpha : p_1^*\Lc \fleche p_2^*\Lc$
est un isomorphisme tel que, à des isomorphismes canoniques près,
$(p_{23}^*\alpha) \circ (p_{12}^*\alpha) = p_{13}^*\alpha$. Un morphisme
de $(\Lc,\alpha)$ vers $(\Mc, \beta)$ est un morphisme $\gamma :
\Lc \fleche \Mc$ tel que $(p_2^*\gamma)\circ \alpha= \beta\circ
(p_1^*\gamma)$.
\end{sousprop}

\begin{sousprop}
\label{pic_egal_h1}
Soit $\X$ un $S$-champ alg\'ebrique. Alors :
$$\Pic \X \simeq H^1(\X,\gm),$$
o\`u $H^1(\X, \gm)$ est le premier groupe de cohomologie du faisceau $\gm$ sur
$\X$ muni de la topologie lisse-\'etale calcul\'e au sens des foncteurs
d\'eriv\'es.
\end{sousprop}
\begin{demo}
Soit $T\fleche \X$ une pr\'esentation de $\X$. On reprend ici les notations
de l'annexe (\ref{annexe_desc_coh}) qui pr\'ec\`edent le th\'eor\`eme
(\ref{suite_spectrale_de_descente}). On consid\`ere en particulier le premier
groupe de cohomologie \og \`a la \v{C}ech \fg\ associ\'e \`a cette
pr\'esentation :
$$\check{H}^1(H^0(T^{\bullet},\gm))=
\frac{\Ker(p_{23}^*-p_{13}^*+p_{12}^*)}{\Im(p_1^*-p_2^*)}.$$
Se donner un 1-cocycle de
\v Cech \`a valeurs dans $\gm$ revient \`a se donner une donn\'ee de descente
sur $\Oc_T$, et deux telles donn\'ees de descente $g_1, g_2$ d\'efinissent le
m\^eme \'el\'ement dans $\check{H}^1(H^0(T^{\bullet},\gm))$ si et seulement si
$(\Oc_T,g_1)$ et $(\Oc_T,g_2)$ sont
isomorphes dans la cat\'egorie des faisceaux inversibles sur $T$ munis d'une
donn\'ee de descente relativement \`a $T\fleche \X$. Compte tenu de
l'\'equivalence (\ref{descente_fidelement_plate}) entre cette cat\'egorie et la
cat\'egorie des faisceaux inversibles sur $\X$, le groupe
$\check{H}^1(H^0(T^{\bullet},\gm))$ s'identifie \`a l'ensemble des classes
d'isomorphie de faisceaux inversibles sur $\X$ dont l'image inverse sur $T$ est
isomorphe \`a $\Oc_T$, donc au groupe $\Ker(\Pic(\X) \fleche \Pic(T))$. De plus
cet isomorphisme est fonctoriel en $T$. En passant \`a la limite inductive, on
obtient un isomorphisme
$$\lind \check{H}^1(H^0(T^{\bullet},\gm)) \simeq \Pic \X,$$
o\`u la limite inductive est prise sur l'ensemble des pr\'esentations lisses
$T\fleche \X$ de $\X$. On dispose par ailleurs de la suite exacte en bas
degr\'es associ\'ee \`a la suite spectrale de descente
(\ref{suite_spectrale_de_descente}) :
$$0\fleche \check{H}^1(H^0(T^{\bullet},\gm)) \fleche H^1(\X,\gm) \fleche
\check{H}^0(H^1(T^{\bullet},\gm)).$$
Cette suite exacte \'etant fonctorielle en $T$, on voit que les injections
des $\check{H}^1(H^0(T^{\bullet},\gm))$ dans $H^1(\X,\gm)$ induisent
un morphisme
injectif de $\lind \check{H}^1(H^0(T^{\bullet},\gm))$ dans $H^1(\X,\gm)$.
Pour conclure, il ne reste plus qu'\`a montrer, vu la suite exacte en bas
degr\'es ci-dessus, que pour tout $x\in H^1(\X,\gm)$, il existe une
pr\'esentation $T\fleche \X$ telle que l'image de $x$ dans $H^1(T,\gm)$ soit
nulle.
C'est un exercice facile laissé au lecteur.
\end{demo}

On peut calculer les groupes de Picard relatifs à l'aide des images directes supérieures de $\gm$.

\begin{sousprop}
\label{description_coh_de_picet_et_picfppf}
Soit $\X$ un $S$-champ alg\'ebrique. Alors pour tout $S$-sch\'ema $T$:
$$\begin{array}{l}
\piceth(T)\simeq H^0(T,R^1{f_T}_* \gm)\\
\picfppf(T)\simeq H_{\text{\rm pl}}^0(T,R^1f_{T*}^{\text{\rm pl}} \gm).
\end{array}$$
\end{sousprop}
\begin{demo}
D'après la proposition~\ref{prop_image_directe_sup},
la restriction au site étale de $T$ du faisceau $R^1\hspace{-0.56pt} f_{T*}\gm$ est le
faisceau \'etale associ\'e \`a
$$\xymatrix{U \ar@{|->}[r] & H^1(\X_T\times_T U, \gm)=\Pic_{\X_T}(U)}.$$
Donc $(R^1f_{T*}\gm)_{\text{ét}}=\Pic_{\X_T/T \, \text{(\'Et)}}$ et en particulier
$$H^0(T, R^1f_{T*}\gm)=\Pic_{\X_T/T \, \text{(\'Et)}}(T)=\picet(T).$$
On obtient la seconde assertion en appliquant exactement le même raisonnement pour la topologie \emph{fppf}. On utilise~\ref{images_directes_superieures_fppf}
au lieu de~\ref{prop_image_directe_sup}.
\end{demo}

\subsection{Comparaison entre les différents foncteurs de Picard}

\begin{sousremarque}\rm
\label{P_isom_Pet_si_k_alg_clos}
On peut d'ores et déjà comparer les groupes de sections des différents
foncteurs de Picard dans quelques cas particuliers grâce à la remarque de topologie
suivante. Supposons que $U$ soit un $S$-schéma affine tel que, pour toute
famille couvrante $(U_i\fleche U)_{i\in I}$ pour une certaine topologie
\emph{(Top)}, il existe un indice $i\in I$ tel que le morphisme $U_i\fleche U$
ait une section. Alors le morphisme $\Pic(\X\times_S U) \fleche
\Pic_{\X/S \, \rm{(Top)}}(U)$ est un isomorphisme. En particulier, si $U$ est le spectre d'un anneau
local (resp. d'un anneau local strictement hensélien, d'un corps algébriquement
clos) alors $\piczar(U)$ (resp. $\picet(U)$, $\picfppf(U)$) s'identifie à $\Pic(\X\times_S U)$.
\end{sousremarque}

La notion de morphisme cohomologiquement plat en
dimension z\'ero se généralise sans difficultés aux champs algébriques. 
\begin{sousdefi}
On dit qu'un morphisme de champs alg\'ebriques $f : \X \fleche \Y$ est
cohomologiquement plat en dimension z\'ero s'il est plat et si le
morphisme
$$f^{\sharp} : \Oc_{\Y} \flechelongue f_*\Oc_{\X}$$
est un isomorphisme universellement (en particulier $\X$ et $\Y$ ont alors les mêmes sections globales).
\end{sousdefi}

\begin{sousremarque} \rm
\label{remarque_section}
Un morphisme cohomologiquement plat en dimension z\'ero est surjectif. En effet, pour tout point $s : \Spec k \fleche \Y$ de $\Y$, la fibre $\X_s$
de $\X$ au-dessus de $s$ est non vide puisque
$\Gamma(\X_s,\Oc_{\X_s})$ est isomorphe à $k$.
On en déduit que si $S$ est un schéma et
$f~:~\X \fleche S$ un morphisme localement de présentation finie et
cohomologiquement plat en dimension zéro de $S$-champs algébriques
alors $f$ a une section localement pour la topologie \emph{(fppf)} sur $S$.
\end{sousremarque}

\begin{sousremarque}\rm
La notion de platitude cohomologique est intimement liée à la conne\-xité des fibres géométriques. En effet si $\X$ est cohomologiquement plat en dimension zéro sur $\Y$, ses fibres sont géométriquement connexes. 
Réciproquement, dans le cas où le morphisme $f$ est propre et plat, il suffit que les fibres de $\X$ soient géométriquement connexes et géométriquement réduites pour que $\X$ soit cohomologiquement plat en dimension zéro sur $\Y$ (au moins lorsque $\Y$ est noethérien et réduit).
\end{sousremarque}

Le théorème suivant permet de comparer les différents foncteurs de Picard d'un champ algébrique. Il généralise \cite{poly_Kleiman} 2.5.
\begin{sousthm}
\label{comparaison_des_foncteurs_de_Picard}
Soit $f : \X\fleche S$ un $S$-champ alg\'ebrique cohomologiquement plat en
dimension z\'ero. Alors les morphismes naturels
$$\xymatrix{P_{\X/S} \ar[r]^-{i_1} &\piczar \ar[r]^{i_2} &
\piceth}$$
sont injectifs, et le morphisme naturel
$$\flechen{\piceth}{}{\picfppf}$$
est un isomorphisme. Si de plus $f$ a une section localement pour la topologie
de Zariski, alors $i_2$ est un isomorphisme. Enfin si $f$ a une section, $i_1$
est lui aussi un isomorphisme.
\end{sousthm}
\begin{demo}
Soit $T$ un $S$-sch\'ema.
La suite spectrale de Leray
$$H^p(T, R^qf_{T*}\gm) \Rightarrow H^{p+q}(\X_T,\gm)$$
induit la suite exacte longue en bas degr\'es suivante :
\begin{multline} 
0\fleche H^1(T,f_{T*}\gm) \fleche H^1(\X_T,\gm)\fleche H^0(T, R^1f_{T*}\gm)
\fleche  \\
\fleche H^2(T,f_{T*}\gm) \fleche H^2(\X_T,\gm) \nonumber
\end{multline}
 (o\`u tous les calculs sont effectu\'es pour la topologie lisse-\'etale sur le
champ alg\'ebrique consid\'er\'e. En particulier lorsqu'il s'agit d'un sch\'ema,
cela revient d'après~(\ref{coh_et_egale_coh_liset}) \`a calculer sa cohomologie
pour la topologie \'etale.) Le morphisme
$f$ \'etant cohomologiquement plat en dimension z\'ero, on a $f_{T*}\gm=\gm$.
D'apr\`es les propri\'et\'es pr\'ec\'edentes, le d\'ebut de la suite exacte
ci-dessus fournit la suite exacte :
$$0\fleche \Pic(T)\fleche \Pic(\X_T)\fleche \picet(T)$$
ce qui montre que le morphisme naturel $P_{\X/S} \fleche \picet$
est injectif. Il en r\'esulte, d'une part, que le morphisme naturel
$P_{\X/S} \fleche \piczar$
est injectif, et d'autre part, en appliquant le foncteur \og faisceau associ\'e
pour la topologie de Zariski \fg, que le morphisme naturel
$\piczar \fleche \picet$
est injectif.

Dans le cas o\`u $f$ a une section, le morphisme induit par $f$ de
$H^2(T,f_{T*}\gm)$ vers $H^2(\X_T,\gm)$ a une r\'etraction, donc il est
injectif. Mais alors la fl\`eche de
$H^0(T, R^1f_{T*}\gm)$ vers $H^2(T,f_{T*}\gm)$
est nulle, et donc $H^1(\X_T,\gm)\fleche H^0(T, R^1f_{T*}\gm)$
est surjectif, d'o\`u une suite exacte courte :
$$0\fleche \Pic(T)\fleche \Pic(\X_T)\fleche \picet(T) \fleche 0$$
ce qui montre que $i_1$ et $i_2$ sont des isomorphismes.

L'avant-dernière assertion en découle immédiatement puisque la question de savoir si $i_2$ est un isomorphisme est locale pour la topologie de Zariski.

Il nous reste \`a montrer que $\picet \fleche \picfppf$ est un
isomorphisme.
On a un diagramme commutatif dans lequel les lignes sont les suites exactes de bas degré associées aux suites
spectrales de Leray pour la cohomologie lisse-étale et pour la cohomologie \emph{fppf} :
$$\hskip -5pt \shorthandoff{!;:?}
\xymatrix@!0 @R=3.5pc @C=7pc{ 
H^1(T,\gm) \ar[r] \ar[d] & H^1(\X_T,\gm)\ar[r]\ar[d] & H^0(T, R^1f_{T*}\gm)
\ar[r]\ar[d] & H^2(T,\gm) \ar[r]\ar[d] & H^2(\X_T,\gm)\ar[d] \\
H^1_{\text{pl}}(T,\gm) \ar[r] & H^1_{\text{pl}}(\X_T,\gm)\ar[r] & H^0_{\text{pl}}(T, R^1f_{T*}^{\text{pl}}\gm)
\ar[r] & H^2_{\text{pl}}(T,\gm) \ar[r] & H^2_{\text{pl}}(\X_T,\gm).
}$$
Le diagramme ci-dessus tient compte du fait que, par platitude cohomologique, $f_{T*}\gm=\gm$ et
$f_{T*}^{\text{pl}}\gm=\gm$. Comme $\gm$ est un groupe lisse
sur la base, le théorème (\ref{coh_fppf_groupe_lisse}) nous dit que les deux flèches verticales de gauche
et les deux flèches verticales de droite sont des isomorphismes. D'après le lemme des 5, celle du milieu en est un aussi.
Or d'après (\ref{description_coh_de_picet_et_picfppf}) cette flèche est précisément le morphisme
$\picet(T)\fleche \picfppf(T)$.
\end{demo}

\subsection{Le champ de Picard d'un champ algébrique}

\label{champ_de_Picard}

Avec l'apparition des champs, un nouvel objet \og classifiant \fg\ pour les faisceaux inversibles est venu s'ajouter aux cinq foncteurs précédents : le champ de Picard. Il r\'esulte ais\'ement de (\ref{descente_fidelement_plate}) que le
$S$-groupoïde $\champic(\X/S)$ défini en introduction est bien un
$S$-champ. De plus la formation de $\champic(\X/S)$ commute au changement de base.
Par ailleurs, en notant $\fHom(\X,\Y)$ le $S$-champ dont la fibre en $U$ est la cat\'egorie
$\Hom(\X\times_SU,\Y\times_SU)$, il est clair que
$\champic(\X/S)$ est isomorphe à $\fHom(\X,\bgm)$. La proposition~(\ref{comparaison_champ_foncteur}) et le corollaire~(\ref{foncteur_pic_repres_implique_champic_alg}) montrent que dans le cas cohomologiquement plat, les questions de l'algébricité du champ de Picard et de la représentabilité du foncteur de Picard sont essentiellement équivalentes.

\begin{sousremarque}\rm
La terminologie adoptée ici, qui semble s'imposer naturellement, n'entre pas en conflit avec celle de \cite{LMB}~(14.4.2), où l'on appelle \og champ de Picard \fg\ un champ muni d'un morphisme d'addition qui en fait une sorte de \og champ en groupes\fg. En effet, si $\X$ est un $S$-champ, \emph{le} champ de Picard de $\X$ défini ci-dessus est bien \emph{un} champ de Picard au
sens de~\cite{LMB}~(14.4.2), le 1-morphisme d'addition
$$+ : \champic(\X/S)\times_S \champic(\X/S) \flechelongue \champic(\X/S)$$
étant donné par le produit tensoriel de faisceaux inversibles.
\end{sousremarque}

Le champ $\champic(\X/S)$ est muni d'un morphisme naturel vers le
$S$-groupo\"\i de
associ\'e au pr\'efaisceau $P_{\X/S}$, d\'efini sur $\champic(\X/S)_U$
par le foncteur qui envoie un faisceau inversible sur sa classe d'isomorphie
dans $P_{\X/S}(U)$. On en d\'eduit par composition un
morphisme naturel $\pi_{\text{ét}}$ (resp. $\pi_{\text{\emph{fppf}}}$) de $\champic(\X/S)$ vers le $S$-espace $\picet$ (resp. $\picfppf$). Il sera noté $\pi$ lorsque $\picet$ et $\picfppf$ coïncident (par exemple dans le cas cohomologiquement plat).

\begin{sousprop}
Le morphisme de $S$-champs $\pi_{\text{\emph{ét}}}$ (resp. $\pi_{\text{fppf}}$) est une gerbe (resp. une gerbe
\emph{(fppf)}), autrement dit c'est un \'epimorphisme étale (resp. \emph{fppf}) et son morphisme diagonal aussi.
\end{sousprop}
\begin{demo}
Notons $P$ le $S$-espace $\picet$ (resp. $\picfppf$).
Soient $U\in \ob\aff$ et $x\in P(U)$. Il existe par définition de $P$ une famille couvrante (pour la topologie considérée) dans $\aff$, que l'on peut
supposer r\'eduite \`a un \'el\'ement $(U'\fleche U)$, et un \'el\'ement $l\in
P_{\X/S}(U')$, tels que l'image de $l$ dans $P(U')$ soit \'egale \`a $x_{|U'}$.
En d'autres termes il existe un faisceau inversible $\Lc$ dans
$\champic(\X/S)_{U'}$ dont l'image dans $P(U')$ est $x_{|U'}$, ce qui montre que
$\pi_{\text{ét}}$ (resp. $\pi_{\text{\emph{fppf}}}$) est un \'epimorphisme.
Ce n'est pas plus difficile pour le morphisme diagonal.
\end{demo}

\begin{sousprop}
\label{comparaison_champ_foncteur}
Si
$\champic(\X/S)$ est un champ alg\'ebrique (resp. algébrique au sens d'Artin), et si f est cohomologiquement plat
en dimension z\'ero, alors $\pic$ est repr\'esentable par un $S$-espace
alg\'ebrique (resp. algébrique au sens d'Artin), et le 1-morphisme $\pi : \champic(\X/S)\fleche \pic$ est
fid\`element plat et localement de pr\'esentation finie.
\end{sousprop}
\begin{demo}
D'apr\`es le corollaire (10.8) de \cite{LMB}, il suffit de v\'erifier que pour
tout $U\in\ob\aff$ et pour tout faisceau inversible $\Lc$ sur $\X\times_S U$, le
$U$-espace
alg\'ebrique en groupes $\fIsom(\Lc,\Lc)$ est plat et localement de
pr\'esentation finie. Or, par
platitude cohomologique,
il est isomorphe \`a $\gm$.
\end{demo}

Réciproquement, on peut essayer de déduire l'algébricité du champ de Picard de celle du foncteur de Picard. Voyons comment.
Le champ des faisceaux inversibles sur $S$, c'est-à-dire le champ dont la catégorie fibre au-dessus d'un $S$-schéma $U$ est la catégorie des faisceaux inversibles sur $U$, s'identifie naturellement au champ $\bgm$ sur $S$. Le foncteur \og image inverse par $f$\fg\ induit donc un morphisme de $S$-champs algébriques
$$\flechen{\bgm}{f^*}{\champic(\X/S).}$$
De plus ce morphisme est un monomorphisme (\emph{i.e.} les foncteurs $f^*_U$ sont pleinement fidèles) dès que $f$ est cohomologiquement plat en dimension zéro. On a donc dans ce cas une \og suite exacte \fg\ de champs de Picard :
$$\xymatrix{0 \ar[r]& \bgm \ar[r]^-{f^*} &\champic(\X/S) \ar[r]^-{\pi}& \pic \ar[r]&0}.$$
Autrement dit $f^*$ est un monomorphisme, $\pi$ est un épimorphisme, un objet de $\champic(\X/S)$ est envoyé sur $0$ par $\pi$ si et seulement s'il provient de $\bgm$, et tous les morphismes sont compatibles aux morphismes d'addition des champs de Picard considérés\footnote{Si $f$ n'est pas cohomologiquement plat en dimension zéro, $f^*$ n'est plus un monomorphisme, mais le reste est presque vrai : il faut juste remplacer \og provient \fg\ par \og provient localement pour la topologie étale\fg.}. La proposition ci-dessous nous permet de déterminer quand cette suite exacte est scindée.

\begin{sousprop}
\label{comparaison_champ_foncteur2}
Soit $\X$ un $S$-champ algébrique
cohomologiquement plat en dimension zéro.
\begin{itemize}
\item[(i)] Les propositions suivantes sont équivalentes :
\begin{itemize}
\item[a)] La suite exacte ci-dessus est scindée, autrement dit il existe un morphisme $s$ de $\pic$ dans $\champic(\X/S)$ tel que $\pi\circ s$ soit égal à l'identité.
\item[b)] Il existe un isomorphisme
de $\bgm\times_S \pic$ dans $\champic(\X/S)$
compatible avec les projections sur $\pic$.
\item[c)] Il existe sur le $S$-champ (non nécessairement algébrique) $\X\times_S\pic$ un faisceau inversible \og universel\fg\ $\Pc$  qui représente le foncteur $\pic$, i.e. tel que pour tout $U\in\ob\aff$ et tout élément $l$ de $\pic(U)$ on ait $l=[\Pc_{|_{\X\times_S U}}]=\pi(\Pc_{|_{\X\times_S U}})$.
$$\xymatrix{\X\times_S U \ar[r] \ar[d] \cartesien &
\X\times_S \pic \ar[r] \ar[d] \cartesien &\X \ar[d] \\
U \ar[r]_l & \pic \ar[r]& S}$$
\end{itemize}
\item[(ii)] Les propositions suivantes sont équivalentes :
\begin{itemize}
\item[d)] Le morphisme naturel de $P_{\X/S}$ dans $\pic$ est un isomorphisme.
\item[e)] Il existe un élément de $P_{\X/S}(\pic)$ (c'est-à-dire $\Hom(\pic,P_{\X/S})$ lorsque $\pic$ n'est pas représentable) qui a pour image l'identité dans $\pic(\pic)$ (ou, ce qui revient au même, $P_{\X/S}\fleche \pic$ a une section).
\end{itemize}
\item[(iii)] Les conditions de (i) impliquent celles de (ii), et la réciproque est vraie si $\pic$ est représentable par un schéma. Toutes ces conditions sont vérifiées si le morphisme structural $\X \fleche S$ a une section.
\end{itemize}
\end{sousprop}
\begin{demo}
Les implications $b\Rightarrow a \Rightarrow e$, $c \Rightarrow a$ et $d \Leftrightarrow e$ sont évidentes.

Montrons que a) implique b). Soit $s$ une section de $\pi$. On vérifie facilement que le morphisme $(\Mc,l) \mapsto (f^*\Mc)\otimes s(l)$ de $\bgm\times_S \pic$ dans $\champic(\X/S)$
est un isomorphisme. (On utilise ici la platitude cohomologique en dimension zéro.) Un quasi-inverse est donné par le foncteur de $\champic(\X/S)$ vers $\bgm \times_S \pic$ qui à un faisceau inversible $\Mc$ associe le couple
$(\Mc\otimes s(\pi(\Mc))^{-1},\pi(\Mc))$, en identifiant $\bgm$ à son image essentielle dans $\champic(\X/S)$. 

Montrons maintenant que a) implique c)\footnote{C'est évident si $\pic$ est représentable par un schéma ! Ce que l'on ne suppose pas ici\dots}. On se donne une section $s$ de $\pi$ et l'on va construire un faisceau inversible $\Pc$ (au sens de \cite{Mumford_Picard_groups}) sur le champ $\X\times_S \pic$. Il faut donc construire pour tout $S$-schéma affine $U$ et tout objet $x$ de $(\X\times_S \pic)_U$ un faisceau inversible $\Pc(x)$ sur $U$, et des isomorphismes de transition entre les $\Pc(x)$. Soit $x$ un objet de $(\X\times_S \pic)_U$. On note $l=f_P\circ x$ et $t_x$ la section de $f_U$ induite par $x$, comme dans le diagramme ci-dessous.
$$\xymatrix{\X\times_S U \ar[r] \ar[d]^{f_U} &
\X\times_S \pic \ar[r] \ar[d]^{f_P} \cartesien &\X \ar[d]^f \\
U \ar[r]_l \ar[ru]_x \ar@/^6mm/[u]^{t_x}& \pic \ar[r]& S}$$
On a un faisceau inversible $s(l)$ sur $\X\times_S U$. On pose $\Pc(x)=t_x^*s(l)$. Les isomorphismes de transition sont définis de manière évidente et l'on obtient ainsi un faisceau inversible $\Pc$ sur $\X\times_S \pic$ dont il ne reste plus qu'à montrer qu'il représente $\pic$. Soit $l$ un élément de $\pic(U)$. Il faut montrer que
$l=\pi(\Pc_{|_{\X\times_S U}})$. Il est évident, vu la construction de $\Pc$, que $\Pc_{|_{\X\times_S U}}$ est isomorphe à $s(l)$. Comme $s$ est une section de $\pi$, on a bien le résultat attendu.

Supposons que $\pic$ soit représentable par un schéma et montrons que e) implique c). Sous l'hypothèse e), il existe par définition du foncteur $P_{\X/S}$ un faisceau inversible $\Pc$ sur $\X\times_S \pic$ qui induit l'élément identité de $\pic(\pic)$. Il est clair que $\Pc$ est universel.

Enfin dans le cas où $f$ a une section $\sigma$, montrons que la condition b) est vérifiée. Le morphisme $f^*$ admet dans ce cas une rétraction $\sigma^*$.
$$\xymatrix{0 \ar[r]& \bgm \ar[r]^-{f^*} &\champic(\X/S) \ar@/^5mm/[l]^{\sigma^*}\ar[r]^-{\pi}& \pic \ar[r]&0}.$$
On vérifie alors facilement que le 1-morphisme
${\Mc} \mapsto {(\sigma^*\Mc,\pi(\Mc))}$
est un isomorphisme de $\champic(\X/S)$ dans $\bgm\times_S \pic$.
\end{demo}

\begin{sousremarque}\rm
On se donne à la fois une section $\sigma$ de $f$ et un faisceau inversible universel $\Pc$ sur $\X\times_S \pic$ (on ne suppose pas $\pic$ représentable). Alors l'isomorphisme $\bgm\times_S \pic \fleche \champic(\X/S)$ induit par $\Pc$ est un quasi-inverse de $(\sigma^*, \pi)$ si et seulement si $\sigma^*\Pc$ est trivial. En particulier on obtient toujours un quasi-inverse en \og rigidifiant \fg\ $\Pc$ le long de $\sigma$, c'est-à-dire en remplaçant $\Pc$ par $\Pc\otimes (f^*\sigma^*\Pc)^{-1}$.
\end{sousremarque}

\begin{souscor}
\label{foncteur_pic_repres_implique_champic_alg}
Soit $\X$ un $S$-champ algébrique
localement de présentation finie et cohomologiquement plat en dimension zéro. On suppose que le foncteur de Picard $\pic$ est représentable par un espace algébrique. Alors le champ de Picard $\champic(\X/S)$ est algébrique.
\end{souscor}
\begin{demo}
Pour un $S$-champ, être algébrique est une condition locale pour la topologie \emph{fppf}. On peut donc supposer (\ref{remarque_section}) que $f$ a une section et il suffit d'appliquer la proposition précédente.
\end{demo}

On peut donner une autre interprétation de l'équivalence entre b) et c).

\begin{souscor}
Soit $\omega$ l'élément du groupe $H^2(\Pic_{\X/S}, \gm)$ défini par la $\gm$-gerbe $\champic(\X/S)$. Alors $\omega$ est une classe d'obstruction à l'existence d'un faisceau inversible universel sur $\X\times_S \pic$ (au sens où un tel faisceau existe si et seulement si $\omega$ est nul).
\end{souscor}

\section{Représentabilité}

\subsection{Un champ algébrique au sens d'Artin}
\label{algebricite}

Nous donnons dans cette partie une démonstration directe, en termes de faisceaux
inversibles, du théorème d'Aoki cité en introduction.

\smallskip
\noindent
{\sc Déformations de faisceaux inversibles}

Commen\c cons par une petite remarque d'alg\`ebre que nous utiliserons
abondamment par la suite et que, par commodit\'e, nous \'enon\c cons sous la
forme d'un lemme.

\begin{souslem}
\label{defm_lemme_alg}
Soient $A$ un anneau et $I$ un id\'eal de carr\'e nul de $A$. On note $\pi$ la
projection canonique $\pi : A \fleche A/I$.\\
1) Le morphisme de groupes ab\'eliens $\pi^{\times} : A^{\times} \fleche
(A/I)^{\times}$ induit par $\pi$ est surjectif.\\
2) L'application $x\mapsto 1+x$ induit un isomorphisme de groupes ab\'eliens de
$I$ sur $\Ker \pi^{\times}$. $\square$
\end{souslem}

Soit $\X$ un champ alg\'ebrique sur un sch\'ema $T$ et soit $\Lc$ un faisceau
inversible sur $\X$. On consid\`ere une immersion ferm\'ee
$$\xymatrix{i : \X \ar[r] & \Xt}$$
d\'efinie par un id\'eal quasi-coh\'erent $I$ de $\Xt$ de carr\'e nul.

\begin{sousremarque} \rm
Si $\X$ et $\Xt$ sont des champs de Deligne-Mumford, le morphisme $i$ induit une équivalence de sites entre les sites étales de $\X$ et de $\Xt$, ce qui permet d'identifier les topos étales de $\X$ et de $\Xt$. Il est alors évident que la catégorie des $\Oc_{\X}$-modules quasi-cohérents est équivalente à la catégorie des $\Oc_{\Xt}$-modules quasi-cohérents annulés par $I$. Lorsque l'on travaille avec des champs d'Artin (donc avec leurs sites lisses-étales) il faut être plus prudent. En effet le foncteur naturel du site lisse-étale de $\Xt$ vers celui de $\X$ n'est même plus fidèle si bien que ces derniers ne sont \emph{a priori} pas équivalents. Cependant, la descente fidèlement plate des modules quasi-cohérents (\emph{cf.} par exemple \cite{LMB}~(13.5)) nous permet encore d'identifier la catégorie des $\Oc_{\X}$-modules à la catégorie des $\Oc_{\Xt}$-modules annulés par $I$. En particulier, l'idéal $I$ peut être vu comme un $\Oc_{\X}$-module. 
\end{sousremarque}

On note $\defm(\Lc)$ la cat\'egorie des d\'eformations de $\Lc$ \`a $\Xt$
d\'efinie de la mani\`ere suivante. Un objet de $\defm(\Lc)$ est un couple
$(\Lt,\lambda)$ o\`u $\Lt$ est un faisceau inversible sur $\Xt$ et $\lambda$ est
un isomorphisme $\xymatrix@C=1pc{\lambda : i^*\Lt\ar[r]^-{\sim}&\Lc}$. Un morphisme de
$(\Lt,\lambda)$ vers $(\widetilde{\Mc},\mu)$ est un isomorphisme $\alpha :
\raisebox{.7ex}{\xymatrix@C=1pc{\Lt\ar[r]^-{\sim}&\widetilde{\Mc}}}$ tel que
$\mu \circ i^*\alpha=\lambda$. On
note $\ov{\defm(\Lc)}$ l'ensemble des classes d'isomorphie de $\defm(\Lc)$. Le théorème ci-dessous donne une description de la catégorie des déformations de $\Lc$ à $\Xt$. Un petit calcul de complexe cotangent (voir \cite{Brochard_these} lemme 3.2.3) permet de le déduire du théorème analogue démontré par Aoki pour les déformations de morphismes de champs algébriques (cf. \cite{Aoki_Hom}~2.1). Il est cependant nettement plus facile à obtenir que ce dernier, comme nous pouvons le voir ci-dessous.

\begin{sousthm}
\label{thm_defm_fi}
\begin{enumerate}
\item[(1)] Il existe un \'el\'ement $\omega\in H^2(\X,I)$
dont l'annulation \'equivaut \`a l'existence d'une d\'eformation de $\Lc$ \`a
$\Xt$.
\item[(2)] Si $\omega=0$, alors $\ov{\defm(\Lc)}$ est un torseur sous
$H^1(\X,I)$.
\item[(3)] Si $(\Lt,\lambda)$ est une d\'eformation de $\Lc$, son groupe
d'automorphismes est isomorphe \`a $H^0(\X,I)$.
\end{enumerate}
\end{sousthm}
\begin{demo}
Commençons par le troisième point.
Par d\'efinition un automorphisme de $(\Lt,\lambda)$ est un automorphisme
$\varphi$ de $\Lt$ qui induit l'identité de $i^*\Lt$. Autrement dit,
$\Aut_{\defm(\Lc)}(\Lt,\lambda)$ est le noyau du morphisme $\Aut(\Lt) \fleche
\Aut(\Lc)$ induit par $i^*$, c'est-à-dire le noyau du morphisme
$$H^0(\Xt,\Oc_{\Xt})^{\times} \flechelongue
\left(\frac{H^0(\Xt,\Oc_{\Xt})}{H^0(\Xt,I)}\right)^{\times}.$$
D'apr\`es le lemme (\ref{defm_lemme_alg}) 2), il est donc isomorphe \`a $H^0(\Xt,I)$ via
l'application $x\mapsto 1+x$.

Montrons maintenant que $\ov{\defm(\Lc)}$ est isomorphe \`a l'ensemble 
$\Pic_{[\Lc]}(\Xt)$ des \'el\'ements de $\Pic(\Xt)$ qui
sont envoy\'es sur la classe $[\Lc]$ de $\Lc$ dans $\Pic(\X)$. On a une application naturelle
$\phi : \ov{\defm(\Lc)} \fleche \Pic_{[\Lc]}(\Xt)$ qui
\`a une d\'eformation $(\Lt,\lambda)$ associe la classe de $\Lt$ dans
$\Pic(\Xt)$. Elle est clairement surjective. De plus si $(\Lt,\lambda)$ et
$(\Mt,\mu)$ sont tels que $[\Lt]=[\Mt]$, il existe un isomorphisme
$\alpha : \Lt \fleche \Mt$. On va montrer que l'on peut choisir $\alpha$ de
telle sorte que $i^*\alpha$ soit égal à $\mu^{-1} \circ \lambda$. Il suffit pour cela de voir
que $\Aut(\Lt) \fleche \Aut(i^*\Lt)$ est surjectif : il
n'y aura plus alors qu'\`a corriger $\alpha$ par un automorphisme convenablement
choisi de $\Lt$. Or ce morphisme de groupes s'identifie \`a
$\Aut(\Oc_{\Xt})\fleche\Aut(\Oc_{\X})$, qui est surjectif en vertu de
(\ref{defm_lemme_alg}), 1).

Donc $\ov{\defm(\Lc)}$ est
naturellement un torseur sous le noyau du morphisme $\Pic(\Xt)\fleche \Pic(\X)$.
Calculons ce noyau. On a une suite exacte de faisceaux quasi-coh\'erents sur
$\Xt$ :
$$0 \flechelongue I \flechelongue \Oc_{\Xt} \flechelongue i_*\Oc_{\X} \flechelongue 0.$$
Elle induit \emph{via} l'application exponentielle
une suite exacte de faisceaux ab\'eliens sur $\Xt$:
$$0 \flechelongue I \flechelongue \Oc_{\Xt}^{\times} \flechelongue i_*\Oc_{\X}^{\times} \flechelongue 0.$$
La suite exacte longue de cohomologie associ\'ee nous donne :
$$H^1(\Xt, I) \flechelongue H^1(\Xt,\Oc_{\Xt}^{\times}) \flechelongue
H^1(\Xt, i_*\Oc_{\X}^{\times}) \flechelongue H^2(\Xt, I).$$
Or d'apr\`es \ref{pic_egal_h1}, le groupe $H^1(\Xt,\Oc_{\Xt}^{\times})$ est
isomorphe \`a $\Pic(\Xt)$. De plus le lemme \ref{coh_et_ext_inf}
fourni en annexe
montre que $H^1(\Xt, i_*\Oc_{\X}^{\times})\simeq H^1(\X,\Oc_{\X}^{\times})\simeq
\Pic(\X)$. Montrons enfin que la premi\`ere fl\`eche de la suite exacte longue
ci-dessus est injective. Il suffit pour cela de voir que le morphisme
$H^0(\Xt,\Oc_{\Xt}^{\times}) \fleche
H^0(\Xt, i_*\Oc_{\X}^{\times})$ est surjectif, ce qui r\'esulte encore du lemme
(\ref{defm_lemme_alg}). On obtient donc la suite exacte 
$$0\flechelongue H^1(\Xt, I) \flechelongue \Pic(\Xt) \flechelongue
\Pic(\X) \flechelongue H^2(\Xt, I),$$
ce qui ach\`eve notre d\'emonstration.
\end{demo}

\begin{sousremarque} \rm
\label{rem_noyau}
Nous avons montr\'e au passage que le noyau de $\Pic(\Xt)\fleche \Pic(\X)$
s'identifie \`a $H^1(\Xt, I)$.
\end{sousremarque}

Nous pouvons maintenant utiliser le théorème~\cite{Artin_Versal_defm}~(5.3) d'Artin pour (re)démontrer le théorème \ref{algebricite_champ_Picard}.

%

\medskip
\noindent
{\bf Nouvelle démonstration du théorème \ref{algebricite_champ_Picard}}
On note $\Pc$ le champ $\champic(\X/S)$. Si $A$ est un anneau sur $S$, on note $\Pc(A)$ la catégorie fibre de $\Pc$ au-dessus de $\Spec A$. On peut supposer que $S$ est le spectre d'un anneau $R$. On peut même supposer que $R$ est de type fini sur $\Z$ grâce à la proposition~(4.18)~(ii) de~\cite{LMB}. On vérifie une à une les conditions de~\cite{Artin_Versal_defm}~(5.3).

\begin{etape}{Présentation finie}
Il faut montrer que si $A$ est une limite inductive d'anneaux noethériens $\lind A_i$, alors le foncteur $$\disp \lind \Pc(A_i) \flechelongue \Pc(A)$$
est une équivalence de catégories. Ceci résulte immédiatement de la proposition~(4.18) de~\cite{LMB} compte tenu du fait que $\Pc(A)$ s'identifie à $\Hom_A(\X_A, \bgm)$ et que $\bgm$ est de présentation finie. (Ici un indice comme $\X_A$ d\'esigne le changement de base $\X\times_S \Spec A$.)
\end{etape}

\begin{etape}{Première condition de Schlessinger (appelée (S1') dans \cite{Artin_Versal_defm})}
On se donne un carré cartésien d'anneaux noethériens
$$\xymatrix{B' \ar[r]\ar[d] \cartesien &A' \ar[d]^p\\ B\ar[r] & A}$$
où $A'$ est une extension infinitésimale de $A$. Il faut montrer que pour tout objet $a$ de $\Pc(A)$ le foncteur
$$\Pc_a(B') \flechelongue \Pc_a(A')\times\Pc_a(B)$$
est une équivalence de catégories, où pour toute $A$-algèbre $R'$, $\Pc_a(R')$ désigne la catégorie des objets de $\Pc(R')$ au-dessus de $a$ avec les morphismes au-dessus de l'identité de $a$. C'est une conséquence immédiate du théorème suivant.

\begin{sousthm}
\label{equiv_cat_fi_pour_defo}
Si $\X$ est un champ alg\'ebrique
\emph{plat} sur $S$, le foncteur naturel
$$F(\X) : \inv(\X_{B'}) \flechelongue \inv(\X_{A'}) \times_{\inv(\X_{A})} \inv(\X_{B}).$$
est une \'equivalence de cat\'egories (où $\inv(\X)$ est la catégorie des faisceaux inversibles sur $\X$).
\end{sousthm}
\begin{demo}

Commen\c cons par l'observation suivante. Soient $\U$ un $S$-champ
alg\'ebrique, $U^1$ un
$S$-espace alg\'ebrique, et soit $U^1 \fleche \U$ un \'epimorphisme. On note $U^2=U^1\times_{\U} U^1$ et
$U^3=U^1\times_{\U} U^1\times_{\U} U^1$.
Alors la cat\'egorie $\inv(\U)$ est \'equivalente \`a la
cat\'egorie ${\rm Desc}(U^1,\U)$ décrite dans~\ref{descente_fidelement_plate}.
De plus la construction de ${\rm Desc}(U^1,\U)$ commute aux produits fibr\'es de
cat\'egories, de sorte que si les foncteurs naturels
$F(U^i)$ pour $i=1, 2, 3,$
sont des \'equivalences de cat\'egories, il en est de m\^eme du foncteur
$F(\Uc)$.

Revenons \`a la d\'emonstration proprement dite du th\'eor\`eme~\ref{equiv_cat_fi_pour_defo}.
\vskip .2cm
\emph{Premier cas : $\X$ est un $S$-sch\'ema affine $\Spec R'$.}\\
Alors l'assertion r\'esulte du th\'eor\`eme~\cite{Ferrand}~2.2 et du lemme
suivant.

\begin{souslem}
Si $R'$ est une $R$-alg\`ebre plate, alors le carr\'e obtenu par changement de
base
$$\xymatrix{B'\otimes_R R'\ar[r]\ar[d] &
A'\otimes_R R'\ar[d]^{p_{R'}} \\
B\otimes_R R' \ar[r] &A\otimes_R R'}$$
est encore cart\'esien, et $p_{R'}$ est surjectif. $\square$
\end{souslem}

\emph{$2^{\textrm{\`eme}}$ cas : $\X$ est un $S$-sch\'ema qui est union
disjointe de sch\'emas affines.}\\
On a $\X=\disp \coprod_i X_i$, o\`u chaque $X_i$ est affine. Alors
$\inv(\X)=\disp \prod_i \inv(X_i)$ universellement, et l'assertion est vraie
pour chacun des $S$-sch\'emas $X_i$, donc il suffit de montrer que les produits
fibr\'es de cat\'egories commutent aux produits arbitraires, ce qui est purement
formel et \'evident.

\emph{$3^{\textrm{\`eme}}$ cas : $\X$ est un sch\'ema s\'epar\'e.}\\
Soit $X$ l'union disjointe d'une famille d'ouverts affines recouvrant $\X$.
Alors $X$, $X\times_{\X} X$, et $X\times_{\X} X\times_{\X} X$ sont unions
disjointes de sch\'emas affines, donc d'apr\`es la remarque qui amor\c cait
notre d\'emonstration, l'assertion r\'esulte du second cas.

\emph{$4^{\textrm{\`eme}}$ cas : $\X$ est un $S$-sch\'ema quelconque.}\\
De m\^eme soit $X$ la somme des ouverts d'un recouvrement
affine de $\X$. Alors les schémas $X$, $X\times_{\X} X$ et $X\times_{\X}
X\times_{\X} X$ sont s\'epar\'es, donc il suffit d'appliquer le troisi\`eme cas.

\emph{$5^{\textrm{\`eme}}$ cas : $\X$ est un $S$-espace alg\'ebrique.}\\
On choisit un morphisme \'etale surjectif $X\fleche \X$ o\`u $X$ est un
sch\'ema. Alors $X\times_{\X} X$ et $X\times_{\X}
X\times_{\X} X$ sont encore des sch\'emas, de sorte que l'assertion r\'esulte du
cas pr\'ec\'edent.

\emph{$6^{\textrm{\`eme}}$ cas : $\X$ est un $S$-champ alg\'ebrique.}\\
On choisit une pr\'esentation $P:X\fleche \X$ par un espace alg\'ebrique. Alors 
$X\times_{\X} X$ et $X\times_{\X}
X\times_{\X} X$ sont encore des espaces alg\'ebriques, donc le r\'esultat se d\'eduit du cinqui\`eme cas.
\end{demo}
\end{etape}

\begin{etape}{Seconde condition de Schlessinger (appelée (S2) dans \cite{Artin_Versal_defm})}
Soient $A$ une $R$-algèbre noethérienne et $M$ un $A$-module de type fini (on peut supposer $A=R$). On note $A[M]$ l'extension infinitésimale de $A$ par $M$. Si $a$ est un objet de $\Pc(A)$, on note
$$D_a(M)=\Pc_a(A[M])/\sim.$$
Il faut vérifier que $D_a(M)$ est un $A$-module de type fini. Or $D_a(M)$ n'est autre que l'ensemble appelé $\ov{\defm(a)}$ ci-dessus, et l'on a vu que cet ensemble est naturellement un torseur sous $H^1(\Xt,I)$, où $\Xt$ désigne le produit fibré $\X\times_S \Spec A[M]$. Il suffit de vérifier que $H^1(\Xt,I)$ est un $A$-module de type fini. On vérifie facilement que l'idéal $I$ qui définit l'extension infinitésimale $\xymatrix@C=1pc{\X \ar@{^(->}[r]&\Xt}$ s'identifie à $i_*f^*M$. D'après~\ref{coh_et_ext_inf} le module $H^1(\Xt,I)$ est donc isomorphe à $H^1(\X, f^*M)$. L'assertion résulte maintenant de la
finitude de la cohomologie des faisceaux cohérents sur un champ algébrique
propre et localement noethérien (\cite{Olsson_lemme_chow} théorème~(1.2), ou
\cite{Faltings_finiteness_coco}).
\end{etape}

\begin{etape}{Fin de la condition (1) du théorème (5.3)}
Il faut vérifier que si $a$ est un objet de $\Pc(A)$ et si $M$ est un $A$-module de type fini, alors l'ensemble $\Aut_a(M)$ des automorphismes infinitésimaux, défini par
$$\Aut_a(M):=\Ker(\Aut_a(A[M]) \fleche \Aut_a(A))$$
est un $A$-module de type fini. D'après le troisième point du théorème~\ref{thm_defm_fi} ci-dessus, $\Aut_a(M)$ est naturellement isomorphe à $H^0(\X, f^*M)$, qui est bien de type fini d'après le théorème de finitude de Faltings-Olsson.
\end{etape}

\begin{etape}{Commutation aux limites projectives}
Il faut vérifier que si $A$ est un anneau noethérien local complet d'idéal maximal $\mgo$, alors le foncteur naturel
$$\Pc(A) \flechelongue \lpro \Pc(A/\mgo^{n+1})$$
est une équivalence de catégories.
C'est une conséquence immédiate du théorème~\cite{Olsson_lemme_chow}~(1.4) d'Olsson selon lequel le foncteur naturel $\disp \Coh(\X) \fleche \lpro \Coh(\X_{A/\mgo^{n+1}})$
est une équivalence.
\end{etape}

\begin{etape}{Conditions sur la théorie des déformations}
Il faut vérifier que les modules $\Aut_a(M)$, $D_a(M)$ et $\Oc_a(M)$ vérifient les conditions~(4.1) de~\cite{Artin_Versal_defm}, c'est-à-dire :
\begin{itemize}
\item[(i)] qu'ils commutent à la localisation étale, \emph{i.e.} si $A\fleche B$ est étale alors $\Aut_b(M\otimes_A B)$ (où $b$ est l'image de $a$) est isomorphe à $\Aut_a(M)\otimes_A B$, et \emph{idem} pour les autres ;
\item[(ii)] qu'ils sont compatibles aux limites projectives, \emph{i.e.} si
$\mgo$ est un idéal maximal de $A$ et si $\hat{A}$ est le complété de $A$ relativement à $\mgo$, alors le morphisme naturel
$$\Aut_a(M)\otimes_A \hat{A} \flechelongue \Aut_{\hat{a}}(M\otimes_A \hat{A})$$
est un isomorphisme, etc\dots
\item[(iii)] (constructibilité) Il existe un ouvert dense de points de type fini $p\in \Spec A$ tels que
$$\Aut_a(M)\otimes_A \kappa(p) \flechelongue \Aut_{a_p}(M\otimes_A \kappa(p))$$
soit un isomorphisme, etc\dots
\end{itemize}
Or d'après ce qui précède, ces modules s'identifient respectivement à $H^0(\X,f^*M)$, $H^1(\X, f^*M)$ et $H^2(\X,f^*M)$. La condition (i) résulte de la proposition~(\ref{lemme_coh_et_chgt_base}). La condition (ii) est une conséquence de~\cite{Olsson_lemme_chow}.

Enfin pour la condition (iii), nous devons montrer que si $A_0$ est une $R$-algèbre intègre de
type fini, il existe un ouvert non vide $U$ de $\Spec A_0$ tel que pour tout
point fermé $s$ dans $U$ le $\kappa(s)$-espace vectoriel
$D(k,M\otimes_{A_0}\kappa(s),{\xi_0}_s)$ soit isomorphe à
$D(A_0,M,\xi_0)\otimes_{A_0}\kappa(s)$. En d'autres termes on demande que pour
un tel point $s$, le morphisme canonique
\begin{equation}
\label{morph_condition_4c}
H^1(\X_{A_0},f_{A_0}^*M)\otimes_{A_0}\kappa(s) \fleche
H^1(\X_{s},f_{s}^*(M\otimes_{A_0}\kappa(s)))
\end{equation}
soit un isomorphisme. Montrons d'abord qu'il existe un ouvert non vide de $\Spec
A_0$ sur lequel $M$ est libre de rang fini. On note $K$ le corps des fractions
de $A_0$. Comme $M$ est de type fini, $M_K$ est un $K$-espace vectoriel de
dimension finie. Soient $\nlist{x}$ des éléments de $M$ dont les images dans
$M_K$ engendrent $M_K$. Soit $(y_1, \dots, y_n)$ une famille génératrice de $M$.
Il existe $f\in A_0$ non nul tel que pour tout $j$ on puisse écrire $y_j$ comme
combinaison linéaire des $\frac{x_i}{f}$ dans $M_K$. Donc quitte à localiser par
$f$ on peut supposer que la famille $\nlist{x}$ engendre $M$. C'est alors
aussitôt une base de $M$. En effet, on a un morphisme surjectif $\varphi : A^n
\fleche M$ qui induit un isomorphisme $\varphi_K : A^n_K\fleche M_K$, et le
morphisme de localisation $A^n \fleche K^n$ est injectif puisque $A$ est
intègre, donc $\varphi$ est aussi injectif. Maintenant, $M$ est libre de rang
fini sur $\Spec A_0$, donc $f_{A_0}^*M\simeq (\Oc_{\X\times_S \Spec A_0})^n$ est
cohérent et plat sur $\Spec A_0$ (car $f$ est plat). Donc le résultat de Mumford
(\cite{Varietes_abeliennes} paragraphe~5) généralisé par Aoki aux champs
algébriques (\cite{Aoki_Hom}, théorème~(A.1)) s'applique et il existe un ouvert
non vide $U$ de $\Spec A_0$ tel que, pour tout point $s$ de $U$, le morphisme
\eqref{morph_condition_4c} ci-dessus soit un isomorphisme. 
\end{etape}

\begin{etape}{Quasi-séparation du morphisme diagonal}
Soient $A$ une $R$-algèbre de type fini et $\Lc$ un faisceau inversible sur $\X_A$. Soit $\varphi$ un automorphisme de $\Lc$ qui induit l'identité dans $\Pc(k)$ pour un ensemble dense de points $A\fleche k$. Il faut montrer que $\varphi$ est l'identité sur un ouvert non vide de $\Spec A$. Remarquons tout de suite qu'en fait, en regardant bien l'article~\cite{Artin_Versal_defm} d'Artin, on voit qu'on peut supposer $A$ intègre. Soit $X \fleche \X$ une présentation de $\X$, avec $X$ affine d'anneau $B$. On note encore $\varphi$ l'élément de $\Gamma(\X_A, \Oc_{\X_A})^{\times}$ correspondant à l'automorphisme $\varphi$. Il suffit clairement de montrer que $\varphi_{|_X}$ vaut 1. Notons $b$ l'élément $\varphi_{|_X}-1$ de $B$, et $H$ l'ensemble des points $p$ de $\Spec A$ tels que $b_p=0$. Par hypothèse $H$ est dense dans $\Spec A$. Par ailleurs, d'après~\cite{EGA4_3}~(9.4.6), $H$ est constructible. Donc le point générique $\xi$ de $\Spec A$ est dans $H$ et $b_{\xi}=0$. Comme $B$ est plat sur $A$, ceci prouve que $b$ est nul. $\square$
\end{etape}

\subsection{Quasi-séparation}
\label{quasi-separation}

Nous montrons maintenant que, sous les mêmes hypothèses qu'en \ref{algebricite_champ_Picard}, le champ de Picard est algébrique au sens de \cite{LMB}, c'est-à-dire en un sens un peu plus fort que celui d'Artin. La différence porte uniquement sur les conditions de finitude imposées au morphisme diagonal : les champs algébriques de \cite{LMB} sont supposés quasi-séparés, autrement dit leur diagonale est quasi-compacte (et séparée). La question de la quasi-séparation étant parfois délicate, nous proposons ci-dessous un critère général de quasi-séparation, inspiré des techniques employées par Artin dans \cite{Global_Analysis_1}.

\begin{sousprop}
Soient $S$ un sch\'ema localement noeth\'erien et $\X$ un $S$-champ
alg\'ebrique (au sens d'Artin)
localement de pr\'esentation finie. On suppose que les deux conditions suivantes
sont remplies.

(i) Pour tout $U\in\ob\aff$ et tout $x\in\ob \X_U$, le morphisme
$\Auts(x)\fleche U$ est quasi-compact.

(ii) Soit $U\in\ob\aff$ int\`egre, et soient $x, y\in \ob \X_U$. On suppose
qu'il existe un ensemble dense de points $t$ de $U$, tels qu'il existe une
extension $L(t)$ de $\kappa(t)$ telle que $x_{L(t)} \simeq y_{L(t)}$. Alors $x$
et $y$ sont isomorphes sur un ouvert dense de $U$.

Alors $\X$ est quasi-s\'epar\'e (donc est algébrique au sens de \cite{LMB}).
\end{sousprop}
\begin{demo}
Soit $U\fleche \X\times_S\X$ une pr\'esentation de
$\X\times_S\X$, o\`u $U$ est un sch\'ema localement noeth\'erien. Alors par
descente fidèlement plate, pour montrer que $\Delta: \X \fleche \X\times_S
\X$ est quasi-compact, il suffit de montrer que le morphisme obtenu par
changement de base
$$\xymatrix{I\ar[r] \ar[d] \cartesien &U\ar[d]^{(x,y)}\\
  \X\ar[r]& \X\times_S \X}$$
est quasi-compact. Comme c'est une question locale sur $U$, on peut supposer $U$
affine.

Soit $G=\Auts(y)$. L'action \`a gauche naturelle de $G$ sur $I$ fait clairement
de $I\fleche U$ un pseudo-torseur\footnote{Autrement dit la fl\`eche
$I\fleche U$ est invariante sous $G$, et le morphisme naturel $G\times I \fleche
I\times_U I$ est un isomorphisme.} sous $G$. Par r\'ecurrence noeth\'erienne sur
les ferm\'es de $U$, on peut supposer que pour tout ferm\'e strict $F$ de $U$,
$I\times_U F$ est quasi-compact. Il suffit maintenant de trouver un ouvert non
vide $V$ de $U$ tel que $I\times_U V$ soit quasi-compact. En effet, si $V$ est
un tel ouvert, alors en notant $F=U\setminus V$, l'espace alg\'ebrique
$I'=(I\times_U V )\coprod (I\times_U F)$ est quasi-compact, et le morphisme
$I'\fleche I$ obtenu par changement de base \`a partir de $V\coprod F \fleche U$
est surjectif, donc $I$ est quasi-compact.

Posons $Z=\overline{f(I)}$. Soit $W$ un ouvert affine int\`egre de
$Z_{\text{r\'ed}}$. Alors le morphisme $f_W: I\times_U W \fleche W$ est
dominant, et la condition (ii) signifie pr\'ecis\'ement que dans ce cas, quitte
\`a choisir $W$ plus petit, $f_W$ a une section. Ceci prouve que $I\times_U W$
est un torseur sous $G$, donc est quasi-compact. Ceci \'etant, on peut \'ecrire
$W=Z_{\text{r\'ed}}\cap V=Z_{\text{r\'ed}}\times_U V$ o\`u $V$ est un ouvert
(non vide) de $U$. On v\'erifie sans peine que le morphisme induit $I\times_U W
\fleche I\times_U V$ est surjectif, ce qui montre que $I\times_U V$ est
quasi-compact, et ach\`eve la d\'emonstration.
\end{demo}

Nous sommes maintenant en mesure de montrer que le champ de Picard est quasi-séparé. Ceci résulte immédiatement du critère ci-dessus et du lemme suivant.

\begin{souslem}
\label{lemme_locale_separation_du_foncteur_de_Picard}
Soient $A$ un anneau int\`egre, $T=\Spec A$, et $\X \fleche T$ un $T$-champ
alg\'ebrique propre, plat, de pr\'esentation finie, et cohomologiquement
plat en dimension z\'ero. Soit $\Lc$ un faisceau inversible sur $\X$.
\begin{itemize}
\item[(i)] Si $A$ est un anneau local (en particulier si $A$ est un anneau de
valuation discr\`ete), de corps
r\'esiduel $k$ et de corps des fractions $K$, alors $\Lc$
est trivial si et seulement si $\Lc_k$ et $\Lc_K$ le sont.
\item[(ii)] S'il existe un ensemble dense $\Sc$ de points de $T$ tel que pour tout
$t\in \Sc$, $\Lc_t$ soit trivial, alors il existe un ouvert non vide $U$ de $T$
tel que $\Lc_U$ soit trivial.
\end{itemize}
\end{souslem}
\begin{demo}
D'apr\`es l'annexe A de~\cite{Aoki_Hom} les corollaires~1 et 2 du paragraphe 5
de~\cite{Varietes_abeliennes} sont encore valables. Ils nous apprennent que la
fonction $$t\longmapsto d(t):=\dim_{\kappa(t)} H^0(\X_t, \Lc_t)$$ est semi-continue
sup\'erieurement, et que si elle est constante sur
$T$, alors le morphisme naturel
$$(f_*\Lc)\otimes_{\Oc_T} \kappa(t) \flechelongue H^0(\X_t, \Lc_t)$$ est un isomorphisme
pour tout $t\in T$. De plus en un point o\`u $\Lc_t$ est trivial, on a $H^0(\X_t,
\Lc_t)=H^0(\X_t, \Oc_{\X_t})\simeq \kappa(t)$ par platitude cohomologique, donc
$d(t)=1$.

\medskip
(i) Ici comme $\Lc_K$ et $\Lc_k$ sont triviaux, la fonction $d$ vaut 1 au point
g\'en\'erique et au point ferm\'e, donc elle est constante
\'egale \`a 1 sur $T$. D'apr\`es ce qui pr\'ec\`ede, on a en particulier un
isomorphisme $H^0(T, f_*\Lc)\otimes_A k=H^0(\X,\Lc)\otimes_A k\fleche
H^0(\X_k,\Lc_k)$. Le faisceau $\Lc_k$ \'etant trivial, il a une section globale
partout non nulle sur $\X_k$. D'apr\`es l'isomorphisme pr\'ec\'edent,
cette section provient d'une
section globale $s$ de $\Lc$ (comme $k$ est un quotient de $A$, les \'el\'ements
de $H^0(\X,\Lc)\otimes_A k$ peuvent tous s'\'ecrire sous la forme $s\otimes 1$).
Alors $s$ est non nulle en tout point $x$ de $\X$ d'image $\Spec k$ dans $T$.
Or l'ensemble $C$ des points de
$\X$ o\`u $s$ est nulle est un ferm\'e de $|\X|$ (voir par exemple~\cite{Brochard_these}~2.1.2.2). Comme $f$ est propre, l'image
de $C$ est un ferm\'e de $T$, qui de plus ne contient pas le point ferm\'e.
N\'ecessairement $C$ est vide, donc $s$ est partout non nulle et $\Lc$ est trivial (\cite{Brochard_these}~2.1.2).

(ii) La fonction $d$ \'etant semi-continue sup\'erieurement, il existe un ouvert
non vide sur lequel elle est constante et quitte \`a remplacer $T$ par cet
ouvert, on peut supposer qu'elle est constante. D\`es lors pour tout $t\in T$ le
morphisme naturel
$$H^0(\X,\Lc)\otimes_A \kappa(t) \flechelongue H^0(\X_t, \Lc_t)$$
est un isomorphisme. Soit $t\in \Sc$. Alors $\Lc_t$ est trivial donc il a une
section globale $s_t\in H^0(\X_t, \Lc_t)$ partout non nulle sur $\X_t$. D'apr\`es
l'isomorphisme pr\'ec\'edent, elle provient d'un \'el\'ement du groupe
$H^0(\X,\Lc)\otimes_A \kappa(t)$ que l'on peut \'ecrire sous la forme
$s\otimes(\frac1f)$ o\`u $f\in A$ n'est pas dans l'id\'eal premier correspondant
\`a $t$. Quitte \`a remplacer $T$ par $D(f)$, on peut supposer $f=1$ et on a
donc trouv\'e une section $s\in H^0(\X,\Lc)$ dont la r\'eduction $s_t$ \`a $\X_t$
est partout non nulle.

Maintenant l'ensemble $C$ des points de $\X$ o\`u $s$ est nulle est un ferm\'e
de $|\X|$, et son image est
un ferm\'e $F$ de $T$, qui ne contient pas $t$. Son compl\'ementaire
$U=T\setminus F$ est alors un ouvert non vide de $T$, et la section $s_U$
induite par $s$ sur $\X\times_T U$ est partout non nulle sur $\X\times_T U$, de
sorte que $\Lc_U$ est trivial.
\end{demo}


\begin{souscor}
Si $\X$ est propre, plat et cohomologiquement plat sur $S$, alors le foncteur de Picard $\pic$ est un espace algébrique localement séparé. $\square$
\end{souscor}

\section{Autour de l'élément neutre}

\subsection{Lissité et dimension}
Soient $k$ un corps, $P$ un $k$-espace alg\'ebrique, et $y : \Spec k \fleche P$
un $k$-point de $P$. Soit $k[\eps]=k[X]/(X^2)$. On note $i$ l'immersion ferm\'ee de $\Spec k$ dans $\Spec k[\eps]$, et  $T_yP$ l'espace tangent \`a $P$ en $y$.
$$T_yP = \{\varphi : \Spec k[\eps] \fleche P \ | \ \varphi\circ i = y\}$$
Si $f : X\fleche P$ est un morphisme de $k$-espaces alg\'ebriques, et si $x$ est
un $k$-point de $X$, on note $T_xf : T_xX \fleche T_{f(x)}P$ le 
morphisme naturel ${\varphi} \mapsto {f\circ \varphi.}$
Nous rappelons ci-dessous quelques r\'esultats bien utiles.

\begin{souslem}
\label{lemme_espace_tangent}
\begin{itemize}
\item[1)] Soit $f : X \fleche P$ un morphisme formellement lisse (resp.
formellement non ramifi\'e, formellement \'etale) de $k$-espaces alg\'ebriques
et soit $x : \Spec k \fleche X$ un point $k$-rationnel de $X$. Alors
$T_xf : T_xX \fleche T_{f(x)}P$ est surjective (resp. injective, resp. un
isomorphisme).
\item[2)] Soient $X$ un $k$-espace alg\'ebrique localement de type fini, $x$ un
$k$-point de $X$, et $L$ une
extension de $k$. On note $X_L$ le $L$-espace alg\'ebrique obtenu par changement
de base, et $x_L$ le point de $X_L$ induit par $x$.
Alors le morphisme naturel
$$(T_xX)\otimes_k L \flechelongue T_{x_L}(X_L)$$
est un isomorphisme.
\end{itemize}
\end{souslem}
\begin{demo}
Le point 1) résulte immédiatement des définitions. Pour le second point, on se ramène par~\cite{Knutson}~II~(6.4) au cas où $X$ est un schéma affine d'anneau $A$. L'assertion n'est maintenant plus qu'un exercice facile d'algèbre commutative (rédigé en détail dans~\cite{Brochard_these}).
\end{demo}

Nous rappelons également le r\'esultat suivant, valable aussi pour un espace algébrique puisqu'un espace algébrique en groupes sur un corps est toujours un schéma.  

\begin{sousprop}[\cite{poly_Kleiman} 5.13 et 5.14]
\label{prop_espaces_alg_en_groupes}
Soient $P$ un $k$-schéma en groupes localement de type fini, et $e$
le $k$-point neutre. Alors
$P$ a la m\^eme dimension en tout point. De plus cette dimension est
inf\'erieure \`a $\dim_k(T_eP)$, et les propri\'et\'es suivantes sont
\'equivalentes :
\begin{itemize}
\item[(i)] $\dim P =\dim_k(T_eP)$;
\item[(ii)] $P$ est lisse en 0;
\item[(iii)] $P$ est lisse.
\end{itemize}
Elles sont v\'erifi\'ees lorsque $k$ est de caract\'eristique nulle.
\end{sousprop}

\begin{sousthm}
Soient $k$ un corps, $S$ son spectre, et $\X$ un $S$-champ alg\'ebrique.
On note $\Pic_{\X/k}$ le foncteur de Picard relatif $\piceth$ et
on suppose qu'il est repr\'esentable par un $S$-espace
alg\'ebrique localement de type fini.

a) Alors l'espace tangent \`a l'origine est $$T_0\Pic_{\X/k}
=H^1(\X,\Oc_{\X}).$$

b) L'espace alg\'ebrique $\Pic_{\X/k}$ a la m\^eme dimension en tout point. De
plus cette dimension est inf\'erieure \`a $\dim_k
H^1(\X,\Oc_{\X})$, et il y a \'egalit\'e si et seulement si $\Pic_{\X/k}$ est
lisse \`a l'origine. Dans ce cas, $\Pic_{\X/k}$ est lisse de dimension $\dim_k
H^1(\X,\Oc_{\X})$ partout. Il en est toujours ainsi lorsque $k$ est de caract\'eristique nulle.
\end{sousthm}
\begin{demo}
a) On note encore $P=\Pic_{\X/k}$. L'espace tangent \`a
$P$ en 0 est par définition le noyau de
$P(k[\eps]) \fleche P(k)$. On note $\Xt=\X\otimes_k k[\eps]$, et $i$
l'injection canonique de $\Xt$ dans $\X$. On voit facilement que l'id\'eal $I$ sur $\Xt$ de carr\'e
nul qui d\'efinit $\X$ comme sous-champ ferm\'e de $\Xt$ est isomorphe \`a
$i_*\Oc_{\X}$. La remarque~(\ref{rem_noyau}) montre alors que le noyau de $\Pic(\Xt)\rightarrow \Pic(\X)$ s'identifie à $H^1(\Xt,i_*\Oc_{\X}).$
Or d'apr\`es le lemme~\ref{coh_et_ext_inf}, $H^1(\Xt,i_*\Oc_{\X})$ est
isomorphe \`a $H^1(\X,\Oc_{\X})$. On a par ailleurs un carr\'e commutatif :
$$\xymatrix{ \Pic(\Xt)\ar[r] \ar[d] & \Pic(\X)\ar[d] \\
P(k[\eps]) \ar[r]& P(k),}$$
qui induit un morphisme entre les noyaux des fl\`eches horizontales, et donc,
d'apr\`es ce qui pr\'ec\`ede, un morphisme $v: H^1(\X,\Oc_{\X}) \fleche
T_0P$. Montrons que $v$ est un isomorphisme. Remarquons tout d'abord que
si $k$ est alg\'ebriquement clos, les fl\`eches verticales du carr\'e ci-dessus
sont des isomorphismes en vertu de (\ref{P_isom_Pet_si_k_alg_clos}), de sorte
que $v$ en est un aussi.
Passons maintenant au cas g\'en\'eral. Soit $\ov{k}$ une cl\^oture alg\'ebrique
de $k$. Le carr\'e ci-dessus s'envoie alors sur le carr\'e correspondant obtenu
apr\`es extension du corps de base \`a $\ov{k}$. On en d\'eduit en regardant les
noyaux un carr\'e commutatif de $k$-espaces vectoriels :
$$\xymatrix{H^1(\X,\Oc_{\X})\ar[r] \ar[d]_v &
H^1(\X_{\ov{k}},\Oc_{\X_{\ov{k}}})\ar[d]^{\ov{v}} \\
T_0P \ar[r] &T_0(P_{\ov{k}}).}$$
Il résulte du cas précédent, de~\ref{lemme_coh_et_chgt_base} et de~\ref{lemme_espace_tangent} 2) que $v\otimes_k \ov{k}$ est un isomorphisme, donc $v$ aussi. Pour b) il suffit d'appliquer~\ref{prop_espaces_alg_en_groupes}.
\end{demo}

\subsection{Construction de la composante neutre}

On rappelle le résultat suivant concernant la composante neutre d'un $k$-schéma
en groupes localement de type fini.

\begin{sousthm}[\cite{poly_Kleiman}~5.1]
Soient $k$ un corps et $G$ un $k$-schéma en groupes localement de type fini. Alors $G$ est séparé. Soit
$G^0$ la composante connexe de l'élément neutre de $G$. Alors $G^0$ est un sous-schéma en groupes ouvert et
fermé de $G$. Il est de type fini et géométriquement irréductible. De plus, la formation de $G^0$ commute
aux extensions du corps $k$.
\end{sousthm}

Naturellement, tout ceci est valable pour le foncteur de Picard $\Pic_{\X/k}$ dès
qu'il est représentable par un espace algébrique localement de type fini sur
$k$. En effet, c'est alors automatiquement un schéma en vertu d'un lemme
d'Artin (\cite{Global_Analysis_1} lemme~4.2 p.~43). En particulier, c'est le cas
dès que $\X$ est un champ algébrique propre et cohomologiquement plat en dimension zéro sur $k$.

La notion de faisceaux inversibles algébriquement équivalents se généralise très bien aux champs algébriques, et comme
dans le cas des schémas elle permet de caractériser les $k$-points du foncteur de Picard qui sont dans la composante neutre.

\begin{sousdefi}[\cite{poly_Kleiman}~5.9]
Soit $\X$ un champ algébrique sur un corps $k$. Soient $\Lc$ et $\Nc$ deux faisceaux inversibles sur $\X$. On dit que $\Lc$ et $\Nc$ sont
\emph{algébriquement équivalents} s'il existe une suite de $k$-schémas connexes de type fini $T_1, \dots, T_n$, des points
géométriques $s_i, t_i$ de $T_i$ (pour tout $i$) ayant tous le même corps, et un faisceau inversible $\Mc_i$ sur
$\X\times_k T_i$ tels que
$$\Lc_{s_1} \simeq \Mc_{1,s_1}, \quad\Mc_{1,t_1}\simeq \Mc_{2,s_2},\quad \dots, \quad\Mc_{n-1,t_{n-1}}\simeq \Mc_{n,s_n},
\quad\Mc_{n,t_n}\simeq \Nc_{t_n}$$
\end{sousdefi}

\begin{sousthm}[\cite{poly_Kleiman}~5.10]
Soit $\X$ un champ algébrique sur un corps $k$. On suppose que $\Pic_{\X/k}$ est représentable par un schéma localement
de type fini. Soit $\Lc$ un faisceau inversible sur $\X$ et soit $\lambda$ le point correspondant de $\Pic_{\X/k}$. Alors
$\Lc$ est algébriquement équivalent à $\Oc_{\X}$ si et seulement si $\lambda$ est dans la composante neutre $\Pic^0_{\X/k}$.
\end{sousthm}
\begin{demo} On peut recopier telle quelle la démonstration qui accompagne l'énoncé référencé ci-dessus.
\end{demo}

Comme toujours, la définition de la composante neutre est plus délicate lorsque la base n'est plus un corps, mais un schéma quelconque. En effet, on dispose, dans chaque fibre $\Pic_{\X_s/\kappa(s)}$, d'un ouvert $\Pic_{\X_s/\kappa(s)}^0$. Mais il se pourrait très bien que la réunion des $\Pic_{\X_s/\kappa(s)}^0$ ne soit pas un ouvert de $\Pic_{\X/S}$. On peut tout de même définir \emph{a priori} la composante neutre comme sous-foncteur de $\Pic_{\X/S}$ (voir la remarque~\ref{rem_caractérisation_fonctorielle_CC_neutre}).

\begin{sousdefi} Soient $S$ un schéma et $\X$ un $S$-champ algébrique. On suppose que le foncteur de Picard
$\pic$ est représentable
par un espace algébrique localement de type fini. On désigne alors par $\piczero$ le sous-foncteur de $\pic$ défini de la
manière suivante. Pour tout $S$-schéma $S'$ et pour tout $\xi \in \pic(S')$, on dit que $\xi$ appartient à
$\piczero(S')$ si et seulement si pour tout point $s'$ de $S'$, l'élément $\xi_{|_{s'}}$
de $\Pic_{\X_{s'}/\kappa(s')}(\kappa(s'))$ est dans la composante neutre $\Pic^0_{\X_{s'}/\kappa(s')}(\kappa(s'))$.
\end{sousdefi}

\begin{sousremarque}\rm
Il est clair que $\piczero$ est bien un sous-foncteur en groupes de $\pic$. De plus, la formation de $\piczero$ commute
au changement de base. Par ailleurs, si $S$ est le spectre d'un corps $k$ et si $\Pic_{\X/k}$ est représentable par un schéma
en groupes localement de type fini, alors le sous-foncteur $\Pic^0_{\X/k}$ ainsi défini coïncide avec le sous-foncteur
ouvert défini par la composante connexe de l'élément neutre dans $\Pic_{\X/k}$. 
\end{sousremarque}

Les lemmes techniques qui suivent ont pour but de généraliser au cas des espaces algébriques un résultat de Grothendieck sur les composantes connexes des fibres le long d'une section (voir~\ref{CC_des_fibres_lemmeB_version2}).

\begin{souslem}
\label{CC_des_fibres_lemme_projection}
Soient $k$ un corps, $X$ un $k$-schéma connexe localement de type fini, $L$ une extension de $k$ et $p$ la projection de $X_L$ dans $X$.
Alors les fibres de $p$ rencontrent toutes les composantes connexes de $X_L$. Autrement dit, pour toute composante
connexe $U$ de $X_L$, le morphisme induit $p_{|_U} : U \fleche X$ est surjectif.
\end{souslem}
\begin{demo} En dévissant l'extension, il suffit clairement de traiter le cas d'une extension algébrique et le cas d'une extension
transcendante pure. Dans le premier cas, le morphisme $\Spec L \fleche \Spec k$ est universellement ouvert (\cite{EGA4_2}~2.4.9)
et universellement fermé (\cite{EGA2}~6.1.10). Si $U$ est une composante connexe de $X_L$, son image $p(U)$ dans $X$ est ouverte, fermée
et non vide, donc c'est $X$ tout entier.

Supposons maintenant l'extension $L/k$ transcendante pure. Si $\Omega$ est une autre extension de $k$, alors l'anneau
$L\otimes_k \Omega$ est intègre.
Il en résulte que les fibres de $p$ sont géométriquement intègres.
Notons $(X_i)_{i\in I}$ la famille des composantes connexes de $X_L$. Chaque
$X_i$ est une réunion de fibres, donc les $p(X_i)$ forment une partition de $X$, et ils sont ouverts puisque $p$
est ouvert. Par connexité un seul d'entre eux est non vide, donc $X_L$ est connexe.
\end{demo}

\begin{souslem}
 \label{CC_des_fibres_lemmeA_version2}
Soient $\xymatrix{V \ar[r]^e &X \ar[r]^f &S}$ deux morphismes de schémas, avec $f$ localement de présentation finie. On note $g=f\circ e$. Si $s\in S$, on note $X_s$ la fibre $f^{-1}(s)$ et $X_s^0(e)$ la réunion des composantes connexes de $X_s$ qui rencontrent $e(V)$, et l'on note $X^0(e)$ la réunion des $X_s^0(e)$. (On notera $X^0_s$ et $X^0$ s'il n'y a pas d'ambiguïté.)
\begin{itemize}
 \item[(i)] Soit $\varphi : S'\fleche S$ un morphisme de changement de base. On adopte les notations du diagramme suivant :
$$\xymatrix{V' \ar[r]^{e'} \ar[d] \cartesien & X' \ar[r]^{f'} \ar[d]_{\varphi'} \cartesien & S'\ar[d]_{\varphi}\\
V \ar[r]^{e} & X \ar[r]^{f}  & S}$$
Alors $X'^0=\varphi'^{-1}(X^0)$.

Dans le cas où $\varphi=g$, le morphisme $f' : X' \fleche S'$ a une section $\sigma$. On note $C_{s'}'^0$ la composante connexe de $\sigma(s')$ dans $X'_{s'}$ pour tout $s'\in S'$, et $C'^0$ la réunion des $C'^0_{s'}$. Alors $X^0=\varphi'(C'^0)$.
\item[(ii)] On suppose que $f$ est universellement ouvert en tout point de $X^0$, et que ses fibres sont géométriquement réduites. On suppose de plus que l'on est dans l'un des deux cas suivants :
\begin{itemize}
\item[a)] $g$ est universellement ouvert;
\item[b)] $e$ est une immersion ouverte, et $S$ est localement noethérien.
\end{itemize}
Alors $X^0$ est un ouvert de $X$.
\end{itemize}
\end{souslem}
\begin{demo}
 (i) Pour la première assertion, on peut supposer que $S$ et $S'$ sont des spectres de corps. L'inclusion $X'^0 \subset \varphi'^{-1}(X^0)$ est évidente. Pour l'inclusion réciproque et pour la seconde assertion, on utilise le lemme~\ref{CC_des_fibres_lemme_projection}. Les détails sont laissés au lecteur.

(ii) On effectue le changement de base par $g : V \fleche S$. Avec les notations de (i), $C'^0$ est un ouvert de $X'$ par~\cite{EGA4_3}~(15.6.4). D'après (i), $X^0$ est égal à $\varphi'(C'^0)$. On en déduit immédiatement que $X^0$ est ouvert dans le cas a).

Dans le cas b), la question est clairement locale sur $S$ donc on peut supposer $S$ noethérien. Montrons que l'on peut supposer $X$ quasi-compact. Si $U$ est un ouvert quasi-compact de $X$, on note $e_U$ l'immersion de $V\cap U$ dans $U$. Notons $\Uc$ la réunion des $U^0(e_U)$ pour $U$ ouvert quasi-compact de $X$ et montrons que $\Uc=X^0$. Il est clair que $\Uc$ est inclus dans $X^0$. Réciproquement, soit $x$ un point de $X^0$, et $s=f(x)$. Alors il existe un point $v$ de $V_s$ tel que $x$ et $v$ soient dans la même composante connexe de $X_s$. En considérant une chaîne (finie) de composantes irréductibles qui relient $x$ et $v$ dans $X_s$, on peut trouver un ouvert quasi-compact $U$ de $X$ qui contient $x$ et $v$ et tel que $x$ et $v$ soient dans la même composante connexe de $U_s$ (procéder comme dans~\cite{EGA4_3}~(15.6.4)). Ceci prouve que $x$ appartient à $U^0_s(e_U)$, donc à $\Uc$.

On suppose donc $X$ noethérien. Alors $X^0$ est constructible puisque c'est l'image par $\varphi'$ de $C'^0$, qui est ouvert. Donc pour montrer qu'il est ouvert il suffit de montrer qu'il est stable par générisation. Soit $x\in X^0$ et soit $x'$ une générisation de $x$. D'après~\cite{EGA3_2}~(7.1.9) il existe un anneau de valuation discrète A et un morphisme $u$ de $S'=\Spec A$ dans $X$ qui envoie le point fermé $\xi$ sur $x$ et le point générique $\eta$ sur $x'$. On effectue le changement de base par $f\circ u : S' \fleche S$. Le morphisme $u$ induit une section $u'$ de $f'$. Le point $u'(\eta)$ est une générisation de $u'(\xi)$. D'après (i), $u'(\xi)$ est un point de $X'^0$, et il suffit de prouver que $u'(\eta)$ appartient à $X'^0$ puisque $\varphi'(X'^0) \subset X^0$ (avec les notations de (i)).

On est donc ramené au cas où $S$ est le spectre d'un anneau de valuation discrète. Le complémentaire de $X^0_{\xi}$ dans la fibre spéciale étant un fermé de $X$, quitte à le supprimer on peut supposer que $X^0_{\xi}=X_{\xi}$. Maintenant $f$ est universellement ouvert en tout point de la fibre spéciale (par hypothèse), mais aussi en tout point de la fibre générique (évident). Donc il est universellement ouvert et l'on est ramené au cas a).
\end{demo}

\begin{souslem}
\label{CC_des_fibres_lemmeB_version2}
Soient $S$ un schéma et $f : X \fleche S$ un $S$-espace algébrique localement de présentation finie, à fibres géométriquement réduites, muni d'une section $e : S \fleche X$. Pour tout $s$ appartenant à $S$, on note $X^0_s$ la composante connexe de $e(s)$ dans la fibre $X_s$. On note $X^0$ la réunion des $X^0_s$.
\begin{itemize}
\item[(i)]
On suppose que $S$ est localement noethérien, que $f$ est universellement ouvert en tout point de $X^0$, et que pour tout $s$, $X^0_s$ est irréductible. Alors $X^0$ est un ouvert de $X$.
\item[(ii)]
On suppose que pour tout $s$, $X^0_s$ est géométriquement irréductible, et que la fonction $s\mapsto \dim(X_s^0)$ est localement constante sur $S$. Alors $f$ est universellement ouvert en tout point de $X^0$.
\end{itemize}
\end{souslem}
\begin{demo}
(i) Soit $\pi : X_1 \fleche X$ une présentation de $X$. On note $S_1$ le produit fibré $X_1\times_X S$, $e_1$ la projection de $S_1$ sur $X_1$ et $f_1$ le composé $f\circ \pi$. On applique alors le lemme~\ref{CC_des_fibres_lemmeA_version2} au diagramme
$$\xymatrix{S_1\ar[r]^{e_1} & X_1 \ar[r]^{f_1} & S}.$$
On vérifie facilement que $\pi$ envoie $X_1^0(e_1)$ dans $X^0$. On en déduit, vu que $\pi$ est étale, que $f_1$ est universellement ouvert en tout point de $X_1^0(e_1)$. Il résulte alors de~\ref{CC_des_fibres_lemmeA_version2}~(ii), cas~a) que $X_1^0(e_1)$ est ouvert dans $X_1$.

On note maintenant $U=\pi(X_1^0(e_1))$ (c'est un ouvert de $X$ inclus dans $X^0$) et $V_1=\pi^{-1}(U)$. Autrement dit $V_1$ est le saturé de $X_1^0(e_1)$ pour la relation d'équivalence définie par $\pi$. On note alors $W_1= X_1^0(V_1 \inj X_1)$. De même que précedemment on vérifie facilement que $\pi(W_1)$ est inclus dans $X^0$, si bien que $f_1$ est universellement ouvert en tout point de $W_1$. En appliquant le lemme~\ref{CC_des_fibres_lemmeA_version2}~(ii), cas~b) au diagramme
$$\xymatrix{V_1\ar@{^(->}[r] & X_1 \ar[r]^{f_1} & S},$$
on en déduit que $W_1$ est un ouvert de $X_1$. Il suffit donc de montrer que $X^0$ est égal à $\pi(W_1)$ pour conclure.

On peut pour cela supposer que $S$ est le spectre d'un corps. De plus, $\pi$ étant surjectif, il suffit de montrer que $\pi^{-1}(X^0)$ est inclus dans $W_1$. Maintenant $X^0$ est un ouvert irréductible de $X$. Notons $\eta$ son point générique. Vu la construction de $W_1$ à partir de $V_1$, il suffit de montrer que $V_1$ contient le point générique de toute composante irréductible de $\pi^{-1}(X^0)$. Or, $\pi$ étant ouvert, il envoie nécéssairement tous ces points génériques sur $\eta$ d'après~\cite{EGA4_1}~(1.10.4). Comme $V_1$ est saturé pour la relation d'équivalence définie par $\pi$, il suffit de montrer qu'il contient le point générique d'\emph{une} composante irréductible de $\pi^{-1}(X^0)$, ce qui est évident puisque c'en est un ouvert non vide.

(ii) La formation de $X^0$ commutant au changement de base, il suffit de montrer que $f$ est ouvert en tout point de $X^0$. Comme $f$ est localement de présentation finie, il suffit de montrer que $f$ est générisant en tout point de $X^0$. Soient $x\in X^0$ et $s'$ une générisation de $s=f(x)$. On note $x'$ le point générique de $X^0_{s'}$. Il suffit maintenant de montrer que $X^0_s$ est inclus dans $\ov{X^0_{s'}}$.

Or on a les inégalités :
$$\dim(X^0_{s'}) \le \dim_{e(s)}(\ov{X^0_{s'}}\cap f^{-1}(s)) \le \dim(X^0_s).$$
La première résulte du théorème de semi-continuité de Chevalley (\cite{EGA4_3}~(13.1.3)), compte tenu du fait que $X^0_{s'}$ est égal à $\ov{X^0_{s'}}\cap f^{-1}(s')$. La seconde résulte du fait que toute composante irréductible de $\ov{X^0_{s'}}\cap f^{-1}(s)$ contenant $e(s)$ est nécessairement incluse dans $X^0_s$. Si l'on suppose de plus que la fonction $t\mapsto \dim(X^0_t)$ est localement constante, alors les membres de gauche et de droite sont égaux, d'où $\dim_{e(s)}(\ov{X^0_{s'}}\cap f^{-1}(s)) = \dim(X^0_s)$. Maintenant soit $Z$  une composante irréductible de $\ov{X^0_{s'}}\cap f^{-1}(s)$ contenant $e(s)$ et de dimension $\dim(X^0_s)$. Alors $Z$ est nécessairement inclus dans $X^0_s$, donc $Z=X^0_s$ puisqu'ils sont de irréductibles et de même dimension, d'où $X^0_s \subset \ov{X^0_{s'}}$.
\end{demo}

\begin{sousremarque}\rm
\label{rem_caractérisation_fonctorielle_CC_neutre}
Comme sous-foncteur de $X$, l'ouvert $X^0$ est caractérisé par la propriété suivante.
Pour tout $S$-schéma $T$ et pour tout $\xi \in X(T)$, $\xi$ appartient à $X^0(T)$ si et seulement si 
pour tout $s\in S$, le point $\xi_s\in X_s(T_s)$ obtenu par changement de base appartient à $X^0_s(T_s)$, autrement
dit le morphisme correspondant $T_s \fleche X_s$ se factorise par $X^0_s$.
Notons aussi que la formation de $X^0$ commute au changement de base. \end{sousremarque}


Nous pouvons maintenant énoncer le résultat principal de cette section.

\begin{sousprop} Soient $S$ un schéma localement noethérien et $\X$ un $S$-champ algébrique. On suppose que le foncteur de Picard $\pic$ 
est représentable par un espace algébrique localement de présentation finie sur $S$ et que les fibres $\Pic_{\X_s/\kappa(s)}$ sont géométriquement réduites. On suppose de plus que l'une ou l'autre des hypothèses suivantes est vérifiée :
\begin{itemize}
\item[a)] le morphisme $\pic \fleche S$ est universellement ouvert en tout point de $\piczero$ (vrai par exemple s'il est plat) ;
\item[b)] la fonction $s\mapsto \dim(\Pic_{\X_s/\kappa(s)})$ est localement constante sur $S$ (hypothèse vérifiée par exemple si $\pic$ est lisse le long de la section unité).
\end{itemize}
 Alors le morphisme naturel 
$$\xymatrix{\piczero  \ar[r]& \pic}$$
est une immersion ouverte. De plus, $\piczero$ est de type fini sur $S$.
\end{sousprop}
\begin{demo} Que $\piczero$ soit un ouvert de $\pic$ résulte immédiatement du lemme~(\ref{CC_des_fibres_lemmeB_version2}).
Comme il est par hypothèse localement de type fini sur $S$, il ne reste plus qu'à montrer qu'il est quasi-compact
lorsque $S$ est localement noethérien. Nous adoptons pour cela la même démarche que celle de Kleiman dans~\cite{poly_Kleiman}~5.20. C'est une question locale sur $S$, donc on peut supposer que
$S$ est affine noethérien. Notons $\sigma$ le morphisme naturel de $\piczero$ vers $S$. Par récurrence noethérienne
sur les fermés de $S$, on peut supposer que pour tout fermé strict $Z$ de $S$, $\sigma^{-1}(Z)=\piczero\times_S Z$ est 
quasi-compact. Il suffit alors de construire un ouvert non vide $U$ de $S$ tel que $\sigma^{-1}(U)$ soit quasi-compact.

Soit $\pi : V \fleche \piczero$ un morphisme étale, où $V$ est un schéma affine non vide, et soit $U=\sigma(\pi(V))$.
Comme $\sigma \circ \pi$ est un morphisme
lisse, $U$ est un ouvert de $S$. Pour montrer que cet ouvert $U$ convient, on va construire un morphisme
surjectif de $V\times_S V$ vers $\piczero \times_S U$. Ce sera suffisant, puisque $V\times_S V$ est quasi-compact.
On a un diagramme commutatif de $S$-espaces algébriques :
$$\xymatrix{V\times_S V \ar[r]^-{\pi\times \pi} \ar[d] & \piczero\times_S \piczero \ar[r]^-{\alpha} & \piczero \ar[d]^{\sigma}\\
U \ar[rr] && S}$$
où $\alpha$ est défini fonctoriellement par $\alpha(g,h)=g.h^{-1}$.
On en déduit un morphisme $\alpha' : V\times_S V \fleche \piczero\times_S U$, dont il ne nous reste plus qu'à montrer
la surjectivité. On peut supposer pour cela que $S$ est le spectre d'un corps algébriquement clos.
Mais dans ce cas on a $U=S$ et il faut montrer que le morphisme $\alpha \circ (\pi \times \pi)$ est surjectif.
De plus, vu que $S$ est le spectre d'un corps, $\piczero$ est un schéma en groupes et on peut supposer que $V$ est un
ouvert de $\piczero$. Le résultat est alors classique et facile (\emph{cf.}~\cite{poly_Kleiman}).
\end{demo}


\subsection{Propreté de la composante neutre}

\begin{sousthm}
\label{composante_neutre_propre_sur_k}
On suppose que $\X$ est un champ algébrique propre, géométriquement normal et cohomologiquement plat en
dimension zéro sur $\Spec k$. Alors la composante neutre $\Pic^0_{\X/k}$
du schéma de Picard est propre sur $k$.
\end{sousthm}
\begin{demo}
Si $\ov{k}$ est une clôture algébrique de $k$ alors le champ $\X_{\ov{k}}$ obtenu par changement de base est normal d'après~\cite{EGA4_2}~(6.7.7). Il vérifie clairement les autres hypothèses du théorème donc par descente fidèlement plate on peut supposer le corps $k$ algébriquement clos.

Il suffit (cf. argumentaire de Kleiman au cours de la démonstration du théorème~5.4 de~\cite{poly_Kleiman}) de montrer que tout morphisme de schémas de $G$ vers $\Pic_{\X/k}$ est constant, avec $G=\mathbb{G}_a$ ou $G=\gm$. (En fait il suffirait même de le faire pour $\gm$ puisqu'alors c'est aussi vrai pour $G=\mathbb{G}_a$, mais
cette restriction n'apporte pas grand-chose.) En effet, il suffit de montrer que le
réduit $(\Pic^0_{\X/k})_{\textrm{réd}}$ est propre, or ce dernier est lisse donc on peut
lui appliquer le théorème de structure de Chevalley et Rosenlicht (cf. par
exemple~\cite{Conrad_thm_Chevalley}, théorème~1.1). On en déduit que
$(\Pic^0_{\X/k})_{\textrm{réd}}$ a un sous-groupe algébrique linéaire $H$, fermé et
distingué dans $(\Pic^0_{\X/k})_{\textrm{réd}}$, tel que le quotient
$(\Pic^0_{\X/k})_{\textrm{réd}}/H$ soit une variété abélienne. Il suffit de montrer
que $H$ est trivial. Le groupe $H$ est commutatif donc résoluble. Il est dès lors triangulable d'après le théorème de Lie-Kolchin. On en déduit que s'il était non trivial, il contiendrait un sous-groupe isomorphe à $\gm$ ou $\mathbb{G}_a$ (voir par exemple le livre de Springer~\cite{Springer_Linear_Algebraic_Groups}, lemme~6.3.4).

Vu que $f$ a une section, et vu que les groupes de Picard de $\Spec k$ et
de $G$ sont triviaux, dire que tout morphisme de schémas de $G$ vers $\Pic_{\X/k}$ est constant revient à dire, grâce au théorème~(\ref{comparaison_des_foncteurs_de_Picard}), que
le morphisme naturel
$$\xymatrix{\Pic(\X) \ar[r]& \Pic(\X\times G)}$$
est un isomorphisme.
La suite exacte des termes de bas degré associée à la suite spectrale de Leray du morphisme $p : \X\times G \fleche \X$
s'écrit :
$$\xymatrix{0 \ar[r]& H^1(\X,p_*\gm) \ar[r] & H^1(\X\times G, \gm) \ar[r] & H^0(\X, R^1p_*\gm).}$$
Commençons par montrer que le faisceau $R^1p_* \gm$ est nul. On sait d'après le calcul des images directes supérieures
effectué en annexe, que c'est le faisceau associé au préfaisceau qui à tout ouvert lisse-étale $(U,u)$ de $\X$
associe $H^1(U\times G,\gm)$. Donc il suffit de montrer que
pour tout schéma affine  $U$ lisse sur $\X$ et pour tout $\xi\in H^1(U\times G,\gm)$,
il existe une famille couvrante étale $V\fleche U$ telle
que l'élément $\xi_{|_V}$ de $H^1(V\times G,\gm)$ soit nul. Mais pour démontrer ceci il suffit clairement de savoir que 
le morphisme $\Pic(U) \fleche \Pic(U\times G)$ est surjectif. Nous sommes donc ramenés à montrer que si $U$ est un schéma affine normal
(que l'on peut aussi supposer intègre, vu que $U$ est de toute manière somme disjointe finie
de schémas intègres) sur $\Spec k$, alors $\Pic(U\times G)$ s'identifie à $\Pic(U)$.
Ce fait est démontré par Kleiman dans~\cite{poly_Kleiman}, au cours de la démonstration du théorème~(5.4).

Calculons maintenant le faisceau $p_*\gm$. Considérons d'abord le cas où $G=\gm$. Nous allons montrer que
$p_*\gm$ s'identifie à $\gm\times \Z$. Il suffit bien évidemment de vérifier que ces deux faisceaux coïncident
sur le site $\lisets(\X)$. Si $U=\Spec A$ est un schéma affine lisse sur $\X$, il est en particulier normal, donc somme disjointe
finie de schémas affines intègres, si bien que l'on peut supposer $U$ intègre. Alors $\gm\times \Z (U)=A^{\times}\times \Z$.
Par ailleurs, vu que $p$ est représentable, on a
$$p_*\gm(U)=\gm((\X\times G) \times_{\X} U)=\gm(U\times G)=A[X,X^{-1}]^{\times}.$$
Or lorsque $A$ est intègre, il est clair que les éléments inversibles de l'anneau $A[X,X^{-1}]$ sont les
éléments de la forme $aX^n$ avec $a\in A^{\times}$ et $n\in \Z$.
Dans le cas où $G=\mathbb{G}_a$, le lecteur vérifiera facilement que l'on trouve $p_*\gm=\gm$.

Or le groupe $H^1(\X,\Z)$ est réduit à zéro d'après le théorème~(\ref{H1_normal}) ci-dessous. On a alors, que $G$ soit égal à $\mathbb{G}_a$ ou $\gm$,
$$H^1(\X,p_*\gm)=H^1(\X,\gm)=\Pic(\X),$$
ce qui, vu la suite exacte évoquée ci-dessus, fournit l'isomorphisme désiré.
\end{demo}

\begin{sousthm}
\label{H1_normal}
Soit $\X$ un champ algébrique localement noethérien et normal. Alors le groupe $H^1(\X,\Z)$ est nul.
\end{sousthm}

Nous allons pour cela montrer que tout $\Z$-torseur sur $\X$ est trivial. Notre démarche est
fortement inspirée de l'étude des préschémas constants tordus quasi-isotriviaux proposée dans SGA3
(\cite{SGA3} exposé~X, paragraphe~5).

\begin{souslem} Soit $f : \X \fleche \Y$ un morphisme représentable de $S$-champs algébriques. Soit $y$ un point
de l'espace topologique $|\Y|$ sous-jacent à $\Y$. Les propositions suivantes sont équivalentes :
\begin{itemize}
\item[(i)] Pour un représentant $\Spec K \fleche \Y$ de $y$, le morphisme $f_K$ induit par changement de base
de $\X_K=\X\times_{\Y} \Spec K$ vers $\Spec K$ est fini.
\item[(ii)] Pour tout représentant $\Spec K \fleche \Y$ de $y$, le morphisme $f_K$ est fini.
\end{itemize}
Lorsqu'elles sont vérifiées, on dit que $f_y : \X_y \fleche y$ est fini, ou encore que $\X_y$ est fini.
\end{souslem}
\begin{demo}
Il suffit clairement de montrer que si $L$ est une extension de $K$ et si $f_L$ est fini, alors $f_K$ est fini. Or si
$f_L$ est fini, $\X_L$ est un schéma affine, et $\X_K$ aussi par descente \emph{fpqc} pour les morphismes affines.
Maintenant la descente fidèlement plate pour les morphismes finis de schémas (\cite{EGA4_2}~(2.7.1)) assure que $f_K$ est fini.
\end{demo}

Le lemme suivant généralise le lemme~5.13 de~\cite{SGA3}, exposé~X.

\begin{souslem}
\label{morphisme_fini_condition_ouverte_et_fermee}
Soit $p :\Pc \fleche \X$ un morphisme représentable de $S$-champs algébriques, avec $\X$ localement noethérien.
On suppose qu'il existe une présentation $X\fleche \X$ de $\X$ telle que $P=\Pc \times_{\X} X$ soit une union
disjointe de copies de $X$. Soit $\Cc$ un sous-champ ouvert et fermé de $\Pc$. Alors l'ensemble
des points $x$ de $|\X|$ tels que $\Cc_x$ soit fini est ouvert et fermé dans $|\X|$. Si on note $\Uc$ le sous-champ
ouvert et fermé que cet ensemble définit, le champ $\Cc_{\Uc}=\Cc\times_{\X}\Uc$ est fini sur $\Uc$.
\end{souslem}
\begin{demo}
Notons $C$ le sous-espace algébrique ouvert et fermé de $P$ obtenu par changement de base à partir de $\Cc$.
Les propriétés que l'on veut montrer sont clairement de nature locale pour la topologie lisse sur $\X$, 
donc il suffit de montrer qu'elles sont vérifiées par $C \fleche P \fleche X$. Comme $X$ est lui aussi localement
noethérien, ses composantes connexes sont ouvertes et fermées donc on peut supposer que $X$ est un schéma connexe.
Dans ce cas, vu que $P$ est une union disjointe de copies de $X$, le sous-schéma ouvert et fermé $C$ est
lui-même l'union disjointe de certaines de ces copies. Si elles sont en nombre fini alors $C$ est fini sur $X$, sinon
l'ensemble des points $x$ de $X$ où $C_x$ est fini est vide.
\end{demo}
\begin{democ}{du théorème \ref{H1_normal}}
Nous allons utiliser la description du premier groupe de cohomologie en termes de torseurs
(cf. paragraphe~\ref{Cohomologie_et_torseurs}).
Il suffit en vertu de~(\ref{thm_coh_torseur}) de montrer que tout $\Z$-torseur sur $\X$ au sens de~(\ref{def_torseur}) est trivial.
Le cas où $\X$ est le spectre d'un corps est bien connu et nous nous en servirons par la suite.
Remarquons tout d'abord que dans tous les champs algébriques (\emph{a fortiori} tous les schémas ou espaces algébriques) qui vont
intervenir au cours de la démonstration, les composantes connexes sont irréductibles.
En effet, ils seront tous normaux et localement noethériens car ce sont là des propriétés de nature locale pour la topologie lisse.
Notre affirmation résulte alors de la proposition~4.13 de~\cite{LMB}.

Donnons-nous donc un $\Z$-torseur sur $\X$, c'est-à-dire un morphisme $p : \Pc \fleche \X$ représentable et lisse muni d'une
action de $\Z$ qui en fait un torseur (\emph{cf.} \ref{def_torseur}). Pour montrer que $\Pc$ est trivial, on peut supposer $\X$ connexe donc irréductible.
Soit $\Cc$ une composante connexe de $\Pc$. Notons $\eta$ le point générique de
l'espace $|\X|$ sous-jacent à $\X$.

\begin{etape}{La fibre générique $\Cc_{\eta}$ est finie.}
Soit $u : U\fleche \X$ un morphisme lisse, où $U$ est un schéma affine irréductible (donc intègre, puisque $U$ est normal).
On note $s$ le point générique de $U$ et on adopte encore les notations du diagramme suivant :
$$\xymatrix{\Cc_s \ar[r] \ar[d] \cartesien & \Cc_U \ar[r] \ar[d] \cartesien & \Cc \ar[d]\\
\Pc_s \ar[r] \ar[d] \cartesien & \Pc_U \ar[r] \ar[d] \cartesien & \Pc \ar[d]^p \\
\Spec \kappa(s) \ar[r] & U\ar[r]^u & \X.}$$
Notons $(\Cc_i)_{i\in I}$ les composantes connexes de $\Cc_U$ et pour chaque $i\in I$ notons $\eta_i$ le point
générique de $|\Cc_i|$. Notons enfin $\xi$ le point générique de $|\Cc|$ et $\Spec K \fleche \Cc$ l'un de ses
représentants. Pour tout $i$ le morphisme de $\Cc_i$ vers $\Cc$ est lisse donc générisant (\cite{LMB}~(5.8)) si bien qu'il envoie
le point $\eta_i$ sur $\xi$. On en déduit (cf. par exemple \cite{LMB}~(5.4)~(iv)) que le champ $\Cc_{i,K}=\Cc_i \times_{\Cc} \Spec K$
est non vide. Maintenant, les $|\Cc_{i,K}|$ forment une partition ouverte de $|\Cc_{U,K}|$, où $\Cc_{U,K}$ est le produit
fibré $\Cc_U\times_{\Cc} \Spec K$. Or le champ
$$\Cc_{U,K} = U\times_{\X} \Spec K$$
est quasi-compact car $\X$ est quasi-séparé donc les $|\Cc_{i,K}|$ sont en nombre fini et finalement $\Cc_U$ n'a qu'un
nombre fini de composantes connexes.

Par ailleurs, pour chaque $i$ le morphisme naturel de $\Cc_i$ vers $U$ est lui aussi générisant donc il envoie $\eta_i$
sur $s$. En particulier sa fibre générique est irréductible donc $\Cc_s$ est une union finie d'irréductibles. Comme
$\Pc_s$ est un $\Z$-torseur sur un corps, il est nécessairement trivial, donc $\Cc_s$ est une union disjointe finie
de copies de $\Spec \kappa(s)$. Or le morphisme de $\Spec \kappa(s)$  vers $\X$ est un représentant
de $\eta$ donc $\Cc_{\eta}$ est fini.
\end{etape}

\begin{etape}{Montrons que $\Cc$ est fini.}
D'après le lemme~(\ref{morphisme_fini_condition_ouverte_et_fermee}), l'ensemble des points $x$ de $|\X|$ où
$\Cc_x$ est fini est ouvert et fermé dans $|\X|$. Or il est non vide puisqu'il contient $\eta$, donc par
connexité c'est $|\X|$ tout entier. Le même lemme prouve alors que le morphisme de $\Cc$ vers $\X$ est fini.
\end{etape}

\begin{etape}{Montrons que $\Cc \fleche \X$ est étale.}
Il s'agit d'un morphisme fini. En particulier il est schématique. Notre assertion résulte alors du
fait qu'un morphisme fini et lisse de schémas est étale (\cite{SGA1}~II~1.4).
\end{etape}

\begin{etape}{Montrons que $\Cc \fleche \X$ est radiciel.}
Supposons qu'il existe un corps $K$ et un morphisme $\Spec K \fleche \X$ tel que le schéma $\Cc_K$ obtenu par changement
de base contienne au moins deux points $x_1$ et $x_2$. Notons $c_1$ et $c_2$ leurs images
dans $|\Cc| \subset |\Pc|$. Soit $n$ l'unique élément de $\Z\setminus\{0\}$ tel que l'automorphisme $\tau_n$ correspondant
envoie $x_1$ sur $x_2$. On note encore $\tau_n$ l'automorphisme de $|\Pc|$ qui correspond à $n$. Il est clair qu'il envoie
$c_1$ sur $c_2$. Or $\tau_n(|\Cc|)$ est un connexe qui contient $c_2$ donc il est inclus dans la composante connexe de
$c_2$, à savoir $|\Cc|$. On montre de même que $\tau_n^{-1}(|\Cc|)$ est inclus dans $|\Cc|$ donc $\tau_n$ induit
un automorphisme de $|\Cc|$ et tous les $\tau_n^k(c_1)$, $k\in \Z$, sont dans $|\Cc|$. Donc tous les $\tau_n^k(x_1)$
sont dans l'ouvert $\Cc_{K}$ de $\Pc_K$, ce qui contredit le fait qu'il est de type fini.
\end{etape}

\begin{etapefinale}{Conclusion.}
Le morphisme de $\Cc$ vers $\X$ est schématique, étale, radiciel et de type fini donc par~\cite{SGA1}~I~5.1 c'est une
immersion ouverte. De plus il est aussi fermé puisqu'il est fini donc par connexité de $\X$ il est surjectif.
Ceci prouve que c'est un isomorphisme, donc le torseur $\Pc \fleche \X$ est trivial.
\end{etapefinale}
\end{democ}

\begin{sousthm}
Soit $\X$ un $S$-champ algébrique. On suppose que $\pic$ est représentable par un espace algébrique localement de présentation finie, que $\piczero$ est un ouvert de $\pic$ et qu'il est de plus séparé et de type fini sur $S$. Alors l'ensemble $U$ des points $s$ de $S$ tels que  $\Pic^0_{\X_s/\kappa(s)}$ soit propre sur $\kappa(s)$ est un ouvert de $S$. De plus $\piczero\times_S U$ est propre sur $U$. 
\end{sousthm}
\begin{demo}
On note $P$ l'espace algébrique $\piczero$ et $f$ son morphisme structural vers $S$. Des arguments standard de passage à la limite permettent de supposer $S$ noethérien.
Il faut montrer que tout point $s$ de $S$ tel que  $P_s$ soit propre sur $\kappa(s)$ admet un voisinage ouvert $U$ tel que le morphisme induit $f^{-1}(U) \fleche U$ soit propre. D'après~\cite{EGA4_3}~8.10.5~(xii) (combiné avec le lemme de Chow pour les espaces algébriques), on peut supposer que $S$ est le spectre d'un anneau local noethérien $A$, et il faut montrer que $f$ est propre dès que la fibre du point fermé l'est. Par un argument de descente fidèlement plate, on peut même supposer $A$ complet. Nous montrons d'abord l'assertion suivante.

\begin{etape}{Assertion : Si $f :P \fleche S$ est un morphisme séparé et de type fini d'espaces algébriques (où $S$ est le spectre d'un anneau local noethérien complet) et si la fibre spéciale de $f$ est propre, alors $P$ admet un sous-espace ouvert et fermé $P_0$ qui est propre sur $S$ et qui contient la fibre spéciale.}
\indent
On peut clairement supposer la fibre spéciale non vide. Commençons par traiter le cas où $P$ est irréductible. En utilisant le lemme de Chow pour les espaces algébriques (\cite{Knutson}~IV~3.1), on peut trouver un schéma irréductible $V$ et un morphisme surjectif et projectif $g$ de $V$ vers $P$. Maintenant, vu que la fibre spéciale de $P$ est propre, celle de $V$ l'est aussi. D'après \cite{EGA3_1}~5.5.2, on en déduit que $V$ lui-même est propre sur $S$, ce qui prouve que $P$ est propre puisqu'il est déjà séparé et de type fini et que $g$ est surjectif.

Supposons maintenant $P$ connexe. On le décompose alors en composantes irréductibles $P_1, \dots, P_n$. Soit $Z$ la réunion des composantes qui rencontrent la fibre spéciale. Chacune de ces composantes est propre d'après le cas précédent, si bien que $Z$ est propre. De plus, $Z$ est non vide puisque la fibre spéciale est supposée non vide. Soit $F$ la réunion des autres composantes. Supposons $F$ non vide. Alors $F$ rencontre $Z$ puisque $P$ est connexe, et le produit fibré $Z\times_P F$ est propre sur $S$. En particulier, son image est un fermé non vide de $S$, donc elle contient le point fermé. Ceci prouve que $F$ rencontre la fibre spéciale, ce qui est absurde. Donc $F$ est vide et finalement $P$ est propre.

Pour prouver l'assertion dans le cas général, il suffit d'appliquer le résultat à chacune des composantes connexes de $P$, puis de regrouper celles qui rencontrent la fibre spéciale d'une part, et les autres d'autre part. 
\end{etape}

\begin{etapefinale}{Conclusion}
D'après l'assertion précédente $P$ s'écrit comme l'union disjointe d'un sous-espace propre $P_0$ et d'un sous-espace $P_1$ qui s'envoie dans le complémentaire du point fermé. Vu que $S$ est connexe, l'image de la section neutre est entièrement contenue dans $P_0$. Mais comme par ailleurs les fibres de $f$ sont connexes, le fait qu'elles rencontrent $P_0$ montre qu'elles sont incluses dans $P_0$, si bien que $P_1$ est vide et que $P$ est propre.
\end{etapefinale}
\end{demo}

\section{Quelques exemples}

\subsection{Espace de module des courbes elliptiques}

Soient $S$ un schéma et $\Mc_{1,1,S}$ le $S$-champ qui classifie les courbes elliptiques. Mumford a calculé en 1965 le groupe de Picard de $\Mc_{1,1,S}$ lorsque $S$ est le spectre d'un corps de caractéristique différente de 2 et 3. Ce groupe est isomorphe à $\Z/12\Z$, engendré par le \emph{fibré de Hodge} $\lambda$ défini de la manière suivante. Si $t : T \fleche \Mc_{1,1,S}$ est un $T$-point de $\Mc_{1,1,S}$ correspondant à une courbe elliptique $f : E \fleche T$, la restriction $t^*\lambda$ est le faisceau inversible $f_*\Omega_{E/T}$. Fulton et Olsson généralisent ce résultat dans \cite{Fulton_Olsson}. Ils montrent que si la base $S$ est réduite ou est un $\Z[\frac12]$-schéma, alors le groupe $\Pic(\Mc_{1,1,S})$ s'identifie à
$$\Pic(\A_S^1)\times \Z/12\Z(S).$$
On en déduit que le foncteur de Picard $\Pic_{\Mc_{1,1,S}/S}$ n'est pas représentable. En effet, dans le cas contraire, le foncteur de Picard de la droite affine serait représentable\footnote{au moins dans le cas où $S$ est un $\Z[\frac12]$-schéma. Pour le cas général, on obtient aussi une contradiction car pour montrer que $\Pic_{\A^1_S/S}$ n'est pas représentable, il suffit de considérer des schémas de base réduits.}, ce que l'on sait être faux.

Fulton et Olsson calculent aussi le groupe de Picard de la compactification standard $\overline{\Mc_{1,1,S}}$ de $\Mc_{1,1,S}$. Il est isomorphe à $\Z(S)\times \Pic(S)$. On en déduit immédiatement que le foncteur de Picard de $\overline{\Mc_{1,1,S}}$ est représentable, isomorphe à $\Z$, et que son champ de Picard est isomorphe à $\Z \times \bgm$. 

\subsection{Espaces projectifs à poids}

\begin{sousdefi}
Soient $n\ge 1$ et $a_0, \dots, a_n$ des entiers strictement positifs. Le groupe $\gm$ agit sur $\A^{n+1}_{\Z}\setminus \{0\}$ par
$$\lambda \cdot (x_0, \dots, x_n) = (\lambda^{a_0}x_0, \dots, \lambda^{a_n}x_n).$$
On note $\P(a_0, \dots, a_n)$ le champ quotient $[\A^{n+1}_{\Z}\setminus \{0\}/\gm]$ de $\A^{n+1}_{\Z}\setminus \{0\}$ par cette action, et on l'appelle espace projectif de poids $(a_0, \dots, a_n)$.
\end{sousdefi}

Nous allons voir que la théorie développée ci-dessus permet de retrouver très rapidement le résultat de~\cite{Noohi_weighted_stacks}, à savoir que le champ de Picard de $\P(a_0, \dots, a_n)$ est isomorphe à $\Z\times \bgm$. On note $\X=\P(a_0, \dots, a_n)$ et $\X_T=\X\times T$ pour tout schéma $T$. Il est bien connu que $\X$ est propre, plat et cohomologiquement plat sur $S= \Spec \Z$, donc on sait déjà que $\pic$ est représentable par un espace algébrique localement séparé. De plus $\X \fleche S$ a une section donc son champ de Picard est isomorphe à $\pic \times \bgm$ et il suffit de voir que $\pic$ est isomorphe à $\Z$.

\begin{souslem}
\label{lemme_espaces_projectifs_a_poids}
Le morphisme $\pic \fleche S$ est non ramifié.
\end{souslem}
\begin{demo}
Il est déjà localement de présentation finie. Il suffit donc de montrer que pour tout schéma affine $S'$, tout sous-schéma fermé $S'_0$ de $S'$ défini par un idéal
$I$ de carré nul, et tout morphisme de $S'$ dans $S$, l'application canonique
de $\Hom_S(S', \pic)$ dans $\Hom_S(S_0', \pic)$
est injective. On se donne un faisceau inversible $\Lc$ sur $\X_{S'}$ dont la restriction $\Lc_0$ à $\X_{S'_0}$ provient de la base $S'_0$. Il faut montrer que $\Lc$ provient de la base $S'$. Par hypothèse, $\Lc_0$ provient d'un faisceau $\Bc_0$ sur $S'_0$, qui lui-même provient d'un faisceau inversible $\Bc$ sur $S'$ puisque $H^2(S', I)$ est nul (on utilise~(\ref{thm_defm_fi})). Maintenant $\Lc$ et $\Bc_{|_{\X_{S'}}}$ sont deux déformations de $\Lc_0$ à $\X_{S'}$. D'après le théorème~(\ref{thm_defm_fi}) elles sont nécessairement isomorphes puisque $H^1(\X_{S'}, I)$ est nul.
\end{demo}

Se donner un faisceau inversible sur $\X_S$ revient à se donner un faisceau inversible $\gm$-équivariant sur $\A^{n+1}_S\setminus \{0\}$, c'est-à-dire un faisceau inversible $\Lc$ sur $\A^{n+1}_S\setminus \{0\}$ et un élément de $\gm(\A^{n+1}_S\setminus \{0\}\times\gm)$ qui vérifie une certaine condition de cocycle. Si $d$ est un entier, l'élément $\lambda^d$ de $\gm(\gm)=\Z[\lambda,\lambda^{-1}]^{\times}$ définit un élément de $\gm(\A^{n+1}_S\setminus \{0\}\times\gm)$ qui vérifie bien cette condition de cocycle. On note $\Oc(d)$ le faisceau inversible ainsi construit sur $\X_S$. Alors l'application $d\mapsto \Oc(d)$ définit un morphisme d'espaces algébriques
$$\varphi : \Z \flechelongue \pic.$$

On va construire son inverse \emph{via} le degré. Si $S$ est le spectre d'un corps $k$ un calcul élémentaire montre que l'application de $\Z$ dans $\Pic(\X_k)$ décrite ci-dessus est un isomorphisme. Si $\Lc$ est un faisceau inversible sur $\X_k$, on note $\deg(\Lc)$ l'entier correspondant à $\Lc$. Maintenant si $S$ est une base quelconque et $\Lc$ un faisceau inversible sur $\X_S$, on définit une fonction $\deg_{\Lc}$ sur $S$ à valeurs dans $\Z$ en associant à chaque point $s$ de $S$ le degré du faisceau $\Lc_s$ sur $\X_{\kappa(s)}$. D'après~\ref{lemme_espaces_projectifs_a_poids}, la diagonale $\Delta : \pic \fleche \pic\times_S \pic$ est une immersion ouverte. On en déduit facilement que $\deg_{\Lc}$ est une fonction localement constante sur $S$. L'application $\Lc \mapsto \deg_{\Lc}$ définit donc un morphisme de faisceaux de $\pic$ dans $\Z$ dont il est clair qu'il est un inverse de $\varphi$.

\subsection{Racine \iem{n} d'un faisceau inversible}

Soient $X$ un $S$-schéma et $\Lc$ un faisceau inversible sur $X$.
Soit $n$ un entier strictement positif. On fabrique un
champ $[\Lc^{\frac1n}]$ à partir de ces données de la manière suivante.
Si $U$ est un objet de $\aff$, $[\Lc^{\frac1n}]_U$ est la catégorie des triplets
$(x,\Mc,\varphi)$ où
$$\left\{ \begin{array}{l}
x : U\fleche X \text{ est un élément de }X(U) \\
\Mc \text{ est un faisceau inversible sur }$U$\\
\varphi : \Mc^{\otimes n} \fleche x^*\Lc \text{ est un isomorphisme de faisceaux inversibles.}
\end{array} \right.$$

On note $\pi$ le morphisme canonique de $[\Lc^{\frac1n}]$ dans $X$. 
Si $U\in\ob\aff$ et si $\alpha$ est un objet de $[\Lc^{\frac1n}]_U$, le foncteur $\fAut_U(\alpha)$ est représentable par $\mun$. Plus généralement, si $\alpha_1,\alpha_2$ sont deux objets de $[\Lc^{\frac1n}]_U$, alors le foncteur $\fIsom(\alpha_1,\alpha_2)$
est représentable par un schéma fini sur $U$ (localement ce schéma est de la forme $\Spec(A[X]/(X^n-\gamma))$ où $\Spec A$ est un ouvert de $U$ qui trivialise les objets $\alpha_1$ et $\alpha_2$ et $\gamma$ est un élément de $A^{\times}$). En d'autres termes le morphisme diagonal
$\Delta_{\pi}$
est schématique et fini. En particulier $[\Lc^{\frac1n}]$ est un $S$-préchamp. Il est clair que c'est même un $S$-champ \emph{fppf}.
Si $f : Y \fleche X$ est un morphisme de $S$-schémas, alors le $S$-champ $[\Lc^{\frac1n}]\times_X Y$ est canoniquement
1-isomorphe à $[(f^*\Lc)^{\frac1n}]$.

\begin{sousprop}
\begin{itemize}
 \item[1)]
Le champ $[\Lc^{\frac1n}]$ est une gerbe pour la topologie $\emph{fppf}$ sur $X$. Si $S$ est un $\Z[\frac1n]$-schéma, alors $[\Lc^{\frac1n}]$
est même une gerbe pour la topologie étale.
\item[2)] Si $\Lc$ a une racine $n^{\text{ième}}$ sur $X$ alors $[\Lc^{\frac1n}]$ est canoniquement (une fois qu'on a fixé un faisceau inversible $\Mc$ et un isomorphisme entre
$\Mc^{\otimes n}$ et $\Lc$) 1-isomorphe au champ classifiant $\bmunfppf$ du groupe $\mun$ pour la topologie \emph{fppf}.
En particulier c'est un champ algébrique (\cite{LMB}~(10.6) et~(10.13.1)).
\item[3)] Le champ $[\Lc^{\frac1n}]$ est algébrique.
\end{itemize}
\end{sousprop}
\begin{demo}
1) Il est clair que $[\Lc^{\frac1n}]$ a des objets partout localement pour la topologie de Zariski, donc $\pi$ est un épimorphisme.
Pour montrer que le morphisme diagonal est un épimorphisme on est ramené à montrer que pour tout schéma affine $\Spec A$, tout élément de $A^{\times}$ admet une racine
$n^{\text{ième}}$ localement pour la topologie considérée. Si $\gamma\in A^{\times}$,
le morphisme de $\Spec A[X]/(X^n-\gamma)$ vers $\Spec A$
est une famille couvrante pour la topologie \emph{fppf} qui répond au problème posé.
Si de plus $n$ est inversible alors c'est même une famille couvrante pour la topologie étale.

2) La donnée d'un faisceau inversible $\Mc$ sur $X$ et d'un isomorphisme
$\varphi$ de $\Mc^{\otimes n}$ vers $\Lc$
définit une section $s : X \fleche [\Lc^{\frac1n}]$ du morphisme structural $\pi$. Donc la gerbe \emph{fppf}
$[\Lc^{\frac1n}]$ sur $X$ est une gerbe neutre (au sens de \cite{LMB}~(3.20)). Le résultat découle donc de l'analogue \emph{fppf}
de \cite{LMB}~(3.21).

3) résulte
de 2) et du fait que $\Lc$ est localement trivial pour la topologie de Zariski.
\end{demo}

\begin{sousremarque}\rm
Si $n$ est inversible, alors $\mun$ est étale et les champs $\bmunfppf$ et $\bmun$ coïncident (\cite{LMB}~(9.6)). Dans
ce cas ce sont des champs de Deligne-Mumford. Notons que si $n$ n'est pas inversible, alors $\mun$ n'est pas
lisse, et le champ $\bmun$ qui classifie les $\mun$-torseurs étales n'a aucune raison \emph{a priori} d'être algébrique.
\end{sousremarque}

\begin{sousremarque}\rm
On suppose que $S$ est un $\Z[\frac1n]$-schéma et que $X$ est noethérien. Alors $\pi$ est propre, lisse, de
présentation finie, et cohomologiquement plat en dimension zéro. En particulier si $X/S$ vérifie ces
propriétés, le morphisme $[\Lc^{\frac1n}] \fleche S$ les vérifie aussi.
\end{sousremarque}

\bigskip
\noindent
{\sc Calcul du groupe de Picard de $[\Lc^{\frac1n}]$}

\begin{souslem}
Soient $X$ un schéma et $A$ un schéma en groupes commutatifs sur $X$. Soit $\Fc$ un
faisceau inversible sur une $A$-gerbe (\emph{fppf}) $\pi : \X \fleche X$.
Il existe un unique $X$-morphisme de schémas en groupes
$$\chi_{\Fc} : A \flechelongue \gm$$
tel que l'action naturelle de $A$ sur $\Fc$ soit induite par $\chi_{\Fc}$ et par la multiplication $\Fc\times \gm \fleche \Fc$
induite par la structure de $\Oc_{\X}$-module de $\Fc$, autrement dit tel que le diagramme suivant soit commutatif :
$$\shorthandoff{!;:?}
\xymatrix@!0 @R=1.5pc @C=3pc{A\times \Fc \ar[rd] \ar[dd] &\\ &\Fc\\ \gm\times\Fc \ar[ur]&}$$
\end{souslem}
\begin{demo} Un faisceau inversible sur la gerbe $\X$ est la donnée, pour tout $x\in\ob\X_U$, d'un faisceau inversible $\Fc_x$
sur $U$, et pour tout morphisme $\varphi : x\fleche x'$ dans $\X$, d'un isomorphisme
$$L_{\Fc}(\varphi) : \Fc_x \fleche \pi(\varphi)^*\Fc_{x'}$$
ces isomorphismes vérifiant de plus une condition de compatibilité évidente.

Construisons d'abord $\chi_{\Fc}(U)$ pour un $U\in\ob{\rm (Aff/}X)$ sur lequel $\X$ a des objets.
Soit $x\in\ob\X_U$ et soit $g\in A(U)$. Via l'identification entre $A(U)$ et $\Aut(x)$, $g$ correspond à un
automorphisme $\varphi$ de $x$, et induit de ce fait un automorphisme $L_{\Fc}(\varphi)$ de $\Fc_x$. Cet
automorphisme correspond à la multiplication par un unique élément de $\gm(U)$, dont on vérifie facilement qu'il ne dépend pas du choix de $x$ (en utilisant la condition de compatibilité entre les $L_{\Fc}(\varphi)$ et le fait que deux objets de $\X_U$ sont localement isomorphes
pour la topologie \emph{fppf}). D'où le morphisme
$\chi_{\Fc}(U) : A(U) \fleche \gm(U)$. Il est clair, vu sa construction, que ce morphisme est déterminé de manière
unique par les actions naturelles de $A$ et de $\gm$ sur $\Fc$.
\'Etant donné que $\X$ a des objets partout localement pour la topologie \emph{fppf}, cette collection de morphismes
se prolonge de manière unique en un caractère $\chi_{\Fc}$ de $A$ vérifiant les propriétés annoncées.
\end{demo}

\begin{sousremarque}\rm
Si $\chi$ est un caractère fixé de $A$, un faisceau inversible $\Fc$ sur $\X$ est un faisceau $\chi$-tordu
de degré $d$ (au sens de Lieblich, \cite{Lieblich}~2.1.2.2) si et seulement si $\chi_{\Fc}=\chi^d$.
\end{sousremarque}

\begin{souspptes}
\label{racine_pptes_caractere}
\begin{itemize}
\item[(1)] La construction de $\chi_{\Fc}$ est compatible au changement de base.
\item[(2)] Si $\Fc$ et $\Gc$ sont des faisceaux inversibles sur $\X$, alors $$\chi_{\Fc\otimes \Gc} = \chi_{\Fc} . \chi_{\Gc}\, .$$
\item[(3)] Un faisceau inversible $\Fc$ sur $\X$ provient de $X$ si et seulement si $\chi_{\Fc}$ est trivial.
\end{itemize}
\end{souspptes}
\begin{demo}
Les deux premières propriétés sont évidentes. Montrons le dernier
point.
Supposons tout d'abord que $\Fc$ soit isomorphe à un faisceau de la forme $\pi^*\Mc$, où $\Mc$ est un faisceau
inversible sur $X$. Il est clair que $\chi_{\Fc}$ est égal à $\chi_{\pi^*\Mc}$ donc il suffit de montrer que
$\chi_{\pi^*\Mc}$ est trivial, c'est-à-dire que pour tout objet $x$ de $\X$ et pour tout automorphisme
$\varphi$ de $x$, l'automorphisme $L_{\pi^*\Mc}(\varphi)$ de $(\pi^*\Mc)_x$ est l'identité. C'est évident par
construction de l'image inverse.

Réciproquement, supposons $\chi_{\Fc}$ trivial, \emph{i.e.} supposons que pour tout objet $x$ de $\X$ et tout
automorphisme $\varphi$ de $x$, $L_{\Fc}(\varphi)$ soit l'identité de $\Fc_x$. Provenir de la base est une question locale sur $X$ pour la topologie \emph{fppf} (voir par exemple le lemme~1.2.2.7 de~\cite{Brochard_these}).
On peut donc supposer que
le morphisme structural $\pi : \X \fleche X$ a une section $s : X \fleche \X$.

On va alors montrer que $\Fc$ est isomorphe à $\pi^*s^*\Fc$. Le faisceau $\pi^*s^*\Fc$ est celui qui à tout objet $x$ de $\X$
associe $\Fc_{s(\pi(x))}$, les isomorphismes de changement de base étant simplement les isomorphismes canoniques.
Nous allons construire une collection d'isomorphismes $\rho_x$ de $\Fc_x$ dans $\Fc_{s(\pi(x))}$ (pour chaque objet $x$ de $\X$), 
compatibles avec les $L_{\Fc}(\varphi)$ et les $L_{\pi^*s^*\Fc}(\varphi)$.

Si $x$ et $s(\pi(x))$ sont isomorphes dans $\X_U$ on choisit un isomorphisme $\varphi$ de $x$ dans $s(\pi(x))$ et on pose
$\rho_x=L_{\Fc}(\varphi)$. Si $\varphi_1$ et $\varphi_2$ sont deux tels isomorphismes, alors $(\varphi_2)^{-1}\circ \varphi_1$
est un automorphisme de $x$, donc d'après l'hypothèse sur $\Fc$ on a $L_{\Fc}((\varphi_2)^{-1}\circ \varphi_1)=\Id_{\Fc_x}$
de sorte que $\rho_x$ est bien défini et ne dépend pas du choix de $\varphi$.

Dans le cas général, on sait que $x$ et $s(\pi(x))$ sont localement isomorphes pour la topologie $\emph{fppf}$. Vu l'unicité dans la construction de $\rho_x$ lorsque $x$ est isomorphe à $s(\pi(x))$, il est clair
qu'il existe un unique isomorphisme $\rho_x : \Fc_x \fleche \Fc_{s(\pi(x))}$ compatible avec ceux construits dans le cas
précédent. La collection de tous les $\rho_x$ ainsi construits répond au problème posé.
\end{demo}

\begin{sousexemple}[groupe de Picard de B$G$] \rm
En utilisant cette construction, on retrouve facilement le groupe de Picard du champ classifiant B$G$, où $G$ est un
$X$-schéma en groupes abéliens. En effet, le morphisme structural $\pi : \text{B}G \fleche X$ a une section, donc
$\pi^* : \Pic(X) \fleche \Pic(\text{B}G)$ a une rétraction et en particulier il est injectif. D'après~\ref{racine_pptes_caractere}, l'application $\Fc \mapsto \chi_{\Fc}$ induit un morphisme de groupes de $\Pic(\text{B}G)$ dans $\widehat{G}$
dont $\Pic(X)$ est le noyau. Ce morphisme est naturellement scindé : si $\chi : G \fleche \gm$ est un caractère
de $G$, on lui associe la classe du faisceau inversible $\Lc(\chi)$ construit de la manière suivante. Pour tout
$U\in\ob\aff$ et tout $G$-torseur $\widetilde{U}$ on définit $\Lc(\chi)_{\widetilde{U}}$ comme étant le
faisceau inversible correspondant au $\gm$-torseur sur $U$ obtenu à partir de $\widetilde{U}$ par extension
du groupe structural via le caractère $\chi$. On a donc une suite exacte courte scindée :
$$\xymatrix{1 \ar[r]& \Pic(X) \ar[r] &\Pic(\text{B}G)\ar[r] & \widehat{G}\ar[r] &1}$$
de sorte que $\Pic(\text{B}G)$ est naturellement isomorphe au produit $\Pic(X)\times \widehat{G}$.
\end{sousexemple}

Dans le cas du champ $[\Lc^{\frac1n}]$, on a un faisceau inversible \og canonique\fg, que nous noterons $\Omega$, défini par $\Mc_{\alpha}=\Mc$ pour tout $U\in\ob\aff$ et tout objet $\alpha=(x,\Mc,\varphi)$ de
$[\Lc^{\frac1n}]_U$. Il est clair que le caractère $\chi_{\Omega}$ associé à $\Omega$ est simplement l'injection canonique
$$\chi : \mun \flechelongue \gm.$$
En particulier, vu que $\chi$ n'est pas le caractère trivial, on peut en déduire que le faisceau $\Omega$ ne
provient pas de la base $X$ ! On a enfin un isomorphisme canonique $\Phi : \flechen{\Omega^{\otimes n}}{\sim}{\pi^*\Lc}.$

\begin{sousprop} On note $l$ la classe du faisceau $\Lc$ dans $\Pic(X)$ et $\omega$ celle de $\Omega$ dans $\Pic([\Lc^{\frac1n}])$.
\begin{itemize}
\item[(1)] Le morphisme $\pi^* : \Pic(X) \fleche \Pic([\Lc^{\frac1n}])$ est injectif, et l'on a une suite exacte courte :
$$\xymatrix{1 \ar[r]& \Pic(X) \ar[r]& \Pic([\Lc^{\frac1n}]) \ar[r]& \widehat{\mun} \ar[r]& 1.}$$
\item[(2)] Le groupe $\Pic([\Lc^{\frac1n}])$ est isomorphe au quotient du groupe $\Pic(X) \times H^0(X,\Z)$ par le sous $H^0(X,\Z)$-module engendré par $(l^{-1},n)$ (autrement dit par la relation $\omega^n=l$).
\end{itemize}
\end{sousprop}
\begin{demo}
La propriété~(\ref{racine_pptes_caractere})~(2) montre que l'application $\Fc \mapsto \chi_{\Fc}$ induit un morphisme
de groupes de $\Pic([\Lc^{\frac1n}])$ dans $\widehat{\mun}$. Il est clair que ce morphisme est surjectif, vu que
$\widehat{\mun}$ est isomorphe au groupe $H^0(X,\Z/n\Z)$,
engendré par l'injection canonique $\chi : \mun \fleche \gm$, et que $\chi=\chi_{\Omega}$.
La propriété~(\ref{racine_pptes_caractere})~(3) montre que la suite ci-dessus est exacte en $\Pic([\Lc^{\frac1n}])$.
Pour en finir avec le premier point il nous reste donc juste à montrer l'injectivité de $\pi^*$.

Soit $\Nc$ un faisceau inversible sur $X$ et soit $f$ un isomorphisme de $\pi^*\Nc$ dans $\Oc_{[\Lc^{\frac1n}]}$. Il s'agit
de montrer que $\Nc$ est trivial. L'isomorphisme $f$ est donné par une collection d'isomorphismes
$$(f_{\alpha} : \flechen{(\pi^*\Nc)_{\alpha}=x^*\Nc}{\sim}{(\Oc_{[\Lc^{\frac1n}]})_{\alpha}=\Oc_U})_{\alpha=(x,\Mc,\varphi)\in [\Lc^{\frac1n}]_U}. $$
En utilisant la condition de compatibilité vérifiée par les $f_{\alpha}$ et la structure de gerbe de $[\Lc^{\frac1n}]$, on construit une famille compatible d'isomorphismes $f_x : \flechen{x^*\Nc}{\sim}{\Oc_U}$ indexée par les $x\in X(U)$, $U\in\ob\aff$, qui définit donc un
isomorphisme de $\Nc$ dans $\Oc_X$.

Pour le point (2), notons $G$ le quotient du groupe $\Pic(X) \times H^0(X,\Z)$ par le sous $H^0(X,\Z)$-module engendré par la relation $\omega^n=l$. On a clairement un morphisme de $G$ dans $\Pic([\Lc^{\frac1n}])$ qui envoie $(0,1)$ sur
$\omega$. En utilisant les propriétés précédentes, et le fait que $\widehat{\mun}$ est isomorphe à $H^0(X,\Z/n\Z)$
et engendré par $\chi_{\Omega}$, on vérifie facilement que ce morphisme est un isomorphisme.
\end{demo}

\begin{sousexemple}\rm
Prenons pour $X$ l'espace projectif $\P^k$ sur $\Spec \Z$. On fixe un entier relatif $l$ et on pose $\Lc=\Oc(l)$. La proposition précédente permet de calculer $\Pic([\Lc^{\frac1n}])$ pour tout $n$ appartenant à $\N^*$. Par exemple si $l=1$, on trouve $\frac1n\Z$. Si $l$ est un multiple de $n$, on est dans le cas où $\Lc$ a une racine \iem{n} et l'on trouve $\Z\times \Z/n\Z$. Dans le cas général, le groupe $\Pic([\Lc^{\frac1n}])$ est isomorphe (de manière non canonique) à $\frac{d}{n}\Z\times \Z/d\Z$ où $d$ est le pgcd de $n$ et $l$.
\end{sousexemple}

\noindent
{\sc Foncteur de Picard relatif de $[\Lc^{\frac1n}]/S$}

Notons $\X=[\Lc^{\frac1n}]$. Pour tout schéma $U$ sur $S$, on a une suite exacte courte :
$$\xymatrix{1 \ar[r]& \Pic(X\times_S U) \ar[r]& \Pic(\X\times_S U) \ar[r]& H^0(X\times_S U, \Z/n\Z) \ar[r]& 1.}$$
En appliquant le foncteur \og faisceau \'etale associé\fg\ elle induit une suite exacte de faisceaux étales :
\begin{equation}
\label{sec_racine_nieme}
\xymatrix{1 \ar[r]& \Pic_{X/S} \ar[r]^{\varphi_0}& \Pic_{\X/S} \ar[r]^{\chi}& f_*\Z/n\Z \ar[r]& 1.}
\end{equation}

\begin{sousremarque}\rm
Si $\Lc$ a une racine ${n}^{\textrm{ième}}$ $\Rc$, la suite exacte (\ref{sec_racine_nieme}) est scindée par $i \mapsto (\omega r^{-1})^i$ où $r$ est la classe de $\Rc$ dans $\Pic(X)$, si bien que $\pic$ s'identifie au produit $\Pic_{X/S}\times_S f_*\Z/n\Z$.
\end{sousremarque}

\begin{sousremarque}\rm
Le faisceau $f_*\Z/n\Z$ n'est \emph{a priori} pas représentable. En conséquence, dans le cas général, il ne suffit pas que $\Pic_{X/S}$ soit représentable pour que $\pic$ le soit, même lorsque la gerbe $\X$ est triviale. Cependant si $f$ est ouvert, dominant et à fibres géométriquement connexes (par exemple s'il est localement de type fini et cohomologiquement plat), alors $f_*\Z/n\Z= \Z/n\Z$. 
\end{sousremarque}

On suppose maintenant que $f_*\Z/n\Z$ coïncide avec $\Z/n\Z$ et l'on considère le produit $\Pic_{X/S}\times_S \Z$ de $\Pic_{X/S}$ par
le groupe constant $\Z$. On va voir que $\pic$ est isomorphe au quotient de
$\Pic_{X/S}\times_S \Z$ par la relation $\omega^n=l$. On note $H$ le sous-groupe engendré par $(l^{-1},n)$. En notant $\varphi_i$ le morphisme composé
$\xymatrix{\Pic_{X/S} \ar[r]^{\varphi_0}& \pic \ar[r]^{\mu_{\omega^i}}& \pic.}$
(où $\mu_{\omega^i}$ est la multiplication par $\omega^i$)
on définit donc un morphisme
$$\xymatrix{\varphi : \Pic_{X/S}\times_S \Z \ar[r]& \pic}$$
dont il est clair qu'il est invariant sous $H$ et universel pour cette propriété. Le foncteur $\pic$ s'identifie donc bien au quotient évoqué ci-dessus. On peut construire ce quotient \og à la main \fg\ comme suit (voir figure ci-dessous). Pour tout couple d'entiers $(i,k)$ on identifie les copies de $\Pic_{X/S}$ numéro $i$ et $i+nk$ via l'isomorphisme de translation 
$\mu_{l^k} : (\Pic_{X/S})_{i+nk} \fleche (\Pic_{X/S})_i.$
La loi de groupe est induite naturellement par celle de $\Pic_{X/S}$ et par la relation $\omega^n=l$. Ceci montre en particulier que si $\Pic_{X/S}$ est représentable, alors $\pic$ l'est aussi\footnote{Mais bien sûr, on le savait déjà dans le cas où $f$ est propre, plat et cohomologiquement plat en dimension zéro.}. On en déduit aussi que $\Pic_{X/S}$ et $\pic$ ont la même composante neutre.

%
%
%
\ifx\figfortexisloaded\relax \else\let\figfortexisloaded=\relax\fi
\message{version 1.8}
\newif\iftextures
\catcode`\@=11
\newdimen\epsil@n\epsil@n=0.00005pt
\newdimen\Cepsil@n\Cepsil@n=0.005pt
\newdimen\dcq@\dcq@=254pt
\newdimen\PI@\PI@=3.141592pt
\newdimen\DemiPI@deg\DemiPI@deg=90pt
\newdimen\PI@deg\PI@deg=180pt
\newdimen\DePI@deg\DePI@deg=360pt
\chardef\t@n=10
\chardef\c@nt=100
\chardef\@lxxiv=74
\chardef\@xci=91
\mathchardef\@nMnCQn=9949
\chardef\@vi=6
\chardef\@xxx=30
\chardef\@lvi=56
\chardef\@lxxi=71
\chardef\@lxxxv=85
\mathchardef\@mmmmlxviii=4068
\mathchardef\@ccclx=360
\mathchardef\@dccxx=720
\newcount\p@rtent \newcount\f@ctech \newcount\result@tent
\newdimen\v@lmin \newdimen\v@lmax \newdimen\v@leur
\newdimen\result@t\newdimen\result@@t
\newdimen\mili@u \newdimen\c@rre \newdimen\delt@
\def\degT@rd{0.017453 }  
\def\rdT@deg{57.295779 } 
{\catcode`p=12 \catcode`t=12 \gdef\v@leurseule#1pt{#1}}
\def\repdecn@mb#1{\expandafter\v@leurseule\the#1\space}
\def\arct@n#1(#2,#3){{\v@lmin=#2\v@lmax=#3%
    \maxim@m{\mili@u}{-\v@lmin}{\v@lmin}\maxim@m{\c@rre}{-\v@lmax}{\v@lmax}%
    \delt@=\mili@u\m@ech\mili@u%
    \ifdim\c@rre>\@nMnCQn\mili@u\divide\v@lmax\tw@\c@lATAN\v@leur(\z@,\v@lmax)
    \else%
    \maxim@m{\mili@u}{-\v@lmin}{\v@lmin}\maxim@m{\c@rre}{-\v@lmax}{\v@lmax}%
    \m@ech\c@rre%
    \ifdim\mili@u>\@nMnCQn\c@rre\divide\v@lmin\tw@
    \maxim@m{\mili@u}{-\v@lmin}{\v@lmin}\c@lATAN\v@leur(\mili@u,\z@)%
    \else\c@lATAN\v@leur(\delt@,\v@lmax)\fi\fi%
    \ifdim\v@lmin<\z@\v@leur=-\v@leur\ifdim\v@lmax<\z@\advance\v@leur-\PI@%
    \else\advance\v@leur\PI@\fi\fi%
    \global\result@t=\v@leur}#1=\result@t}
\def\m@ech#1{\ifdim#1>1.646pt\divide\mili@u\t@n\divide\c@rre\t@n\m@ech#1\fi}
\def\c@lATAN#1(#2,#3){{\v@lmin=#2\v@lmax=#3\v@leur=\z@\delt@=\tw@ pt%
    \un@iter{0.785398}{\v@lmax<}%
    \un@iter{0.463648}{\v@lmax<}%
    \un@iter{0.244979}{\v@lmax<}%
    \un@iter{0.124355}{\v@lmax<}%
    \un@iter{0.062419}{\v@lmax<}%
    \un@iter{0.031240}{\v@lmax<}%
    \un@iter{0.015624}{\v@lmax<}%
    \un@iter{0.007812}{\v@lmax<}%
    \un@iter{0.003906}{\v@lmax<}%
    \un@iter{0.001953}{\v@lmax<}%
    \un@iter{0.000976}{\v@lmax<}%
    \un@iter{0.000488}{\v@lmax<}%
    \un@iter{0.000244}{\v@lmax<}%
    \un@iter{0.000122}{\v@lmax<}%
    \un@iter{0.000061}{\v@lmax<}%
    \un@iter{0.000030}{\v@lmax<}%
    \un@iter{0.000015}{\v@lmax<}%
    \global\result@t=\v@leur}#1=\result@t}
\def\un@iter#1#2{%
    \divide\delt@\tw@\edef\dpmn@{\repdecn@mb{\delt@}}%
    \mili@u=\v@lmin%
    \ifdim#2\z@%
      \advance\v@lmin-\dpmn@\v@lmax\advance\v@lmax\dpmn@\mili@u%
      \advance\v@leur-#1pt%
    \else%
      \advance\v@lmin\dpmn@\v@lmax\advance\v@lmax-\dpmn@\mili@u%
      \advance\v@leur#1pt%
    \fi}
\def\c@ssin#1#2#3{\expandafter\ifx\csname COS@\number#3\endcsname\relax\c@lCS{#3pt}%
    \expandafter\xdef\csname COS@\number#3\endcsname{\repdecn@mb\result@t}%
    \expandafter\xdef\csname SIN@\number#3\endcsname{\repdecn@mb\result@@t}\fi%
    \edef#1{\csname COS@\number#3\endcsname}\edef#2{\csname SIN@\number#3\endcsname}}
\def\c@lCS#1{{\mili@u=#1\p@rtent=\@ne%
    \relax\ifdim\mili@u<\z@\red@ng<-\else\red@ng>+\fi\f@ctech=\p@rtent%
    \relax\ifdim\mili@u<\z@\mili@u=-\mili@u\f@ctech=-\f@ctech\fi\c@@lCS}}
\def\c@@lCS{\v@lmin=\mili@u\c@rre=-\mili@u\advance\c@rre\DemiPI@deg\v@lmax=\c@rre%
    \mili@u\@lxxi\mili@u\divide\mili@u\@mmmmlxviii%
    \edef\v@larg{\repdecn@mb{\mili@u}}\mili@u=-\v@larg\mili@u%
    \edef\v@lmxde{\repdecn@mb{\mili@u}}%
    \c@rre\@lxxi\c@rre\divide\c@rre\@mmmmlxviii%
    \edef\v@largC{\repdecn@mb{\c@rre}}\c@rre=-\v@largC\c@rre%
    \edef\v@lmxdeC{\repdecn@mb{\c@rre}}%
    \fctc@s\mili@u\v@lmin\global\result@t\p@rtent\v@leur%
    \let\t@mp=\v@larg\let\v@larg=\v@largC\let\v@largC=\t@mp%
    \let\t@mp=\v@lmxde\let\v@lmxde=\v@lmxdeC\let\v@lmxdeC=\t@mp%
    \fctc@s\c@rre\v@lmax\global\result@@t\f@ctech\v@leur}
\def\fctc@s#1#2{\v@leur=#1\relax\ifdim#2<\@lxxxv\p@\cosser@h\else\sinser@t\fi}
\def\cosser@h{\advance\v@leur\@lvi\p@\divide\v@leur\@lvi%
    \v@leur=\v@lmxde\v@leur\advance\v@leur\@xxx\p@%
    \v@leur=\v@lmxde\v@leur\advance\v@leur\@ccclx\p@%
    \v@leur=\v@lmxde\v@leur\advance\v@leur\@dccxx\p@\divide\v@leur\@dccxx}
\def\sinser@t{\v@leur=\v@lmxdeC\p@\advance\v@leur\@vi\p@%
    \v@leur=\v@largC\v@leur\divide\v@leur\@vi}
\def\red@ng#1#2{\relax\ifdim\mili@u#1#2\DemiPI@deg\advance\mili@u#2-\PI@deg%
    \p@rtent=-\p@rtent\red@ng#1#2\fi}
\def\invers@#1#2{{\v@leur=#2\maxim@m{\v@lmax}{-\v@leur}{\v@leur}%
    \f@ctech=\@ne\m@inv@rs%
    \multiply\v@leur\f@ctech\edef\v@lv@leur{\repdecn@mb{\v@leur}}%
    \p@rtentiere{\p@rtent}{\v@leur}\v@lmin=\p@\divide\v@lmin\p@rtent%
    \inv@rs@\multiply\v@lmax\f@ctech\global\result@t=\v@lmax}#1=\result@t}
\def\m@inv@rs{\ifdim\v@lmax<\p@\multiply\v@lmax\t@n\multiply\f@ctech\t@n\m@inv@rs\fi}
\def\inv@rs@{\v@lmax=-\v@lmin\v@lmax=\v@lv@leur\v@lmax%
    \advance\v@lmax\tw@ pt\v@lmax=\repdecn@mb{\v@lmin}\v@lmax%
    \delt@=\v@lmax\advance\delt@-\v@lmin\ifdim\delt@<\z@\delt@=-\delt@\fi%
    \ifdim\delt@>\epsil@n\v@lmin=\v@lmax\inv@rs@\fi}
\def\minim@m#1#2#3{\relax\ifdim#2<#3#1=#2\else#1=#3\fi}
\def\maxim@m#1#2#3{\relax\ifdim#2>#3#1=#2\else#1=#3\fi}
\def\p@rtentiere#1#2{#1=#2\divide#1by65536 }
\def\r@undint#1#2{{\v@leur=#2\divide\v@leur\t@n\p@rtentiere{\p@rtent}{\v@leur}%
    \v@leur=\p@rtent pt\global\result@t=\t@n\v@leur}#1=\result@t}
\def\sqrt@#1#2{{\v@leur=#2%
    \minim@m{\v@lmin}{\p@}{\v@leur}\maxim@m{\v@lmax}{\p@}{\v@leur}%
    \f@ctech=\@ne\m@sqrt@\sqrt@@%
    \mili@u=\v@lmin\advance\mili@u\v@lmax\divide\mili@u\tw@\multiply\mili@u\f@ctech%
    \global\result@t=\mili@u}#1=\result@t}
\def\m@sqrt@{\ifdim\v@leur>\dcq@\divide\v@leur\c@nt\v@lmax=\v@leur%
    \multiply\f@ctech\t@n\m@sqrt@\fi}
\def\sqrt@@{\mili@u=\v@lmin\advance\mili@u\v@lmax\divide\mili@u\tw@%
    \c@rre=\repdecn@mb{\mili@u}\mili@u%
    \ifdim\c@rre<\v@leur\v@lmin=\mili@u\else\v@lmax=\mili@u\fi%
    \delt@=\v@lmax\advance\delt@-\v@lmin\ifdim\delt@>\epsil@n\sqrt@@\fi}
\def\extrairelepremi@r#1\de#2{\expandafter\lepremi@r#2@#1#2}
\def\lepremi@r#1,#2@#3#4{\def#3{#1}\def#4{#2}\ignorespaces}
\def\@cfor#1:=#2\do#3{%
  \edef\@fortemp{#2}%
  \ifx\@fortemp\empty\else\@cforloop#2,\@nil,\@nil\@@#1{#3}\fi}
\def\@cforloop#1,#2\@@#3#4{%
  \def#3{#1}%
  \ifx#3\Fig@nnil\let\@nextwhile=\Fig@fornoop\else#4\relax\let\@nextwhile=\@cforloop\fi%
  \@nextwhile#2\@@#3{#4}}

\def\@ecfor#1:=#2\do#3{%
  \def\@@cfor{\@cfor#1:=}%
  \edef\@@@cfor{#2}%
  \expandafter\@@cfor\@@@cfor\do{#3}}
\def\Fig@nnil{\@nil}
\def\Fig@fornoop#1\@@#2#3{}
\def\trtlis@rg#1#2{\def\list@@rg{#1}%
    \@ecfor\p@rv@l:=\list@@rg\do{\expandafter#2\p@rv@l|}}
\newbox\b@xvisu
\newtoks\let@xte
\newif\ifitis@K
\newcount\s@mme
\newcount\l@mbd@un \newcount\l@mbd@de
\newcount\superc@ntr@l\superc@ntr@l=\@ne        
\newcount\typec@ntr@l\typec@ntr@l=\superc@ntr@l 
\newdimen\v@lX  \newdimen\v@lY  \newdimen\v@lZ
\newdimen\v@lXa \newdimen\v@lYa \newdimen\v@lZa
\newdimen\unit@\unit@=\p@ 
\def\unit@util{pt}
\def\ptT@ptps{0.996264 }
\def\ptpsT@pt{1.00375 }
\def\ptT@unit@{1} 
\def\setunit@#1{\def\unit@util{#1}\setunit@@#1:\invers@{\result@t}{\unit@}%
    \edef\ptT@unit@{\repdecn@mb\result@t}}
\def\setunit@@#1#2:{\ifcat#1a\unit@=\@ne#1#2\else\unit@=#1#2\fi}
\def\d@fm@cdim#1#2{{\v@leur=#2\v@leur=\ptT@unit@\v@leur\xdef#1{\repdecn@mb\v@leur}}}
\newif\ifBdingB@x\BdingB@xtrue
\newdimen\c@@rdXmin \newdimen\c@@rdYmin  
\newdimen\c@@rdXmax \newdimen\c@@rdYmax
\def\b@undb@x#1#2{\ifBdingB@x%
    \relax\ifdim#1<\c@@rdXmin\global\c@@rdXmin=#1\fi%
    \relax\ifdim#2<\c@@rdYmin\global\c@@rdYmin=#2\fi%
    \relax\ifdim#1>\c@@rdXmax\global\c@@rdXmax=#1\fi%
    \relax\ifdim#2>\c@@rdYmax\global\c@@rdYmax=#2\fi\fi}
\def\b@undb@xP#1{{\Figg@tXY{#1}\b@undb@x{\v@lX}{\v@lY}}}
\def\ellBB@x#1;#2,#3(#4,#5,#6){{\s@uvc@ntr@l\et@tellBB@x%
    \setc@ntr@l{2}\figptell-2::#1;#2,#3(#4,#6)\b@undb@xP{-2}%
    \figptell-2::#1;#2,#3(#5,#6)\b@undb@xP{-2}%
    \c@ssin{\C@}{\S@}{#6}\v@lmin=\C@ pt\v@lmax=\S@ pt%
    \mili@u=#3\v@lmin\delt@=#2\v@lmax\arct@n\v@leur(\delt@,\mili@u)%
    \mili@u=-#3\v@lmax\delt@=#2\v@lmin\arct@n\c@rre(\delt@,\mili@u)%
    \v@leur=\rdT@deg\v@leur\advance\v@leur-\DePI@deg%
    \c@rre=\rdT@deg\c@rre\advance\c@rre-\DePI@deg%
    \v@lmin=#4pt\v@lmax=#5pt%
    \loop\ifdim\v@leur<\v@lmax\ifdim\v@leur>\v@lmin%
    \edef\@ngle{\repdecn@mb\v@leur}\figptell-2::#1;#2,#3(\@ngle,#6)%
    \b@undb@xP{-2}\fi\advance\v@leur\PI@deg\repeat%
    \loop\ifdim\c@rre<\v@lmax\ifdim\c@rre>\v@lmin%
    \edef\@ngle{\repdecn@mb\c@rre}\figptell-2::#1;#2,#3(\@ngle,#6)%
    \b@undb@xP{-2}\fi\advance\c@rre\PI@deg\repeat%
    \resetc@ntr@l\et@tellBB@x}\ignorespaces}
\def\initb@undb@x{\c@@rdXmin=\maxdimen\c@@rdYmin=\maxdimen%
    \c@@rdXmax=-\maxdimen\c@@rdYmax=-\maxdimen}
\def\c@ntr@lnum#1{%
    \relax\ifnum\typec@ntr@l=\@ne%
    \ifnum#1<\z@%
    \immediate\write16{*** Forbidden point number (#1). Abort.}\end\fi\fi%
    \set@bjc@de{#1}}
\def\set@bjc@de#1{\edef\objc@de{@BJ\ifnum#1<\z@ M\romannumeral-#1\else\romannumeral#1\fi}}
\def\setc@ntr@l#1{\ifnum\superc@ntr@l>#1\typec@ntr@l=\superc@ntr@l%
    \else\typec@ntr@l=#1\fi}
\def\resetc@ntr@l#1{\global\superc@ntr@l=#1\setc@ntr@l{#1}}
\def\s@uvc@ntr@l#1{\edef#1{\the\superc@ntr@l}}
\def\c@lproscalDD#1[#2,#3]{{\Figg@tXY{#2}%
    \edef\Xu@{\repdecn@mb{\v@lX}}\edef\Yu@{\repdecn@mb{\v@lY}}\Figg@tXY{#3}%
    \global\result@t=\Xu@\v@lX\global\advance\result@t\Yu@\v@lY}#1=\result@t}
\def\c@lproscalTD#1[#2,#3]{{\Figg@tXY{#2}\edef\Xu@{\repdecn@mb{\v@lX}}%
    \edef\Yu@{\repdecn@mb{\v@lY}}\edef\Zu@{\repdecn@mb{\v@lZ}}%
    \Figg@tXY{#3}\global\result@t=\Xu@\v@lX\global\advance\result@t\Yu@\v@lY%
    \global\advance\result@t\Zu@\v@lZ}#1=\result@t}
\def\c@lprovec#1{%
    \det@rmC\v@lZa(\v@lX,\v@lY,\v@lmin,\v@lmax)%
    \det@rmC\v@lXa(\v@lY,\v@lZ,\v@lmax,\v@leur)%
    \det@rmC\v@lYa(\v@lZ,\v@lX,\v@leur,\v@lmin)%
    \Figv@ctCreg#1(\v@lXa,\v@lYa,\v@lZa)}
\def\det@rm#1[#2,#3]{{\Figg@tXY{#2}\Figg@tXYa{#3}%
    \delt@=\repdecn@mb{\v@lX}\v@lYa\advance\delt@-\repdecn@mb{\v@lY}\v@lXa%
    \global\result@t=\delt@}#1=\result@t}
\def\det@rmC#1(#2,#3,#4,#5){{\global\result@t=\repdecn@mb{#2}#5%
    \global\advance\result@t-\repdecn@mb{#3}#4}#1=\result@t}
\def\getredf@ctDD#1(#2,#3){{\maxim@m{\v@lXa}{-#2}{#2}\maxim@m{\v@lYa}{-#3}{#3}%
    \maxim@m{\v@lXa}{\v@lXa}{\v@lYa}
    \ifdim\v@lXa>\@xci pt\divide\v@lXa\@xci%
    \p@rtentiere{\p@rtent}{\v@lXa}\advance\p@rtent\@ne\else\p@rtent=\@ne\fi%
    \global\result@tent=\p@rtent}#1=\result@tent\ignorespaces}
\def\getredf@ctTD#1(#2,#3,#4){{\maxim@m{\v@lXa}{-#2}{#2}\maxim@m{\v@lYa}{-#3}{#3}%
    \maxim@m{\v@lZa}{-#4}{#4}\maxim@m{\v@lXa}{\v@lXa}{\v@lYa}%
    \maxim@m{\v@lXa}{\v@lXa}{\v@lZa}
    \ifdim\v@lXa>\@lxxiv pt\divide\v@lXa\@lxxiv%
    \p@rtentiere{\p@rtent}{\v@lXa}\advance\p@rtent\@ne\else\p@rtent=\@ne\fi%
    \global\result@tent=\p@rtent}#1=\result@tent\ignorespaces}
\def\FigptintercircB@zDD#1:#2:#3,#4[#5,#6,#7,#8]{{\s@uvc@ntr@l\et@tfigptintercircB@zDD%
    \setc@ntr@l{2}\figvectPDD-1[#5,#8]\Figg@tXY{-1}\getredf@ctDD\f@ctech(\v@lX,\v@lY)%
    \mili@u=#4\unit@\divide\mili@u\f@ctech\c@rre=\repdecn@mb{\mili@u}\mili@u%
    \figptBezierDD-5::#3[#5,#6,#7,#8]%
    \v@lmin=#3\p@\v@lmax=\v@lmin\advance\v@lmax0.1\p@%
    \loop\edef\T@{\repdecn@mb{\v@lmax}}\figptBezierDD-2::\T@[#5,#6,#7,#8]%
    \figvectPDD-1[-5,-2]\n@rmeucCDD{\delt@}{-1}\ifdim\delt@<\c@rre\v@lmin=\v@lmax%
    \advance\v@lmax0.1\p@\repeat%
    \loop\mili@u=\v@lmin\advance\mili@u\v@lmax%
    \divide\mili@u\tw@\edef\T@{\repdecn@mb{\mili@u}}\figptBezierDD-2::\T@[#5,#6,#7,#8]%
    \figvectPDD-1[-5,-2]\n@rmeucCDD{\delt@}{-1}\ifdim\delt@>\c@rre\v@lmax=\mili@u%
    \else\v@lmin=\mili@u\fi\v@leur=\v@lmax\advance\v@leur-\v@lmin%
    \ifdim\v@leur>\epsil@n\repeat\figptcopyDD#1:#2/-2/%
    \resetc@ntr@l\et@tfigptintercircB@zDD}\ignorespaces}
\def\inters@cDD#1:#2[#3,#4;#5,#6]{{\s@uvc@ntr@l\et@tinters@cDD%
    \setc@ntr@l{2}\vecunit@{-1}{#4}\vecunit@{-2}{#6}%
    \Figg@tXY{-1}\setc@ntr@l{1}\Figg@tXYa{#3}%
    \edef\A@{\repdecn@mb{\v@lX}}\edef\B@{\repdecn@mb{\v@lY}}%
    \v@lmin=\B@\v@lXa\advance\v@lmin-\A@\v@lYa%
    \Figg@tXYa{#5}\setc@ntr@l{2}\Figg@tXY{-2}%
    \edef\C@{\repdecn@mb{\v@lX}}\edef\D@{\repdecn@mb{\v@lY}}%
    \v@lmax=\D@\v@lXa\advance\v@lmax-\C@\v@lYa%
    \delt@=\A@\v@lY\advance\delt@-\B@\v@lX%
    \invers@{\v@leur}{\delt@}\edef\v@ldelta{\repdecn@mb{\v@leur}}%
    \v@lXa=\A@\v@lmax\advance\v@lXa-\C@\v@lmin%
    \v@lYa=\B@\v@lmax\advance\v@lYa-\D@\v@lmin%
    \v@lXa=\v@ldelta\v@lXa\v@lYa=\v@ldelta\v@lYa%
    \setc@ntr@l{1}\Figp@intregDD#1:{#2}(\v@lXa,\v@lYa)%
    \resetc@ntr@l\et@tinters@cDD}\ignorespaces}
\def\inters@cTD#1:#2[#3,#4;#5,#6]{{\s@uvc@ntr@l\et@tinters@cTD%
    \setc@ntr@l{2}\figvectNVTD-1[#4,#6]\figvectNVTD-2[#6,-1]\figvectPTD-1[#3,#5]%
    \r@pPSTD\v@leur[-2,-1,#4]\edef\v@lcoef{\repdecn@mb{\v@leur}}%
    \figpttraTD#1:{#2}=#3/\v@lcoef,#4/\resetc@ntr@l\et@tinters@cTD}\ignorespaces}
\def\r@pPSTD#1[#2,#3,#4]{{\Figg@tXY{#2}\edef\Xu@{\repdecn@mb{\v@lX}}%
    \edef\Yu@{\repdecn@mb{\v@lY}}\edef\Zu@{\repdecn@mb{\v@lZ}}%
    \Figg@tXY{#3}\v@lmin=\Xu@\v@lX\advance\v@lmin\Yu@\v@lY\advance\v@lmin\Zu@\v@lZ%
    \Figg@tXY{#4}\v@lmax=\Xu@\v@lX\advance\v@lmax\Yu@\v@lY\advance\v@lmax\Zu@\v@lZ%
    \invers@{\v@leur}{\v@lmax}\global\result@t=\repdecn@mb{\v@leur}\v@lmin}%
    #1=\result@t}
\def\n@rminfDD#1#2{{\Figg@tXY{#2}\maxim@m{\v@lX}{\v@lX}{-\v@lX}%
    \maxim@m{\v@lY}{\v@lY}{-\v@lY}\maxim@m{\global\result@t}{\v@lX}{\v@lY}}%
    #1=\result@t}
\def\n@rminfTD#1#2{{\Figg@tXY{#2}\maxim@m{\v@lX}{\v@lX}{-\v@lX}%
    \maxim@m{\v@lY}{\v@lY}{-\v@lY}\maxim@m{\v@lZ}{\v@lZ}{-\v@lZ}%
    \maxim@m{\v@lX}{\v@lX}{\v@lY}\maxim@m{\global\result@t}{\v@lX}{\v@lZ}}%
    #1=\result@t}
\def\n@rmeucCDD#1#2{\Figg@tXY{#2}\divide\v@lX\f@ctech\divide\v@lY\f@ctech%
    #1=\repdecn@mb{\v@lX}\v@lX\v@lX=\repdecn@mb{\v@lY}\v@lY\advance#1\v@lX}
\def\n@rmeucCTD#1#2{\Figg@tXY{#2}%
    \divide\v@lX\f@ctech\divide\v@lY\f@ctech\divide\v@lZ\f@ctech%
    #1=\repdecn@mb{\v@lX}\v@lX\v@lX=\repdecn@mb{\v@lY}\v@lY\advance#1\v@lX%
    \v@lX=\repdecn@mb{\v@lZ}\v@lZ\advance#1\v@lX}
\def\n@rmeucSVDD#1#2{{\Figg@tXY{#2}%
    \v@lXa=\repdecn@mb{\v@lX}\v@lX\v@lYa=\repdecn@mb{\v@lY}\v@lY%
    \advance\v@lXa\v@lYa\sqrt@{\global\result@t}{\v@lXa}}#1=\result@t}
\def\n@rmeucSVTD#1#2{{\Figg@tXY{#2}\v@lXa=\repdecn@mb{\v@lX}\v@lX%
    \v@lYa=\repdecn@mb{\v@lY}\v@lY\v@lZa=\repdecn@mb{\v@lZ}\v@lZ%
    \advance\v@lXa\v@lYa\advance\v@lXa\v@lZa\sqrt@{\global\result@t}{\v@lXa}}#1=\result@t}
\def\n@rmeucDD#1#2{{\Figg@tXY{#2}\getredf@ctDD\f@ctech(\v@lX,\v@lY)%
    \divide\v@lX\f@ctech\divide\v@lY\f@ctech%
    \v@lXa=\repdecn@mb{\v@lX}\v@lX\v@lYa=\repdecn@mb{\v@lY}\v@lY%
    \advance\v@lXa\v@lYa\sqrt@{\global\result@t}{\v@lXa}%
    \global\multiply\result@t\f@ctech}#1=\result@t}
\def\n@rmeucTD#1#2{{\Figg@tXY{#2}\getredf@ctTD\f@ctech(\v@lX,\v@lY,\v@lZ)%
    \divide\v@lX\f@ctech\divide\v@lY\f@ctech\divide\v@lZ\f@ctech%
    \v@lXa=\repdecn@mb{\v@lX}\v@lX%
    \v@lYa=\repdecn@mb{\v@lY}\v@lY\v@lZa=\repdecn@mb{\v@lZ}\v@lZ%
    \advance\v@lXa\v@lYa\advance\v@lXa\v@lZa\sqrt@{\global\result@t}{\v@lXa}%
    \global\multiply\result@t\f@ctech}#1=\result@t}
\def\vecunit@DD#1#2{{\Figg@tXY{#2}\getredf@ctDD\f@ctech(\v@lX,\v@lY)%
    \divide\v@lX\f@ctech\divide\v@lY\f@ctech%
    \Figv@ctCreg#1(\v@lX,\v@lY)\n@rmeucSV{\v@lYa}{#1}%
    \invers@{\v@lXa}{\v@lYa}\edef\v@lv@lXa{\repdecn@mb{\v@lXa}}%
    \v@lX=\v@lv@lXa\v@lX\v@lY=\v@lv@lXa\v@lY%
    \Figv@ctCreg#1(\v@lX,\v@lY)\multiply\v@lYa\f@ctech\global\result@t=\v@lYa}}
\def\vecunit@TD#1#2{{\Figg@tXY{#2}\getredf@ctTD\f@ctech(\v@lX,\v@lY,\v@lZ)%
    \divide\v@lX\f@ctech\divide\v@lY\f@ctech\divide\v@lZ\f@ctech%
    \Figv@ctCreg#1(\v@lX,\v@lY,\v@lZ)\n@rmeucSV{\v@lYa}{#1}%
    \invers@{\v@lXa}{\v@lYa}\edef\v@lv@lXa{\repdecn@mb{\v@lXa}}%
    \v@lX=\v@lv@lXa\v@lX\v@lY=\v@lv@lXa\v@lY\v@lZ=\v@lv@lXa\v@lZ%
    \Figv@ctCreg#1(\v@lX,\v@lY,\v@lZ)\multiply\v@lYa\f@ctech\global\result@t=\v@lYa}}
\def\vecunitC@TD[#1,#2]{\Figg@tXYa{#1}\Figg@tXY{#2}%
    \advance\v@lX-\v@lXa\advance\v@lY-\v@lYa\advance\v@lZ-\v@lZa\c@lvecunitTD}
\def\vecunitCV@TD#1{\Figg@tXY{#1}\c@lvecunitTD}
\def\c@lvecunitTD{\getredf@ctTD\f@ctech(\v@lX,\v@lY,\v@lZ)%
    \divide\v@lX\f@ctech\divide\v@lY\f@ctech\divide\v@lZ\f@ctech%
    \v@lXa=\repdecn@mb{\v@lX}\v@lX%
    \v@lYa=\repdecn@mb{\v@lY}\v@lY\v@lZa=\repdecn@mb{\v@lZ}\v@lZ%
    \advance\v@lXa\v@lYa\advance\v@lXa\v@lZa\sqrt@{\v@lYa}{\v@lXa}%
    \invers@{\v@lXa}{\v@lYa}\edef\v@lv@lXa{\repdecn@mb{\v@lXa}}%
    \v@lX=\v@lv@lXa\v@lX\v@lY=\v@lv@lXa\v@lY\v@lZ=\v@lv@lXa\v@lZ}
\def\figgetangleDD#1[#2,#3,#4]{\ifps@cri{\s@uvc@ntr@l\et@tfiggetangleDD\setc@ntr@l{2}%
    \figvectPDD-1[#2,#3]\figvectPDD-2[#2,#4]\vecunit@{-1}{-1}%
    \c@lproscalDD\delt@[-2,-1]\figvectNVDD-1[-1]\c@lproscalDD\v@leur[-2,-1]%
    \arct@n\v@lmax(\delt@,\v@leur)\v@lmax=\rdT@deg\v@lmax%
    \ifdim\v@lmax<\z@\advance\v@lmax\DePI@deg\fi\xdef#1{\repdecn@mb{\v@lmax}}%
    \resetc@ntr@l\et@tfiggetangleDD}\ignorespaces\fi}
\def\figgetangleTD#1[#2,#3,#4,#5]{\ifps@cri{\s@uvc@ntr@l\et@tfiggetangleTD\setc@ntr@l{2}%
    \figvectPTD-1[#2,#3]\figvectPTD-2[#2,#5]\figvectNVTD-3[-1,-2]%
    \figvectPTD-2[#2,#4]\figvectNVTD-4[-3,-1]%
    \vecunit@{-1}{-1}\c@lproscalTD\delt@[-2,-1]\c@lproscalTD\v@leur[-2,-4]%
    \arct@n\v@lmax(\delt@,\v@leur)\v@lmax=\rdT@deg\v@lmax%
    \ifdim\v@lmax<\z@\advance\v@lmax\DePI@deg\fi\xdef#1{\repdecn@mb{\v@lmax}}%
    \resetc@ntr@l\et@tfiggetangleTD}\ignorespaces\fi}    
\def\figgetdist#1[#2,#3]{\ifps@cri{\s@uvc@ntr@l\et@tfiggetdist\setc@ntr@l{2}%
    \figvectP-1[#2,#3]\n@rmeuc{\v@lX}{-1}\v@lX=\ptT@unit@\v@lX\xdef#1{\repdecn@mb{\v@lX}}%
    \resetc@ntr@l\et@tfiggetdist}\ignorespaces\fi}
\def\Figg@tT#1{\c@ntr@lnum{#1}%
    {\expandafter\expandafter\expandafter\extr@ctT\csname\objc@de\endcsname:%
     \ifnum\B@@ltxt=\z@\ptn@me{#1}\else\csname\objc@de T\endcsname\fi}}
\def\extr@ctT#1,#2,#3/#4:{\def\B@@ltxt{#3}}
\def\Figg@tXY#1{\c@ntr@lnum{#1}%
    \expandafter\expandafter\expandafter\extr@ctC\csname\objc@de\endcsname:}
\def\extr@ctCDD#1/#2,#3,#4:{\v@lX=#2\v@lY=#3}
\def\extr@ctCTD#1/#2,#3,#4:{\v@lX=#2\v@lY=#3\v@lZ=#4}
\def\Figg@tXYa#1{\c@ntr@lnum{#1}%
    \expandafter\expandafter\expandafter\extr@ctCa\csname\objc@de\endcsname:}
\def\extr@ctCaDD#1/#2,#3,#4:{\v@lXa=#2\v@lYa=#3}
\def\extr@ctCaTD#1/#2,#3,#4:{\v@lXa=#2\v@lYa=#3\v@lZa=#4}
\def\figinit#1{\initpr@lim\Figinit@#1,:\initpss@ttings\ignorespaces}
\def\Figinit@#1,#2:{\setunit@{#1}\def\t@xt@{#2}\ifx\t@xt@\empty\else\Figinit@@#2:\fi}
\def\Figinit@@#1#2:{\if#12 \else\Figs@tproj{#1}\initTD@\fi}
\newif\ifTr@isDim
\def\UnD@fined{UNDEFINED}
\def\ifundefined#1{\expandafter\ifx\csname#1\endcsname\relax}
\def\initpr@lim{\initb@undb@x\figsetmark{}\figsetptname{$A_{##1}$}\def\Sc@leFact{1}%
    \initDD@\figsetroundcoord{yes}\ps@critrue\expandafter\setupd@te\defaultupdate:%
    \edef\disob@unit{\UnD@fined}\edef\t@rgetpt{\UnD@fined}}
\def\initDD@{\Tr@isDimfalse%
    \ifPDFm@ke%
     \let\Ps@rcerc=\Ps@rcercBz%
     \let\Ps@rell=\Ps@rellBz%
    \fi
    \let\c@lDCUn=\c@lDCUnDD%
    \let\c@lDCDeux=\c@lDCDeuxDD%
    \let\c@ldefproj=\relax%
    \let\c@lproscal=\c@lproscalDD%
    \let\c@lprojSP=\relax%
    \let\extr@ctC=\extr@ctCDD%
    \let\extr@ctCa=\extr@ctCaDD%
    \let\extr@ctCF=\extr@ctCFDD%
    \let\Figp@intreg=\Figp@intregDD%
    \let\Figpts@xes=\Figpts@xesDD%
    \let\n@rmeucSV=\n@rmeucSVDD\let\n@rmeuc=\n@rmeucDD\let\n@rminf=\n@rminfDD%
    \let\pr@dMatV=\pr@dMatVDD%
    \let\ps@xes=\ps@xesDD%
    \let\vecunit@=\vecunit@DD%
    \let\figcoord=\figcoordDD%
    \let\figgetangle=\figgetangleDD%
    \let\figpt=\figptDD%
    \let\figptBezier=\figptBezierDD%
    \let\figptbary=\figptbaryDD%
    \let\figptcirc=\figptcircDD%
    \let\figptcircumcenter=\figptcircumcenterDD%
    \let\figptcopy=\figptcopyDD%
    \let\figptcurvcenter=\figptcurvcenterDD%
    \let\figptell=\figptellDD%
    \let\figptendnormal=\figptendnormalDD%
    \def\figptinterlineplane{\un@v@ilable{figptinterlineplane}}%
    \let\figptinterlines=\inters@cDD%
    \let\figptorthocenter=\figptorthocenterDD%
    \let\figptorthoprojline=\figptorthoprojlineDD%
    \def\figptorthoprojplane{\un@v@ilable{figptorthoprojplane}}%
    \let\figptrot=\figptrotDD%
    \let\figptscontrol=\figptscontrolDD%
    \let\figptsintercirc=\figptsintercircDD%
    \let\figptsinterlinell=\figptsinterlinellDD%
    \let\figptsorthoprojline=\figptsorthoprojlineDD%
    \let\figptsrot=\figptsrotDD%
    \let\figptssym=\figptssymDD%
    \let\figptstra=\figptstraDD%
    \let\figptsym=\figptsymDD%
    \let\figpttraC=\figpttraCDD%
    \let\figpttra=\figpttraDD%
    \def\figptvisilimSL{\un@v@ilable{figptvisilimSL}}%
    \def\figsetobdist{\un@v@ilable{figsetobdist}}%
    \def\figsettarget{\un@v@ilable{figsettarget}}%
    \def\figsetview{\un@v@ilable{figsetview}}%
    \let\figvectDBezier=\figvectDBezierDD%
    \let\figvectN=\figvectNDD%
    \let\figvectNV=\figvectNVDD%
    \let\figvectP=\figvectPDD%
    \let\figvectU=\figvectUDD%
    \let\psarccircP=\psarccircPDD%
    \let\psarccirc=\psarccircDD%
    \let\psarcell=\psarcellDD%
    \let\psarcellPA=\psarcellPADD%
    \let\psarrowBezier=\psarrowBezierDD%
    \let\psarrowcircP=\psarrowcircPDD%
    \let\psarrowcirc=\psarrowcircDD%
    \let\psarrowhead=\psarrowheadDD%
    \let\psarrow=\psarrowDD%
    \let\psBezier=\psBezierDD%
    \let\pscirc=\pscircDD%
    \let\pscurve=\pscurveDD%
    \let\psnormal=\psnormalDD%
    }
\def\initTD@{\Tr@isDimtrue\initb@undb@xTD\newt@rgetptfalse\newdis@bfalse%
    \let\c@lDCUn=\c@lDCUnTD%
    \let\c@lDCDeux=\c@lDCDeuxTD%
    \let\c@ldefproj=\c@ldefprojTD%
    \let\c@lproscal=\c@lproscalTD%
    \let\extr@ctC=\extr@ctCTD%
    \let\extr@ctCa=\extr@ctCaTD%
    \let\extr@ctCF=\extr@ctCFTD%
    \let\Figp@intreg=\Figp@intregTD%
    \let\Figpts@xes=\Figpts@xesTD%
    \let\n@rmeucSV=\n@rmeucSVTD\let\n@rmeuc=\n@rmeucTD\let\n@rminf=\n@rminfTD%
    \let\pr@dMatV=\pr@dMatVTD%
    \let\ps@xes=\ps@xesTD%
    \let\vecunit@=\vecunit@TD%
    \let\figcoord=\figcoordTD%
    \let\figgetangle=\figgetangleTD%
    \let\figpt=\figptTD%
    \let\figptBezier=\figptBezierTD%
    \let\figptbary=\figptbaryTD%
    \let\figptcirc=\figptcircTD%
    \let\figptcircumcenter=\figptcircumcenterTD%
    \let\figptcopy=\figptcopyTD%
    \let\figptcurvcenter=\figptcurvcenterTD%
    \let\figptinterlineplane=\figptinterlineplaneTD%
    \let\figptinterlines=\inters@cTD%
    \let\figptorthocenter=\figptorthocenterTD%
    \let\figptorthoprojline=\figptorthoprojlineTD%
    \let\figptorthoprojplane=\figptorthoprojplaneTD%
    \let\figptrot=\figptrotTD%
    \let\figptscontrol=\figptscontrolTD%
    \let\figptsintercirc=\figptsintercircTD%
    \let\figptsorthoprojline=\figptsorthoprojlineTD%
    \let\figptsorthoprojplane=\figptsorthoprojplaneTD%
    \let\figptsrot=\figptsrotTD%
    \let\figptssym=\figptssymTD%
    \let\figptstra=\figptstraTD%
    \let\figptsym=\figptsymTD%
    \let\figpttraC=\figpttraCTD%
    \let\figpttra=\figpttraTD%
    \let\figptvisilimSL=\figptvisilimSLTD%
    \let\figsetobdist=\figsetobdistTD%
    \let\figsettarget=\figsettargetTD%
    \let\figsetview=\figsetviewTD%
    \let\figvectDBezier=\figvectDBezierTD%
    \let\figvectN=\figvectNTD%
    \let\figvectNV=\figvectNVTD%
    \let\figvectP=\figvectPTD%
    \let\figvectU=\figvectUTD%
    \let\psarccircP=\psarccircPTD%
    \let\psarccirc=\psarccircTD%
    \let\psarcell=\psarcellTD%
    \let\psarcellPA=\psarcellPATD%
    \let\psarrowBezier=\psarrowBezierTD%
    \let\psarrowcircP=\psarrowcircPTD%
    \let\psarrowcirc=\psarrowcircTD%
    \let\psarrowhead=\psarrowheadTD%
    \let\psarrow=\psarrowTD%
    \let\psBezier=\psBezierTD%
    \let\pscirc=\pscircTD%
    \let\pscurve=\pscurveTD%
    }
\def\un@v@ilable#1{\immediate\write16{*** The macro #1 is not available in the current context.}}
\def\figinsert#1{{\def\t@xt@{#1}\relax\ifx\t@xt@\empty\Figinsert@\DefGIfilen@me,:%
    \else\expandafter\FiginsertNu@#1 :\fi}\ignorespaces}
\def\FiginsertNu@#1 #2:{\def\t@xt@{#1}\relax\ifx\t@xt@\empty\def\t@xt@{#2}%
    \ifx\t@xt@\empty\Figinsert@\DefGIfilen@me,:\else\FiginsertNu@#2:\fi%
    \else\expandafter\FiginsertNd@#1 #2:\fi}
\def\FiginsertNd@#1#2:{\ifcat#1a\Figinsert@#1#2,:\else\Figinsert@\DefGIfilen@me,#1#2,:\fi}
\def\Figinsert@#1,#2:{\def\t@xt@{#2}\ifx\t@xt@\empty\xdef\Sc@leFact{1}\else%
    \def\Xarg@##1,{\def\@rgdeux{##1}}\Xarg@#2\xdef\Sc@leFact{\@rgdeux}\fi\@psfgetbb{#1}%
    \v@lX=\@psfllx\p@\v@lX=\ptpsT@pt\v@lX\v@lX=\Sc@leFact\v@lX%
    \v@lY=\@psflly\p@\v@lY=\ptpsT@pt\v@lY\v@lY=\Sc@leFact\v@lY%
    \b@undb@x{\v@lX}{\v@lY}%
    \v@lX=\@psfurx\p@\v@lX=\ptpsT@pt\v@lX\v@lX=\Sc@leFact\v@lX%
    \v@lY=\@psfury\p@\v@lY=\ptpsT@pt\v@lY\v@lY=\Sc@leFact\v@lY%
    \b@undb@x{\v@lX}{\v@lY}%
    \ifPDFm@ke\Figinclud@PDF{#1}{\Sc@leFact}\else%
    \v@lX=\c@nt pt\v@lX=\Sc@leFact\v@lX\edef\F@ct{\repdecn@mb{\v@lX}}%
    \iftextures\special{postscriptfile #1 vscale=\F@ct\space hscale=\F@ct}%
    \else\includegraphics{#1}\fi\fi%
    \message{[#1]}\ignorespaces}
\def\figinsertE#1{\FiginsertE@#1,:\ignorespaces}
\def\FiginsertE@#1,#2:{{\def\t@xt@{#2}\ifx\t@xt@\empty\xdef\Sc@leFact{1}\else%
    \def\Xarg@##1,{\def\@rgdeux{##1}}\Xarg@#2\xdef\Sc@leFact{\@rgdeux}\fi%
    \pdfximage{#1}\setbox\Gb@x=\hbox{\pdfrefximage\pdflastximage}%
    \v@lX=\z@\v@lY=-\Sc@leFact\dp\Gb@x\b@undb@x{\v@lX}{\v@lY}%
    \advance\v@lX\Sc@leFact\wd\Gb@x\advance\v@lY\Sc@leFact\dp\Gb@x%
    \advance\v@lY\Sc@leFact\ht\Gb@x\b@undb@x{\v@lX}{\v@lY}%
    \v@lX=\Sc@leFact\wd\Gb@x\pdfximage width \v@lX {#1}%
    \rlap{\pdfrefximage\pdflastximage}\message{[#1]}}\ignorespaces}
\def\figptDD#1:#2(#3,#4){\ifps@cri\c@ntr@lnum{#1}%
    {\v@lX=#3\unit@\v@lY=#4\unit@\Fig@dmpt{#2}{\z@}}\ignorespaces\fi}
\def\Fig@dmpt#1#2{\def\t@xt@{#1}\ifx\t@xt@\empty\def\B@@ltxt{\z@}%
    \else\expandafter\gdef\csname\objc@de T\endcsname{#1}\def\B@@ltxt{\@ne}\fi%
    \expandafter\xdef\csname\objc@de\endcsname{\ifitis@vect@r\C@dCl@svect%
    \else\C@dCl@spt\fi,\z@,\B@@ltxt/\the\v@lX,\the\v@lY,#2}}
\def\C@dCl@spt{P}
\def\C@dCl@svect{V}
\def\figptTD#1:#2(#3,#4){\ifps@cri\c@ntr@lnum{#1}%
    \def\c@@rdYZ{#4,0,0}\extrairelepremi@r\c@@rdY\de\c@@rdYZ%
    \extrairelepremi@r\c@@rdZ\de\c@@rdYZ%
    {\v@lX=#3\unit@\v@lY=\c@@rdY\unit@\v@lZ=\c@@rdZ\unit@\Fig@dmpt{#2}{\the\v@lZ}%
    \b@undb@xTD{\v@lX}{\v@lY}{\v@lZ}}\ignorespaces\fi}
\def\Figp@intregDD#1:#2(#3,#4){\c@ntr@lnum{#1}%
    {\result@t=#4\v@lX=#3\v@lY=\result@t\Fig@dmpt{#2}{\z@}}\ignorespaces}
\def\Figp@intregTD#1:#2(#3,#4){\c@ntr@lnum{#1}%
    \def\c@@rdYZ{#4,\z@,\z@}\extrairelepremi@r\c@@rdY\de\c@@rdYZ%
    \extrairelepremi@r\c@@rdZ\de\c@@rdYZ%
    {\v@lX=#3\v@lY=\c@@rdY\v@lZ=\c@@rdZ\Fig@dmpt{#2}{\the\v@lZ}%
    \b@undb@xTD{\v@lX}{\v@lY}{\v@lZ}}\ignorespaces}
\def\figptBezierDD#1:#2:#3[#4,#5,#6,#7]{\ifps@cri{\s@uvc@ntr@l\et@tfigptBezierDD%
    \FigptBezier@#3[#4,#5,#6,#7]\Figp@intregDD#1:{#2}(\v@lX,\v@lY)%
    \resetc@ntr@l\et@tfigptBezierDD}\ignorespaces\fi}
\def\figptBezierTD#1:#2:#3[#4,#5,#6,#7]{\ifps@cri{\s@uvc@ntr@l\et@tfigptBezierTD%
    \FigptBezier@#3[#4,#5,#6,#7]\Figp@intregTD#1:{#2}(\v@lX,\v@lY,\v@lZ)%
    \resetc@ntr@l\et@tfigptBezierTD}\ignorespaces\fi}
\def\FigptBezier@#1[#2,#3,#4,#5]{\setc@ntr@l{2}%
    \edef\T@{#1}\v@leur=\p@\advance\v@leur-#1pt\edef\UNmT@{\repdecn@mb{\v@leur}}%
    \figptcopy-4:/#2/\figptcopy-3:/#3/\figptcopy-2:/#4/\figptcopy-1:/#5/%
    \l@mbd@un=-4 \l@mbd@de=-\thr@@\p@rtent=\m@ne\c@lDecast%
    \l@mbd@un=-4 \l@mbd@de=-\thr@@\p@rtent=-\tw@\c@lDecast%
    \l@mbd@un=-4 \l@mbd@de=-\thr@@\p@rtent=-\thr@@\c@lDecast\Figg@tXY{-4}}
\def\c@lDCUnDD#1#2{\Figg@tXY{#1}\v@lX=\UNmT@\v@lX\v@lY=\UNmT@\v@lY%
    \Figg@tXYa{#2}\advance\v@lX\T@\v@lXa\advance\v@lY\T@\v@lYa%
    \Figp@intregDD#1:(\v@lX,\v@lY)}
\def\c@lDCUnTD#1#2{\Figg@tXY{#1}\v@lX=\UNmT@\v@lX\v@lY=\UNmT@\v@lY\v@lZ=\UNmT@\v@lZ%
    \Figg@tXYa{#2}\advance\v@lX\T@\v@lXa\advance\v@lY\T@\v@lYa\advance\v@lZ\T@\v@lZa%
    \Figp@intregTD#1:(\v@lX,\v@lY,\v@lZ)}
\def\c@lDecast{\relax\ifnum\l@mbd@un<\p@rtent\c@lDCUn{\l@mbd@un}{\l@mbd@de}%
    \advance\l@mbd@un\@ne\advance\l@mbd@de\@ne\c@lDecast\fi}
\def\figptmap#1:#2=#3/#4/#5/{\ifps@cri{\s@uvc@ntr@l\et@tfigptmap%
    \setc@ntr@l{2}\figvectP-1[#4,#3]\Figg@tXY{-1}%
    \pr@dMatV/#5/\figpttra#1:{#2}=#4/1,-1/%
    \resetc@ntr@l\et@tfigptmap}\ignorespaces\fi}
\def\pr@dMatVDD/#1,#2;#3,#4/{\v@lXa=#1\v@lX\advance\v@lXa#2\v@lY%
    \v@lYa=#3\v@lX\advance\v@lYa#4\v@lY\Figv@ctCreg-1(\v@lXa,\v@lYa)}
\def\pr@dMatVTD/#1,#2,#3;#4,#5,#6;#7,#8,#9/{%
    \v@lXa=#1\v@lX\advance\v@lXa#2\v@lY\advance\v@lXa#3\v@lZ%
    \v@lYa=#4\v@lX\advance\v@lYa#5\v@lY\advance\v@lYa#6\v@lZ%
    \v@lZa=#7\v@lX\advance\v@lZa#8\v@lY\advance\v@lZa#9\v@lZ%
    \Figv@ctCreg-1(\v@lXa,\v@lYa,\v@lZa)}
\def\figptbaryDD#1:#2[#3;#4]{\ifps@cri{\edef\list@num{#3}\extrairelepremi@r\p@int\de\list@num%
    \s@mme=\z@\@ecfor\c@ef:=#4\do{\advance\s@mme\c@ef}%
    \edef\listec@ef{#4,0}\extrairelepremi@r\c@ef\de\listec@ef%
    \Figg@tXY{\p@int}\divide\v@lX\s@mme\divide\v@lY\s@mme%
    \multiply\v@lX\c@ef\multiply\v@lY\c@ef%
    \@ecfor\p@int:=\list@num\do{\extrairelepremi@r\c@ef\de\listec@ef%
           \Figg@tXYa{\p@int}\divide\v@lXa\s@mme\divide\v@lYa\s@mme%
           \multiply\v@lXa\c@ef\multiply\v@lYa\c@ef%
           \advance\v@lX\v@lXa\advance\v@lY\v@lYa}%
    \Figp@intregDD#1:{#2}(\v@lX,\v@lY)}\ignorespaces\fi}
\def\figptbaryTD#1:#2[#3;#4]{\ifps@cri{\edef\list@num{#3}\extrairelepremi@r\p@int\de\list@num%
    \s@mme=\z@\@ecfor\c@ef:=#4\do{\advance\s@mme\c@ef}%
    \edef\listec@ef{#4,0}\extrairelepremi@r\c@ef\de\listec@ef%
    \Figg@tXY{\p@int}\divide\v@lX\s@mme\divide\v@lY\s@mme\divide\v@lZ\s@mme%
    \multiply\v@lX\c@ef\multiply\v@lY\c@ef\multiply\v@lZ\c@ef%
    \@ecfor\p@int:=\list@num\do{\extrairelepremi@r\c@ef\de\listec@ef%
           \Figg@tXYa{\p@int}\divide\v@lXa\s@mme\divide\v@lYa\s@mme\divide\v@lZa\s@mme%
           \multiply\v@lXa\c@ef\multiply\v@lYa\c@ef\multiply\v@lZa\c@ef%
           \advance\v@lX\v@lXa\advance\v@lY\v@lYa\advance\v@lZ\v@lZa}%
    \Figp@intregTD#1:{#2}(\v@lX,\v@lY,\v@lZ)}\ignorespaces\fi}
\def\figptbaryR#1:#2[#3;#4]{\ifps@cri{%
    \v@leur=\z@\@ecfor\c@ef:=#4\do{\maxim@m{\v@lmax}{\c@ef pt}{-\c@ef pt}%
    \ifdim\v@lmax>\v@leur\v@leur=\v@lmax\fi}%
    \ifdim\v@leur<\p@\f@ctech=\@M\else\ifdim\v@leur<\t@n\p@\f@ctech=\@m\else%
    \ifdim\v@leur<\c@nt\p@\f@ctech=\c@nt\else\ifdim\v@leur<\@m\p@\f@ctech=\t@n\else%
    \f@ctech=\@ne\fi\fi\fi\fi%
    \def\listec@ef{0}%
    \@ecfor\c@ef:=#4\do{\sc@lec@nvRI{\c@ef pt}\edef\listec@ef{\listec@ef,\the\s@mme}}%
    \extrairelepremi@r\c@ef\de\listec@ef\figptbary#1:#2[#3;\listec@ef]}\ignorespaces\fi}
\def\sc@lec@nvRI#1{\v@leur=#1\p@rtentiere{\s@mme}{\v@leur}\advance\v@leur-\s@mme\p@%
    \multiply\v@leur\f@ctech\p@rtentiere{\p@rtent}{\v@leur}%
    \multiply\s@mme\f@ctech\advance\s@mme\p@rtent}
\def\figptcircDD#1:#2:#3;#4(#5){\ifps@cri{\s@uvc@ntr@l\et@tfigptcircDD%
    \c@lptellDD#1:{#2}:#3;#4,#4(#5)\resetc@ntr@l\et@tfigptcircDD}\ignorespaces\fi}
\def\figptcircTD#1:#2:#3,#4,#5;#6(#7){\ifps@cri{\s@uvc@ntr@l\et@tfigptcircTD%
    \setc@ntr@l{2}\c@lExtAxes#3,#4,#5(#6)\figptellP#1:{#2}:#3,-4,-5(#7)%
    \resetc@ntr@l\et@tfigptcircTD}\ignorespaces\fi}
\def\figptcircumcenterDD#1:#2[#3,#4,#5]{\ifps@cri{\s@uvc@ntr@l\et@tfigptcircumcenterDD%
    \setc@ntr@l{2}\figvectNDD-5[#3,#4]\figptbaryDD-3:[#3,#4;1,1]%
                  \figvectNDD-6[#4,#5]\figptbaryDD-4:[#4,#5;1,1]%
    \resetc@ntr@l{2}\inters@cDD#1:{#2}[-3,-5;-4,-6]%
    \resetc@ntr@l\et@tfigptcircumcenterDD}\ignorespaces\fi}
\def\figptcircumcenterTD#1:#2[#3,#4,#5]{\ifps@cri{\s@uvc@ntr@l\et@tfigptcircumcenterTD%
    \setc@ntr@l{2}\figvectNTD-1[#3,#4,#5]%
    \figvectPTD-3[#3,#4]\figvectNVTD-5[-1,-3]\figptbaryTD-3:[#3,#4;1,1]%
    \figvectPTD-4[#4,#5]\figvectNVTD-6[-1,-4]\figptbaryTD-4:[#4,#5;1,1]%
    \resetc@ntr@l{2}\inters@cTD#1:{#2}[-3,-5;-4,-6]%
    \resetc@ntr@l\et@tfigptcircumcenterTD}\ignorespaces\fi}
\def\figptcopyDD#1:#2/#3/{\ifps@cri{\Figg@tXY{#3}%
    \Figp@intregDD#1:{#2}(\v@lX,\v@lY)}\ignorespaces\fi}
\def\figptcopyTD#1:#2/#3/{\ifps@cri{\Figg@tXY{#3}%
    \Figp@intregTD#1:{#2}(\v@lX,\v@lY,\v@lZ)}\ignorespaces\fi}
\def\figptcurvcenterDD#1:#2:#3[#4,#5,#6,#7]{\ifps@cri{\s@uvc@ntr@l\et@tfigptcurvcenterDD%
    \setc@ntr@l{2}\c@lcurvradDD#3[#4,#5,#6,#7]\edef\Sprim@{\repdecn@mb{\result@t}}%
    \figptBezierDD-1::#3[#4,#5,#6,#7]\figpttraDD#1:{#2}=-1/\Sprim@,-5/%
    \resetc@ntr@l\et@tfigptcurvcenterDD}\ignorespaces\fi}
\def\figptcurvcenterTD#1:#2:#3[#4,#5,#6,#7]{\ifps@cri{\s@uvc@ntr@l\et@tfigptcurvcenterTD%
    \setc@ntr@l{2}\figvectDBezierTD -5:1,#3[#4,#5,#6,#7]%
    \figvectDBezierTD -6:2,#3[#4,#5,#6,#7]\vecunit@TD{-5}{-5}%
    \edef\Sprim@{\repdecn@mb{\result@t}}\figvectNVTD-1[-6,-5]%
    \figvectNVTD-5[-5,-1]\c@lproscalTD\v@leur[-6,-5]%
    \invers@{\v@leur}{\v@leur}\v@leur=\Sprim@\v@leur\v@leur=\Sprim@\v@leur%
    \figptBezierTD-1::#3[#4,#5,#6,#7]\edef\Sprim@{\repdecn@mb{\v@leur}}%
    \figpttraTD#1:{#2}=-1/\Sprim@,-5/\resetc@ntr@l\et@tfigptcurvcenterTD}\ignorespaces\fi}
\def\c@lcurvradDD#1[#2,#3,#4,#5]{{\figvectDBezierDD -5:1,#1[#2,#3,#4,#5]%
    \figvectDBezierDD -6:2,#1[#2,#3,#4,#5]\vecunit@DD{-5}{-5}%
    \edef\Sprim@{\repdecn@mb{\result@t}}\figvectNVDD-5[-5]\c@lproscalDD\v@leur[-6,-5]%
    \invers@{\v@leur}{\v@leur}\v@leur=\Sprim@\v@leur\v@leur=\Sprim@\v@leur%
    \global\result@t=\v@leur}}
\def\figptellDD#1:#2:#3;#4,#5(#6,#7){\ifps@cri{\s@uvc@ntr@l\et@tfigptell%
    \c@lptellDD#1::#3;#4,#5(#6)\figptrotDD#1:{#2}=#1/#3,#7/%
    \resetc@ntr@l\et@tfigptell}\ignorespaces\fi}
\def\c@lptellDD#1:#2:#3;#4,#5(#6){\c@ssin{\C@}{\S@}{#6}\v@lmin=\C@ pt\v@lmax=\S@ pt%
    \v@lmin=#4\v@lmin\v@lmax=#5\v@lmax%
    \edef\Xc@mp{\repdecn@mb{\v@lmin}}\edef\Yc@mp{\repdecn@mb{\v@lmax}}%
    \setc@ntr@l{2}\figvectC-1(\Xc@mp,\Yc@mp)\figpttraDD#1:{#2}=#3/1,-1/}
\def\figptellP#1:#2:#3,#4,#5(#6){\ifps@cri{\s@uvc@ntr@l\et@tfigptellP%
    \setc@ntr@l{2}\figvectP-1[#3,#4]\figvectP-2[#3,#5]%
    \v@leur=#6pt\c@lptellP{#3}{-1}{-2}\figptcopy#1:{#2}/-3/%
    \resetc@ntr@l\et@tfigptellP}\ignorespaces\fi}
\def\c@lptellP#1#2#3{\edef\@ngle{\repdecn@mb\v@leur}\c@ssin{\C@}{\S@}{\@ngle}%
    \figpttra-3:=#1/\C@,#2/\figpttra-3:=-3/\S@,#3/}
\def\figptendnormalDD#1:#2:#3,#4[#5,#6]{\ifps@cri{\s@uvc@ntr@l\et@tfigptendnormal%
    \Figg@tXYa{#5}\Figg@tXY{#6}%
    \advance\v@lX-\v@lXa\advance\v@lY-\v@lYa%
    \setc@ntr@l{2}\Figv@ctCreg-1(\v@lX,\v@lY)\vecunit@{-1}{-1}\Figg@tXY{-1}%
    \delt@=#3\unit@\maxim@m{\delt@}{\delt@}{-\delt@}\edef\l@ngueur{\repdecn@mb{\delt@}}%
    \v@lX=\l@ngueur\v@lX\v@lY=\l@ngueur\v@lY%
    \delt@=\p@\advance\delt@-#4pt\edef\l@ngueur{\repdecn@mb{\delt@}}%
    \figptbaryR-1:[#5,#6;#4,\l@ngueur]\Figg@tXYa{-1}%
    \advance\v@lXa\v@lY\advance\v@lYa-\v@lX%
    \setc@ntr@l{1}\Figp@intregDD#1:{#2}(\v@lXa,\v@lYa)\resetc@ntr@l\et@tfigptendnormal}%
    \ignorespaces\fi}
\def\figptexcenter#1:#2[#3,#4,#5]{\ifps@cri{\let@xte={-}%
    \Figptexinsc@nter#1:#2[#3,#4,#5]}\ignorespaces\fi}
\def\figptincenter#1:#2[#3,#4,#5]{\ifps@cri{\let@xte={}%
    \Figptexinsc@nter#1:#2[#3,#4,#5]}\ignorespaces\fi}
\let\figptinscribedcenter=\figptincenter
\def\Figptexinsc@nter#1:#2[#3,#4,#5]{%
    \figgetdist\LA@[#4,#5]\figgetdist\LB@[#3,#5]\figgetdist\LC@[#3,#4]%
    \figptbaryR#1:{#2}[#3,#4,#5;\the\let@xte\LA@,\LB@,\LC@]}
\def\figptinterlineplaneTD#1:#2[#3,#4;#5,#6]{\ifps@cri{\s@uvc@ntr@l\et@tfigptinterlineplane%
    \setc@ntr@l{2}\figvectPTD-1[#3,#5]\vecunit@TD{-2}{#6}%
    \r@pPSTD\v@leur[-2,-1,#4]\edef\v@lcoef{\repdecn@mb{\v@leur}}%
    \figpttraTD#1:{#2}=#3/\v@lcoef,#4/\resetc@ntr@l\et@tfigptinterlineplane}\ignorespaces\fi}
\def\figptorthocenterDD#1:#2[#3,#4,#5]{\ifps@cri{\s@uvc@ntr@l\et@tfigptorthocenterDD%
    \setc@ntr@l{2}\figvectNDD-3[#3,#4]\figvectNDD-4[#4,#5]%
    \resetc@ntr@l{2}\inters@cDD#1:{#2}[#5,-3;#3,-4]%
    \resetc@ntr@l\et@tfigptorthocenterDD}\ignorespaces\fi}
\def\figptorthocenterTD#1:#2[#3,#4,#5]{\ifps@cri{\s@uvc@ntr@l\et@tfigptorthocenterTD%
    \setc@ntr@l{2}\figvectNTD-1[#3,#4,#5]%
    \figvectPTD-2[#3,#4]\figvectNVTD-3[-1,-2]%
    \figvectPTD-2[#4,#5]\figvectNVTD-4[-1,-2]%
    \resetc@ntr@l{2}\inters@cTD#1:{#2}[#5,-3;#3,-4]%
    \resetc@ntr@l\et@tfigptorthocenterTD}\ignorespaces\fi}
\def\figptorthoprojlineDD#1:#2=#3/#4,#5/{\ifps@cri{\s@uvc@ntr@l\et@tfigptorthoprojlineDD%
    \setc@ntr@l{2}\figvectPDD-3[#4,#5]\figvectNVDD-4[-3]\resetc@ntr@l{2}%
    \inters@cDD#1:{#2}[#3,-4;#4,-3]\resetc@ntr@l\et@tfigptorthoprojlineDD}\ignorespaces\fi}
\def\figptorthoprojlineTD#1:#2=#3/#4,#5/{\ifps@cri{\s@uvc@ntr@l\et@tfigptorthoprojlineTD%
    \setc@ntr@l{2}\figvectPTD-1[#4,#3]\figvectPTD-2[#4,#5]\vecunit@TD{-2}{-2}%
    \c@lproscalTD\v@leur[-1,-2]\edef\v@lcoef{\repdecn@mb{\v@leur}}%
    \figpttraTD#1:{#2}=#4/\v@lcoef,-2/\resetc@ntr@l\et@tfigptorthoprojlineTD}\ignorespaces\fi}
\def\figptorthoprojplaneTD#1:#2=#3/#4,#5/{\ifps@cri{\s@uvc@ntr@l\et@tfigptorthoprojplane%
    \setc@ntr@l{2}\figvectPTD-1[#3,#4]\vecunit@TD{-2}{#5}%
    \c@lproscalTD\v@leur[-1,-2]\edef\v@lcoef{\repdecn@mb{\v@leur}}%
    \figpttraTD#1:{#2}=#3/\v@lcoef,-2/\resetc@ntr@l\et@tfigptorthoprojplane}\ignorespaces\fi}
\def\figpthom#1:#2=#3/#4,#5/{\ifps@cri{\s@uvc@ntr@l\et@tfigpthom%
    \setc@ntr@l{2}\figvectP-1[#4,#3]\figpttra#1:{#2}=#4/#5,-1/%
    \resetc@ntr@l\et@tfigpthom}\ignorespaces\fi}
\def\figptrotDD#1:#2=#3/#4,#5/{\ifps@cri{\s@uvc@ntr@l\et@tfigptrotDD%
    \c@ssin{\C@}{\S@}{#5}\setc@ntr@l{2}\figvectPDD-1[#4,#3]\Figg@tXY{-1}%
    \v@lXa=\C@\v@lX\advance\v@lXa-\S@\v@lY%
    \v@lYa=\S@\v@lX\advance\v@lYa\C@\v@lY%
    \Figv@ctCreg-1(\v@lXa,\v@lYa)\figpttraDD#1:{#2}=#4/1,-1/%
    \resetc@ntr@l\et@tfigptrotDD}\ignorespaces\fi}
\def\figptrotTD#1:#2=#3/#4,#5,#6/{\ifps@cri{\s@uvc@ntr@l\et@tfigptrotTD%
    \c@ssin{\C@}{\S@}{#5}%
    \setc@ntr@l{2}\figptorthoprojplaneTD-3:=#4/#3,#6/\figvectPTD-2[-3,#3]%
    \n@rmeucTD\v@leur{-2}\ifdim\v@leur<\Cepsil@n\Figg@tXYa{#3}\else%
    \edef\v@lcoef{\repdecn@mb{\v@leur}}\figvectNVTD-1[#6,-2]%
    \Figg@tXYa{-1}\v@lXa=\v@lcoef\v@lXa\v@lYa=\v@lcoef\v@lYa\v@lZa=\v@lcoef\v@lZa%
    \v@lXa=\S@\v@lXa\v@lYa=\S@\v@lYa\v@lZa=\S@\v@lZa\Figg@tXY{-2}%
    \advance\v@lXa\C@\v@lX\advance\v@lYa\C@\v@lY\advance\v@lZa\C@\v@lZ%
    \Figg@tXY{-3}\advance\v@lXa\v@lX\advance\v@lYa\v@lY\advance\v@lZa\v@lZ\fi%
    \Figp@intregTD#1:{#2}(\v@lXa,\v@lYa,\v@lZa)\resetc@ntr@l\et@tfigptrotTD}\ignorespaces\fi}
\def\figptsymDD#1:#2=#3/#4,#5/{\ifps@cri{\s@uvc@ntr@l\et@tfigptsymDD%
    \resetc@ntr@l{2}\figptorthoprojlineDD-5:=#3/#4,#5/\figvectPDD-2[#3,-5]%
    \figpttraDD#1:{#2}=#3/2,-2/\resetc@ntr@l\et@tfigptsymDD}\ignorespaces\fi}
\def\figptsymTD#1:#2=#3/#4,#5/{\ifps@cri{\s@uvc@ntr@l\et@tfigptsymTD%
    \resetc@ntr@l{2}\figptorthoprojplaneTD-3:=#3/#4,#5/\figvectPTD-2[#3,-3]%
    \figpttraTD#1:{#2}=#3/2,-2/\resetc@ntr@l\et@tfigptsymTD}\ignorespaces\fi}
\def\figpttraDD#1:#2=#3/#4,#5/{\ifps@cri{\Figg@tXYa{#5}\v@lXa=#4\v@lXa\v@lYa=#4\v@lYa%
    \Figg@tXY{#3}\advance\v@lX\v@lXa\advance\v@lY\v@lYa%
    \Figp@intregDD#1:{#2}(\v@lX,\v@lY)}\ignorespaces\fi}
\def\figpttraTD#1:#2=#3/#4,#5/{\ifps@cri{\Figg@tXYa{#5}\v@lXa=#4\v@lXa\v@lYa=#4\v@lYa%
    \v@lZa=#4\v@lZa\Figg@tXY{#3}\advance\v@lX\v@lXa\advance\v@lY\v@lYa%
    \advance\v@lZ\v@lZa\Figp@intregTD#1:{#2}(\v@lX,\v@lY,\v@lZ)}\ignorespaces\fi}
\def\figpttraCDD#1:#2=#3/#4,#5/{\ifps@cri{\v@lXa=#4\unit@\v@lYa=#5\unit@%
    \Figg@tXY{#3}\advance\v@lX\v@lXa\advance\v@lY\v@lYa%
    \Figp@intregDD#1:{#2}(\v@lX,\v@lY)}\ignorespaces\fi}
\def\figpttraCTD#1:#2=#3/#4,#5,#6/{\ifps@cri{\v@lXa=#4\unit@\v@lYa=#5\unit@\v@lZa=#6\unit@%
    \Figg@tXY{#3}\advance\v@lX\v@lXa\advance\v@lY\v@lYa\advance\v@lZ\v@lZa%
    \Figp@intregTD#1:{#2}(\v@lX,\v@lY,\v@lZ)}\ignorespaces\fi}
\def\figptsaxes#1:#2(#3){\ifps@cri{\an@lys@xes#3,:\ifx\t@xt@\empty%
    \ifTr@isDim\Figpts@xes#1:#2(0,#3,0,#3,0,#3)\else\Figpts@xes#1:#2(0,#3,0,#3)\fi%
    \else\Figpts@xes#1:#2(#3)\fi}\ignorespaces\fi}
\def\Figpts@xesDD#1:#2(#3,#4,#5,#6){%
    \s@mme=#1\figpttraC\the\s@mme:$x$=#2/#4,0/%
    \advance\s@mme\@ne\figpttraC\the\s@mme:$y$=#2/0,#6/}
\def\Figpts@xesTD#1:#2(#3,#4,#5,#6,#7,#8){%
    \s@mme=#1\figpttraC\the\s@mme:$x$=#2/#4,0,0/%
    \advance\s@mme\@ne\figpttraC\the\s@mme:$y$=#2/0,#6,0/%
    \advance\s@mme\@ne\figpttraC\the\s@mme:$z$=#2/0,0,#8/}
\def\figptsmap#1=#2/#3/#4/{\ifps@cri{\s@uvc@ntr@l\et@tfigptsmap%
    \setc@ntr@l{2}\def\list@num{#2}\s@mme=#1%
    \@ecfor\p@int:=\list@num\do{\figvectP-1[#3,\p@int]\Figg@tXY{-1}%
    \pr@dMatV/#4/\figpttra\the\s@mme:=#3/1,-1/\advance\s@mme\@ne}%
    \resetc@ntr@l\et@tfigptsmap}\ignorespaces\fi}
\def\figptscontrolDD#1[#2,#3,#4,#5]{\ifps@cri{\s@uvc@ntr@l\et@tfigptscontrolDD\setc@ntr@l{2}%
    \v@lX=\z@\v@lY=\z@\Figtr@nptDD{-5}{#2}\Figtr@nptDD{2}{#5}%
    \divide\v@lX\@vi\divide\v@lY\@vi%
    \Figtr@nptDD{3}{#3}\Figtr@nptDD{-1.5}{#4}\Figp@intregDD-1:(\v@lX,\v@lY)%
    \v@lX=\z@\v@lY=\z@\Figtr@nptDD{2}{#2}\Figtr@nptDD{-5}{#5}%
    \divide\v@lX\@vi\divide\v@lY\@vi\Figtr@nptDD{-1.5}{#3}\Figtr@nptDD{3}{#4}%
    \s@mme=#1\advance\s@mme\@ne\Figp@intregDD\the\s@mme:(\v@lX,\v@lY)%
    \figptcopyDD#1:/-1/\resetc@ntr@l\et@tfigptscontrolDD}\ignorespaces\fi}
\def\figptscontrolTD#1[#2,#3,#4,#5]{\ifps@cri{\s@uvc@ntr@l\et@tfigptscontrolTD\setc@ntr@l{2}%
    \v@lX=\z@\v@lY=\z@\v@lZ=\z@\Figtr@nptTD{-5}{#2}\Figtr@nptTD{2}{#5}%
    \divide\v@lX\@vi\divide\v@lY\@vi\divide\v@lZ\@vi%
    \Figtr@nptTD{3}{#3}\Figtr@nptTD{-1.5}{#4}\Figp@intregTD-1:(\v@lX,\v@lY,\v@lZ)%
    \v@lX=\z@\v@lY=\z@\v@lZ=\z@\Figtr@nptTD{2}{#2}\Figtr@nptTD{-5}{#5}%
    \divide\v@lX\@vi\divide\v@lY\@vi\divide\v@lZ\@vi\Figtr@nptTD{-1.5}{#3}\Figtr@nptTD{3}{#4}%
    \s@mme=#1\advance\s@mme\@ne\Figp@intregTD\the\s@mme:(\v@lX,\v@lY,\v@lZ)%
    \figptcopyTD#1:/-1/\resetc@ntr@l\et@tfigptscontrolTD}\ignorespaces\fi}
\def\Figtr@nptDD#1#2{\Figg@tXYa{#2}\v@lXa=#1\v@lXa\v@lYa=#1\v@lYa%
    \advance\v@lX\v@lXa\advance\v@lY\v@lYa}
\def\Figtr@nptTD#1#2{\Figg@tXYa{#2}\v@lXa=#1\v@lXa\v@lYa=#1\v@lYa\v@lZa=#1\v@lZa%
    \advance\v@lX\v@lXa\advance\v@lY\v@lYa\advance\v@lZ\v@lZa}
\def\figptscontrolcurve#1,#2[#3]{\ifps@cri{\s@uvc@ntr@l\et@tfigptscontrolcurve%
    \def\list@num{#3}\extrairelepremi@r\Ak@\de\list@num%
    \extrairelepremi@r\Ai@\de\list@num\extrairelepremi@r\Aj@\de\list@num%
    \s@mme=#1\figptcopy\the\s@mme:/\Ai@/%
    \setc@ntr@l{2}\figvectP -1[\Ak@,\Aj@]%
    \@ecfor\Ak@:=\list@num\do{\advance\s@mme\@ne\figpttra\the\s@mme:=\Ai@/\curv@roundness,-1/%
       \figvectP -1[\Ai@,\Ak@]\advance\s@mme\@ne\figpttra\the\s@mme:=\Aj@/-\curv@roundness,-1/%
       \advance\s@mme\@ne\figptcopy\the\s@mme:/\Aj@/%
       \edef\Ai@{\Aj@}\edef\Aj@{\Ak@}}\advance\s@mme-#1\divide\s@mme\thr@@%
       \xdef#2{\the\s@mme}%
    \resetc@ntr@l\et@tfigptscontrolcurve}\ignorespaces\fi}
\def\figptsintercircDD#1[#2,#3;#4,#5]{\ifps@cri{\s@uvc@ntr@l\et@tfigptsintercircDD%
    \setc@ntr@l{2}\let\c@lNVintc=\c@lNVintcDD\Figptsintercirc@#1[#2,#3;#4,#5]%
    \resetc@ntr@l\et@tfigptsintercircDD}\ignorespaces\fi}
\def\figptsintercircTD#1[#2,#3;#4,#5;#6]{\ifps@cri{\s@uvc@ntr@l\et@tfigptsintercircTD%
    \setc@ntr@l{2}\let\c@lNVintc=\c@lNVintcTD\vecunitC@TD[#2,#6]%
    \Figv@ctCreg-3(\v@lX,\v@lY,\v@lZ)\Figptsintercirc@#1[#2,#3;#4,#5]%
    \resetc@ntr@l\et@tfigptsintercircTD}\ignorespaces\fi}
\def\Figptsintercirc@#1[#2,#3;#4,#5]{\figvectP-1[#2,#4]%
    \vecunit@{-1}{-1}\delt@=\result@t\f@ctech=\result@tent%
    \s@mme=#1\advance\s@mme\@ne\figptcopy#1:/#2/\figptcopy\the\s@mme:/#4/%
    \ifdim\delt@=\z@\else%
    \v@lmin=#3\unit@\v@lmax=#5\unit@\v@leur=\v@lmin\advance\v@leur\v@lmax%
    \ifdim\v@leur>\delt@%
    \v@leur=\v@lmin\advance\v@leur-\v@lmax\maxim@m{\v@leur}{\v@leur}{-\v@leur}%
    \ifdim\v@leur<\delt@%
    \divide\v@lmin\f@ctech\divide\v@lmax\f@ctech\divide\delt@\f@ctech%
    \v@lmin=\repdecn@mb{\v@lmin}\v@lmin\v@lmax=\repdecn@mb{\v@lmax}\v@lmax%
    \invers@{\v@leur}{\delt@}\advance\v@lmax-\v@lmin%
    \v@lmax=-\repdecn@mb{\v@leur}\v@lmax\advance\delt@\v@lmax\delt@=.5\delt@%
    \v@lmax=\delt@\multiply\v@lmax\f@ctech%
    \edef\t@ille{\repdecn@mb{\v@lmax}}\figpttra-2:=#2/\t@ille,-1/%
    \delt@=\repdecn@mb{\delt@}\delt@\advance\v@lmin-\delt@%
    \sqrt@{\v@leur}{\v@lmin}\multiply\v@leur\f@ctech\edef\t@ille{\repdecn@mb{\v@leur}}%
    \c@lNVintc\figpttra#1:=-2/-\t@ille,-1/\figpttra\the\s@mme:=-2/\t@ille,-1/\fi\fi\fi}
\def\c@lNVintcDD{\Figg@tXY{-1}\Figv@ctCreg-1(-\v@lY,\v@lX)} 
\def\c@lNVintcTD{{\Figg@tXY{-3}\v@lmin=\v@lX\v@lmax=\v@lY\v@leur=\v@lZ%
    \Figg@tXY{-1}\c@lprovec{-3}\vecunit@{-3}{-3}
    \Figg@tXY{-1}\v@lmin=\v@lX\v@lmax=\v@lY%
    \v@leur=\v@lZ\Figg@tXY{-3}\c@lprovec{-1}}} 
\def\figptsinterlinellDD#1[#2,#3,#4,#5;#6,#7]{\ifps@cri{\s@uvc@ntr@l\et@tfigptsinterlinellDD%
    \figptcopy#1:/#6/\s@mme=#1\advance\s@mme\@ne\figptcopy\the\s@mme:/#7/%
    \v@lmin=#3\unit@\v@lmax=#4\unit@
    \setc@ntr@l{2}\figptbaryDD-4:[#6,#7;1,1]\figptsrotDD-3=-4,#7/#2,-#5/
    \Figg@tXY{-3}\Figg@tXYa{#2}\advance\v@lX-\v@lXa\advance\v@lY-\v@lYa
    \figvectP-1[-3,-2]\Figg@tXYa{-1}\figvectP-3[-4,#7]\Figptsint@rLE{#1}
    \resetc@ntr@l\et@tfigptsinterlinellDD}\ignorespaces\fi}
\def\figptsinterlinellP#1[#2,#3,#4;#5,#6]{\ifps@cri{\s@uvc@ntr@l\et@tfigptsinterlinellP%
    \figptcopy#1:/#5/\s@mme=#1\advance\s@mme\@ne\figptcopy\the\s@mme:/#6/\setc@ntr@l{2}%
    \figvectP-1[#2,#3]\vecunit@{-1}{-1}\v@lmin=\result@t
    \figvectP-2[#2,#4]\vecunit@{-2}{-2}\v@lmax=\result@t
    \figptbary-4:[#5,#6;1,1]
    \figvectP-3[#2,-4]\c@lproscal\v@lX[-3,-1]\c@lproscal\v@lY[-3,-2]
    \figvectP-3[-4,#6]\c@lproscal\v@lXa[-3,-1]\c@lproscal\v@lYa[-3,-2]
    \Figptsint@rLE{#1}\resetc@ntr@l\et@tfigptsinterlinellP}\ignorespaces\fi}
\def\Figptsint@rLE#1{%
    \getredf@ctDD\f@ctech(\v@lmin,\v@lmax)%
    \getredf@ctDD\p@rtent(\v@lX,\v@lY)\ifnum\p@rtent>\f@ctech\f@ctech=\p@rtent\fi%
    \getredf@ctDD\p@rtent(\v@lXa,\v@lYa)\ifnum\p@rtent>\f@ctech\f@ctech=\p@rtent\fi%
    \divide\v@lmin\f@ctech\divide\v@lmax\f@ctech\divide\v@lX\f@ctech\divide\v@lY\f@ctech%
    \divide\v@lXa\f@ctech\divide\v@lYa\f@ctech%
    \c@rre=\repdecn@mb\v@lXa\v@lmax\mili@u=\repdecn@mb\v@lYa\v@lmin%
    \getredf@ctDD\f@ctech(\c@rre,\mili@u)%
    \c@rre=\repdecn@mb\v@lX\v@lmax\mili@u=\repdecn@mb\v@lY\v@lmin%
    \getredf@ctDD\p@rtent(\c@rre,\mili@u)\ifnum\p@rtent>\f@ctech\f@ctech=\p@rtent\fi%
    \divide\v@lmin\f@ctech\divide\v@lmax\f@ctech\divide\v@lX\f@ctech\divide\v@lY\f@ctech%
    \divide\v@lXa\f@ctech\divide\v@lYa\f@ctech%
    \v@lmin=\repdecn@mb{\v@lmin}\v@lmin\v@lmax=\repdecn@mb{\v@lmax}\v@lmax%
    \edef\G@xde{\repdecn@mb\v@lmin}\edef\P@xde{\repdecn@mb\v@lmax}%
    \c@rre=-\v@lmax\v@leur=\repdecn@mb\v@lY\v@lY\advance\c@rre\v@leur\c@rre=\G@xde\c@rre%
    \v@leur=\repdecn@mb\v@lX\v@lX\v@leur=\P@xde\v@leur\advance\c@rre\v@leur
    \v@lmin=\repdecn@mb\v@lYa\v@lmin\v@lmax=\repdecn@mb\v@lXa\v@lmax%
    \mili@u=\repdecn@mb\v@lX\v@lmax\advance\mili@u\repdecn@mb\v@lY\v@lmin
    \v@lmax=\repdecn@mb\v@lXa\v@lmax\advance\v@lmax\repdecn@mb\v@lYa\v@lmin
    \ifdim\v@lmax>\epsil@n%
    \maxim@m{\v@leur}{\c@rre}{-\c@rre}\maxim@m{\v@lmin}{\mili@u}{-\mili@u}%
    \maxim@m{\v@leur}{\v@leur}{\v@lmin}\maxim@m{\v@lmin}{\v@lmax}{-\v@lmax}%
    \maxim@m{\v@leur}{\v@leur}{\v@lmin}\p@rtentiere{\p@rtent}{\v@leur}\advance\p@rtent\@ne%
    \divide\c@rre\p@rtent\divide\mili@u\p@rtent\divide\v@lmax\p@rtent%
    \delt@=\repdecn@mb{\mili@u}\mili@u\v@leur=\repdecn@mb{\v@lmax}\c@rre%
    \advance\delt@-\v@leur\ifdim\delt@<\z@\else\sqrt@\delt@\delt@%
    \invers@\v@lmax\v@lmax\edef\Uns@rAp{\repdecn@mb\v@lmax}%
    \v@leur=-\mili@u\advance\v@leur-\delt@\v@leur=\Uns@rAp\v@leur%
    \edef\t@ille{\repdecn@mb\v@leur}\figpttra#1:=-4/\t@ille,-3/\s@mme=#1\advance\s@mme\@ne%
    \v@leur=-\mili@u\advance\v@leur\delt@\v@leur=\Uns@rAp\v@leur%
    \edef\t@ille{\repdecn@mb\v@leur}\figpttra\the\s@mme:=-4/\t@ille,-3/\fi\fi}
\def\figptsorthoprojlineDD#1=#2/#3,#4/{\ifps@cri{\s@uvc@ntr@l\et@tfigptsorthoprojlineDD%
    \setc@ntr@l{2}\figvectPDD-3[#3,#4]\figvectNVDD-4[-3]\resetc@ntr@l{2}%
    \def\list@num{#2}\s@mme=#1\@ecfor\p@int:=\list@num\do{%
    \inters@cDD\the\s@mme:[\p@int,-4;#3,-3]\advance\s@mme\@ne}%
    \resetc@ntr@l\et@tfigptsorthoprojlineDD}\ignorespaces\fi}
\def\figptsorthoprojlineTD#1=#2/#3,#4/{\ifps@cri{\s@uvc@ntr@l\et@tfigptsorthoprojlineTD%
    \setc@ntr@l{2}\figvectPTD-2[#3,#4]\vecunit@TD{-2}{-2}%
    \def\list@num{#2}\s@mme=#1\@ecfor\p@int:=\list@num\do{%
    \figvectPTD-1[#3,\p@int]\c@lproscalTD\v@leur[-1,-2]%
    \edef\v@lcoef{\repdecn@mb{\v@leur}}\figpttraTD\the\s@mme:=#3/\v@lcoef,-2/%
    \advance\s@mme\@ne}\resetc@ntr@l\et@tfigptsorthoprojlineTD}\ignorespaces\fi}
\def\figptsorthoprojplaneTD#1=#2/#3,#4/{\ifps@cri{\s@uvc@ntr@l\et@tfigptsorthoprojplane%
    \setc@ntr@l{2}\vecunit@TD{-2}{#4}%
    \def\list@num{#2}\s@mme=#1\@ecfor\p@int:=\list@num\do{\figvectPTD-1[\p@int,#3]%
    \c@lproscalTD\v@leur[-1,-2]\edef\v@lcoef{\repdecn@mb{\v@leur}}%
    \figpttraTD\the\s@mme:=\p@int/\v@lcoef,-2/\advance\s@mme\@ne}%
    \resetc@ntr@l\et@tfigptsorthoprojplane}\ignorespaces\fi}
\def\figptshom#1=#2/#3,#4/{\ifps@cri{\s@uvc@ntr@l\et@tfigptshom%
    \setc@ntr@l{2}\def\list@num{#2}\s@mme=#1%
    \@ecfor\p@int:=\list@num\do{\figvectP-1[#3,\p@int]%
    \figpttra\the\s@mme:=#3/#4,-1/\advance\s@mme\@ne}%
    \resetc@ntr@l\et@tfigptshom}\ignorespaces\fi}
\def\figptsrotDD#1=#2/#3,#4/{\ifps@cri{\s@uvc@ntr@l\et@tfigptsrotDD%
    \c@ssin{\C@}{\S@}{#4}\setc@ntr@l{2}\def\list@num{#2}\s@mme=#1%
    \@ecfor\p@int:=\list@num\do{\figvectPDD-1[#3,\p@int]\Figg@tXY{-1}%
    \v@lXa=\C@\v@lX\advance\v@lXa-\S@\v@lY%
    \v@lYa=\S@\v@lX\advance\v@lYa\C@\v@lY%
    \Figv@ctCreg-1(\v@lXa,\v@lYa)\figpttraDD\the\s@mme:=#3/1,-1/\advance\s@mme\@ne}%
    \resetc@ntr@l\et@tfigptsrotDD}\ignorespaces\fi}
\def\figptsrotTD#1=#2/#3,#4,#5/{\ifps@cri{\s@uvc@ntr@l\et@tfigptsrotTD%
    \c@ssin{\C@}{\S@}{#4}%
    \setc@ntr@l{2}\def\list@num{#2}\s@mme=#1%
    \@ecfor\p@int:=\list@num\do{\figptorthoprojplaneTD-3:=#3/\p@int,#5/%
    \figvectPTD-2[-3,\p@int]%
    \figvectNVTD-1[#5,-2]\n@rmeucTD\v@leur{-2}\edef\v@lcoef{\repdecn@mb{\v@leur}}%
    \Figg@tXYa{-1}\v@lXa=\v@lcoef\v@lXa\v@lYa=\v@lcoef\v@lYa\v@lZa=\v@lcoef\v@lZa%
    \v@lXa=\S@\v@lXa\v@lYa=\S@\v@lYa\v@lZa=\S@\v@lZa\Figg@tXY{-2}%
    \advance\v@lXa\C@\v@lX\advance\v@lYa\C@\v@lY\advance\v@lZa\C@\v@lZ%
    \Figg@tXY{-3}\advance\v@lXa\v@lX\advance\v@lYa\v@lY\advance\v@lZa\v@lZ%
    \Figp@intregTD\the\s@mme:(\v@lXa,\v@lYa,\v@lZa)\advance\s@mme\@ne}%
    \resetc@ntr@l\et@tfigptsrotTD}\ignorespaces\fi}
\def\figptssymDD#1=#2/#3,#4/{\ifps@cri{\s@uvc@ntr@l\et@tfigptssymDD%
    \setc@ntr@l{2}\figvectPDD-3[#3,#4]\Figg@tXY{-3}\Figv@ctCreg-4(-\v@lY,\v@lX)%
    \resetc@ntr@l{2}\def\list@num{#2}\s@mme=#1%
    \@ecfor\p@int:=\list@num\do{\inters@cDD-5:[#3,-3;\p@int,-4]\figvectPDD-2[\p@int,-5]%
    \figpttraDD\the\s@mme:=\p@int/2,-2/\advance\s@mme\@ne}%
    \resetc@ntr@l\et@tfigptssymDD}\ignorespaces\fi}
\def\figptssymTD#1=#2/#3,#4/{\ifps@cri{\s@uvc@ntr@l\et@tfigptssymTD%
    \setc@ntr@l{2}\vecunit@TD{-2}{#4}\def\list@num{#2}\s@mme=#1%
    \@ecfor\p@int:=\list@num\do{\figvectPTD-1[\p@int,#3]%
    \c@lproscalTD\v@leur[-1,-2]\v@leur=2\v@leur\edef\v@lcoef{\repdecn@mb{\v@leur}}%
    \figpttraTD\the\s@mme:=\p@int/\v@lcoef,-2/\advance\s@mme\@ne}%
    \resetc@ntr@l\et@tfigptssymTD}\ignorespaces\fi}
\def\figptstraDD#1=#2/#3,#4/{\ifps@cri{\Figg@tXYa{#4}\v@lXa=#3\v@lXa\v@lYa=#3\v@lYa%
    \def\list@num{#2}\s@mme=#1\@ecfor\p@int:=\list@num\do{\Figg@tXY{\p@int}%
    \advance\v@lX\v@lXa\advance\v@lY\v@lYa%
    \Figp@intregDD\the\s@mme:(\v@lX,\v@lY)\advance\s@mme\@ne}}\ignorespaces\fi}
\def\figptstraTD#1=#2/#3,#4/{\ifps@cri{\Figg@tXYa{#4}\v@lXa=#3\v@lXa\v@lYa=#3\v@lYa%
    \v@lZa=#3\v@lZa\def\list@num{#2}\s@mme=#1\@ecfor\p@int:=\list@num\do{\Figg@tXY{\p@int}%
    \advance\v@lX\v@lXa\advance\v@lY\v@lYa\advance\v@lZ\v@lZa%
    \Figp@intregTD\the\s@mme:(\v@lX,\v@lY,\v@lZ)\advance\s@mme\@ne}}\ignorespaces\fi}
\def\figptvisilimSLTD#1:#2[#3,#4;#5,#6]{\ifps@cri{\s@uvc@ntr@l\et@tfigptvisilimSLTD%
    \setc@ntr@l{2}\figvectP-1[#3,#4]\n@rminf{\delt@}{-1}%
    \ifcase\curr@ntproj\v@lX=\cxa@\p@\v@lY=-\p@\v@lZ=\cxb@\p@
    \Figv@ctCreg-2(\v@lX,\v@lY,\v@lZ)\figvectP-3[#5,#6]\figvectNV-1[-2,-3]%
    \or\figvectP-1[#5,#6]\vecunitCV@TD{-1}\v@lmin=\v@lX\v@lmax=\v@lY
    \v@leur=\v@lZ\v@lX=\cza@\p@\v@lY=\czb@\p@\v@lZ=\czc@\p@\c@lprovec{-1}%
    \or\c@ley@pt{-2}\figvectN-1[#5,#6,-2]\fi
    \edef\Ai@{#3}\edef\Aj@{#4}\figvectP-2[#5,\Ai@]\c@lproscal\v@leur[-1,-2]%
    \ifdim\v@leur>\z@\p@rtent=\@ne\else\p@rtent=\m@ne\fi%
    \figvectP-2[#5,\Aj@]\c@lproscal\v@leur[-1,-2]%
    \ifdim\p@rtent\v@leur>\z@\figptcopy#1:#2/#3/%
    \message{*** \BS@ figptvisilimSL: points are on the same side.}\else%
    \figptcopy-3:/#3/\figptcopy-4:/#4/%
    \loop\figptbary-5:[-3,-4;1,1]\figvectP-2[#5,-5]\c@lproscal\v@leur[-1,-2]%
    \ifdim\p@rtent\v@leur>\z@\figptcopy-3:/-5/\else\figptcopy-4:/-5/\fi%
    \divide\delt@\tw@\ifdim\delt@>\epsil@n\repeat%
    \figptbary#1:#2[-3,-4;1,1]\fi\resetc@ntr@l\et@tfigptvisilimSLTD}\ignorespaces\fi}
\def\c@ley@pt#1{\t@stp@r\ifitis@K\v@lX=\cza@\p@\v@lY=\czb@\p@\v@lZ=\czc@\p@%
    \Figv@ctCreg-1(\v@lX,\v@lY,\v@lZ)\Figp@intreg-2:(\wd\Bt@rget,\ht\Bt@rget,\dp\Bt@rget)%
    \figpttra#1:=-2/-\disob@intern,-1/\else\end\fi}
\def\t@stp@r{\itis@Ktrue\ifnewt@rgetpt\else\itis@Kfalse%
    \message{*** \BS@ figptvisilimXX: target point undefined.}\fi\ifnewdis@b\else%
    \itis@Kfalse\message{*** \BS@ figptvisilimXX: observation distance undefined.}\fi%
    \ifitis@K\else\message{*** This macro must be called after \BS@ psbeginfig or after
    having set the missing parameter(s) with \BS@ figset proj()}\fi}
\def\figscan#1(#2,#3){{\s@uvc@ntr@l\et@tfigscan\@psfgetbb{#1}\if@psfbbfound\else%
    \def\@psfllx{0}\def\@psflly{20}\def\@psfurx{540}\def\@psfury{640}\fi%
    \unit@=\@ne bp\setc@ntr@l{2}\figsetmark{}%
    \def\minst@p{20pt}%
    \v@lX=\@psfllx\p@\v@lX=\Sc@leFact\v@lX\r@undint\v@lX\v@lX%
    \v@lY=\@psflly\p@\v@lY=\Sc@leFact\v@lY\ifdim\v@lY>\z@\r@undint\v@lY\v@lY\fi%
    \delt@=\@psfury\p@\delt@=\Sc@leFact\delt@%
    \advance\delt@-\v@lY\v@lXa=\@psfurx\p@\v@lXa=\Sc@leFact\v@lXa\v@leur=\minst@p%
    \edef\valv@lY{\repdecn@mb{\v@lY}}\edef\LgTr@it{\the\delt@}%
    \loop\ifdim\v@lX<\v@lXa\edef\valv@lX{\repdecn@mb{\v@lX}}%
    \figptDD -1:(\valv@lX,\valv@lY)\figwriten -1:\hbox{\vrule height\LgTr@it}(0)%
    \ifdim\v@leur<\minst@p\else\figsetmark{\raise-8bp\hbox{$\scriptscriptstyle\triangle$}}%
    \figwrites -1:\@ffichnb{0}{\valv@lX}(6)\v@leur=\z@\figsetmark{}\fi%
    \advance\v@leur#2pt\advance\v@lX#2pt\repeat%
    \def\minst@p{10pt}%
    \v@lX=\@psfllx\p@\v@lX=\Sc@leFact\v@lX\ifdim\v@lX>\z@\r@undint\v@lX\v@lX\fi%
    \v@lY=\@psflly\p@\v@lY=\Sc@leFact\v@lY\r@undint\v@lY\v@lY%
    \delt@=\@psfurx\p@\delt@=\Sc@leFact\delt@%
    \advance\delt@-\v@lX\v@lYa=\@psfury\p@\v@lYa=\Sc@leFact\v@lYa\v@leur=\minst@p%
    \edef\valv@lX{\repdecn@mb{\v@lX}}\edef\LgTr@it{\the\delt@}%
    \loop\ifdim\v@lY<\v@lYa\edef\valv@lY{\repdecn@mb{\v@lY}}%
    \figptDD -1:(\valv@lX,\valv@lY)\figwritee -1:\vbox{\hrule width\LgTr@it}(0)%
    \ifdim\v@leur<\minst@p\else\figsetmark{$\triangleright$\kern4bp}%
    \figwritew -1:\@ffichnb{0}{\valv@lY}(6)\v@leur=\z@\figsetmark{}\fi%
    \advance\v@leur#3pt\advance\v@lY#3pt\repeat%
    \resetc@ntr@l\et@tfigscan}\ignorespaces}
\def\figshowpts[#1,#2]{{\figsetmark{$\bullet$}\figsetptname{\bf ##1}%
    \p@rtent=#2\relax\ifnum\p@rtent<\z@\p@rtent=\z@\fi%
    \s@mme=#1\relax\ifnum\s@mme<\z@\s@mme=\z@\fi%
    \loop\ifnum\s@mme<\p@rtent\pt@rvect{\s@mme}%
    \ifitis@K\figwriten{\the\s@mme}:(4pt)\fi\advance\s@mme\@ne\repeat%
    \pt@rvect{\s@mme}\ifitis@K\figwriten{\the\s@mme}:(4pt)\fi}\ignorespaces}
\def\pt@rvect#1{\set@bjc@de{#1}%
    \expandafter\expandafter\expandafter\inqpt@rvec\csname\objc@de\endcsname:}
\def\inqpt@rvec#1#2:{\if#1\C@dCl@spt\itis@Ktrue\else\itis@Kfalse\fi}
\def\figshowsettings{{%
    \immediate\write16{====================================================================}%
    \immediate\write16{ Current settings about:}%
    \immediate\write16{ --- GENERAL ---}%
    \immediate\write16{Scale factor and Unit = \unit@util\space (\the\unit@)
     \space -> \BS@ figinit{ScaleFactorUnit}}%
    \immediate\write16{Update mode = \ifpstestm@de yes\else no\fi
     \space-> \BS@ psset(update=yes/no) or \BS@ pssetdefault(update=yes/no)}%
    \immediate\write16{ --- PRINTING ---}%
    \immediate\write16{Implicit point name = \ptn@me{i} \space-> \BS@ figsetptname{Name}}%
    \immediate\write16{Point marker = \the\c@nsymb \space -> \BS@ figsetmark{Mark}}%
    \immediate\write16{Print rounded coordinates = \ifr@undcoord yes\else no\fi
     \space-> \BS@ figsetroundcoord{yes/no}}%
    \immediate\write16{ --- GRAPHICAL (general) ---}%
    \immediate\write16{First-level (or primary) settings:}%
    \immediate\write16{ Color = \curr@ntcolor \space-> \BS@ psset(color=ColorDefinition)}%
    \immediate\write16{ Filling mode = \iffillm@de yes\else no\fi
     \space-> \BS@ psset(fillmode=yes/no)}%
    \immediate\write16{ Line join = \curr@ntjoin \space-> \BS@ psset(join=miter/round/bevel)}%
    \immediate\write16{ Line style = \curr@ntdash \space-> \BS@ psset(dash=Index/Pattern)}%
    \immediate\write16{ Line width = \curr@ntwidth
     \space-> \BS@ psset(width=real in PostScript units)}%
    \immediate\write16{Second-level (or secondary) settings:}%
    \immediate\write16{ Color = \sec@ndcolor \space-> \BS@ psset second(color=ColorDefinition)}%
    \immediate\write16{ Line style = \curr@ntseconddash
     \space-> \BS@ psset second(dash=Index/Pattern)}%
    \immediate\write16{ Line width = \curr@ntsecondwidth
     \space-> \BS@ psset second(width=real in PostScript units)}%
    \immediate\write16{Third-level (or ternary) settings:}%
    \immediate\write16{ Color = \th@rdcolor \space-> \BS@ psset third(color=ColorDefinition)}%
    \immediate\write16{ --- GRAPHICAL (specific) ---}%
    \immediate\write16{Arrow-head:}%
    \immediate\write16{ (half-)Angle = \@rrowheadangle
     \space-> \BS@ psset arrowhead(angle=real in degrees)}%
    \immediate\write16{ Filling mode = \if@rrowhfill yes\else no\fi
     \space-> \BS@ psset arrowhead(fillmode=yes/no)}%
    \immediate\write16{ "Outside" = \if@rrowhout yes\else no\fi
     \space-> \BS@ psset arrowhead(out=yes/no)}%
    \immediate\write16{ Length = \@rrowheadlength
     \if@rrowratio\space(not active)\else\space(active)\fi
     \space-> \BS@ psset arrowhead(length=real in user coord.)}%
    \immediate\write16{ Ratio = \@rrowheadratio
     \if@rrowratio\space(active)\else\space(not active)\fi
     \space-> \BS@ psset arrowhead(ratio=real in [0,1])}%
    \immediate\write16{Curve: Roundness = \curv@roundness
     \space-> \BS@ psset curve(roundness=real in [0,0.5])}%
    \immediate\write16{Mesh: Diagonal = \the\c@ntrolmesh
     \space-> \BS@ psset mesh(diag=integer in {-1,0,1})}%
    \immediate\write16{Flow chart:}%
    \immediate\write16{ Arrow position = \@rrowp@s
     \space-> \BS@ psset flowchart(arrowposition=real in [0,1])}%
    \immediate\write16{ Arrow reference point = \ifcase\@rrowr@fpt start\else end\fi
     \space-> \BS@ psset flowchart(arrowrefpt = start/end)}%
    \immediate\write16{ Line type = \ifcase\fclin@typ@ curve\else polygon\fi
     \space-> \BS@ psset flowchart(line=polygon/curve)}%
    \immediate\write16{ Padding = (\Xp@dd, \Yp@dd)
     \space-> \BS@ psset flowchart(padding = real in user coord.)}%
    \immediate\write16{\space\space\space\space(or
     \BS@ psset flowchart(xpadding=real, ypadding=real) )}%
    \immediate\write16{ Radius = \fclin@r@d
     \space-> \BS@ psset flowchart(radius=positive real in user coord.)}%
    \immediate\write16{ Shape = \fcsh@pe
     \space-> \BS@ psset flowchart(shape = rectangle, ellipse or lozenge)}%
    \immediate\write16{ Thickness = \thickn@ss
     \space-> \BS@ psset flowchart(thickness = real in user coord.)}%
    \ifTr@isDim%
    \immediate\write16{ --- 3D to 2D PROJECTION ---}%
    \immediate\write16{Projection : \typ@proj \space-> \BS@ figinit{ScaleFactorUnit, ProjType}}%
    \immediate\write16{Longitude (psi) = \v@lPsi \space-> \BS@ figset proj(psi=real in degrees)}%
    \ifcase\curr@ntproj\immediate\write16{Depth coeff. (Lambda)
     \space = \v@lTheta \space-> \BS@ figset proj(lambda=real in [0,1])}%
    \else\immediate\write16{Latitude (theta)
     \space = \v@lTheta \space-> \BS@ figset proj(theta=real in degrees)}%
    \fi%
    \ifnum\curr@ntproj=\tw@%
    \immediate\write16{Observation distance = \disob@unit
     \space-> \BS@ figset proj(dist=real in user coord.)}%
    \immediate\write16{Target point = \t@rgetpt \space-> \BS@ figset proj(targetpt=pt number)}%
     \v@lX=\ptT@unit@\wd\Bt@rget\v@lY=\ptT@unit@\ht\Bt@rget\v@lZ=\ptT@unit@\dp\Bt@rget%
    \immediate\write16{ Its coordinates are
     (\repdecn@mb{\v@lX}, \repdecn@mb{\v@lY}, \repdecn@mb{\v@lZ})}%
    \fi%
    \fi%
    \immediate\write16{====================================================================}%
    \ignorespaces}}
{\catcode`\/=0 \catcode`/\=12 /gdef/BS@{\}}
\newif\ifitis@vect@r
\def\figvectC#1(#2,#3){{\itis@vect@rtrue\figpt#1:(#2,#3)}\ignorespaces}
\def\Figv@ctCreg#1(#2,#3){{\itis@vect@rtrue\Figp@intreg#1:(#2,#3)}\ignorespaces}
\def\figvectDBezierDD#1:#2,#3[#4,#5,#6,#7]{\ifps@cri{\s@uvc@ntr@l\et@tfigvectDBezierDD%
    \FigvectDBezier@#2,#3[#4,#5,#6,#7]\v@lX=\c@ef\v@lX\v@lY=\c@ef\v@lY%
    \Figv@ctCreg#1(\v@lX,\v@lY)\resetc@ntr@l\et@tfigvectDBezierDD}\ignorespaces\fi}
\def\figvectDBezierTD#1:#2,#3[#4,#5,#6,#7]{\ifps@cri{\s@uvc@ntr@l\et@tfigvectDBezierTD%
    \FigvectDBezier@#2,#3[#4,#5,#6,#7]\v@lX=\c@ef\v@lX\v@lY=\c@ef\v@lY\v@lZ=\c@ef\v@lZ%
    \Figv@ctCreg#1(\v@lX,\v@lY,\v@lZ)\resetc@ntr@l\et@tfigvectDBezierTD}\ignorespaces\fi}
\def\FigvectDBezier@#1,#2[#3,#4,#5,#6]{\setc@ntr@l{2}%
    \edef\T@{#2}\v@leur=\p@\advance\v@leur-#2pt\edef\UNmT@{\repdecn@mb{\v@leur}}%
    \ifnum#1=\tw@\def\c@ef{6}\else\def\c@ef{3}\fi%
    \figptcopy-4:/#3/\figptcopy-3:/#4/\figptcopy-2:/#5/\figptcopy-1:/#6/%
    \l@mbd@un=-4 \l@mbd@de=-\thr@@\p@rtent=\m@ne\c@lDecast%
    \ifnum#1=\tw@\c@lDCDeux{-4}{-3}\c@lDCDeux{-3}{-2}\c@lDCDeux{-4}{-3}\else%
    \l@mbd@un=-4 \l@mbd@de=-\thr@@\p@rtent=-\tw@\c@lDecast%
    \c@lDCDeux{-4}{-3}\fi\Figg@tXY{-4}}
\def\c@lDCDeuxDD#1#2{\Figg@tXY{#2}\Figg@tXYa{#1}%
    \advance\v@lX-\v@lXa\advance\v@lY-\v@lYa\Figp@intregDD#1:(\v@lX,\v@lY)}
\def\c@lDCDeuxTD#1#2{\Figg@tXY{#2}\Figg@tXYa{#1}\advance\v@lX-\v@lXa%
    \advance\v@lY-\v@lYa\advance\v@lZ-\v@lZa\Figp@intregTD#1:(\v@lX,\v@lY,\v@lZ)}
\def\figvectNDD#1[#2,#3]{\ifps@cri{\Figg@tXYa{#2}\Figg@tXY{#3}%
    \advance\v@lX-\v@lXa\advance\v@lY-\v@lYa%
    \Figv@ctCreg#1(-\v@lY,\v@lX)}\ignorespaces\fi}
\def\figvectNTD#1[#2,#3,#4]{\ifps@cri{\vecunitC@TD[#2,#4]\v@lmin=\v@lX\v@lmax=\v@lY%
    \v@leur=\v@lZ\vecunitC@TD[#2,#3]\c@lprovec{#1}}\ignorespaces\fi}
\def\figvectNVDD#1[#2]{\ifps@cri{\Figg@tXY{#2}\Figv@ctCreg#1(-\v@lY,\v@lX)}\ignorespaces\fi}
\def\figvectNVTD#1[#2,#3]{\ifps@cri{\vecunitCV@TD{#3}\v@lmin=\v@lX\v@lmax=\v@lY%
    \v@leur=\v@lZ\vecunitCV@TD{#2}\c@lprovec{#1}}\ignorespaces\fi}
\def\figvectPDD#1[#2,#3]{\ifps@cri{\Figg@tXYa{#2}\Figg@tXY{#3}%
    \advance\v@lX-\v@lXa\advance\v@lY-\v@lYa%
    \Figv@ctCreg#1(\v@lX,\v@lY)}\ignorespaces\fi}
\def\figvectPTD#1[#2,#3]{\ifps@cri{\Figg@tXYa{#2}\Figg@tXY{#3}%
    \advance\v@lX-\v@lXa\advance\v@lY-\v@lYa\advance\v@lZ-\v@lZa%
    \Figv@ctCreg#1(\v@lX,\v@lY,\v@lZ)}\ignorespaces\fi}
\def\figvectUDD#1[#2]{\ifps@cri{\n@rmeuc\v@leur{#2}\invers@\v@leur\v@leur%
    \delt@=\repdecn@mb{\v@leur}\unit@\edef\v@ldelt@{\repdecn@mb{\delt@}}%
    \Figg@tXY{#2}\v@lX=\v@ldelt@\v@lX\v@lY=\v@ldelt@\v@lY%
    \Figv@ctCreg#1(\v@lX,\v@lY)}\ignorespaces\fi}
\def\figvectUTD#1[#2]{\ifps@cri{\n@rmeuc\v@leur{#2}\invers@\v@leur\v@leur%
    \delt@=\repdecn@mb{\v@leur}\unit@\edef\v@ldelt@{\repdecn@mb{\delt@}}%
    \Figg@tXY{#2}\v@lX=\v@ldelt@\v@lX\v@lY=\v@ldelt@\v@lY\v@lZ=\v@ldelt@\v@lZ%
    \Figv@ctCreg#1(\v@lX,\v@lY,\v@lZ)}\ignorespaces\fi}
\def\figvisu#1#2#3{\c@ldefproj\initb@undb@x\setbox\b@xvisu=\hbox{\ignorespaces#3}%
    \v@lXa=-\c@@rdYmin\v@lYa=\c@@rdYmax\advance\v@lYa-\c@@rdYmin%
    \v@lX=\c@@rdXmax\advance\v@lX-\c@@rdXmin%
    \setbox#1=\hbox{#2}\v@lY=-\v@lX\maxim@m{\v@lX}{\v@lX}{\wd#1}%
    \advance\v@lY\v@lX\divide\v@lY\tw@\advance\v@lY-\c@@rdXmin%
    \setbox#1=\vbox{\parindent0mm\hsize=\v@lX\vskip\v@lYa%
    \rlap{\hskip\v@lY\smash{\raise\v@lXa\box\b@xvisu}}%
    \def\t@xt@{#2}\ifx\t@xt@\empty\else\medskip\centerline{#2}\fi}\wd#1=\v@lX}
\newbox\Bt@rget\setbox\Bt@rget=\null
\newbox\BminTD@\setbox\BminTD@=\null
\newbox\BmaxTD@\setbox\BmaxTD@=\null
\newif\ifnewt@rgetpt\newif\ifnewdis@b
\def\b@undb@xTD#1#2#3{%
    \relax\ifdim#1<\wd\BminTD@\global\wd\BminTD@=#1\fi%
    \relax\ifdim#2<\ht\BminTD@\global\ht\BminTD@=#2\fi%
    \relax\ifdim#3<\dp\BminTD@\global\dp\BminTD@=#3\fi%
    \relax\ifdim#1>\wd\BmaxTD@\global\wd\BmaxTD@=#1\fi%
    \relax\ifdim#2>\ht\BmaxTD@\global\ht\BmaxTD@=#2\fi%
    \relax\ifdim#3>\dp\BmaxTD@\global\dp\BmaxTD@=#3\fi}
\def\c@ldefdisob{{\ifdim\wd\BminTD@<\maxdimen\v@leur=\wd\BmaxTD@\advance\v@leur-\wd\BminTD@%
    \delt@=\ht\BmaxTD@\advance\delt@-\ht\BminTD@\maxim@m{\v@leur}{\v@leur}{\delt@}%
    \delt@=\dp\BmaxTD@\advance\delt@-\dp\BminTD@\maxim@m{\v@leur}{\v@leur}{\delt@}%
    \v@leur=5\v@leur\else\v@leur=800pt\fi\c@ldefdisob@{\v@leur}}}
\def\c@ldefdisob@#1{{\v@leur=#1\ifdim\v@leur<\p@\v@leur=800pt\fi%
    \xdef\disob@intern{\repdecn@mb{\v@leur}}%
    \delt@=\ptT@unit@\v@leur\xdef\disob@unit{\repdecn@mb{\delt@}}%
    \f@ctech=\@ne\loop\ifdim\v@leur>\t@n pt\divide\v@leur\t@n\multiply\f@ctech\t@n\repeat%
    \xdef\disob@{\repdecn@mb{\v@leur}}\xdef\divf@ctproj{\the\f@ctech}}%
    \global\newdis@btrue}
\def\c@ldeft@rgetpt{\newt@rgetpttrue\def\t@rgetpt{CenterBoundBox}{%
    \delt@=\wd\BmaxTD@\advance\delt@-\wd\BminTD@\divide\delt@\tw@%
    \v@leur=\wd\BminTD@\advance\v@leur\delt@\global\wd\Bt@rget=\v@leur%
    \delt@=\ht\BmaxTD@\advance\delt@-\ht\BminTD@\divide\delt@\tw@%
    \v@leur=\ht\BminTD@\advance\v@leur\delt@\global\ht\Bt@rget=\v@leur%
    \delt@=\dp\BmaxTD@\advance\delt@-\dp\BminTD@\divide\delt@\tw@%
    \v@leur=\dp\BminTD@\advance\v@leur\delt@\global\dp\Bt@rget=\v@leur}}
\def\c@ldefprojTD{\ifnewt@rgetpt\else\c@ldeft@rgetpt\fi\ifnewdis@b\else\c@ldefdisob\fi}
\def\c@lprojcav{
    \v@lZa=\cxa@\v@lY\advance\v@lX\v@lZa%
    \v@lZa=\cxb@\v@lY\v@lY=\v@lZ\advance\v@lY\v@lZa\ignorespaces}
\def\c@lprojrea{
    \advance\v@lX-\wd\Bt@rget\advance\v@lY-\ht\Bt@rget\advance\v@lZ-\dp\Bt@rget%
    \v@lZa=\cza@\v@lX\advance\v@lZa\czb@\v@lY\advance\v@lZa\czc@\v@lZ%
    \divide\v@lZa\divf@ctproj\advance\v@lZa\disob@ pt\invers@{\v@lZa}{\v@lZa}%
    \v@lZa=\disob@\v@lZa\edef\v@lcoef{\repdecn@mb{\v@lZa}}%
    \v@lXa=\cxa@\v@lX\advance\v@lXa\cxb@\v@lY\v@lXa=\v@lcoef\v@lXa%
    \v@lY=\cyb@\v@lY\advance\v@lY\cya@\v@lX\advance\v@lY\cyc@\v@lZ%
    \v@lY=\v@lcoef\v@lY\v@lX=\v@lXa\ignorespaces}
\def\c@lprojort{
    \v@lXa=\cxa@\v@lX\advance\v@lXa\cxb@\v@lY%
    \v@lY=\cyb@\v@lY\advance\v@lY\cya@\v@lX\advance\v@lY\cyc@\v@lZ%
    \v@lX=\v@lXa\ignorespaces}
\def\Figptpr@j#1:#2/#3/{{\Figg@tXY{#3}\superc@lprojSP%
    \Figp@intregDD#1:{#2}(\v@lX,\v@lY)}\ignorespaces}
\def\figsetobdistTD(#1){{\ifcurr@ntPS%
    \immediate\write16{*** \BS@ figsetobdist is ignored inside a
     \BS@ psbeginfig-\BS@ psendfig block.}%
    \else\v@leur=#1\unit@\c@ldefdisob@{\v@leur}\fi}\ignorespaces}
\def\Figs@tproj#1{%
    \if#13 \d@faultproj\else\if#1c\d@faultproj%
    \else\if#1o\xdef\curr@ntproj{1}\xdef\typ@proj{orthogonal}%
         \figsetviewTD(\def@ultpsi,\def@ulttheta)%
         \global\let\c@lprojSP=\c@lprojort\global\let\superc@lprojSP=\c@lprojort%
    \else\if#1r\xdef\curr@ntproj{2}\xdef\typ@proj{realistic}%
         \figsetviewTD(\def@ultpsi,\def@ulttheta)%
         \global\let\c@lprojSP=\c@lprojrea\global\let\superc@lprojSP=\c@lprojrea%
    \else\d@faultproj\message{*** Unknown projection. Cavalier projection assumed.}%
    \fi\fi\fi\fi}
\def\d@faultproj{\xdef\curr@ntproj{0}\xdef\typ@proj{cavalier}\figsetviewTD(\def@ultpsi,0.5)%
         \global\let\c@lprojSP=\c@lprojcav\global\let\superc@lprojSP=\c@lprojcav}
\def\figsettargetTD[#1]{{\ifcurr@ntPS%
    \immediate\write16{*** \BS@ figsettarget is ignored inside a
     \BS@ psbeginfig-\BS@ psendfig block.}%
    \else\global\newt@rgetpttrue\xdef\t@rgetpt{#1}\Figg@tXY{#1}\global\wd\Bt@rget=\v@lX%
    \global\ht\Bt@rget=\v@lY\global\dp\Bt@rget=\v@lZ\fi}\ignorespaces}
\def\figsetviewTD(#1){\ifcurr@ntPS%
     \immediate\write16{*** \BS@ figsetview is ignored inside a
     \BS@ psbeginfig-\BS@ psendfig block.}\else\Figsetview@#1,:\fi\ignorespaces}
\def\Figsetview@#1,#2:{{\xdef\v@lPsi{#1}\def\t@xt@{#2}%
    \ifx\t@xt@\empty\def\@rgdeux{\v@lTheta}\else%
    \def\Xarg@##1,{\edef\@rgdeux{##1}}\Xarg@#2\fi%
    \c@ssin{\costhet@}{\sinthet@}{#1}\v@lmin=\costhet@ pt\v@lmax=\sinthet@ pt%
    \ifcase\curr@ntproj%
    \v@leur=\@rgdeux\v@lmin\xdef\cxa@{\repdecn@mb{\v@leur}}%
    \v@leur=\@rgdeux\v@lmax\xdef\cxb@{\repdecn@mb{\v@leur}}\v@leur=\@rgdeux pt%
    \relax\ifdim\v@leur>\p@\message{*** Lambda too large ! See \BS@ figset proj() !}\fi%
    \else%
    \v@lmax=-\v@lmax\xdef\cxa@{\repdecn@mb{\v@lmax}}\xdef\cxb@{\costhet@}%
    \ifx\t@xt@\empty\edef\@rgdeux{\def@ulttheta}\fi\c@ssin{\C@}{\S@}{\@rgdeux}%
    \v@lmax=-\S@ pt%
    \v@leur=\v@lmax\v@leur=\costhet@\v@leur\xdef\cya@{\repdecn@mb{\v@leur}}%
    \v@leur=\v@lmax\v@leur=\sinthet@\v@leur\xdef\cyb@{\repdecn@mb{\v@leur}}%
    \xdef\cyc@{\C@}\v@lmin=-\C@ pt%
    \v@leur=\v@lmin\v@leur=\costhet@\v@leur\xdef\cza@{\repdecn@mb{\v@leur}}%
    \v@leur=\v@lmin\v@leur=\sinthet@\v@leur\xdef\czb@{\repdecn@mb{\v@leur}}%
    \xdef\czc@{\repdecn@mb{\v@lmax}}\fi%
    \xdef\v@lTheta{\@rgdeux}}}
\def\def@ultpsi{40}
\def\def@ulttheta{25}
\def\figset#1(#2){\def\t@xt@{#1}\ifx\t@xt@\empty\trtlis@rg{#2}{\Figsetwr@te}
    \else\keln@mde#1|%
    \def\n@mref{pr}\ifx\l@debut\n@mref\ifcurr@ntPS
     \immediate\write16{*** \BS@ figset proj(...) is ignored inside a
     \BS@ psbeginfig-\BS@ psendfig block.}\else\trtlis@rg{#2}{\Figsetpr@j}\fi\else%
    \def\n@mref{wr}\ifx\l@debut\n@mref\trtlis@rg{#2}{\Figsetwr@te}\else
    \immediate\write16{*** Unknown keyword: \BS@ figset #1(...)}%
    \fi\fi\fi\ignorespaces}
\def\Figsetpr@j#1=#2|{\keln@mtr#1|%
    \def\n@mref{dep}\ifx\l@debut\n@mref\Figsetd@p{#2}\else
    \def\n@mref{dis}\ifx\l@debut\n@mref%
     \ifnum\curr@ntproj=\tw@\figsetobdist(#2)\else\Figset@rr\fi\else
    \def\n@mref{lam}\ifx\l@debut\n@mref\Figsetd@p{#2}\else
    \def\n@mref{lat}\ifx\l@debut\n@mref\Figsetth@{#2}\else
    \def\n@mref{lon}\ifx\l@debut\n@mref\figsetview(#2)\else
    \def\n@mref{psi}\ifx\l@debut\n@mref\figsetview(#2)\else
    \def\n@mref{tar}\ifx\l@debut\n@mref%
     \ifnum\curr@ntproj=\tw@\figsettarget[#2]\else\Figset@rr\fi\else
    \def\n@mref{the}\ifx\l@debut\n@mref\Figsetth@{#2}\else
    \immediate\write16{*** Unknown attribute: \BS@ figset proj(..., #1=...).}%
    \fi\fi\fi\fi\fi\fi\fi\fi}
\def\Figsetd@p#1{\ifnum\curr@ntproj=\z@\figsetview(\v@lPsi,#1)\else\Figset@rr\fi}
\def\Figsetth@#1{\ifnum\curr@ntproj=\z@\Figset@rr\else\figsetview(\v@lPsi,#1)\fi}
\def\Figset@rr{\message{*** \BS@ figset proj(): Attribute "\n@mref" ignored, incompatible
    with current projection}}
\def\initb@undb@xTD{\wd\BminTD@=\maxdimen\ht\BminTD@=\maxdimen\dp\BminTD@=\maxdimen%
    \wd\BmaxTD@=-\maxdimen\ht\BmaxTD@=-\maxdimen\dp\BmaxTD@=-\maxdimen}
\newbox\Gb@x      
\newbox\Gb@xSC    
\newtoks\c@nsymb  
\newif\ifr@undcoord\newif\ifunitpr@sent
\def\unssqrttw@{0.707106 }
\def\figAst{\raise-1.15ex\hbox{$\ast$}}
\def\figBullet{\raise-1.15ex\hbox{$\bullet$}}
\def\figCirc{\raise-1.15ex\hbox{$\circ$}}
\def\figDiamond{\raise-1.15ex\hbox{$\diamond$}}%
\def\boxit#1#2{\leavevmode\hbox{\vrule\vbox{\hrule\vglue#1%
    \vtop{\hbox{\kern#1{#2}\kern#1}\vglue#1\hrule}}\vrule}}
\def\centertext#1#2{\vbox{\hsize#1\parindent0cm%
    \leftskip=0pt plus 1fil\rightskip=0pt plus 1fil\parfillskip=0pt{#2}}}
\def\lefttext#1#2{\vbox{\hsize#1\parindent0cm\rightskip=0pt plus 1fil#2}}
\def\c@nterpt{\ignorespaces%
    \kern-.5\wd\Gb@xSC%
    \raise-.5\ht\Gb@xSC\rlap{\hbox{\raise.5\dp\Gb@xSC\hbox{\copy\Gb@xSC}}}%
    \kern .5\wd\Gb@xSC\ignorespaces}
\def\b@undb@xSC#1#2{{\v@lXa=#1\v@lYa=#2%
    \v@leur=\ht\Gb@xSC\advance\v@leur\dp\Gb@xSC%
    \advance\v@lXa-.5\wd\Gb@xSC\advance\v@lYa-.5\v@leur\b@undb@x{\v@lXa}{\v@lYa}%
    \advance\v@lXa\wd\Gb@xSC\advance\v@lYa\v@leur\b@undb@x{\v@lXa}{\v@lYa}}}
\def\@keldist#1#2{\edef\Dist@n{#2}\y@tiunit{\Dist@n}%
    \ifunitpr@sent#1=\Dist@n\else#1=\Dist@n\unit@\fi}
\def\y@tiunit#1{\unitpr@sentfalse\expandafter\y@tiunit@#1:}
\def\y@tiunit@#1#2:{\ifcat#1a\unitpr@senttrue\else\def\l@suite{#2}%
    \ifx\l@suite\empty\else\y@tiunit@#2:\fi\fi}
\def\figcoordDD#1{{\v@lX=\ptT@unit@\v@lX\v@lY=\ptT@unit@\v@lY%
    \ifr@undcoord\ifcase#1\v@leur=0.5pt\or\v@leur=0.05pt\or\v@leur=0.005pt%
    \or\v@leur=0.0005pt\else\v@leur=\z@\fi%
    \ifdim\v@lX<\z@\advance\v@lX-\v@leur\else\advance\v@lX\v@leur\fi%
    \ifdim\v@lY<\z@\advance\v@lY-\v@leur\else\advance\v@lY\v@leur\fi\fi%
    (\@ffichnb{#1}{\repdecn@mb{\v@lX}},\ifmmode\else\thinspace\fi%
    \@ffichnb{#1}{\repdecn@mb{\v@lY}})}}
\def\@ffichnb#1#2{\def\@@ffich{\@ffich#1(}\edef\n@mbre{#2}%
    \expandafter\@@ffich\n@mbre)}
\def\@ffich#1(#2.#3){{#2\ifnum#1>\z@.\fi\def\dig@ts{#3}\s@mme=\z@%
    \loop\ifnum\s@mme<#1\expandafter\@ffichdec\dig@ts:\advance\s@mme\@ne\repeat}}
\def\@ffichdec#1#2:{\relax#1\def\dig@ts{#20}}
\def\figcoordTD#1{{\v@lX=\ptT@unit@\v@lX\v@lY=\ptT@unit@\v@lY\v@lZ=\ptT@unit@\v@lZ%
    \ifr@undcoord\ifcase#1\v@leur=0.5pt\or\v@leur=0.05pt\or\v@leur=0.005pt%
    \or\v@leur=0.0005pt\else\v@leur=\z@\fi%
    \ifdim\v@lX<\z@\advance\v@lX-\v@leur\else\advance\v@lX\v@leur\fi%
    \ifdim\v@lY<\z@\advance\v@lY-\v@leur\else\advance\v@lY\v@leur\fi%
    \ifdim\v@lZ<\z@\advance\v@lZ-\v@leur\else\advance\v@lZ\v@leur\fi\fi%
    (\@ffichnb{#1}{\repdecn@mb{\v@lX}},\ifmmode\else\thinspace\fi%
     \@ffichnb{#1}{\repdecn@mb{\v@lY}},\ifmmode\else\thinspace\fi%
     \@ffichnb{#1}{\repdecn@mb{\v@lZ}})}}
\def\figsetroundcoord#1{\expandafter\Figsetr@undcoord#1:\ignorespaces}
\def\Figsetr@undcoord#1#2:{\if#1n\r@undcoordfalse\else\r@undcoordtrue\fi}
\def\Figsetwr@te#1=#2|{\keln@mun#1|%
    \def\n@mref{m}\ifx\l@debut\n@mref\figsetmark{#2}\else
    \immediate\write16{*** Unknown attribute: \BS@ figset (..., #1=...)}%
    \fi}
\def\figsetmark#1{\c@nsymb={#1}\setbox\Gb@xSC=\hbox{\the\c@nsymb}\ignorespaces}
\def\figsetptname#1{\def\ptn@me##1{#1}\ignorespaces}
\def\FigWrit@L#1:#2(#3,#4){\ignorespaces\@keldist\v@leur{#3}\@keldist\delt@{#4}%
    \C@rp@r@m\def\list@num{#1}\@ecfor\p@int:=\list@num\do{\FigWrit@pt{\p@int}{#2}}}
\def\FigWrit@pt#1#2{\FigWp@r@m{#1}{#2}\Vc@rrect\figWp@si%
    \ifdim\wd\Gb@xSC>\z@\b@undb@xSC{\v@lX}{\v@lY}\fi\figWBB@x}
\def\FigWp@r@m#1#2{\Figg@tXY{#1}%
    \setbox\Gb@x=\hbox{\def\t@xt@{#2}\ifx\t@xt@\empty\Figg@tT{#1}\else#2\fi}\c@lprojSP}
\let\Vc@rrect=\relax
\let\C@rp@r@m=\relax
\def\figwrite[#1]#2{{\ignorespaces\def\list@num{#1}\@ecfor\p@int:=\list@num\do{%
    \setbox\Gb@x=\hbox{\def\t@xt@{#2}\ifx\t@xt@\empty\Figg@tT{\p@int}\else#2\fi}%
    \Figwrit@{\p@int}}}\ignorespaces}
\def\Figwrit@#1{\Figg@tXY{#1}\c@lprojSP%
    \rlap{\kern\v@lX\raise\v@lY\hbox{\unhcopy\Gb@x}}\v@leur=\v@lY%
    \advance\v@lY\ht\Gb@x\b@undb@x{\v@lX}{\v@lY}\advance\v@lX\wd\Gb@x%
    \v@lY=\v@leur\advance\v@lY-\dp\Gb@x\b@undb@x{\v@lX}{\v@lY}}
\def\figwritec[#1]#2{{\ignorespaces\def\list@num{#1}%
    \@ecfor\p@int:=\list@num\do{\Figwrit@c{\p@int}{#2}}}\ignorespaces}
\def\Figwrit@c#1#2{\FigWp@r@m{#1}{#2}%
    \rlap{\kern\v@lX\raise\v@lY\hbox{\rlap{\kern-.5\wd\Gb@x%
    \raise-.5\ht\Gb@x\hbox{\raise.5\dp\Gb@x\hbox{\unhcopy\Gb@x}}}}}%
    \v@leur=\ht\Gb@x\advance\v@leur\dp\Gb@x%
    \advance\v@lX-.5\wd\Gb@x\advance\v@lY-.5\v@leur\b@undb@x{\v@lX}{\v@lY}%
    \advance\v@lX\wd\Gb@x\advance\v@lY\v@leur\b@undb@x{\v@lX}{\v@lY}}
\def\figwritep[#1]{{\ignorespaces\def\list@num{#1}\setbox\Gb@x=\hbox{\c@nterpt}%
    \@ecfor\p@int:=\list@num\do{\Figwrit@{\p@int}}}\ignorespaces}
\def\figwritew#1:#2(#3){\figwritegcw#1:{#2}(#3,0pt)}
\def\figwritee#1:#2(#3){\figwritegce#1:{#2}(#3,0pt)}
\def\figwriten#1:#2(#3){{\def\Vc@rrect{\v@lZ=\v@leur\advance\v@lZ\dp\Gb@x}%
    \Figwrit@NS#1:{#2}(#3)}\ignorespaces}
\def\figwrites#1:#2(#3){{\def\Vc@rrect{\v@lZ=-\v@leur\advance\v@lZ-\ht\Gb@x}%
    \Figwrit@NS#1:{#2}(#3)}\ignorespaces}
\def\Figwrit@NS#1:#2(#3){\let\figWp@si=\FigWp@siNS\let\figWBB@x=\FigWBB@xNS%
    \FigWrit@L#1:{#2}(#3,0pt)}
\def\FigWp@siNS{\rlap{\kern\v@lX\raise\v@lY\hbox{\rlap{\kern-.5\wd\Gb@x%
    \raise\v@lZ\hbox{\unhcopy\Gb@x}}\c@nterpt}}}
\def\FigWBB@xNS{\advance\v@lY\v@lZ%
    \advance\v@lY-\dp\Gb@x\advance\v@lX-.5\wd\Gb@x\b@undb@x{\v@lX}{\v@lY}%
    \advance\v@lY\ht\Gb@x\advance\v@lY\dp\Gb@x%
    \advance\v@lX\wd\Gb@x\b@undb@x{\v@lX}{\v@lY}}
\def\figwritenw#1:#2(#3){{\let\figWp@si=\FigWp@sigW\let\figWBB@x=\FigWBB@xgWE%
    \def\C@rp@r@m{\v@leur=\unssqrttw@\v@leur\delt@=\v@leur%
    \ifdim\delt@=\z@\delt@=\epsil@n\fi}\let@xte={-}\FigWrit@L#1:{#2}(#3,0pt)}\ignorespaces}
\def\figwritesw#1:#2(#3){{\let\figWp@si=\FigWp@sigW\let\figWBB@x=\FigWBB@xgWE%
    \def\C@rp@r@m{\v@leur=\unssqrttw@\v@leur\delt@=-\v@leur%
    \ifdim\delt@=\z@\delt@=-\epsil@n\fi}\let@xte={-}\FigWrit@L#1:{#2}(#3,0pt)}\ignorespaces}
\def\figwritene#1:#2(#3){{\let\figWp@si=\FigWp@sigE\let\figWBB@x=\FigWBB@xgWE%
    \def\C@rp@r@m{\v@leur=\unssqrttw@\v@leur\delt@=\v@leur%
    \ifdim\delt@=\z@\delt@=\epsil@n\fi}\let@xte={}\FigWrit@L#1:{#2}(#3,0pt)}\ignorespaces}
\def\figwritese#1:#2(#3){{\let\figWp@si=\FigWp@sigE\let\figWBB@x=\FigWBB@xgWE%
    \def\C@rp@r@m{\v@leur=\unssqrttw@\v@leur\delt@=-\v@leur%
    \ifdim\delt@=\z@\delt@=-\epsil@n\fi}\let@xte={}\FigWrit@L#1:{#2}(#3,0pt)}\ignorespaces}
\def\figwritegw#1:#2(#3,#4){{\let\figWp@si=\FigWp@sigW\let\figWBB@x=\FigWBB@xgWE%
    \let@xte={-}\FigWrit@L#1:{#2}(#3,#4)}\ignorespaces}
\def\figwritege#1:#2(#3,#4){{\let\figWp@si=\FigWp@sigE\let\figWBB@x=\FigWBB@xgWE%
    \let@xte={}\FigWrit@L#1:{#2}(#3,#4)}\ignorespaces}
\def\FigWp@sigW{\v@lXa=\z@\v@lYa=\ht\Gb@x\advance\v@lYa\dp\Gb@x%
    \ifdim\delt@>\z@\relax%
    \rlap{\kern\v@lX\raise\v@lY\hbox{\rlap{\kern-\wd\Gb@x\kern-\v@leur%
          \raise\delt@\hbox{\raise\dp\Gb@x\hbox{\unhcopy\Gb@x}}}\c@nterpt}}%
    \else\ifdim\delt@<\z@\relax\v@lYa=-\v@lYa%
    \rlap{\kern\v@lX\raise\v@lY\hbox{\rlap{\kern-\wd\Gb@x\kern-\v@leur%
          \raise\delt@\hbox{\raise-\ht\Gb@x\hbox{\unhcopy\Gb@x}}}\c@nterpt}}%
    \else\v@lXa=-.5\v@lYa%
    \rlap{\kern\v@lX\raise\v@lY\hbox{\rlap{\kern-\wd\Gb@x\kern-\v@leur%
          \raise-.5\ht\Gb@x\hbox{\raise.5\dp\Gb@x\hbox{\unhcopy\Gb@x}}}\c@nterpt}}%
    \fi\fi}
\def\FigWp@sigE{\v@lXa=\z@\v@lYa=\ht\Gb@x\advance\v@lYa\dp\Gb@x%
    \ifdim\delt@>\z@\relax%
    \rlap{\kern\v@lX\raise\v@lY\hbox{\c@nterpt\kern\v@leur%
          \raise\delt@\hbox{\raise\dp\Gb@x\hbox{\unhcopy\Gb@x}}}}%
    \else\ifdim\delt@<\z@\relax\v@lYa=-\v@lYa%
    \rlap{\kern\v@lX\raise\v@lY\hbox{\c@nterpt\kern\v@leur%
          \raise\delt@\hbox{\raise-\ht\Gb@x\hbox{\unhcopy\Gb@x}}}}%
    \else\v@lXa=-.5\v@lYa%
    \rlap{\kern\v@lX\raise\v@lY\hbox{\c@nterpt\kern\v@leur%
          \raise-.5\ht\Gb@x\hbox{\raise.5\dp\Gb@x\hbox{\unhcopy\Gb@x}}}}%
    \fi\fi}
\def\FigWBB@xgWE{\advance\v@lY\delt@%
    \advance\v@lX\the\let@xte\v@leur\advance\v@lY\v@lXa\b@undb@x{\v@lX}{\v@lY}%
    \advance\v@lX\the\let@xte\wd\Gb@x\advance\v@lY\v@lYa\b@undb@x{\v@lX}{\v@lY}}
\def\figwritegcw#1:#2(#3,#4){{\let\figWp@si=\FigWp@sigcW\let\figWBB@x=\FigWBB@xgcWE%
    \let@xte={-}\FigWrit@L#1:{#2}(#3,#4)}\ignorespaces}
\def\figwritegce#1:#2(#3,#4){{\let\figWp@si=\FigWp@sigcE\let\figWBB@x=\FigWBB@xgcWE%
    \let@xte={}\FigWrit@L#1:{#2}(#3,#4)}\ignorespaces}
\def\FigWp@sigcW{\rlap{\kern\v@lX\raise\v@lY\hbox{\rlap{\kern-\wd\Gb@x\kern-\v@leur%
     \raise-.5\ht\Gb@x\hbox{\raise\delt@\hbox{\raise.5\dp\Gb@x\hbox{\unhcopy\Gb@x}}}}%
     \c@nterpt}}}
\def\FigWp@sigcE{\rlap{\kern\v@lX\raise\v@lY\hbox{\c@nterpt\kern\v@leur%
    \raise-.5\ht\Gb@x\hbox{\raise\delt@\hbox{\raise.5\dp\Gb@x\hbox{\unhcopy\Gb@x}}}}}}
\def\FigWBB@xgcWE{\v@lZ=\ht\Gb@x\advance\v@lZ\dp\Gb@x%
    \advance\v@lX\the\let@xte\v@leur\advance\v@lY\delt@\advance\v@lY.5\v@lZ%
    \b@undb@x{\v@lX}{\v@lY}%
    \advance\v@lX\the\let@xte\wd\Gb@x\advance\v@lY-\v@lZ\b@undb@x{\v@lX}{\v@lY}}
\def\figwritebn#1:#2(#3){{\def\Vc@rrect{\v@lZ=\v@leur}\Figwrit@NS#1:{#2}(#3)}\ignorespaces}
\def\figwritebs#1:#2(#3){{\def\Vc@rrect{\v@lZ=-\v@leur}\Figwrit@NS#1:{#2}(#3)}\ignorespaces}
\def\figwritebw#1:#2(#3){{\let\figWp@si=\FigWp@sibW\let\figWBB@x=\FigWBB@xbWE%
    \let@xte={-}\FigWrit@L#1:{#2}(#3,0pt)}\ignorespaces}
\def\figwritebe#1:#2(#3){{\let\figWp@si=\FigWp@sibE\let\figWBB@x=\FigWBB@xbWE%
    \let@xte={}\FigWrit@L#1:{#2}(#3,0pt)}\ignorespaces}
\def\FigWp@sibW{\rlap{\kern\v@lX\raise\v@lY\hbox{\rlap{\kern-\wd\Gb@x\kern-\v@leur%
          \hbox{\unhcopy\Gb@x}}\c@nterpt}}}
\def\FigWp@sibE{\rlap{\kern\v@lX\raise\v@lY\hbox{\c@nterpt\kern\v@leur%
          \hbox{\unhcopy\Gb@x}}}}
\def\FigWBB@xbWE{\v@lZ=\ht\Gb@x\advance\v@lZ\dp\Gb@x%
    \advance\v@lX\the\let@xte\v@leur\advance\v@lY\ht\Gb@x\b@undb@x{\v@lX}{\v@lY}%
    \advance\v@lX\the\let@xte\wd\Gb@x\advance\v@lY-\v@lZ\b@undb@x{\v@lX}{\v@lY}}
\newread\frf@g  \newwrite\fwf@g
\newif\ifcurr@ntPS
\newif\ifps@cri
\newif\ifUse@llipse
\newif\ifpsdebugmode \psdebugmodefalse 
\newif\ifPDFm@ke
\ifx\pdfliteral\undefined\else\ifnum\pdfoutput>\z@\PDFm@ketrue\fi\fi
\ifPDFm@ke
 \def\c@mcurveto{c}
 \def\c@mfill{f}
 \def\c@mgsave{q}
 \def\c@mgrestore{Q}
 \def\c@mlineto{l}
 \def\c@mmoveto{m}
 \def\c@msetgray{g}     \def\c@msetgrayStroke{G}
 \def\c@msetcmykcolor{k}\def\c@msetcmykcolorStroke{K}
 \def\c@msetrgbcolor{rg}\def\c@msetrgbcolorStroke{RG}
 \def\d@fprimarC@lor{\curr@ntcolor\space\curr@ntcolorc@md%
               \space\curr@ntcolor\space\curr@ntcolorc@mdStroke}
 \def\d@fsecondC@lor{\sec@ndcolor\space\sec@ndcolorc@md%
               \space\sec@ndcolor\space\sec@ndcolorc@mdStroke}
 \def\d@fthirdC@lor{\th@rdcolor\space\th@rdcolorc@md%
              \space\th@rdcolor\space\th@rdcolorc@mdStroke}
 \def\c@msetdash{d}
 \def\c@msetlinejoin{j}
 \def\c@msetlinewidth{w}
 \def\f@gclosestroke{\immediate\write\fwf@g{s}}
 \def\f@gfill{\immediate\write\fwf@g{\fillc@md}}
 \def\f@gnewpath{}
 \def\f@gstroke{\immediate\write\fwf@g{S}}
\else
 \let\figinsertE=\figinsert
 \def\c@mcurveto{curveto}
 \def\c@mfill{fill}
 \def\c@mgsave{gsave}
 \def\c@mgrestore{grestore}
 \def\c@mlineto{lineto}
 \def\c@mmoveto{moveto}
 \def\c@msetgray{setgray}          \def\c@msetgrayStroke{}
 \def\c@msetcmykcolor{setcmykcolor}\def\c@msetcmykcolorStroke{}
 \def\c@msetrgbcolor{setrgbcolor}  \def\c@msetrgbcolorStroke{}
 \def\d@fprimarC@lor{\curr@ntcolor\space\curr@ntcolorc@md}
 \def\d@fsecondC@lor{\sec@ndcolor\space\sec@ndcolorc@md}
 \def\d@fthirdC@lor{\th@rdcolor\space\th@rdcolorc@md}
 \def\c@msetdash{setdash}
 \def\c@msetlinejoin{setlinejoin}
 \def\c@msetlinewidth{setlinewidth}
 \def\f@gclosestroke{\immediate\write\fwf@g{closepath\space stroke}}
 \def\f@gfill{\immediate\write\fwf@g{\fillc@md}}
 \def\f@gnewpath{\immediate\write\fwf@g{newpath}}
 \def\f@gstroke{\immediate\write\fwf@g{stroke}}
\fi
\def\c@pypsfile#1#2{\c@pyfil@{\immediate\write#1}{#2}}
\def\Figinclud@PDF#1#2{\openin\frf@g=#1\pdfliteral{q #2 0 0 #2 0 0 cm}%
    \c@pyfil@{\pdfliteral}{\frf@g}\pdfliteral{Q}\closein\frf@g}
\newif\ifmored@ta
\def\c@pyfil@#1#2{\def\blankline{\par}{\catcode`\%=12
    \loop\ifeof#2\mored@tafalse\else\mored@tatrue\immediate\read#2 to\tr@c
    \ifx\tr@c\blankline\else#1{\tr@c}\fi\fi\ifmored@ta\repeat}}
\def\keln@mun#1#2|{\def\l@debut{#1}\def\l@suite{#2}}
\def\keln@mde#1#2#3|{\def\l@debut{#1#2}\def\l@suite{#3}}
\def\keln@mtr#1#2#3#4|{\def\l@debut{#1#2#3}\def\l@suite{#4}}
\def\keln@mqu#1#2#3#4#5|{\def\l@debut{#1#2#3#4}\def\l@suite{#5}}
\let\@psffilein=\frf@g 
\newif\if@psffileok    
\newif\if@psfbbfound   
\newif\if@psfverbose   
\@psfverbosetrue
\def\@psfgetbb#1{\global\@psfbbfoundfalse%
\global\def\@psfllx{0}\global\def\@psflly{0}%
\global\def\@psfurx{30}\global\def\@psfury{30}%
\openin\@psffilein=#1
\ifeof\@psffilein\errmessage{I couldn't open #1, will ignore it}\else
   \edef\setcolonc@tcode{\catcode`\noexpand\:\the\catcode`\:\relax}%
   {\@psffileoktrue \chardef\other=12
    \def\do##1{\catcode`##1=\other}\dospecials \catcode`\ =10 \setcolonc@tcode
    \loop
       \read\@psffilein to \@psffileline
       \ifeof\@psffilein\@psffileokfalse\else
          \expandafter\@psfaux\@psffileline:. \\%
       \fi
   \if@psffileok\repeat
   \if@psfbbfound\else
    \if@psfverbose\message{No bounding box comment in #1; using defaults}\fi\fi
   }\closein\@psffilein\fi}%
{\catcode`\%=12 \global\let\@psfpercent=
\long\def\@psfaux#1#2:#3\\{\ifx#1\@psfpercent
   \def\testit{#2}\ifx\testit\@psfbblit
      \@psfgrab #3 . . . \\%
      \@psffileokfalse
      \global\@psfbbfoundtrue
   \fi\else\ifx#1\par\else\@psffileokfalse\fi\fi}%
\def\@psfempty{}%
\def\@psfgrab #1 #2 #3 #4 #5\\{%
\global\def\@psfllx{#1}\ifx\@psfllx\@psfempty
      \@psfgrab #2 #3 #4 #5 .\\\else
   \global\def\@psflly{#2}%
   \global\def\@psfurx{#3}\global\def\@psfury{#4}\fi}%
\def\PSwrit@cmd#1#2#3{{\Figg@tXY{#1}\c@lprojSP\b@undb@x{\v@lX}{\v@lY}%
    \v@lX=\ptT@ptps\v@lX\v@lY=\ptT@ptps\v@lY%
    \immediate\write#3{\repdecn@mb{\v@lX}\space\repdecn@mb{\v@lY}\space#2}}}
\def\PSwrit@cmdS#1#2#3#4#5{{\Figg@tXY{#1}\c@lprojSP\b@undb@x{\v@lX}{\v@lY}%
    \global\result@t=\v@lX\global\result@@t=\v@lY%
    \v@lX=\ptT@ptps\v@lX\v@lY=\ptT@ptps\v@lY%
    \immediate\write#3{\repdecn@mb{\v@lX}\space\repdecn@mb{\v@lY}\space#2}}%
    \edef#4{\the\result@t}\edef#5{\the\result@@t}}
\def\psaltitude#1[#2,#3,#4]{{\ifcurr@ntPS\ifps@cri%
    \PSc@mment{psaltitude Square Dim=#1, Triangle=[#2 / #3,#4]}%
    \s@uvc@ntr@l\et@tpsaltitude\resetc@ntr@l{2}\figptorthoprojline-5:=#2/#3,#4/%
    \figvectP -1[#3,#4]\n@rminf{\v@leur}{-1}\vecunit@{-3}{-1}%
    \figvectP -1[-5,#3]\n@rminf{\v@lmin}{-1}\figvectP -2[-5,#4]\n@rminf{\v@lmax}{-2}%
    \ifdim\v@lmin<\v@lmax\s@mme=#3\else\v@lmax=\v@lmin\s@mme=#4\fi%
    \figvectP -4[-5,#2]\vecunit@{-4}{-4}\delt@=#1\unit@%
    \edef\t@ille{\repdecn@mb{\delt@}}\figpttra-1:=-5/\t@ille,-3/%
    \figptstra-3=-5,-1/\t@ille,-4/\psline[#2,-5]\psline[-1,-2,-3]%
    \ifdim\v@leur<\v@lmax\Pss@tsecondSt\psline[-5,\the\s@mme]\Psrest@reSt\fi%
    \PSc@mment{End psaltitude}\resetc@ntr@l\et@tpsaltitude\fi\fi}}
\def\Ps@rcerc#1;#2(#3,#4){\ellBB@x#1;#2,#2(#3,#4,0)%
    \f@gnewpath{\delt@=#2\unit@\delt@=\ptT@ptps\delt@%
    \BdingB@xfalse%
    \PSwrit@cmd{#1}{\repdecn@mb{\delt@}\space #3\space #4\space arc}{\fwf@g}}}
\def\psarccircDD#1;#2(#3,#4){\ifcurr@ntPS\ifps@cri%
    \PSc@mment{psarccircDD Center=#1 ; Radius=#2 (Ang1=#3, Ang2=#4)}%
    \iffillm@de\Ps@rcerc#1;#2(#3,#4)%
    \f@gfill%
    \else\Ps@rcerc#1;#2(#3,#4)\f@gstroke\fi%
    \PSc@mment{End psarccircDD}\fi\fi}
\def\psarccircTD#1,#2,#3;#4(#5,#6){{\ifcurr@ntPS\ifps@cri\s@uvc@ntr@l\et@tpsarccircTD%
    \PSc@mment{psarccircTD Center=#1,P1=#2,P2=#3 ; Radius=#4 (Ang1=#5, Ang2=#6)}%
    \setc@ntr@l{2}\c@lExtAxes#1,#2,#3(#4)\psarcellPATD#1,-4,-5(#5,#6)%
    \PSc@mment{End psarccircTD}\resetc@ntr@l\et@tpsarccircTD\fi\fi}}
\def\c@lExtAxes#1,#2,#3(#4){%
    \figvectPTD-5[#1,#2]\vecunit@{-5}{-5}\figvectNTD-4[#1,#2,#3]\vecunit@{-4}{-4}%
    \figvectNVTD-3[-4,-5]\delt@=#4\unit@\edef\r@yon{\repdecn@mb{\delt@}}%
    \figpttra-4:=#1/\r@yon,-5/\figpttra-5:=#1/\r@yon,-3/}
\def\psarccircPDD#1;#2[#3,#4]{{\ifcurr@ntPS\ifps@cri\s@uvc@ntr@l\et@tpsarccircPDD%
    \PSc@mment{psarccircPDD Center=#1; Radius=#2, [P1=#3, P2=#4]}%
    \Ps@ngleparam#1;#2[#3,#4]\ifdim\v@lmin>\v@lmax\advance\v@lmax\DePI@deg\fi%
    \edef\@ngdeb{\repdecn@mb{\v@lmin}}\edef\@ngfin{\repdecn@mb{\v@lmax}}%
    \psarccirc#1;\r@dius(\@ngdeb,\@ngfin)%
    \PSc@mment{End psarccircPDD}\resetc@ntr@l\et@tpsarccircPDD\fi\fi}}
\def\psarccircPTD#1;#2[#3,#4,#5]{{\ifcurr@ntPS\ifps@cri\s@uvc@ntr@l\et@tpsarccircPTD%
    \PSc@mment{psarccircPTD Center=#1; Radius=#2, [P1=#3, P2=#4, P3=#5]}%
    \setc@ntr@l{2}\c@lExtAxes#1,#3,#5(#2)\psarcellPP#1,-4,-5[#3,#4]%
    \PSc@mment{End psarccircPTD}\resetc@ntr@l\et@tpsarccircPTD\fi\fi}}
\def\Ps@ngleparam#1;#2[#3,#4]{\setc@ntr@l{2}%
    \figvectPDD-1[#1,#3]\vecunit@{-1}{-1}\Figg@tXY{-1}\arct@n\v@lmin(\v@lX,\v@lY)%
    \figvectPDD-2[#1,#4]\vecunit@{-2}{-2}\Figg@tXY{-2}\arct@n\v@lmax(\v@lX,\v@lY)%
    \v@lmin=\rdT@deg\v@lmin\v@lmax=\rdT@deg\v@lmax%
    \v@leur=#2pt\maxim@m{\mili@u}{-\v@leur}{\v@leur}%
    \edef\r@dius{\repdecn@mb{\mili@u}}}
\def\Ps@rcercBz#1;#2(#3,#4){\Ps@rellBz#1;#2,#2(#3,#4,0)}
\def\Ps@rellBz#1;#2,#3(#4,#5,#6){%
    \ellBB@x#1;#2,#3(#4,#5,#6)\BdingB@xfalse%
    \c@lNbarcs{#4}{#5}\v@leur=#4pt\setc@ntr@l{2}\figptell-13::#1;#2,#3(#4,#6)%
    \f@gnewpath\PSwrit@cmd{-13}{\c@mmoveto}{\fwf@g}%
    \s@mme=\z@\bcl@rellBz#1;#2,#3(#6)\BdingB@xtrue}
\def\bcl@rellBz#1;#2,#3(#4){\relax%
    \ifnum\s@mme<\p@rtent\advance\s@mme\@ne%
    \advance\v@leur\delt@\edef\@ngle{\repdecn@mb\v@leur}\figptell-14::#1;#2,#3(\@ngle,#4)%
    \advance\v@leur\delt@\edef\@ngle{\repdecn@mb\v@leur}\figptell-15::#1;#2,#3(\@ngle,#4)%
    \advance\v@leur\delt@\edef\@ngle{\repdecn@mb\v@leur}\figptell-16::#1;#2,#3(\@ngle,#4)%
    \figptscontrolDD-18[-13,-14,-15,-16]%
    \PSwrit@cmd{-18}{}{\fwf@g}\PSwrit@cmd{-17}{}{\fwf@g}%
    \PSwrit@cmd{-16}{\c@mcurveto}{\fwf@g}%
    \figptcopyDD-13:/-16/\bcl@rellBz#1;#2,#3(#4)\fi}
\def\Ps@rell#1;#2,#3(#4,#5,#6){\ellBB@x#1;#2,#3(#4,#5,#6)%
    \f@gnewpath{\v@lmin=#2\unit@\v@lmin=\ptT@ptps\v@lmin%
    \v@lmax=#3\unit@\v@lmax=\ptT@ptps\v@lmax\BdingB@xfalse%
    \PSwrit@cmd{#1}%
    {#6\space\repdecn@mb{\v@lmin}\space\repdecn@mb{\v@lmax}\space #4\space #5\space ellipse}{\fwf@g}}%
    \global\Use@llipsetrue}
\def\psarcellDD#1;#2,#3(#4,#5,#6){{\ifcurr@ntPS\ifps@cri%
    \PSc@mment{psarcellDD Center=#1 ; XRad=#2, YRad=#3 (Ang1=#4, Ang2=#5, Inclination=#6)}%
    \iffillm@de\Ps@rell#1;#2,#3(#4,#5,#6)%
    \f@gfill%
    \else\Ps@rell#1;#2,#3(#4,#5,#6)\f@gstroke\fi%
    \PSc@mment{End psarcellDD}\fi\fi}}
\def\psarcellTD#1;#2,#3(#4,#5,#6){{\ifcurr@ntPS\ifps@cri\s@uvc@ntr@l\et@tpsarcellTD%
    \PSc@mment{psarcellTD Center=#1 ; XRad=#2, YRad=#3 (Ang1=#4, Ang2=#5, Inclination=#6)}%
    \setc@ntr@l{2}\figpttraC -8:=#1/#2,0,0/\figpttraC -7:=#1/0,#3,0/%
    \figvectC -4(0,0,1)\figptsrot -8=-8,-7/#1,#6,-4/\psarcellPATD#1,-8,-7(#4,#5)%
    \PSc@mment{End psarcellTD}\resetc@ntr@l\et@tpsarcellTD\fi\fi}}
\def\psarcellPADD#1,#2,#3(#4,#5){{\ifcurr@ntPS\ifps@cri\s@uvc@ntr@l\et@tpsarcellPADD%
    \PSc@mment{psarcellPADD Center=#1,PtAxis1=#2,PtAxis2=#3 (Ang1=#4, Ang2=#5)}%
    \setc@ntr@l{2}\figvectPDD-1[#1,#2]\vecunit@DD{-1}{-1}\v@lX=\ptT@unit@\result@t%
    \edef\XR@d{\repdecn@mb{\v@lX}}\Figg@tXY{-1}\arct@n\v@lmin(\v@lX,\v@lY)%
    \v@lmin=\rdT@deg\v@lmin\edef\Inclin@{\repdecn@mb{\v@lmin}}%
    \figgetdist\YR@d[#1,#3]\psarcellDD#1;\XR@d,\YR@d(#4,#5,\Inclin@)%
    \PSc@mment{End psarcellPADD}\resetc@ntr@l\et@tpsarcellPADD\fi\fi}}
\def\psarcellPATD#1,#2,#3(#4,#5){{\ifcurr@ntPS\ifps@cri\s@uvc@ntr@l\et@tpsarcellPATD%
    \PSc@mment{psarcellPATD Center=#1,PtAxis1=#2,PtAxis2=#3 (Ang1=#4, Ang2=#5)}%
    \iffillm@de\Ps@rellPATD#1,#2,#3(#4,#5)%
    \f@gfill%
    \else\Ps@rellPATD#1,#2,#3(#4,#5)\f@gstroke\fi%
    \PSc@mment{End psarcellPATD}\resetc@ntr@l\et@tpsarcellPATD\fi\fi}}
\def\Ps@rellPATD#1,#2,#3(#4,#5){\let\c@lprojSP=\relax%
    \setc@ntr@l{2}\figvectPTD-1[#1,#2]\figvectPTD-2[#1,#3]\c@lNbarcs{#4}{#5}%
    \v@leur=#4pt\c@lptellP{#1}{-1}{-2}\Figptpr@j-5:/-3/%
    \f@gnewpath\PSwrit@cmdS{-5}{\c@mmoveto}{\fwf@g}{\X@un}{\Y@un}%
    \edef\C@nt@r{#1}\s@mme=\z@\bcl@rellPATD}
\def\bcl@rellPATD{\relax%
    \ifnum\s@mme<\p@rtent\advance\s@mme\@ne%
    \advance\v@leur\delt@\c@lptellP{\C@nt@r}{-1}{-2}\Figptpr@j-4:/-3/%
    \advance\v@leur\delt@\c@lptellP{\C@nt@r}{-1}{-2}\Figptpr@j-6:/-3/%
    \advance\v@leur\delt@\c@lptellP{\C@nt@r}{-1}{-2}\Figptpr@j-3:/-3/%
    \v@lX=\z@\v@lY=\z@\Figtr@nptDD{-5}{-5}\Figtr@nptDD{2}{-3}%
    \divide\v@lX\@vi\divide\v@lY\@vi%
    \Figtr@nptDD{3}{-4}\Figtr@nptDD{-1.5}{-6}\v@lmin=\v@lX\v@lmax=\v@lY%
    \v@lX=\z@\v@lY=\z@\Figtr@nptDD{2}{-5}\Figtr@nptDD{-5}{-3}%
    \divide\v@lX\@vi\divide\v@lY\@vi\Figtr@nptDD{-1.5}{-4}\Figtr@nptDD{3}{-6}%
    \BdingB@xfalse%
    \Figp@intregDD-4:(\v@lmin,\v@lmax)\PSwrit@cmdS{-4}{}{\fwf@g}{\X@de}{\Y@de}%
    \Figp@intregDD-4:(\v@lX,\v@lY)\PSwrit@cmdS{-4}{}{\fwf@g}{\X@tr}{\Y@tr}%
    \BdingB@xtrue\PSwrit@cmdS{-3}{\c@mcurveto}{\fwf@g}{\X@qu}{\Y@qu}%
    \B@zierBB@x{1}{\Y@un}(\X@un,\X@de,\X@tr,\X@qu)%
    \B@zierBB@x{2}{\X@un}(\Y@un,\Y@de,\Y@tr,\Y@qu)%
    \edef\X@un{\X@qu}\edef\Y@un{\Y@qu}\figptcopyDD-5:/-3/\bcl@rellPATD\fi}
\def\c@lNbarcs#1#2{%
    \delt@=#2pt\advance\delt@-#1pt\maxim@m{\v@lmax}{\delt@}{-\delt@}%
    \v@leur=\v@lmax\divide\v@leur45 \p@rtentiere{\p@rtent}{\v@leur}\advance\p@rtent\@ne%
    \s@mme=\p@rtent\multiply\s@mme\thr@@\divide\delt@\s@mme}
\def\psarcellPP#1,#2,#3[#4,#5]{{\ifcurr@ntPS\ifps@cri\s@uvc@ntr@l\et@tpsarcellPP%
    \PSc@mment{psarcellPP Center=#1,PtAxis1=#2,PtAxis2=#3 [Point1=#4, Point2=#5]}%
    \setc@ntr@l{2}\figvectP-2[#1,#3]\vecunit@{-2}{-2}\v@lmin=\result@t%
    \invers@{\v@lmax}{\v@lmin}%
    \figvectP-1[#1,#2]\vecunit@{-1}{-1}\v@leur=\result@t%
    \v@leur=\repdecn@mb{\v@lmax}\v@leur\edef\AsB@{\repdecn@mb{\v@leur}}
    \c@lAngle{#1}{#4}{\v@lmin}\edef\@ngdeb{\repdecn@mb{\v@lmin}}%
    \c@lAngle{#1}{#5}{\v@lmax}\ifdim\v@lmin>\v@lmax\advance\v@lmax\DePI@deg\fi%
    \edef\@ngfin{\repdecn@mb{\v@lmax}}\psarcellPA#1,#2,#3(\@ngdeb,\@ngfin)%
    \PSc@mment{End psarcellPP}\resetc@ntr@l\et@tpsarcellPP\fi\fi}}
\def\c@lAngle#1#2#3{\figvectP-3[#1,#2]%
    \c@lproscal\delt@[-3,-1]\c@lproscal\v@leur[-3,-2]%
    \v@leur=\AsB@\v@leur\arct@n#3(\delt@,\v@leur)#3=\rdT@deg#3}
\newif\if@rrowratio\@rrowratiotrue
\newif\if@rrowhfill
\newif\if@rrowhout
\def\Psset@rrowhe@d#1=#2|{\keln@mun#1|%
    \def\n@mref{a}\ifx\l@debut\n@mref\pssetarrowheadangle{#2}\else
    \def\n@mref{f}\ifx\l@debut\n@mref\pssetarrowheadfill{#2}\else
    \def\n@mref{l}\ifx\l@debut\n@mref\pssetarrowheadlength{#2}\else
    \def\n@mref{o}\ifx\l@debut\n@mref\pssetarrowheadout{#2}\else
    \def\n@mref{r}\ifx\l@debut\n@mref\pssetarrowheadratio{#2}\else
    \immediate\write16{*** Unknown attribute: \BS@ psset arrowhead(..., #1=...)}%
    \fi\fi\fi\fi\fi}
\def\pssetarrowheadangle#1{\edef\@rrowheadangle{#1}{\c@ssin{\C@}{\S@}{#1}%
    \xdef\C@AHANG{\C@}\xdef\S@AHANG{\S@}\v@lmax=\S@ pt%
    \invers@{\v@leur}{\v@lmax}\maxim@m{\v@leur}{\v@leur}{-\v@leur}%
    \xdef\UNSS@N{\the\v@leur}}}
\def\pssetarrowheadfill#1{\expandafter\set@rrowhfill#1:}
\def\set@rrowhfill#1#2:{\if#1n\@rrowhfillfalse\else\@rrowhfilltrue\fi}
\def\pssetarrowheadout#1{\expandafter\set@rrowhout#1:}
\def\set@rrowhout#1#2:{\if#1n\@rrowhoutfalse\else\@rrowhouttrue\fi}
\def\pssetarrowheadlength#1{\edef\@rrowheadlength{#1}\@rrowratiofalse}
\def\pssetarrowheadratio#1{\edef\@rrowheadratio{#1}\@rrowratiotrue}
\def\psresetarrowhead{%
    \pssetarrowheadangle{\defaultarrowheadangle}%
    \pssetarrowheadfill{\defaultarrowheadfill}%
    \pssetarrowheadout{\defaultarrowheadout}%
    \pssetarrowheadratio{\defaultarrowheadratio}%
    \d@fm@cdim\defaultarrowheadlength{\defaulth@rdahlength}
    \pssetarrowheadlength{\defaultarrowheadlength}}
\def\defaultarrowheadratio{0.1}
\def\defaultarrowheadangle{20}
\def\defaultarrowheadfill{no}
\def\defaultarrowheadout{no}
\def\defaulth@rdahlength{8pt}
\def\psarrowDD[#1,#2]{{\ifcurr@ntPS\ifps@cri\s@uvc@ntr@l\et@tpsarrow%
    \PSc@mment{psarrowDD [Pt1,Pt2]=[#1,#2]}\pssetfillmode{no}%
    \psarrowheadDD[#1,#2]\setc@ntr@l{2}\psline[#1,-3]%
    \PSc@mment{End psarrowDD}\resetc@ntr@l\et@tpsarrow\fi\fi}}
\def\psarrowTD[#1,#2]{{\ifcurr@ntPS\ifps@cri\s@uvc@ntr@l\et@tpsarrowTD%
    \PSc@mment{psarrowTD [Pt1,Pt2]=[#1,#2]}\resetc@ntr@l{2}%
    \Figptpr@j-5:/#1/\Figptpr@j-6:/#2/\let\c@lprojSP=\relax\psarrowDD[-5,-6]%
    \PSc@mment{End psarrowTD}\resetc@ntr@l\et@tpsarrowTD\fi\fi}}
\def\psarrowheadDD[#1,#2]{{\ifcurr@ntPS\ifps@cri\s@uvc@ntr@l\et@tpsarrowheadDD%
    \if@rrowhfill\def\@hangle{-\@rrowheadangle}\else\def\@hangle{\@rrowheadangle}\fi%
    \if@rrowratio%
    \if@rrowhout\def\@hratio{-\@rrowheadratio}\else\def\@hratio{\@rrowheadratio}\fi%
    \PSc@mment{psarrowheadDD Ratio=\@hratio, Angle=\@hangle, [Pt1,Pt2]=[#1,#2]}%
    \Ps@rrowhead\@hratio,\@hangle[#1,#2]%
    \else%
    \if@rrowhout\def\@hlength{-\@rrowheadlength}\else\def\@hlength{\@rrowheadlength}\fi%
    \PSc@mment{psarrowheadDD Length=\@hlength, Angle=\@hangle, [Pt1,Pt2]=[#1,#2]}%
    \Ps@rrowheadfd\@hlength,\@hangle[#1,#2]%
    \fi%
    \PSc@mment{End psarrowheadDD}\resetc@ntr@l\et@tpsarrowheadDD\fi\fi}}
\def\psarrowheadTD[#1,#2]{{\ifcurr@ntPS\ifps@cri\s@uvc@ntr@l\et@tpsarrowheadTD%
    \PSc@mment{psarrowheadTD [Pt1,Pt2]=[#1,#2]}\resetc@ntr@l{2}%
    \Figptpr@j-5:/#1/\Figptpr@j-6:/#2/\let\c@lprojSP=\relax\psarrowheadDD[-5,-6]%
    \PSc@mment{End psarrowheadTD}\resetc@ntr@l\et@tpsarrowheadTD\fi\fi}}
\def\Ps@rrowhead#1,#2[#3,#4]{\v@leur=#1\p@\maxim@m{\v@leur}{\v@leur}{-\v@leur}%
    \ifdim\v@leur>\Cepsil@n{
    \PSc@mment{ps@rrowhead Ratio=#1, Angle=#2, [Pt1,Pt2]=[#3,#4]}\v@leur=\UNSS@N%
    \v@leur=\curr@ntwidth\v@leur\v@leur=\ptpsT@pt\v@leur\delt@=.5\v@leur
    \setc@ntr@l{2}\figvectPDD-3[#4,#3]%
    \Figg@tXY{-3}\v@lX=#1\v@lX\v@lY=#1\v@lY\Figv@ctCreg-3(\v@lX,\v@lY)%
    \vecunit@{-4}{-3}\mili@u=\result@t%
    \ifdim#2pt>\z@\v@lXa=-\C@AHANG\delt@%
     \edef\c@ef{\repdecn@mb{\v@lXa}}\figpttraDD-3:=-3/\c@ef,-4/\fi%
    \edef\c@ef{\repdecn@mb{\delt@}}%
    \v@lXa=\mili@u\v@lXa=\C@AHANG\v@lXa%
    \v@lYa=\ptpsT@pt\p@\v@lYa=\curr@ntwidth\v@lYa\v@lYa=\sDcc@ngle\v@lYa%
    \advance\v@lXa-\v@lYa\gdef\sDcc@ngle{0}%
    \ifdim\v@lXa>\v@leur\edef\c@efendpt{\repdecn@mb{\v@leur}}%
    \else\edef\c@efendpt{\repdecn@mb{\v@lXa}}\fi%
    \Figg@tXY{-3}\v@lmin=\v@lX\v@lmax=\v@lY%
    \v@lXa=\C@AHANG\v@lmin\v@lYa=\S@AHANG\v@lmax\advance\v@lXa\v@lYa%
    \v@lYa=-\S@AHANG\v@lmin\v@lX=\C@AHANG\v@lmax\advance\v@lYa\v@lX%
    \setc@ntr@l{1}\Figg@tXY{#4}\advance\v@lX\v@lXa\advance\v@lY\v@lYa%
    \setc@ntr@l{2}\Figp@intregDD-2:(\v@lX,\v@lY)%
    \v@lXa=\C@AHANG\v@lmin\v@lYa=-\S@AHANG\v@lmax\advance\v@lXa\v@lYa%
    \v@lYa=\S@AHANG\v@lmin\v@lX=\C@AHANG\v@lmax\advance\v@lYa\v@lX%
    \setc@ntr@l{1}\Figg@tXY{#4}\advance\v@lX\v@lXa\advance\v@lY\v@lYa%
    \setc@ntr@l{2}\Figp@intregDD-1:(\v@lX,\v@lY)%
    \ifdim#2pt<\z@\fillm@detrue\psline[-2,#4,-1]
    \else\figptstraDD-3=#4,-2,-1/\c@ef,-4/\psline[-2,-3,-1]\fi
    \ifdim#1pt>\z@\figpttraDD-3:=#4/\c@efendpt,-4/\else\figptcopyDD-3:/#4/\fi%
    \PSc@mment{End ps@rrowhead}}\fi}
\def\sDcc@ngle{0}
\def\Ps@rrowheadfd#1,#2[#3,#4]{{%
    \PSc@mment{ps@rrowheadfd Length=#1, Angle=#2, [Pt1,Pt2]=[#3,#4]}%
    \setc@ntr@l{2}\figvectPDD-1[#3,#4]\n@rmeucDD{\v@leur}{-1}\v@leur=\ptT@unit@\v@leur%
    \invers@{\v@leur}{\v@leur}\v@leur=#1\v@leur\edef\R@tio{\repdecn@mb{\v@leur}}%
    \Ps@rrowhead\R@tio,#2[#3,#4]\PSc@mment{End ps@rrowheadfd}}}
\def\psarrowBezierDD[#1,#2,#3,#4]{{\ifcurr@ntPS\ifps@cri\s@uvc@ntr@l\et@tpsarrowBezierDD%
    \PSc@mment{psarrowBezierDD Control points=#1,#2,#3,#4}\setc@ntr@l{2}%
    \if@rrowratio\c@larclengthDD\v@leur,10[#1,#2,#3,#4]\else\v@leur=\z@\fi%
    \Ps@rrowB@zDD\v@leur[#1,#2,#3,#4]%
    \PSc@mment{End psarrowBezierDD}\resetc@ntr@l\et@tpsarrowBezierDD\fi\fi}}
\def\psarrowBezierTD[#1,#2,#3,#4]{{\ifcurr@ntPS\ifps@cri\s@uvc@ntr@l\et@tpsarrowBezierTD%
    \PSc@mment{psarrowBezierTD Control points=#1,#2,#3,#4}\resetc@ntr@l{2}%
    \Figptpr@j-7:/#1/\Figptpr@j-8:/#2/\Figptpr@j-9:/#3/\Figptpr@j-10:/#4/%
    \let\c@lprojSP=\relax\ifnum\curr@ntproj<\tw@\psarrowBezierDD[-7,-8,-9,-10]%
    \else\f@gnewpath\PSwrit@cmd{-7}{\c@mmoveto}{\fwf@g}%
    \if@rrowratio\c@larclengthDD\mili@u,10[-7,-8,-9,-10]\else\mili@u=\z@\fi%
    \p@rtent=\NBz@rcs\advance\p@rtent\m@ne\subB@zierTD\p@rtent[#1,#2,#3,#4]%
    \f@gstroke%
    \advance\v@lmin\p@rtent\delt@
    \v@leur=\v@lmin\advance\v@leur0.33333 \delt@\edef\unti@rs{\repdecn@mb{\v@leur}}%
    \v@leur=\v@lmin\advance\v@leur0.66666 \delt@\edef\deti@rs{\repdecn@mb{\v@leur}}%
    \figptcopyDD-8:/-10/\c@lsubBzarc\unti@rs,\deti@rs[#1,#2,#3,#4]%
    \figptcopyDD-8:/-4/\figptcopyDD-9:/-3/\Ps@rrowB@zDD\mili@u[-7,-8,-9,-10]\fi%
    \PSc@mment{End psarrowBezierTD}\resetc@ntr@l\et@tpsarrowBezierTD\fi\fi}}
\def\c@larclengthDD#1,#2[#3,#4,#5,#6]{{\p@rtent=#2\figptcopyDD-5:/#3/%
    \delt@=\p@\divide\delt@\p@rtent\c@rre=\z@\v@leur=\z@\s@mme=\z@%
    \loop\ifnum\s@mme<\p@rtent\advance\s@mme\@ne\advance\v@leur\delt@%
    \edef\T@{\repdecn@mb{\v@leur}}\figptBezierDD-6::\T@[#3,#4,#5,#6]%
    \figvectPDD-1[-5,-6]\n@rmeucDD{\mili@u}{-1}\advance\c@rre\mili@u%
    \figptcopyDD-5:/-6/\repeat\global\result@t=\ptT@unit@\c@rre}#1=\result@t}
\def\Ps@rrowB@zDD#1[#2,#3,#4,#5]{{\pssetfillmode{no}%
    \if@rrowratio\delt@=\@rrowheadratio#1\else\delt@=\@rrowheadlength pt\fi%
    \v@leur=\C@AHANG\delt@\edef\R@dius{\repdecn@mb{\v@leur}}%
    \FigptintercircB@zDD-5::0,\R@dius[#5,#4,#3,#2]%
    \pssetarrowheadlength{\repdecn@mb{\delt@}}\psarrowheadDD[-5,#5]%
    \let\n@rmeuc=\n@rmeucDD\figgetdist\R@dius[#5,-3]%
    \FigptintercircB@zDD-6::0,\R@dius[#5,#4,#3,#2]%
    \figptBezierDD-5::0.33333[#5,#4,#3,#2]\figptBezierDD-3::0.66666[#5,#4,#3,#2]%
    \figptscontrolDD-5[-6,-5,-3,#2]\psBezierDD1[-6,-5,-4,#2]}}
\def\psarrowcircDD#1;#2(#3,#4){{\ifcurr@ntPS\ifps@cri\s@uvc@ntr@l\et@tpsarrowcircDD%
    \PSc@mment{psarrowcircDD Center=#1 ; Radius=#2 (Ang1=#3,Ang2=#4)}%
    \pssetfillmode{no}\Pscirc@rrowhead#1;#2(#3,#4)%
    \setc@ntr@l{2}\figvectPDD -4[#1,-3]\vecunit@{-4}{-4}%
    \Figg@tXY{-4}\arct@n\v@lmin(\v@lX,\v@lY)%
    \v@lmin=\rdT@deg\v@lmin\v@leur=#4pt\advance\v@leur-\v@lmin%
    \maxim@m{\v@leur}{\v@leur}{-\v@leur}%
    \ifdim\v@leur>\DemiPI@deg\relax\ifdim\v@lmin<#4pt\advance\v@lmin\DePI@deg%
    \else\advance\v@lmin-\DePI@deg\fi\fi\edef\ar@ngle{\repdecn@mb{\v@lmin}}%
    \ifdim#3pt<#4pt\psarccirc#1;#2(#3,\ar@ngle)\else\psarccirc#1;#2(\ar@ngle,#3)\fi%
    \PSc@mment{End psarrowcircDD}\resetc@ntr@l\et@tpsarrowcircDD\fi\fi}}
\def\psarrowcircTD#1,#2,#3;#4(#5,#6){{\ifcurr@ntPS\ifps@cri\s@uvc@ntr@l\et@tpsarrowcircTD%
    \PSc@mment{psarrowcircTD Center=#1,P1=#2,P2=#3 ; Radius=#4 (Ang1=#5, Ang2=#6)}%
    \resetc@ntr@l{2}\c@lExtAxes#1,#2,#3(#4)\let\c@lprojSP=\relax%
    \figvectPTD-11[#1,-4]\figvectPTD-12[#1,-5]\c@lNbarcs{#5}{#6}%
    \if@rrowratio\v@lmax=\degT@rd\v@lmax\edef\D@lpha{\repdecn@mb{\v@lmax}}\fi%
    \advance\p@rtent\m@ne\mili@u=\z@%
    \v@leur=#5pt\c@lptellP{#1}{-11}{-12}\Figptpr@j-9:/-3/%
    \f@gnewpath\PSwrit@cmdS{-9}{\c@mmoveto}{\fwf@g}{\X@un}{\Y@un}%
    \edef\C@nt@r{#1}\s@mme=\z@\bcl@rcircTD\f@gstroke%
    \advance\v@leur\delt@\c@lptellP{#1}{-11}{-12}\Figptpr@j-5:/-3/%
    \advance\v@leur\delt@\c@lptellP{#1}{-11}{-12}\Figptpr@j-6:/-3/%
    \advance\v@leur\delt@\c@lptellP{#1}{-11}{-12}\Figptpr@j-10:/-3/%
    \figptscontrolDD-8[-9,-5,-6,-10]%
    \if@rrowratio\c@lcurvradDD0.5[-9,-8,-7,-10]\advance\mili@u\result@t%
    \maxim@m{\mili@u}{\mili@u}{-\mili@u}\mili@u=\ptT@unit@\mili@u%
    \mili@u=\D@lpha\mili@u\advance\p@rtent\@ne\divide\mili@u\p@rtent\fi%
    \Ps@rrowB@zDD\mili@u[-9,-8,-7,-10]%
    \PSc@mment{End psarrowcircTD}\resetc@ntr@l\et@tpsarrowcircTD\fi\fi}}
\def\bcl@rcircTD{\relax%
    \ifnum\s@mme<\p@rtent\advance\s@mme\@ne%
    \advance\v@leur\delt@\c@lptellP{\C@nt@r}{-11}{-12}\Figptpr@j-5:/-3/%
    \advance\v@leur\delt@\c@lptellP{\C@nt@r}{-11}{-12}\Figptpr@j-6:/-3/%
    \advance\v@leur\delt@\c@lptellP{\C@nt@r}{-11}{-12}\Figptpr@j-10:/-3/%
    \figptscontrolDD-8[-9,-5,-6,-10]\BdingB@xfalse%
    \PSwrit@cmdS{-8}{}{\fwf@g}{\X@de}{\Y@de}\PSwrit@cmdS{-7}{}{\fwf@g}{\X@tr}{\Y@tr}%
    \BdingB@xtrue\PSwrit@cmdS{-10}{\c@mcurveto}{\fwf@g}{\X@qu}{\Y@qu}%
    \if@rrowratio\c@lcurvradDD0.5[-9,-8,-7,-10]\advance\mili@u\result@t\fi%
    \B@zierBB@x{1}{\Y@un}(\X@un,\X@de,\X@tr,\X@qu)%
    \B@zierBB@x{2}{\X@un}(\Y@un,\Y@de,\Y@tr,\Y@qu)%
    \edef\X@un{\X@qu}\edef\Y@un{\Y@qu}\figptcopyDD-9:/-10/\bcl@rcircTD\fi}
\def\Pscirc@rrowhead#1;#2(#3,#4){{%
    \PSc@mment{pscirc@rrowhead Center=#1 ; Radius=#2 (Ang1=#3,Ang2=#4)}%
    \v@leur=#2\unit@\edef\s@glen{\repdecn@mb{\v@leur}}\v@lY=\z@\v@lX=\v@leur%
    \resetc@ntr@l{2}\Figv@ctCreg-3(\v@lX,\v@lY)\figpttraDD-5:=#1/1,-3/%
    \figptrotDD-5:=-5/#1,#4/%
    \figvectPDD-3[#1,-5]\Figg@tXY{-3}\v@leur=\v@lX%
    \ifdim#3pt<#4pt\v@lX=\v@lY\v@lY=-\v@leur\else\v@lX=-\v@lY\v@lY=\v@leur\fi%
    \Figv@ctCreg-3(\v@lX,\v@lY)\vecunit@{-3}{-3}%
    \if@rrowratio\v@leur=#4pt\advance\v@leur-#3pt\maxim@m{\mili@u}{-\v@leur}{\v@leur}%
    \mili@u=\degT@rd\mili@u\v@leur=\s@glen\mili@u\edef\s@glen{\repdecn@mb{\v@leur}}%
    \mili@u=#2\mili@u\mili@u=\@rrowheadratio\mili@u\else\mili@u=\@rrowheadlength pt\fi%
    \figpttraDD-6:=-5/\s@glen,-3/\v@leur=#2pt\v@leur=2\v@leur%
    \invers@{\v@leur}{\v@leur}\c@rre=\repdecn@mb{\v@leur}\mili@u
    \mili@u=\c@rre\mili@u=\repdecn@mb{\c@rre}\mili@u%
    \v@leur=\p@\advance\v@leur-\mili@u
    \invers@{\mili@u}{2\v@leur}\delt@=\c@rre\delt@=\repdecn@mb{\mili@u}\delt@%
    \xdef\sDcc@ngle{\repdecn@mb{\delt@}}
    \sqrt@{\mili@u}{\v@leur}\arct@n\v@leur(\mili@u,\c@rre)%
    \v@leur=\rdT@deg\v@leur
    \ifdim#3pt<#4pt\v@leur=-\v@leur\fi%
    \if@rrowhout\v@leur=-\v@leur\fi\edef\cor@ngle{\repdecn@mb{\v@leur}}%
    \figptrotDD-6:=-6/-5,\cor@ngle/\psarrowheadDD[-6,-5]%
    \PSc@mment{End pscirc@rrowhead}}}
\def\psarrowcircPDD#1;#2[#3,#4]{{\ifcurr@ntPS\ifps@cri%
    \PSc@mment{psarrowcircPDD Center=#1; Radius=#2, [P1=#3,P2=#4]}%
    \s@uvc@ntr@l\et@tpsarrowcircPDD\Ps@ngleparam#1;#2[#3,#4]%
    \ifdim\v@leur>\z@\ifdim\v@lmin>\v@lmax\advance\v@lmax\DePI@deg\fi%
    \else\ifdim\v@lmin<\v@lmax\advance\v@lmin\DePI@deg\fi\fi%
    \edef\@ngdeb{\repdecn@mb{\v@lmin}}\edef\@ngfin{\repdecn@mb{\v@lmax}}%
    \psarrowcirc#1;\r@dius(\@ngdeb,\@ngfin)%
    \PSc@mment{End psarrowcircPDD}\resetc@ntr@l\et@tpsarrowcircPDD\fi\fi}}
\def\psarrowcircPTD#1;#2[#3,#4,#5]{{\ifcurr@ntPS\ifps@cri\s@uvc@ntr@l\et@tpsarrowcircPTD%
    \PSc@mment{psarrowcircPTD Center=#1; Radius=#2, [P1=#3,P2=#4,P3=#5]}%
    \figgetangleTD\@ngfin[#1,#3,#4,#5]\v@leur=#2pt%
    \maxim@m{\mili@u}{-\v@leur}{\v@leur}\edef\r@dius{\repdecn@mb{\mili@u}}%
    \ifdim\v@leur<\z@\v@lmax=\@ngfin pt\advance\v@lmax-\DePI@deg%
    \edef\@ngfin{\repdecn@mb{\v@lmax}}\fi\psarrowcircTD#1,#3,#5;\r@dius(0,\@ngfin)%
    \PSc@mment{End psarrowcircPTD}\resetc@ntr@l\et@tpsarrowcircPTD\fi\fi}}
\def\psaxes#1(#2){{\ifcurr@ntPS\ifps@cri\s@uvc@ntr@l\et@tpsaxes%
    \PSc@mment{psaxes Origin=#1 Range=(#2)}\an@lys@xes#2,:\resetc@ntr@l{2}%
    \ifx\t@xt@\empty\ifTr@isDim\ps@xes#1(0,#2,0,#2,0,#2)\else\ps@xes#1(0,#2,0,#2)\fi%
    \else\ps@xes#1(#2)\fi\PSc@mment{End psaxes}\resetc@ntr@l\et@tpsaxes\fi\fi}}
\def\an@lys@xes#1,#2:{\def\t@xt@{#2}}
\def\ps@xesDD#1(#2,#3,#4,#5){%
    \figpttraC-5:=#1/#2,0/\figpttraC-6:=#1/#3,0/\psarrowDD[-5,-6]%
    \figpttraC-5:=#1/0,#4/\figpttraC-6:=#1/0,#5/\psarrowDD[-5,-6]}
\def\ps@xesTD#1(#2,#3,#4,#5,#6,#7){%
    \figpttraC-7:=#1/#2,0,0/\figpttraC-8:=#1/#3,0,0/\psarrowTD[-7,-8]%
    \figpttraC-7:=#1/0,#4,0/\figpttraC-8:=#1/0,#5,0/\psarrowTD[-7,-8]%
    \figpttraC-7:=#1/0,0,#6/\figpttraC-8:=#1/0,0,#7/\psarrowTD[-7,-8]}
\edef\DefGIfilen@me{\jobname GI.anx}
\def\psbeginfig#1{\def\t@xt@{#1}\relax\ifx\t@xt@\empty\Psb@ginfig\DefGIfilen@me%
    \else\expandafter\Psb@ginfigNu@#1 :\fi}
\def\Psb@ginfigNu@#1 #2:{\def\t@xt@{#1}\relax\ifx\t@xt@\empty\def\t@xt@{#2}%
    \ifx\t@xt@\empty\Psb@ginfig\DefGIfilen@me\else\Psb@ginfigNu@#2:\fi%
    \else\Psb@ginfig{#1}\fi}
\def\Psb@ginfig#1{\ifcurr@ntPS\else%
    \edef\PSfilen@me{#1}\edef\auxfilen@me{\jobname.anx}%
    \ifpstestm@de\ps@critrue\else\openin\frf@g=\PSfilen@me\relax%
    \ifeof\frf@g\ps@critrue\else\ps@crifalse\fi\closein\frf@g\fi%
    \curr@ntPStrue\c@ldefproj\expandafter\setupd@te\defaultupdate:%
    \ifps@cri\initb@undb@x%
    \immediate\openout\fwf@g=\auxfilen@me\initpss@ttings\fi%
    \fi}
\def\initpss@ttings{\psreset{arrowhead,curve,first,flowchart,mesh,second,third}%
    \Use@llipsefalse}
\def\B@zierBB@x#1#2(#3,#4,#5,#6){{\c@rre=\t@n\epsil@n
    \v@lmax=#4\advance\v@lmax-#5\v@lmax=\thr@@\v@lmax\advance\v@lmax#6\advance\v@lmax-#3%
    \mili@u=#4\mili@u=-\tw@\mili@u\advance\mili@u#3\advance\mili@u#5%
    \v@lmin=#4\advance\v@lmin-#3\maxim@m{\v@leur}{-\v@lmax}{\v@lmax}%
    \maxim@m{\delt@}{-\mili@u}{\mili@u}\maxim@m{\v@leur}{\v@leur}{\delt@}%
    \maxim@m{\delt@}{-\v@lmin}{\v@lmin}\maxim@m{\v@leur}{\v@leur}{\delt@}%
    \ifdim\v@leur>\c@rre\invers@{\v@leur}{\v@leur}\edef\Uns@rM@x{\repdecn@mb{\v@leur}}%
    \v@lmax=\Uns@rM@x\v@lmax\mili@u=\Uns@rM@x\mili@u\v@lmin=\Uns@rM@x\v@lmin%
    \maxim@m{\v@leur}{-\v@lmax}{\v@lmax}\ifdim\v@leur<\c@rre%
    \maxim@m{\v@leur}{-\mili@u}{\mili@u}\ifdim\v@leur<\c@rre\else%
    \invers@{\mili@u}{\mili@u}\v@leur=-0.5\v@lmin%
    \v@leur=\repdecn@mb{\mili@u}\v@leur\m@jBBB@x{\v@leur}{#1}{#2}(#3,#4,#5,#6)\fi%
    \else\delt@=\repdecn@mb{\mili@u}\mili@u\v@leur=\repdecn@mb{\v@lmax}\v@lmin%
    \advance\delt@-\v@leur\ifdim\delt@<\z@\else\invers@{\v@lmax}{\v@lmax}%
    \edef\Uns@rAp{\repdecn@mb{\v@lmax}}\sqrt@{\delt@}{\delt@}%
    \v@leur=-\mili@u\advance\v@leur\delt@\v@leur=\Uns@rAp\v@leur%
    \m@jBBB@x{\v@leur}{#1}{#2}(#3,#4,#5,#6)%
    \v@leur=-\mili@u\advance\v@leur-\delt@\v@leur=\Uns@rAp\v@leur%
    \m@jBBB@x{\v@leur}{#1}{#2}(#3,#4,#5,#6)\fi\fi\fi}}
\def\m@jBBB@x#1#2#3(#4,#5,#6,#7){{\relax\ifdim#1>\z@\ifdim#1<\p@%
    \edef\T@{\repdecn@mb{#1}}\v@lX=\p@\advance\v@lX-#1\edef\UNmT@{\repdecn@mb{\v@lX}}%
    \v@lX=#4\v@lY=#5\v@lZ=#6\v@lXa=#7\v@lX=\UNmT@\v@lX\advance\v@lX\T@\v@lY%
    \v@lY=\UNmT@\v@lY\advance\v@lY\T@\v@lZ\v@lZ=\UNmT@\v@lZ\advance\v@lZ\T@\v@lXa%
    \v@lX=\UNmT@\v@lX\advance\v@lX\T@\v@lY\v@lY=\UNmT@\v@lY\advance\v@lY\T@\v@lZ%
    \v@lX=\UNmT@\v@lX\advance\v@lX\T@\v@lY%
    \ifcase#2\or\v@lY=#3\or\v@lY=\v@lX\v@lX=#3\fi\b@undb@x{\v@lX}{\v@lY}\fi\fi}}
\def\PsB@zier#1[#2]{{\f@gnewpath%
    \s@mme=\z@\def\list@num{#2,0}\extrairelepremi@r\p@int\de\list@num%
    \PSwrit@cmdS{\p@int}{\c@mmoveto}{\fwf@g}{\X@un}{\Y@un}\p@rtent=#1\bclB@zier}}
\def\bclB@zier{\relax%
    \ifnum\s@mme<\p@rtent\advance\s@mme\@ne\BdingB@xfalse%
    \extrairelepremi@r\p@int\de\list@num\PSwrit@cmdS{\p@int}{}{\fwf@g}{\X@de}{\Y@de}%
    \extrairelepremi@r\p@int\de\list@num\PSwrit@cmdS{\p@int}{}{\fwf@g}{\X@tr}{\Y@tr}%
    \BdingB@xtrue%
    \extrairelepremi@r\p@int\de\list@num\PSwrit@cmdS{\p@int}{\c@mcurveto}{\fwf@g}{\X@qu}{\Y@qu}%
    \B@zierBB@x{1}{\Y@un}(\X@un,\X@de,\X@tr,\X@qu)%
    \B@zierBB@x{2}{\X@un}(\Y@un,\Y@de,\Y@tr,\Y@qu)%
    \edef\X@un{\X@qu}\edef\Y@un{\Y@qu}\bclB@zier\fi}
\def\psBezierDD#1[#2]{\ifcurr@ntPS\ifps@cri%
    \PSc@mment{psBezierDD N arcs=#1, Control points=#2}%
    \iffillm@de\PsB@zier#1[#2]%
    \f@gfill%
    \else\PsB@zier#1[#2]\f@gstroke\fi%
    \PSc@mment{End psBezierDD}\fi\fi}
\def\psBezierTD#1[#2]{\ifcurr@ntPS\ifps@cri\s@uvc@ntr@l\et@tpsBezierTD%
    \PSc@mment{psBezierTD N arcs=#1, Control points=#2}%
    \iffillm@de\PsB@zierTD#1[#2]%
    \f@gfill%
    \else\PsB@zierTD#1[#2]\f@gstroke\fi%
    \PSc@mment{End psBezierTD}\resetc@ntr@l\et@tpsBezierTD\fi\fi}
\def\PsB@zierTD#1[#2]{\ifnum\curr@ntproj<\tw@\PsB@zier#1[#2]\else\PsB@zier@TD#1[#2]\fi}
\def\PsB@zier@TD#1[#2]{{\f@gnewpath%
    \s@mme=\z@\def\list@num{#2,0}\extrairelepremi@r\p@int\de\list@num%
    \let\c@lprojSP=\relax\setc@ntr@l{2}\Figptpr@j-7:/\p@int/%
    \PSwrit@cmd{-7}{\c@mmoveto}{\fwf@g}%
    \loop\ifnum\s@mme<#1\advance\s@mme\@ne\extrairelepremi@r\p@intun\de\list@num%
    \extrairelepremi@r\p@intde\de\list@num\extrairelepremi@r\p@inttr\de\list@num%
    \subB@zierTD\NBz@rcs[\p@int,\p@intun,\p@intde,\p@inttr]\edef\p@int{\p@inttr}\repeat}}
\def\subB@zierTD#1[#2,#3,#4,#5]{\delt@=\p@\divide\delt@\NBz@rcs\v@lmin=\z@%
    {\Figg@tXY{-7}\edef\X@un{\the\v@lX}\edef\Y@un{\the\v@lY}%
    \s@mme=\z@\loop\ifnum\s@mme<#1\advance\s@mme\@ne%
    \v@leur=\v@lmin\advance\v@leur0.33333 \delt@\edef\unti@rs{\repdecn@mb{\v@leur}}%
    \v@leur=\v@lmin\advance\v@leur0.66666 \delt@\edef\deti@rs{\repdecn@mb{\v@leur}}%
    \advance\v@lmin\delt@\edef\trti@rs{\repdecn@mb{\v@lmin}}%
    \figptBezierTD-8::\trti@rs[#2,#3,#4,#5]\Figptpr@j-8:/-8/%
    \c@lsubBzarc\unti@rs,\deti@rs[#2,#3,#4,#5]\BdingB@xfalse%
    \PSwrit@cmdS{-4}{}{\fwf@g}{\X@de}{\Y@de}\PSwrit@cmdS{-3}{}{\fwf@g}{\X@tr}{\Y@tr}%
    \BdingB@xtrue\PSwrit@cmdS{-8}{\c@mcurveto}{\fwf@g}{\X@qu}{\Y@qu}%
    \B@zierBB@x{1}{\Y@un}(\X@un,\X@de,\X@tr,\X@qu)%
    \B@zierBB@x{2}{\X@un}(\Y@un,\Y@de,\Y@tr,\Y@qu)%
    \edef\X@un{\X@qu}\edef\Y@un{\Y@qu}\figptcopyDD-7:/-8/\repeat}}
\def\NBz@rcs{2}
\def\c@lsubBzarc#1,#2[#3,#4,#5,#6]{\figptBezierTD-5::#1[#3,#4,#5,#6]%
    \figptBezierTD-6::#2[#3,#4,#5,#6]\Figptpr@j-4:/-5/\Figptpr@j-5:/-6/%
    \figptscontrolDD-4[-7,-4,-5,-8]}
\def\pscircDD#1(#2){\ifcurr@ntPS\ifps@cri\PSc@mment{pscircDD Center=#1 (Radius=#2)}%
    \psarccircDD#1;#2(0,360)\PSc@mment{End pscircDD}\fi\fi}
\def\pscircTD#1,#2,#3(#4){\ifcurr@ntPS\ifps@cri%
    \PSc@mment{pscircTD Center=#1,P1=#2,P2=#3 (Radius=#4)}%
    \psarccircTD#1,#2,#3;#4(0,360)\PSc@mment{End pscircTD}\fi\fi}
{\catcode`\%=12\gdef\p@urcent{
\def\PSc@mment#1{\ifpsdebugmode\immediate\write\fwf@g{\p@urcent\space#1}\fi}
{\catcode`\[=1\catcode`\{=12\gdef\acc@louv[{}}
{\catcode`\]=2\catcode`\}=12\gdef\acc@lfer{}]]
\def\PSdict@{\ifUse@llipse%
    \immediate\write\fwf@g{/ellipsedict 9 dict def ellipsedict /mtrx matrix put}%
    \immediate\write\fwf@g{/ellipse \acc@louv ellipsedict begin}%
    \immediate\write\fwf@g{ /endangle exch def /startangle exch def}%
    \immediate\write\fwf@g{ /yrad exch def /xrad exch def}%
    \immediate\write\fwf@g{ /rotangle exch def /y exch def /x exch def}%
    \immediate\write\fwf@g{ /savematrix mtrx currentmatrix def}%
    \immediate\write\fwf@g{ x y translate rotangle rotate xrad yrad scale}%
    \immediate\write\fwf@g{ 0 0 1 startangle endangle arc}%
    \immediate\write\fwf@g{ savematrix setmatrix end\acc@lfer def}%
    \fi\PShe@der{EndProlog}}
\def\Pssetc@rve#1=#2|{\keln@mun#1|%
    \def\n@mref{r}\ifx\l@debut\n@mref\pssetroundness{#2}\else
    \immediate\write16{*** Unknown attribute: \BS@ psset curve(..., #1=...)}%
    \fi}
\def\pssetroundness#1{\edef\curv@roundness{#1}}
\def\defaultroundness{0.2} 
\def\pscurveDD[#1]{{\ifcurr@ntPS\ifps@cri\PSc@mment{pscurveDD Points=#1}%
    \s@uvc@ntr@l\et@tpscurveDD%
    \iffillm@de\Psc@rveDD\curv@roundness[#1]%
    \f@gfill%
    \else\Psc@rveDD\curv@roundness[#1]\f@gstroke\fi%
    \PSc@mment{End pscurveDD}\resetc@ntr@l\et@tpscurveDD\fi\fi}}
\def\pscurveTD[#1]{{\ifcurr@ntPS\ifps@cri%
    \PSc@mment{pscurveTD Points=#1}\s@uvc@ntr@l\et@tpscurveTD\let\c@lprojSP=\relax%
    \iffillm@de\Psc@rveTD\curv@roundness[#1]%
    \f@gfill%
    \else\Psc@rveTD\curv@roundness[#1]\f@gstroke\fi%
    \PSc@mment{End pscurveTD}\resetc@ntr@l\et@tpscurveTD\fi\fi}}
\def\Psc@rveDD#1[#2]{%
    \def\list@num{#2}\extrairelepremi@r\Ak@\de\list@num%
    \extrairelepremi@r\Ai@\de\list@num\extrairelepremi@r\Aj@\de\list@num%
    \f@gnewpath\PSwrit@cmdS{\Ai@}{\c@mmoveto}{\fwf@g}{\X@un}{\Y@un}%
    \setc@ntr@l{2}\figvectPDD -1[\Ak@,\Aj@]%
    \@ecfor\Ak@:=\list@num\do{\figpttraDD-2:=\Ai@/#1,-1/\BdingB@xfalse%
       \PSwrit@cmdS{-2}{}{\fwf@g}{\X@de}{\Y@de}%
       \figvectPDD -1[\Ai@,\Ak@]\figpttraDD-2:=\Aj@/-#1,-1/%
       \PSwrit@cmdS{-2}{}{\fwf@g}{\X@tr}{\Y@tr}\BdingB@xtrue%
       \PSwrit@cmdS{\Aj@}{\c@mcurveto}{\fwf@g}{\X@qu}{\Y@qu}%
       \B@zierBB@x{1}{\Y@un}(\X@un,\X@de,\X@tr,\X@qu)%
       \B@zierBB@x{2}{\X@un}(\Y@un,\Y@de,\Y@tr,\Y@qu)%
       \edef\X@un{\X@qu}\edef\Y@un{\Y@qu}\edef\Ai@{\Aj@}\edef\Aj@{\Ak@}}}
\def\Psc@rveTD#1[#2]{\ifnum\curr@ntproj<\tw@\Psc@rvePPTD#1[#2]\else\Psc@rveCPTD#1[#2]\fi}
\def\Psc@rvePPTD#1[#2]{\setc@ntr@l{2}%
    \def\list@num{#2}\extrairelepremi@r\Ak@\de\list@num\Figptpr@j-5:/\Ak@/%
    \extrairelepremi@r\Ai@\de\list@num\Figptpr@j-3:/\Ai@/%
    \extrairelepremi@r\Aj@\de\list@num\Figptpr@j-4:/\Aj@/%
    \f@gnewpath\PSwrit@cmdS{-3}{\c@mmoveto}{\fwf@g}{\X@un}{\Y@un}%
    \figvectPDD -1[-5,-4]%
    \@ecfor\Ak@:=\list@num\do{\Figptpr@j-5:/\Ak@/\figpttraDD-2:=-3/#1,-1/%
       \BdingB@xfalse\PSwrit@cmdS{-2}{}{\fwf@g}{\X@de}{\Y@de}%
       \figvectPDD -1[-3,-5]\figpttraDD-2:=-4/-#1,-1/%
       \PSwrit@cmdS{-2}{}{\fwf@g}{\X@tr}{\Y@tr}\BdingB@xtrue%
       \PSwrit@cmdS{-4}{\c@mcurveto}{\fwf@g}{\X@qu}{\Y@qu}%
       \B@zierBB@x{1}{\Y@un}(\X@un,\X@de,\X@tr,\X@qu)%
       \B@zierBB@x{2}{\X@un}(\Y@un,\Y@de,\Y@tr,\Y@qu)%
       \edef\X@un{\X@qu}\edef\Y@un{\Y@qu}\figptcopyDD-3:/-4/\figptcopyDD-4:/-5/}}
\def\Psc@rveCPTD#1[#2]{\setc@ntr@l{2}%
    \def\list@num{#2}\extrairelepremi@r\Ak@\de\list@num%
    \extrairelepremi@r\Ai@\de\list@num\extrairelepremi@r\Aj@\de\list@num%
    \Figptpr@j-7:/\Ai@/%
    \f@gnewpath\PSwrit@cmd{-7}{\c@mmoveto}{\fwf@g}%
    \figvectPTD -9[\Ak@,\Aj@]%
    \@ecfor\Ak@:=\list@num\do{\figpttraTD-10:=\Ai@/#1,-9/%
       \figvectPTD -9[\Ai@,\Ak@]\figpttraTD-11:=\Aj@/-#1,-9/%
       \subB@zierTD\NBz@rcs[\Ai@,-10,-11,\Aj@]\edef\Ai@{\Aj@}\edef\Aj@{\Ak@}}}
\def\psendfig{\ifcurr@ntPS\ifps@cri\immediate\closeout\fwf@g%
    \immediate\openout\fwf@g=\PSfilen@me\relax%
    \ifPDFm@ke\PSBdingB@x\else%
    \immediate\write\fwf@g{\p@urcent\string!PS-Adobe-2.0 EPSF-2.0}%
    \PShe@der{Creator\string: TeX (fig4tex.tex)}%
    \PShe@der{Title\string: \PSfilen@me}%
    \PShe@der{CreationDate\string: \the\day/\the\month/\the\year}%
    \PSBdingB@x%
    \PShe@der{EndComments}\PSdict@\fi%
    \immediate\write\fwf@g{\c@mgsave}%
    \openin\frf@g=\auxfilen@me\c@pypsfile\fwf@g\frf@g\closein\frf@g%
    \immediate\write\fwf@g{\c@mgrestore}%
    \PSc@mment{End of file.}\immediate\closeout\fwf@g%
    \immediate\openout\fwf@g=\auxfilen@me\immediate\closeout\fwf@g%
    \immediate\write16{File \PSfilen@me\space created.}\fi\fi\curr@ntPSfalse\ps@critrue}
\def\PShe@der#1{\immediate\write\fwf@g{\p@urcent\p@urcent#1}}
\def\PSBdingB@x{{\v@lX=\ptT@ptps\c@@rdXmin\v@lY=\ptT@ptps\c@@rdYmin%
     \v@lXa=\ptT@ptps\c@@rdXmax\v@lYa=\ptT@ptps\c@@rdYmax%
     \PShe@der{BoundingBox\string: \repdecn@mb{\v@lX}\space\repdecn@mb{\v@lY}%
     \space\repdecn@mb{\v@lXa}\space\repdecn@mb{\v@lYa}}}}
\def\psfcconnect[#1]{{\ifcurr@ntPS\ifps@cri\PSc@mment{psfcconnect Points=#1}%
    \pssetfillmode{no}\s@uvc@ntr@l\et@tpsfcconnect\resetc@ntr@l{2}%
    \fcc@nnect@[#1]\resetc@ntr@l\et@tpsfcconnect\PSc@mment{End psfcconnect}\fi\fi}}
\def\fcc@nnect@[#1]{\let\N@rm=\n@rmeucDD\def\list@num{#1}%
    \extrairelepremi@r\Ai@\de\list@num\edef\pr@m{\Ai@}\v@leur=\z@\p@rtent=\@ne\c@llgtot%
    \ifcase\fclin@typ@\edef\list@num{[\pr@m,#1,\Ai@}\expandafter\pscurve\list@num]%
    \else\ifdim\fclin@r@d\p@>\z@\Pslin@conge[#1]\else\psline[#1]\fi\fi%
    \v@leur=\@rrowp@s\v@leur\edef\list@num{#1,\Ai@,0}%
    \extrairelepremi@r\Ai@\de\list@num\mili@u=\epsil@n\c@llgpart%
    \advance\mili@u-\epsil@n\advance\mili@u-\delt@\advance\v@leur-\mili@u%
    \ifcase\fclin@typ@\invers@\mili@u\delt@%
    \ifnum\@rrowr@fpt>\z@\advance\delt@-\v@leur\v@leur=\delt@\fi%
    \v@leur=\repdecn@mb\v@leur\mili@u\edef\v@lt{\repdecn@mb\v@leur}%
    \extrairelepremi@r\Ak@\de\list@num%
    \figvectPDD-1[\pr@m,\Aj@]\figpttraDD-6:=\Ai@/\curv@roundness,-1/%
    \figvectPDD-1[\Ak@,\Ai@]\figpttraDD-7:=\Aj@/\curv@roundness,-1/%
    \delt@=\@rrowheadlength\p@\delt@=\C@AHANG\delt@\edef\R@dius{\repdecn@mb{\delt@}}%
    \ifcase\@rrowr@fpt%
    \FigptintercircB@zDD-8::\v@lt,\R@dius[\Ai@,-6,-7,\Aj@]\psarrowheadDD[-5,-8]\else%
    \FigptintercircB@zDD-8::\v@lt,\R@dius[\Aj@,-7,-6,\Ai@]\psarrowheadDD[-8,-5]\fi%
    \else\advance\delt@-\v@leur%
    \p@rtentiere{\p@rtent}{\delt@}\edef\C@efun{\the\p@rtent}%
    \p@rtentiere{\p@rtent}{\v@leur}\edef\C@efde{\the\p@rtent}%
    \figptbaryDD-5:[\Ai@,\Aj@;\C@efun,\C@efde]\ifcase\@rrowr@fpt%
    \delt@=\@rrowheadlength\unit@\delt@=\C@AHANG\delt@\edef\t@ille{\repdecn@mb{\delt@}}%
    \figvectPDD-2[\Ai@,\Aj@]\vecunit@{-2}{-2}\figpttraDD-5:=-5/\t@ille,-2/\fi%
    \psarrowheadDD[\Ai@,-5]\fi}
\def\c@llgtot{\@ecfor\Aj@:=\list@num\do{\figvectP-1[\Ai@,\Aj@]\N@rm\delt@{-1}%
    \advance\v@leur\delt@\advance\p@rtent\@ne\edef\Ai@{\Aj@}}}
\def\c@llgpart{\extrairelepremi@r\Aj@\de\list@num\figvectP-1[\Ai@,\Aj@]\N@rm\delt@{-1}%
    \advance\mili@u\delt@\ifdim\mili@u<\v@leur\edef\pr@m{\Ai@}\edef\Ai@{\Aj@}\c@llgpart\fi}
\def\Pslin@conge[#1]{\ifnum\p@rtent>\tw@{\def\list@num{#1}%
    \extrairelepremi@r\Ai@\de\list@num\extrairelepremi@r\Aj@\de\list@num%
    \figptcopy-6:/\Ai@/\figvectP-3[\Ai@,\Aj@]\vecunit@{-3}{-3}\v@lmax=\result@t%
    \@ecfor\Ak@:=\list@num\do{\figvectP-4[\Aj@,\Ak@]\vecunit@{-4}{-4}%
    \minim@m\v@lmin\v@lmax\result@t\v@lmax=\result@t%
    \det@rm\delt@[-3,-4]\maxim@m\mili@u{\delt@}{-\delt@}\ifdim\mili@u>\Cepsil@n%
    \ifdim\delt@>\z@\figgetangleDD\Angl@[\Aj@,\Ak@,\Ai@]\else%
    \figgetangleDD\Angl@[\Aj@,\Ai@,\Ak@]\fi%
    \v@leur=\PI@deg\advance\v@leur-\Angl@\p@\divide\v@leur\tw@%
    \edef\Angl@{\repdecn@mb\v@leur}\c@ssin{\C@}{\S@}{\Angl@}\v@leur=\fclin@r@d\unit@%
    \v@leur=\S@\v@leur\mili@u=\C@\p@\invers@\mili@u\mili@u%
    \v@leur=\repdecn@mb{\mili@u}\v@leur%
    \minim@m\v@leur\v@leur\v@lmin\edef\t@ille{\repdecn@mb{\v@leur}}%
    \figpttra-5:=\Aj@/-\t@ille,-3/\psline[-6,-5]\figpttra-6:=\Aj@/\t@ille,-4/%
    \figvectNVDD-3[-3]\figvectNVDD-8[-4]\inters@cDD-7:[-5,-3;-6,-8]%
    \ifdim\delt@>\z@\psarccircP-7;\fclin@r@d[-5,-6]\else\psarccircP-7;\fclin@r@d[-6,-5]\fi%
    \else\psline[-6,\Aj@]\figptcopy-6:/\Aj@/\fi
    \edef\Ai@{\Aj@}\edef\Aj@{\Ak@}\figptcopy-3:/-4/}\psline[-6,\Aj@]}\else\psline[#1]\fi}
\def\psfcnode[#1]#2{{\ifcurr@ntPS\ifps@cri\PSc@mment{psfcnode Points=#1}%
    \s@uvc@ntr@l\et@tpsfcnode\resetc@ntr@l{2}%
    \def\t@xt@{#2}\ifx\t@xt@\empty\def\g@tt@xt{\setbox\Gb@x=\hbox{\Figg@tT{\p@int}}}%
    \else\def\g@tt@xt{\setbox\Gb@x=\hbox{#2}}\fi%
    \v@lmin=\h@rdfcXp@dd\advance\v@lmin\Xp@dd\unit@\multiply\v@lmin\tw@%
    \v@lmax=\h@rdfcYp@dd\advance\v@lmax\Yp@dd\unit@\multiply\v@lmax\tw@%
    \Figv@ctCreg-8(\unit@,-\unit@)\def\list@num{#1}%
    \delt@=\curr@ntwidth bp\divide\delt@\tw@%
    \fcn@de\PSc@mment{End psfcnode}\resetc@ntr@l\et@tpsfcnode\fi\fi}}
\def\d@butn@de{\g@tt@xt\v@lX=\wd\Gb@x%
    \v@lY=\ht\Gb@x\advance\v@lY\dp\Gb@x\advance\v@lX\v@lmin\advance\v@lY\v@lmax}
\def\fcn@deE{%
    \@ecfor\p@int:=\list@num\do{\d@butn@de\v@lX=\unssqrttw@\v@lX\v@lY=\unssqrttw@\v@lY%
    \ifdim\thickn@ss\p@>\z@
    \v@lXa=\v@lX\advance\v@lXa\delt@\v@lXa=\ptT@unit@\v@lXa\edef\XR@d{\repdecn@mb\v@lXa}%
    \v@lYa=\v@lY\advance\v@lYa\delt@\v@lYa=\ptT@unit@\v@lYa\edef\YR@d{\repdecn@mb\v@lYa}%
    \arct@n\v@leur(\v@lXa,\v@lYa)\v@leur=\rdT@deg\v@leur\edef\@nglde{\repdecn@mb\v@leur}%
    {\c@lptellDD-2::\p@int;\XR@d,\YR@d(\@nglde)}
    \advance\v@leur-\PI@deg\edef\@nglun{\repdecn@mb\v@leur}%
    {\c@lptellDD-3::\p@int;\XR@d,\YR@d(\@nglun)}%
    \figptstra-6=-3,-2,\p@int/\thickn@ss,-8/\pssetfillmode{yes}\us@secondC@lor%
    \psline[-2,-3,-6,-5]\psarcell-4;\XR@d,\YR@d(\@nglun,\@nglde,0)\fi
    \v@lX=\ptT@unit@\v@lX\v@lY=\ptT@unit@\v@lY%
    \edef\XR@d{\repdecn@mb\v@lX}\edef\YR@d{\repdecn@mb\v@lY}%
    \pssetfillmode{yes}\us@thirdC@lor\psarcell\p@int;\XR@d,\YR@d(0,360,0)%
    \pssetfillmode{no}\us@primarC@lor\psarcell\p@int;\XR@d,\YR@d(0,360,0)}}
\def\fcn@deL{\delt@=\ptT@unit@\delt@\edef\t@ille{\repdecn@mb\delt@}%
    \@ecfor\p@int:=\list@num\do{\Figg@tXYa{\p@int}\d@butn@de%
    \ifdim\v@lX>\v@lY\itis@Ktrue\else\itis@Kfalse\fi%
    \advance\v@lXa-\v@lX\Figp@intreg-1:(\v@lXa,\v@lYa)%
    \advance\v@lXa\v@lX\advance\v@lYa-\v@lY\Figp@intreg-2:(\v@lXa,\v@lYa)%
    \advance\v@lXa\v@lX\advance\v@lYa\v@lY\Figp@intreg-3:(\v@lXa,\v@lYa)%
    \advance\v@lXa-\v@lX\advance\v@lYa\v@lY\Figp@intreg-4:(\v@lXa,\v@lYa)%
    \ifdim\thickn@ss\p@>\z@\Figg@tXYa{\p@int}\pssetfillmode{yes}\us@secondC@lor
    \c@lpt@xt{-1}{-4}\c@lpt@xt@\v@lXa\v@lYa\v@lX\v@lY\c@rre\delt@%
    \Figp@intregDD-9:(\v@lZ,\v@lYa)\Figp@intregDD-11:(\v@lZa,\v@lYa)%
    \c@lpt@xt{-4}{-3}\c@lpt@xt@\v@lYa\v@lXa\v@lY\v@lX\delt@\c@rre%
    \Figp@intregDD-12:(\v@lXa,\v@lZ)\Figp@intregDD-10:(\v@lXa,\v@lZa)%
    \ifitis@K\figptstra-7=-9,-10,-11/\thickn@ss,-8/\psline[-9,-11,-5,-6,-7]\else%
    \figptstra-7=-10,-11,-12/\thickn@ss,-8/\psline[-10,-12,-5,-6,-7]\fi\fi
    \pssetfillmode{yes}\us@thirdC@lor\psline[-1,-2,-3,-4]%
    \pssetfillmode{no}\us@primarC@lor\psline[-1,-2,-3,-4,-1]}}
\def\c@lpt@xt#1#2{\figvectN-7[#1,#2]\vecunit@{-7}{-7}\figpttra-5:=#1/\t@ille,-7/%
    \figvectP-7[#1,#2]\Figg@tXY{-7}\c@rre=\v@lX\delt@=\v@lY\Figg@tXY{-5}}
\def\c@lpt@xt@#1#2#3#4#5#6{\v@lZ=#6\invers@{\v@lZ}{\v@lZ}\v@leur=\repdecn@mb{#5}\v@lZ%
    \v@lZ=#2\advance\v@lZ-#4\mili@u=\repdecn@mb{\v@leur}\v@lZ%
    \v@lZ=#3\advance\v@lZ\mili@u\v@lZa=-\v@lZ\advance\v@lZa\tw@#1}
\def\fcn@deR{\@ecfor\p@int:=\list@num\do{\Figg@tXYa{\p@int}\d@butn@de%
    \advance\v@lXa-0.5\v@lX\advance\v@lYa-0.5\v@lY\Figp@intreg-1:(\v@lXa,\v@lYa)%
    \advance\v@lXa\v@lX\Figp@intreg-2:(\v@lXa,\v@lYa)%
    \advance\v@lYa\v@lY\Figp@intreg-3:(\v@lXa,\v@lYa)%
    \advance\v@lXa-\v@lX\Figp@intreg-4:(\v@lXa,\v@lYa)%
    \ifdim\thickn@ss\p@>\z@\pssetfillmode{yes}\us@secondC@lor
    \Figv@ctCreg-5(-\delt@,-\delt@)\figpttra-9:=-1/1,-5/%
    \Figv@ctCreg-5(\delt@,-\delt@)\figpttra-10:=-2/1,-5/%
    \Figv@ctCreg-5(\delt@,\delt@)\figpttra-11:=-3/1,-5/%
    \figptstra-7=-9,-10,-11/\thickn@ss,-8/\psline[-9,-11,-5,-6,-7]\fi
    \pssetfillmode{yes}\us@thirdC@lor\psline[-1,-2,-3,-4]%
    \pssetfillmode{no}\us@primarC@lor\psline[-1,-2,-3,-4,-1]}}
\def\Pssetfl@wchart#1=#2|{\keln@mtr#1|%
    \def\n@mref{arr}\ifx\l@debut\n@mref\expandafter\keln@mtr\l@suite|%
     \def\n@mref{owp}\ifx\l@debut\n@mref\edef\@rrowp@s{#2}\else
     \def\n@mref{owr}\ifx\l@debut\n@mref\setfcr@fpt#2|\else
     \immediate\write16{*** Unknown attribute: \BS@ psset flowchart(..., #1=...)}%
     \fi\fi\else%
    \def\n@mref{lin}\ifx\l@debut\n@mref\setfccurv@#2|\else
    \def\n@mref{pad}\ifx\l@debut\n@mref\edef\Xp@dd{#2}\edef\Yp@dd{#2}\else
    \def\n@mref{rad}\ifx\l@debut\n@mref\edef\fclin@r@d{#2}\else
    \def\n@mref{sha}\ifx\l@debut\n@mref\setfcshap@#2|\else
    \def\n@mref{thi}\ifx\l@debut\n@mref\edef\thickn@ss{#2}\else
    \def\n@mref{xpa}\ifx\l@debut\n@mref\edef\Xp@dd{#2}\else
    \def\n@mref{ypa}\ifx\l@debut\n@mref\edef\Yp@dd{#2}\else
    \immediate\write16{*** Unknown attribute: \BS@ psset flowchart(..., #1=...)}%
    \fi\fi\fi\fi\fi\fi\fi\fi}
\def\setfcr@fpt#1#2|{\if#1e\def\@rrowr@fpt{1}\else\def\@rrowr@fpt{0}\fi}
\def\setfccurv@#1#2|{\if#1c\def\fclin@typ@{0}\else\def\fclin@typ@{1}\fi}
\def\setfcshap@#1#2|{%
    \if#1e\let\fcn@de=\fcn@deE\def\h@rdfcXp@dd{4pt}\def\h@rdfcYp@dd{4pt}%
     \edef\fcsh@pe{ellipse}\else%
    \if#1l\let\fcn@de=\fcn@deL\def\h@rdfcXp@dd{4pt}\def\h@rdfcYp@dd{4pt}%
     \edef\fcsh@pe{lozenge}\else%
          \let\fcn@de=\fcn@deR\def\h@rdfcXp@dd{6pt}\def\h@rdfcYp@dd{6pt}%
     \edef\fcsh@pe{rectangle}\fi\fi}
\def\psline[#1]{{\ifcurr@ntPS\ifps@cri\PSc@mment{psline Points=#1}%
    \let\pslign@=\Pslign@P\Pslin@{#1}\PSc@mment{End psline}\fi\fi}}
\def\pslineF#1{{\ifcurr@ntPS\ifps@cri\PSc@mment{pslineF Filename=#1}%
    \let\pslign@=\Pslign@F\Pslin@{#1}\PSc@mment{End pslineF}\fi\fi}}
\def\pslineC(#1){{\ifcurr@ntPS\ifps@cri\PSc@mment{pslineC}%
    \let\pslign@=\Pslign@C\Pslin@{#1}\PSc@mment{End pslineC}\fi\fi}}
\def\Pslin@#1{\iffillm@de\pslign@{#1}%
    \f@gfill%
    \else\pslign@{#1}\ifx\derp@int\premp@int%
    \f@gclosestroke%
    \else\f@gstroke\fi\fi}
\def\Pslign@P#1{\def\list@num{#1}\extrairelepremi@r\p@int\de\list@num%
    \edef\premp@int{\p@int}\f@gnewpath%
    \PSwrit@cmd{\p@int}{\c@mmoveto}{\fwf@g}%
    \@ecfor\p@int:=\list@num\do{\PSwrit@cmd{\p@int}{\c@mlineto}{\fwf@g}%
    \edef\derp@int{\p@int}}}
\def\Pslign@F#1{\s@uvc@ntr@l\et@tPslign@F\setc@ntr@l{2}\openin\frf@g=#1\relax%
    \ifeof\frf@g\message{*** File #1 not found !}\end\else%
    \read\frf@g to\tr@c\edef\premp@int{\tr@c}\expandafter\extr@ctCF\tr@c:%
    \f@gnewpath\PSwrit@cmd{-1}{\c@mmoveto}{\fwf@g}%
    \loop\read\frf@g to\tr@c\ifeof\frf@g\mored@tafalse\else\mored@tatrue\fi%
    \ifmored@ta\expandafter\extr@ctCF\tr@c:\PSwrit@cmd{-1}{\c@mlineto}{\fwf@g}%
    \edef\derp@int{\tr@c}\repeat\fi\closein\frf@g\resetc@ntr@l\et@tPslign@F}
\def\extr@ctCFDD#1 #2:{\v@lX=#1\unit@\v@lY=#2\unit@\Figp@intregDD-1:(\v@lX,\v@lY)}
\def\extr@ctCFTD#1 #2 #3:{\v@lX=#1\unit@\v@lY=#2\unit@\v@lZ=#3\unit@%
    \Figp@intregTD-1:(\v@lX,\v@lY,\v@lZ)}
\def\Pslign@C#1{\s@uvc@ntr@l\et@tPslign@C\setc@ntr@l{2}%
    \def\list@num{#1}\extrairelepremi@r\p@int\de\list@num%
    \edef\premp@int{\p@int}\f@gnewpath%
    \expandafter\Pslign@C@\p@int:\PSwrit@cmd{-1}{\c@mmoveto}{\fwf@g}%
    \@ecfor\p@int:=\list@num\do{\expandafter\Pslign@C@\p@int:%
    \PSwrit@cmd{-1}{\c@mlineto}{\fwf@g}\edef\derp@int{\p@int}}%
    \resetc@ntr@l\et@tPslign@C}
\def\Pslign@C@#1 #2:{{\def\t@xt@{#1}\ifx\t@xt@\empty\Pslign@C@#2:
    \else\extr@ctCF#1 #2:\fi}}
\newcount\c@ntrolmesh
\def\Pssetm@sh#1=#2|{\keln@mun#1|%
    \def\n@mref{d}\ifx\l@debut\n@mref\pssetmeshdiag{#2}\else
    \immediate\write16{*** Unknown attribute: \BS@ psset mesh(..., #1=...)}%
    \fi}
\def\pssetmeshdiag#1{\c@ntrolmesh=#1}
\def\defaultmeshdiag{0}    
\def\psmesh#1,#2[#3,#4,#5,#6]{{\ifcurr@ntPS\ifps@cri%
    \PSc@mment{psmesh N1=#1, N2=#2, Quadrangle=[#3,#4,#5,#6]}%
    \s@uvc@ntr@l\et@tpsmesh\Pss@tsecondSt\setc@ntr@l{2}%
    \ifnum#1>\@ne\Psmeshp@rt#1[#3,#4,#5,#6]\fi%
    \ifnum#2>\@ne\Psmeshp@rt#2[#4,#5,#6,#3]\fi%
    \ifnum\c@ntrolmesh>\z@\Psmeshdi@g#1,#2[#3,#4,#5,#6]\fi%
    \ifnum\c@ntrolmesh<\z@\Psmeshdi@g#2,#1[#4,#5,#6,#3]\fi\Psrest@reSt%
    \psline[#3,#4,#5,#6,#3]\PSc@mment{End psmesh}\resetc@ntr@l\et@tpsmesh\fi\fi}}
\def\Psmeshp@rt#1[#2,#3,#4,#5]{{\l@mbd@un=\@ne\l@mbd@de=#1\loop%
    \ifnum\l@mbd@un<#1\advance\l@mbd@de\m@ne\figptbary-1:[#2,#3;\l@mbd@de,\l@mbd@un]%
    \figptbary-2:[#5,#4;\l@mbd@de,\l@mbd@un]\psline[-1,-2]\advance\l@mbd@un\@ne\repeat}}
\def\Psmeshdi@g#1,#2[#3,#4,#5,#6]{\figptcopy-2:/#3/\figptcopy-3:/#6/%
    \l@mbd@un=\z@\l@mbd@de=#1\loop\ifnum\l@mbd@un<#1%
    \advance\l@mbd@un\@ne\advance\l@mbd@de\m@ne\figptcopy-1:/-2/\figptcopy-4:/-3/%
    \figptbary-2:[#3,#4;\l@mbd@de,\l@mbd@un]%
    \figptbary-3:[#6,#5;\l@mbd@de,\l@mbd@un]\Psmeshdi@gp@rt#2[-1,-2,-3,-4]\repeat}
\def\Psmeshdi@gp@rt#1[#2,#3,#4,#5]{{\l@mbd@un=\z@\l@mbd@de=#1\loop%
    \ifnum\l@mbd@un<#1\figptbary-5:[#2,#5;\l@mbd@de,\l@mbd@un]%
    \advance\l@mbd@de\m@ne\advance\l@mbd@un\@ne%
    \figptbary-6:[#3,#4;\l@mbd@de,\l@mbd@un]\psline[-5,-6]\repeat}}
\def\psnormalDD#1,#2[#3,#4]{{\ifcurr@ntPS\ifps@cri%
    \PSc@mment{psnormal Length=#1, Lambda=#2 [Pt1,Pt2]=[#3,#4]}%
    \s@uvc@ntr@l\et@tpsnormal\resetc@ntr@l{2}\figptendnormal-6::#1,#2[#3,#4]%
    \figptcopyDD-5:/-1/\psarrow[-5,-6]%
    \PSc@mment{End psnormal}\resetc@ntr@l\et@tpsnormal\fi\fi}}
\def\psreset#1{\trtlis@rg{#1}{\Psreset@}}
\def\Psreset@#1|{\keln@mde#1|%
    \def\n@mref{ar}\ifx\l@debut\n@mref\psresetarrowhead\else
    \def\n@mref{cu}\ifx\l@debut\n@mref\psset curve(roundness=\defaultroundness)\else
    \def\n@mref{fi}\ifx\l@debut\n@mref\psset (color=\defaultcolor,dash=\defaultdash,%
         fill=\defaultfill,join=\defaultjoin,width=\defaultwidth)\else
    \def\n@mref{fl}\ifx\l@debut\n@mref\psset flowchart(arrowp=\defaultfcarrowposition,%
	arrowr=\defaultfcarrowrefpt,line=\defaultfcline,xpadd=\defaultfcxpadding,%
	ypadd=\defaultfcypadding,radius=\defaultfcradius,shape=\defaultfcshape,%
	thick=\defaultfcthickness)\else
    \def\n@mref{me}\ifx\l@debut\n@mref\psset mesh(diag=\defaultmeshdiag)\else
    \def\n@mref{se}\ifx\l@debut\n@mref\psresetsecondsettings\else
    \def\n@mref{th}\ifx\l@debut\n@mref\psset third(color=\defaultthirdcolor)\else
    \immediate\write16{*** Unknown keyword #1 (\BS@ psreset).}%
    \fi\fi\fi\fi\fi\fi\fi}
\def\psset#1(#2){\def\t@xt@{#1}\ifx\t@xt@\empty\trtlis@rg{#2}{\Pssetf@rst}
    \else\keln@mde#1|%
    \def\n@mref{ar}\ifx\l@debut\n@mref\trtlis@rg{#2}{\Psset@rrowhe@d}\else
    \def\n@mref{cu}\ifx\l@debut\n@mref\trtlis@rg{#2}{\Pssetc@rve}\else
    \def\n@mref{fi}\ifx\l@debut\n@mref\trtlis@rg{#2}{\Pssetf@rst}\else
    \def\n@mref{fl}\ifx\l@debut\n@mref\trtlis@rg{#2}{\Pssetfl@wchart}\else
    \def\n@mref{me}\ifx\l@debut\n@mref\trtlis@rg{#2}{\Pssetm@sh}\else
    \def\n@mref{se}\ifx\l@debut\n@mref\trtlis@rg{#2}{\Pssets@cond}\else
    \def\n@mref{th}\ifx\l@debut\n@mref\trtlis@rg{#2}{\Pssetth@rd}\else
    \immediate\write16{*** Unknown keyword: \BS@ psset #1(...)}%
    \fi\fi\fi\fi\fi\fi\fi\fi}
\def\pssetdefault#1(#2){\ifcurr@ntPS\immediate\write16{*** \BS@ pssetdefault is ignored
    inside a \BS@ psbeginfig-\BS@ psendfig block.}%
    \immediate\write16{*** It must be called before \BS@ psbeginfig.}\else%
    \def\t@xt@{#1}\ifx\t@xt@\empty\trtlis@rg{#2}{\Pssd@f@rst}\else\keln@mde#1|%
    \def\n@mref{ar}\ifx\l@debut\n@mref\trtlis@rg{#2}{\Pssd@@rrowhe@d}\else
    \def\n@mref{cu}\ifx\l@debut\n@mref\trtlis@rg{#2}{\Pssd@c@rve}\else
    \def\n@mref{fi}\ifx\l@debut\n@mref\trtlis@rg{#2}{\Pssd@f@rst}\else
    \def\n@mref{fl}\ifx\l@debut\n@mref\trtlis@rg{#2}{\Pssd@fl@wchart}\else
    \def\n@mref{me}\ifx\l@debut\n@mref\trtlis@rg{#2}{\Pssd@m@sh}\else
    \def\n@mref{se}\ifx\l@debut\n@mref\trtlis@rg{#2}{\Pssd@s@cond}\else
    \def\n@mref{th}\ifx\l@debut\n@mref\trtlis@rg{#2}{\Pssd@th@rd}\else
    \immediate\write16{*** Unknown keyword: \BS@ pssetdefault #1(...)}%
    \fi\fi\fi\fi\fi\fi\fi\fi\initpss@ttings\fi}
\def\Pssd@f@rst#1=#2|{\keln@mun#1|%
    \def\n@mref{c}\ifx\l@debut\n@mref\edef\defaultcolor{#2}\else
    \def\n@mref{d}\ifx\l@debut\n@mref\edef\defaultdash{#2}\else
    \def\n@mref{f}\ifx\l@debut\n@mref\edef\defaultfill{#2}\else
    \def\n@mref{j}\ifx\l@debut\n@mref\edef\defaultjoin{#2}\else
    \def\n@mref{u}\ifx\l@debut\n@mref\edef\defaultupdate{#2}\pssetupdate{#2}\else
    \def\n@mref{w}\ifx\l@debut\n@mref\edef\defaultwidth{#2}\else
    \immediate\write16{*** Unknown attribute: \BS@ pssetdefault (..., #1=...)}%
    \fi\fi\fi\fi\fi\fi}
\def\Pssd@@rrowhe@d#1=#2|{\keln@mun#1|%
    \def\n@mref{a}\ifx\l@debut\n@mref\edef\defaultarrowheadangle{#2}\else
    \def\n@mref{f}\ifx\l@debut\n@mref\edef\defaultarrowheadangle{#2}\else
    \def\n@mref{l}\ifx\l@debut\n@mref\y@tiunit{#2}\ifunitpr@sent%
     \edef\defaulth@rdahlength{#2}\else\edef\defaulth@rdahlength{#2pt}%
     \message{*** \BS@ pssetdefault (..., #1=#2, ...) : unit is missing, pt is assumed.}%
     \fi\else
    \def\n@mref{o}\ifx\l@debut\n@mref\edef\defaultarrowheadout{#2}\else
    \def\n@mref{r}\ifx\l@debut\n@mref\edef\defaultarrowheadratio{#2}\else
    \immediate\write16{*** Unknown attribute: \BS@ pssetdefault arrowhead(..., #1=...)}%
    \fi\fi\fi\fi\fi}
\def\Pssd@c@rve#1=#2|{\keln@mun#1|%
    \def\n@mref{r}\ifx\l@debut\n@mref\edef\defaultroundness{#2}\else%
    \immediate\write16{*** Unknown attribute: \BS@ pssetdefault curve(..., #1=...)}%
    \fi}
\def\Pssd@fl@wchart#1=#2|{\keln@mtr#1|%
    \def\n@mref{arr}\ifx\l@debut\n@mref\expandafter\keln@mtr\l@suite|%
     \def\n@mref{owp}\ifx\l@debut\n@mref\edef\defaultfcarrowposition{#2}\else
     \def\n@mref{owr}\ifx\l@debut\n@mref\edef\defaultfcarrowrefpt{#2}\else
     \immediate\write16{*** Unknown attribute: \BS@ pssetdefault flowchart(..., #1=...)}%
     \fi\fi\else%
    \def\n@mref{lin}\ifx\l@debut\n@mref\edef\defaultfcline{#2}\else
    \def\n@mref{pad}\ifx\l@debut\n@mref\edef\defaultfcxpadding{#2}%
                    \edef\defaultfcypadding{#2}\else
    \def\n@mref{rad}\ifx\l@debut\n@mref\edef\defaultfcradius{#2}\else
    \def\n@mref{sha}\ifx\l@debut\n@mref\edef\defaultfcshape{#2}\else
    \def\n@mref{thi}\ifx\l@debut\n@mref\edef\defaultfcthickness{#2}\else
    \def\n@mref{xpa}\ifx\l@debut\n@mref\edef\defaultfcxpadding{#2}\else
    \def\n@mref{ypa}\ifx\l@debut\n@mref\edef\defaultfcypadding{#2}\else
    \immediate\write16{*** Unknown attribute: \BS@ pssetdefault flowchart(..., #1=...)}%
    \fi\fi\fi\fi\fi\fi\fi\fi}
\def\defaultfcarrowposition{0.5}\let\defaultfcarrowpos=\defaultfcarrowposition
\def\defaultfcarrowrefpt{start}
\def\defaultfcline{polygon}
\def\defaultfcradius{0}
\def\defaultfcshape{rectangle}
\def\defaultfcthickness{0}\let\defaultfcthick=\defaultfcthickness
\def\defaultfcxpadding{0}\let\defaultfcxpad=\defaultfcxpadding
\def\defaultfcypadding{0}\let\defaultfcypad=\defaultfcypadding
\def\Pssd@m@sh#1=#2|{\keln@mun#1|%
    \def\n@mref{d}\ifx\l@debut\n@mref\edef\defaultmeshdiag{#2}\else%
    \immediate\write16{*** Unknown attribute: \BS@ pssetdefault mesh(..., #1=...)}%
    \fi}
\def\Pssd@s@cond#1=#2|{\keln@mun#1|%
    \def\n@mref{c}\ifx\l@debut\n@mref\edef\defaultsecondcolor{#2}\else%
    \def\n@mref{d}\ifx\l@debut\n@mref\edef\defaultseconddash{#2}\else%
    \def\n@mref{w}\ifx\l@debut\n@mref\edef\defaultsecondwidth{#2}\else%
    \immediate\write16{*** Unknown attribute: \BS@ pssetdefault second(..., #1=...)}%
    \fi\fi\fi}
\def\Pssd@th@rd#1=#2|{\keln@mun#1|%
    \def\n@mref{c}\ifx\l@debut\n@mref\edef\defaultthirdcolor{#2}\else%
    \immediate\write16{*** Unknown attribute: \BS@ pssetdefault third(..., #1=...)}%
    \fi}
\newif\iffillm@de
\def\pssetfillmode#1{\expandafter\setfillm@de#1:}
\def\setfillm@de#1#2:{\if#1n\fillm@defalse\else\fillm@detrue\fi}
\def\defaultfill{no}     
\newif\ifpstestm@de
\def\pssetupdate#1{\ifcurr@ntPS\immediate\write16{*** \BS@ pssetupdate is ignored inside a
     \BS@ psbeginfig-\BS@ psendfig block.}%
    \immediate\write16{*** It must be called before \BS@ psbeginfig.}%
    \else\expandafter\setupd@te#1:\fi}
\def\setupd@te#1#2:{\if#1n\pstestm@defalse\else\pstestm@detrue\fi}
\def\defaultupdate{no}     
\def\Pssetc@lor#1{\ifps@cri\result@tent=\@ne\expandafter\c@lnbV@l#1 :%
    \def\curr@ntcolor{}\def\curr@ntcolorc@md{}%
    \ifcase\result@tent\or\pssetgray{#1}\or\or\pssetrgb{#1}\or\pssetcmyk{#1}\fi\fi}
\def\pssetcmyk#1{\ifps@cri\def\curr@ntcolor{#1}\def\curr@ntcolorc@md{\c@msetcmykcolor}%
    \def\curr@ntcolorc@mdStroke{\c@msetcmykcolorStroke}%
    \ifcurr@ntPS\PSc@mment{pssetcmyk Color=#1}\us@primarC@lor\fi\fi}
\def\pssetrgb#1{\ifps@cri\def\curr@ntcolor{#1}\def\curr@ntcolorc@md{\c@msetrgbcolor}%
    \def\curr@ntcolorc@mdStroke{\c@msetrgbcolorStroke}%
    \ifcurr@ntPS\PSc@mment{pssetrgb Color=#1}\us@primarC@lor\fi\fi}
\def\pssetgray#1{\ifps@cri\def\curr@ntcolor{#1}\def\curr@ntcolorc@md{\c@msetgray}%
    \def\curr@ntcolorc@mdStroke{\c@msetgrayStroke}%
    \ifcurr@ntPS\PSc@mment{pssetgray Gray level=#1}\us@primarC@lor\fi\fi}
\def\us@primarC@lor{\immediate\write\fwf@g{\d@fprimarC@lor}%
    \let\fillc@md=\prfillc@md}
\def\prfillc@md{\d@fprimarC@lor\space\c@mfill}
\def\defaultcolor{0}       
\def\c@lnbV@l#1 #2:{\def\t@xt@{#1}\relax\ifx\t@xt@\empty\c@lnbV@l#2:
    \else\c@lnbV@l@#1 #2:\fi}
\def\c@lnbV@l@#1 #2:{\def\t@xt@{#2}\ifx\t@xt@\empty%
    \def\t@xt@{#1}\ifx\t@xt@\empty\advance\result@tent\m@ne\fi
    \else\advance\result@tent\@ne\c@lnbV@l@#2:\fi}
\def\Blackcmyk{0 0 0 1}
\def\Whitecmyk{0 0 0 0}
\def\Cyancmyk{1 0 0 0}
\def\Magentacmyk{0 1 0 0}
\def\Yellowcmyk{0 0 1 0}
\def\Redcmyk{0 1 1 0}
\def\Greencmyk{1 0 1 0}
\def\Bluecmyk{1 1 0 0}
\def\Graycmyk{0 0 0 0.50}
\def\BrickRedcmyk{0 0.89 0.94 0.28} 
\def\Browncmyk{0 0.81 1 0.60} 
\def\ForestGreencmyk{0.91 0 0.88 0.12} 
\def\Goldenrodcmyk{ 0 0.10 0.84 0} 
\def\Marooncmyk{0 0.87 0.68 0.32} 
\def\Orangecmyk{0 0.61 0.87 0} 
\def\Purplecmyk{0.45 0.86 0 0} 
\def\RoyalBluecmyk{1. 0.50 0 0} 
\def\Violetcmyk{0.79 0.88 0 0} 
\def\Blackrgb{0 0 0}
\def\Whitergb{1 1 1}
\def\Redrgb{1 0 0}
\def\Greenrgb{0 1 0}
\def\Bluergb{0 0 1}
\def\Cyanrgb{0 1 1}
\def\Magentargb{1 0 1}
\def\Yellowrgb{1 1 0}
\def\Grayrgb{0.5 0.5 0.5}
\def\Chocolatergb{0.824 0.412 0.118}
\def\DarkGoldenrodrgb{0.722 0.525 0.043}
\def\DarkOrangergb{1 0.549 0}
\def\Firebrickrgb{0.698 0.133 0.133}
\def\ForestGreenrgb{0.133 0.545 0.133}
\def\Goldrgb{1 0.843 0}
\def\HotPinkrgb{1 0.412 0.706}
\def\Maroonrgb{0.690 0.188 0.376}
\def\Pinkrgb{1 0.753 0.796}
\def\RoyalBluergb{0.255 0.412 0.882}
\def\Pssetf@rst#1=#2|{\keln@mun#1|%
    \def\n@mref{c}\ifx\l@debut\n@mref\Pssetc@lor{#2}\else
    \def\n@mref{d}\ifx\l@debut\n@mref\pssetdash{#2}\else
    \def\n@mref{f}\ifx\l@debut\n@mref\pssetfillmode{#2}\else
    \def\n@mref{j}\ifx\l@debut\n@mref\pssetjoin{#2}\else
    \def\n@mref{u}\ifx\l@debut\n@mref\pssetupdate{#2}\else
    \def\n@mref{w}\ifx\l@debut\n@mref\pssetwidth{#2}\else
    \immediate\write16{*** Unknown attribute: \BS@ psset (..., #1=...)}%
    \fi\fi\fi\fi\fi\fi}
\def\s@uvdash#1{\edef#1{\curr@ntdash}}
\def\defaultdash{1}        
\def\pssetdash#1{\ifps@cri\edef\curr@ntdash{#1}\ifcurr@ntPS\expandafter\Pssetd@sh#1 :\fi\fi}
\def\Pssetd@shI#1{\PSc@mment{pssetdash Index=#1}\ifcase#1%
    \or\immediate\write\fwf@g{[] 0 \c@msetdash}
    \or\immediate\write\fwf@g{[6 2] 0 \c@msetdash}
    \or\immediate\write\fwf@g{[4 2] 0 \c@msetdash}
    \or\immediate\write\fwf@g{[2 2] 0 \c@msetdash}
    \or\immediate\write\fwf@g{[1 2] 0 \c@msetdash}
    \or\immediate\write\fwf@g{[2 4] 0 \c@msetdash}
    \or\immediate\write\fwf@g{[3 5] 0 \c@msetdash}
    \or\immediate\write\fwf@g{[3 3] 0 \c@msetdash}
    \or\immediate\write\fwf@g{[3 5 1 5] 0 \c@msetdash}
    \or\immediate\write\fwf@g{[6 4 2 4] 0 \c@msetdash}
    \fi}
\def\Pssetd@sh#1 #2:{{\def\t@xt@{#1}\ifx\t@xt@\empty\Pssetd@sh#2:
    \else\def\t@xt@{#2}\ifx\t@xt@\empty\Pssetd@shI{#1}\else\s@mme=\@ne\def\debutp@t{#1}%
    \an@lysd@sh#2:\ifodd\s@mme\edef\debutp@t{\debutp@t\space\finp@t}\def\finp@t{0}\fi%
    \PSc@mment{pssetdash Pattern=#1 #2}%
    \immediate\write\fwf@g{[\debutp@t] \finp@t\space\c@msetdash}\fi\fi}}
\def\an@lysd@sh#1 #2:{\def\t@xt@{#2}\ifx\t@xt@\empty\def\finp@t{#1}\else%
    \edef\debutp@t{\debutp@t\space#1}\advance\s@mme\@ne\an@lysd@sh#2:\fi}
\def\s@uvwidth#1{\edef#1{\curr@ntwidth}}
\def\defaultwidth{0.4}     
\def\pssetwidth#1{\ifps@cri\edef\curr@ntwidth{#1}\ifcurr@ntPS%
    \PSc@mment{pssetwidth Width=#1}\immediate\write\fwf@g{#1 \c@msetlinewidth}\fi\fi}
\def\pssetjoin#1{\ifps@cri\edef\curr@ntjoin{#1}\ifcurr@ntPS\expandafter\Pssetj@in#1:\fi\fi}
\def\Pssetj@in#1#2:{\PSc@mment{pssetjoin join=#1}%
    \if#1r\def\t@xt@{1}\else\if#1b\def\t@xt@{2}\else\def\t@xt@{0}\fi\fi%
    \immediate\write\fwf@g{\t@xt@\space\c@msetlinejoin}}
\def\defaultjoin{miter}   
\def\Pssets@cond#1=#2|{\keln@mun#1|%
    \def\n@mref{c}\ifx\l@debut\n@mref\Pssets@condcolor{#2}\else%
    \def\n@mref{d}\ifx\l@debut\n@mref\pssetseconddash{#2}\else%
    \def\n@mref{w}\ifx\l@debut\n@mref\pssetsecondwidth{#2}\else%
    \immediate\write16{*** Unknown attribute: \BS@ psset second(..., #1=...)}%
    \fi\fi\fi}
\def\pssetseconddash#1{\edef\curr@ntseconddash{#1}}
\def\defaultseconddash{4}  
\def\pssetsecondwidth#1{\edef\curr@ntsecondwidth{#1}}
\edef\defaultsecondwidth{\defaultwidth} 
\def\psresetsecondsettings{%
    \pssetseconddash{\defaultseconddash}\pssetsecondwidth{\defaultsecondwidth}%
    \Pssets@condcolor{\defaultsecondcolor}}
\def\Pssets@condcolor#1{\ifps@cri\result@tent=\@ne\expandafter\c@lnbV@l#1 :%
    \def\sec@ndcolor{}\def\sec@ndcolorc@md{}%
    \ifcase\result@tent\or\pssetsecondgray{#1}\or\or\pssetsecondrgb{#1}%
    \or\pssetsecondcmyk{#1}\fi\fi}
\def\pssetsecondcmyk#1{\def\sec@ndcolor{#1}\def\sec@ndcolorc@md{\c@msetcmykcolor}%
    \def\sec@ndcolorc@mdStroke{\c@msetcmykcolorStroke}}
\def\pssetsecondrgb#1{\def\sec@ndcolor{#1}\def\sec@ndcolorc@md{\c@msetrgbcolor}%
    \def\sec@ndcolorc@mdStroke{\c@msetrgbcolorStroke}}
\def\pssetsecondgray#1{\def\sec@ndcolor{#1}\def\sec@ndcolorc@md{\c@msetgray}%
    \def\sec@ndcolorc@mdStroke{\c@msetgrayStroke}}
\def\us@secondC@lor{\immediate\write\fwf@g{\d@fsecondC@lor}%
    \let\fillc@md=\sdfillc@md}
\def\sdfillc@md{\d@fsecondC@lor\space\c@mfill}
\edef\defaultsecondcolor{\defaultcolor} 
\def\Pss@tsecondSt{%
    \s@uvdash{\typ@dash}\pssetdash{\curr@ntseconddash}%
    \s@uvwidth{\typ@width}\pssetwidth{\curr@ntsecondwidth}\us@secondC@lor}
\def\Psrest@reSt{\pssetwidth{\typ@width}\pssetdash{\typ@dash}\us@primarC@lor}
\def\Pssetth@rd#1=#2|{\keln@mun#1|%
    \def\n@mref{c}\ifx\l@debut\n@mref\Pssetth@rdcolor{#2}\else%
    \immediate\write16{*** Unknown attribute: \BS@ psset third(..., #1=...)}%
    \fi}
\def\Pssetth@rdcolor#1{\ifps@cri\result@tent=\@ne\expandafter\c@lnbV@l#1 :%
    \def\th@rdcolor{}\def\th@rdcolorc@md{}%
    \ifcase\result@tent\or\Pssetth@rdgray{#1}\or\or\Pssetth@rdrgb{#1}%
    \or\Pssetth@rdcmyk{#1}\fi\fi}
\def\Pssetth@rdcmyk#1{\def\th@rdcolor{#1}\def\th@rdcolorc@md{\c@msetcmykcolor}%
    \def\th@rdcolorc@mdStroke{\c@msetcmykcolorStroke}}
\def\Pssetth@rdrgb#1{\def\th@rdcolor{#1}\def\th@rdcolorc@md{\c@msetrgbcolor}%
    \def\th@rdcolorc@mdStroke{\c@msetrgbcolorStroke}}
\def\Pssetth@rdgray#1{\def\th@rdcolor{#1}\def\th@rdcolorc@md{\c@msetgray}%
    \def\th@rdcolorc@mdStroke{\c@msetgrayStroke}}
\def\us@thirdC@lor{\immediate\write\fwf@g{\d@fthirdC@lor}%
    \let\fillc@md=\thfillc@md}
\def\thfillc@md{\d@fthirdC@lor\space\c@mfill}
\def\defaultthirdcolor{1}  
\def\pstrimesh#1[#2,#3,#4]{{\ifcurr@ntPS\ifps@cri%
    \PSc@mment{pstrimesh Type=#1, Triangle=[#2,#3,#4]}%
    \s@uvc@ntr@l\et@tpstrimesh\ifnum#1>\@ne\Pss@tsecondSt\setc@ntr@l{2}%
    \Pstrimeshp@rt#1[#2,#3,#4]\Pstrimeshp@rt#1[#3,#4,#2]%
    \Pstrimeshp@rt#1[#4,#2,#3]\Psrest@reSt\fi\psline[#2,#3,#4,#2]%
    \PSc@mment{End pstrimesh}\resetc@ntr@l\et@tpstrimesh\fi\fi}}
\def\Pstrimeshp@rt#1[#2,#3,#4]{{\l@mbd@un=\@ne\l@mbd@de=#1\loop\ifnum\l@mbd@de>\@ne%
    \advance\l@mbd@de\m@ne\figptbary-1:[#2,#3;\l@mbd@de,\l@mbd@un]%
    \figptbary-2:[#2,#4;\l@mbd@de,\l@mbd@un]\psline[-1,-2]%
    \advance\l@mbd@un\@ne\repeat}}
\initpr@lim\initpss@ttings
\catcode`\@=12

\pssetdefault (update=yes)

\newbox\dessin

\def\MyPSfile{FigDom.ps}
\figinit{pt}
\figpt 1:(-120, 60)\figpt 2:(-100, 65)\figpt 3:(-30, 55)\figpt 4:(50, 60)\figpt 5:(70, 57)
\figpt 6:(-120,30) \figpt 7:(-100, 37)\figpt 8:(-30, 23)\figpt 9:(50, 30)\figpt 10:(70, 26)
\figpt 11:(-120,-30)\figpt 12:(-100, -23)\figpt 13:(-30, -34)\figpt 14:(50, -30)\figpt 15:(70, -35)
\figpt 16:(-120,-60) \figpt 17:(-100, -57) \figpt 18:(-30, -67)\figpt 19:(50,-60)
\figpt 20:(70, -65)
\figptbary 21:c.g.[8,13 ; 1,1]
\psbeginfig{\MyPSfile}
\pscurve[1,1,2,3,4,5,5]
\pscurve[6,6,7,8,9,10,10]
\pscurve[11,11,12,13,14,15,15]
\pscurve[16,16,17,18,19,20,20]
\psset arrowhead(fillmode=yes,ratio=0.07)\psarrow [2,13]
\psarrow [7,18]
\psarrow [8,19]
\psarrow [3,14]
\psendfig
\figvisu{\dessin}{Le foncteur de Picard de $[\Lc^{\frac1n}]$}{
\figinsert{\MyPSfile}
\figsetmark{$\figBullet$}
\figwritep[2,3,4,7,8,9,13,14,18,19]
\figwriten 3:$0$(4)
\figwriten 8:$0$(4)
\figwritene 13:$0$(4)
\figwritene 18:$0$(4)
\figwriten 4:$l$(4)
\figwriten 9:$l$(4)
\figwritene 14:$l$(4)
\figwritene 19:$l$(4)
\figwriten 2:$-l$(4)
\figwriten 7:$-l$(4)
\figsetmark{}
\figwritee 5:${(\Pic_{X/S})_n}$(4)
\figwritee 10:${(\Pic_{X/S})_{n-1}}$(4)
\figwritee 15:${(\Pic_{X/S})_0}$(4)
\figwritee 20:${(\Pic_{X/S})_{-1}}$(4)
\figwriten 21:$\vdots$(2)
\figwriten 3:$\vdots$(20)
\figwrites 18:$\vdots$(10)
}

\centerline{\box\dessin}

\medskip
\noindent
{\sc Description du champ de Picard de $[\Lc^{\frac1n}]$}
\smallskip

On a une \og suite exacte \fg\ de champs de Picard :
\begin{equation}
\label{sec_racine_nieme_champ}
\xymatrix{1 \ar[r]& \champic(X/S) \ar[r]^{\pi^*}& \champic(\X/S) \ar[r]^{\chi}& f_*\Z/n\Z \ar[r]& 1.}
\end{equation}

Autrement dit, $\pi^*$ est pleinement fidèle, $\chi$ est un épimorphisme, et si $\Fc$ est un objet de $\champic(\X/S)$, il provient de $\champic(X/S)$ si et seulement si son caractère $\chi_{\Fc}$ est nul. Tout ceci a déjà été prouvé. De même que précédemment, si l'on suppose que $f_*\Z/n\Z=\Z/n\Z$, alors le champ $\champic(\X/S)$ s'identifie au champ obtenu à partir de $\champic(X/S)\times_S \Z$ en recollant les copies numéro $i$ et $i+nk$ le long de l'isomorphisme
$\mu_{l^k} : (\champic(X/S))_{i+nk} \fleche (\champic(X/S))_i$
pour tous $i,k$ appartenant à $\Z$. En particulier il suffit dans ce cas que $\champic(X/S)$ soit algébrique pour que $\champic(\X/S)$ le soit aussi.
Dans le cas où $\Lc $ a une racine \iem{n} $\Rc$ sur $X$, la suite exacte~(\ref{sec_racine_nieme_champ}) est scindée et $\champic(\X/S)$ s'identifie au produit $\champic(X/S) \times_S f_*\Z/n\Z$.

\subsection{Courbes tordues d'Abramovich et Vistoli}
\label{courbes_tordues}

Abramovich et Vistoli ont mis au jour dans \cite{Abramovich_Vistoli_note}, \cite{Abramovich_Vistoli_CMFFS} et \cite{Abramovich_Vistoli_CSSM} une classe de courbes \og tordues\fg\ qui apparaissent naturellement lorsque l'on cherche à compactifier certains espaces de modules. Ces courbes sont des courbes nodales munies d'une \og structure champêtre\fg\ supplémentaire aux points singuliers ou en certains points marqués. Nous nous proposons de décrire le foncteur de Picard des courbes tordues \emph{lisses}. Nous allons voir que la structure supplémentaire modifie le foncteur de Picard de la courbe d'une manière très analogue à ce que nous avons pu observer dans la section précédente. Nous commençons par quelques rappels sur les courbes tordues.

\begin{sousdefi}[\cite{Abramovich_Vistoli_CSSM} 4.1.2 ou \cite{Olsson_log_twisted_curves} 1.2]
Soit $S$ un schéma. Une courbe tordue sur $S$ est un champ de Deligne-Mumford $f : \Cc \fleche S$ modéré, propre, plat et de présentation finie sur $S$ dont les fibres sont purement de dimension~1, géométriquement connexes et ont au plus des singularités nodales, vérifiant de plus les propriétés suivantes :
\begin{itemize}
\item[1)] Si $\pi : \Cc \fleche C$ est l'espace de modules grossier de $\Cc$ et si $C_{\rm{lis}}$ est le lieu lisse de $C$ sur $S$, alors le sous-champ ouvert $\Cc\times_C C_{\rm{lis}}$ est le lieu lisse de $\Cc$ sur $S$.
\item[2)] Pour tout point géométrique $\overline{s} \fleche S$ le morphisme induit $\Cc_{\overline{s}} \fleche C_{\overline{s}}$ est un isomorphisme au-dessus d'un ouvert dense de $C_{\overline{s}}$.
\end{itemize}

Une courbe tordue $n$-pointée est une courbe tordue munie d'une collection $\{\Sigma_i\}_{i=1}^n$ de sous-champs fermés de $\Cc$ deux à deux disjoints tels que :
\begin{itemize}
\item[(i)] Pour tout $i$, le sous-champ fermé $\Sigma_i$ est dans le lieu lisse de $\Cc$.
\item[(ii)] Pour tout $i$, le morphisme $\Sigma_i\fleche \Cc\fleche S$ est une gerbe étale sur $S$.
\item[(iii)] Si $\Cc_{\emph{gén}}$ est l'ouvert complémentaire des $\Sigma_i$ dans
$\Cc_{\rm{lis}}$, alors $\Cc_{\emph{gén}}$ est un schéma.
\end{itemize}
\end{sousdefi}


\begin{sousremarque} \rm 
D'après la proposition~4.1.1 de \cite{Abramovich_Vistoli_CSSM}, l'espace de modules grossier $C$ est une courbe nodale propre et plate sur $S$, de présentation finie et à fibres géométriquement connexes. Si de plus la courbe $\Cc$ est $n$-pointée, alors l'espace de modules grossier $D_i$ du sous-champ fermé $\Sigma_i$ est naturellement un sous-schéma fermé de $C$, et le morphisme composé $D_i \fleche C \fleche S$ est un isomorphisme. Les $D_i$ définissent donc des sections de $C \fleche S$ qui en font une courbe nodale $n$-pointée au sens usuel.
\end{sousremarque}

\begin{sousthm}[\cite{Abramovich_Vistoli_CSSM} 3.2.3 ou \cite{Olsson_log_twisted_curves} 2.2]\ 

Au voisinage (étale) d'un point marqué, la courbe $\Cc \fleche S$ est de la forme $[U/\mu_r]\fleche \Spec A$ où $U=\Spec A[x]$ et où un générateur de $\mu_r$ agit sur $U$ par $x\mapsto \xi.x$ avec $\xi$ une racine primitive \iem{r} de l'unité.

Au voisinage d'un n\oe ud, la courbe $\Cc \fleche S$ est de la forme $[U/\mu_r]\fleche \Spec A$ où $U=\Spec (A[x,y]/(xy-t))$ pour un certain $t\in A$ et où un générateur de $\mu_r$ agit sur $U$ par $(x,y)\mapsto (\xi.x,\xi'.y)$ avec $\xi$ et $\xi'$ des racines primitives ${r}^{\textrm{ièmes}}$ de l'unité.
\end{sousthm}

\begin{sousremarque}\rm
Si $p\in C$ est un n\oe ud fixé, et si $\xi$ et $\xi'$ sont les racines primitives ${r}^{\textrm{ièmes}}$ de l'unité qui apparaissent ci-dessus, on dit que le n\oe ud $p$ est \og balancé\fg\ si le produit $\xi.\xi'$ est égal à 1. On dit que la courbe $\Cc$ est balancée si tous ses n\oe uds sont balancés. Signalons que si $p$ n'est pas balancé, 
on a nécessairement $t=0$ dans la description locale ci-dessus. Autrement dit on ne peut pas faire disparaître le n\oe ud en déformant la courbe (\cite{Olsson_log_twisted_curves}~2.2).
\end{sousremarque}

Signalons qu'Olsson montre dans \cite{Olsson_log_twisted_curves} que se donner une courbe tordue revient à se donner une courbe nodale classique munie d'une certaine \og structure logarithmique\fg. Pour les courbes tordues lisses, Cadman donne une autre description, plus élémentaire.

\begin{sousthm}[\cite{Cadman_USTITCOC} 2.2.4 et 4.1]
\label{Cadman_equiv_AV}
Se donner une courbe tordue $(\Cc, \{\Sigma_i\}_{i=1}^n)$ $n$-pointée lisse sur un schéma $S$ noethérien et connexe est équivalent à se donner une courbe $n$-pointée $(C,\{\sigma_i\}_{i=1}^n)$ lisse sur $S$ et un $n$-uplet $\overrightarrow{r}=(r_1,\dots, r_n)$ d'entiers strictement positifs inversibles sur $S$. La courbe tordue $\Cc$ est alors isomorphe au champ $C_{\D,\overrightarrow{r}}$ défini de la manière suivante.

Chaque section $\sigma_i$ définit un diviseur de Cartier effectif $D_i$ de $C$. On note $s_{D_i}$ la section canonique de $\Oc(D_i)$ qui s'annule sur $D_i$. La collection des $\Oc(D_i)$ et des $s_{D_i}$ correspond à un morphisme $C \fleche [\A^n/\gm^n]$. On note aussi $\theta_{\overrightarrow{r}}$ le morphisme de $[\A^n/\gm^n]$ dans lui-même qui envoie un $n$-uplet de faisceaux inversibles $(L_1, \dots, L_n)$ muni de sections $(t_1, \dots, t_n)$ sur le $n$-uplet $(L_1^{r_1}, \dots, L_n^{r_n})$ muni de $(t_1^{r_1}, \dots, t_n^{r_n})$. On définit alors $C_{\D,\overrightarrow{r}}$ comme étant le produit fibré
$$C\times_{[\A^n/\gm^n],\theta_{\overrightarrow{r}}} [\A^n/\gm^n].$$
\end{sousthm}

Cadman décrit dans \cite{Cadman_USTITCOC} les faisceaux inversibles sur une courbe tordue lisse sur une base connexe et noethérienne (corollaire~3.2.1) : un faisceau inversible sur $C_{\D,\overrightarrow{r}}$ s'écrit de manière unique sous la forme $\pi^*L \otimes \prod_{i=1}^n \Tc_i^{\otimes k_i}$ où $\pi : \Cc \fleche C$ est la projection de $\Cc$ sur son espace de modules grossier, $L$ est un faisceau inversible sur $C$, $\Tc_i$ est le faisceau inversible $\Oc_{\Cc}(\Sigma_i)$ et $k_i$ est un entier compris entre $0$ et $r_i-1$. Il est clair que les hypothèses noethériennes ne sont pas essentielles pour ce résultat. Par ailleurs on peut aussi supprimer l'hypothèse de connexité sur la base : il faut alors remplacer les entiers $k_i$ par des fonctions localement constantes à valeurs dans $\Z$. On obtient ainsi le théorème suivant.

\begin{sousthm}[\cite{Cadman_USTITCOC} corollaire~3.2.1]
Soient $S$ un schéma, $C$ une courbe lisse $n$-pointée sur $S$, $\overrightarrow{r}$ un $n$-uplet d'entiers strictement positifs et $C_{\D,\overrightarrow{r}}$ la courbe tordue associée par la construction de Cadman. Soit $\Lc$ un faisceau inversible sur $C_{\D,\overrightarrow{r}}$. Alors il existe un faisceau inversible $L$ sur $C$ et des fonctions localement constantes $k_i$ appartenant à $H^0(S, \Z)$ prenant leurs valeurs dans $\{0, \dots, r_i-1\}$ tels que
$$\Lc \simeq \pi^* L \otimes \prod_{i=1}^n \Tc_i^{k_i}.$$
De plus les fonctions $k_i$ sont uniques, $L$ est unique à isomorphisme près, et $\Tc_i^{r_i}$ est isomorphe à $\pi^*\Oc(D_i)$. $\square$
\end{sousthm}

En particulier on a une suite exacte courte :
$$\xymatrix{0 \ar[r]& \Pic(C) \ar[r]& \Pic(\Cc) \ar[r]&
\disp \prod_{i=1}^n H^0(S, \Z/r_i\Z) \ar[r]& 0.}$$
Signalons que si $S$ est le spectre d'un corps algébriquement clos, Chiodo (\cite{Chiodo_twisted_curves_spin_structures}) obtient cette suite exacte d'une autre manière pour une courbe tordue quelconque.

Soit $(\Cc,\{\Sigma_i\}_{i=1}^{n})$ une courbe tordue $n$-pointée lisse sur une base $S$ noethérienne et connexe. D'après le théorème~\ref{Cadman_equiv_AV}, $\Cc$ est isomorphe au champ $C_{\D,\overrightarrow{r}}$ où $C$ est l'espace de modules grossier de $\Cc$, $\overrightarrow{r}$ est un $n$-uplet d'entiers positifs et $\D=(D_1,\dots, D_n)$ est le $n$-uplet de diviseurs effectifs de Cartier de $C$ correspondant aux $\Sigma_i$. Pour tout $i$ le morphisme de $D_i$ vers $S$ est un isomorphisme. On note $\Tc_i$ le faisceau $\Oc_{\Cc}(\Sigma_i)$. Alors toutes ces données sont compatibles au changement de base. Plus précisément, si $T\fleche S$ est un morphisme de changement de base, le produit fibré $\Cc\times_S T$ est isomorphe au champ $C'_{\D',\overrightarrow{r}}$ où $C'$ (resp. $D'_i$) est le produit fibré $C\times_S T$ (resp. $D_i\times_S T$) et $\D'=(D'_1,\dots, D'_n)$. De plus le faisceau inversible canonique $\Tc'_i$ n'est autre que $\Phi^*\Tc_i$ où $\Phi$ est la projection de $\Cc'$ sur $\Cc$. En appliquant le théorème précédent à la courbe $\Cc'$, on obtient pour tout $T$ une suite exacte courte
$$\xymatrix{0 \ar[r]& \Pic(C\times_S T) \ar[r]& \Pic(\Cc\times_S T) \ar[r]&
\disp \prod_{i=1}^n H^0(T, \Z/r_i\Z) \ar[r]& 0.}$$
En appliquant le foncteur \og faisceau étale associé\fg\ elle induit une suite exacte courte de faisceaux étales :
$$\xymatrix{0 \ar[r]& \Pic_{C/S} \ar[r]& \Pic_{\Cc/S} \ar[r]&
\disp \prod_{i=1}^n \Z/r_i\Z \ar[r]& 0.}$$

Notons $l_i$ la classe de $\Oc_C(D_i)$ dans $\Pic_{C/S}(S)$ et $t_i$ la classe de $\Tc_i$ dans $\Pic_{\Cc/S}(S)$. On a dans $\Pic_{\Cc/S}$ la relation $t_i^{r_i}=l_i$. En procédant comme pour le cas du champ $[\Lc^{\frac1n}]$, on voit que le foncteur $\Pic_{\Cc/S}$ s'identifie au foncteur quotient de $\Pic_{C/S}\times_S \Z^n$ par les relations $t_i^{r_i}=l_i$ (où par abus les $t_1, \dots, t_n$ désignent aussi les générateurs canoniques de $\Z^n$). On peut construire ce quotient à la main comme suit. Le produit $\Pic_{C/S}\times_S \Z^n$ est une union disjointe de copies de $\Pic_{C/S}$ indexées par les $n$-uplets $\underline{\alpha}=(\alpha_1, \dots, \alpha_n)$ appartenant à $\Z^n$. Alors $\Pic_{\Cc/S}$ est obtenu en identifiant pour tout $\underline{\alpha}$, pour tout entier $k$ appartenant à $\Z$ et pour tout entier $i$ compris entre 1 et $n$, les copies $(\Pic_{C/S})_{\underline{\alpha}}$ et $(\Pic_{C/S})_{(\alpha_1, \dots, \alpha_i+kr_i, \dots,\alpha_n)}$ via l'isomorphisme de multiplication par $l_i$ de
$(\Pic_{C/S})_{(\alpha_1, \dots, \alpha_i+kr_i, \dots,\alpha_n)}$ dans $(\Pic_{C/S})_{\underline{\alpha}}.$
La loi de groupe est évidente.

\begin{sousremarque} \rm
Si l'on ne tient pas compte de la structure de groupe, on voit que $\Pic_{\Cc/S}$ s'identifie à une union disjointe de $r_1\dots r_n$ copies de $\Pic_{C/S}$.
En particulier la composante neutre de $\Pic_{\Cc/S}$ est la même que celle de $\Pic_{C/S}$. Autrement dit, le morphisme de $\Pic_{C/S}$ vers $\Pic_{\Cc/S}$ induit un isomorphisme naturel :
$\xymatrix{\Pic^0_{C/S} \ar[r]^{\sim} & \Pic^0_{\Cc/S}.}$
\end{sousremarque}


\appendix
\setcounter{section}{0}

{\Large \bf Annexes}
\medskip

Ainsi qu'il a été dit en introduction, la présente annexe rassemble les résultats relatifs à la cohomologie des faisceaux sur les champs algébriques nécessaires au texte principal. Pour éviter que cet article ne devienne démesurément long, les résultats sont donnés ici sans démonstration. On pourra trouver ces dernières, ainsi que des commentaires plus fournis, dans la thèse \cite{Brochard_these} de l'auteur.
Voici en résumé la liste des sujets qui y sont abordés.

Dans la première section, nous rappelons les définitions de base et nous vérifions que sur un champ de Deligne-Mumford, les groupes de cohomologie lisse-étale coïncident avec les groupes de cohomologie étale. Puis
nous introduisons un nouveau site, le site lisse-lisse champêtre, dont les objets sont les morphismes représentables et lisses $\Uc \fleche \X$ de \emph{champs algébriques}. Il définit le même topos que le site lisse-étale mais présente au moins deux avantages : il se comporte mieux vis-à-vis des images directes et il a un objet final.

C'est au départ la nécessité de disposer de techniques de descente cohomologique à la Deligne-Saint-Donat qui a motivé le travail du premier paragraphe de la section~\ref{annexe_desc_coh}. Nous avions en particulier besoin d'un analogue pour les champs algébriques de la suite spectrale de descente relative à un morphisme lisse et surjectif de schémas. Il s'est finalement avéré que l'introduction du site lisse-lisse champêtre rendait ce résultat presque trivial (voir proposition~\ref{suite_spectrale_de_descente}). Nous décrivons ensuite une classe de faisceaux acycliques adaptée aux particularités du site lisse-étale. Ces faisceaux \og $\gll$-acycliques\fg\ nous sont surtout utiles pour obtenir la suite spectrale de Leray relative à un morphisme de champs algébriques (thm.~\ref{ss_Leray}). Nous montrons aussi que les images directes supérieures d'un faisceau lisse-étale abélien peuvent être calculées comme l'on imagine (cf. prop.~\ref{prop_image_directe_sup}).

Le premier paragraphe de la section~\ref{coh_et_chgt_de_base} est consacré au changement de base plat. Le second aux extensions infinitésimales. Il est d'usage, lorsque $i : X \fleche \widetilde{X}$ est une extension infinitésimale, d'identifier les catégories de faisceaux Zariski sur $X$ et sur $\widetilde{X}$. Ceci est tout à fait légitime puisque $X$ et $\widetilde{X}$ ont le même espace topologique sous-jacent. Mieux : le foncteur qui à un ouvert étale $\widetilde{U}$ de $\widetilde{X}$ associe l'ouvert étale $\widetilde{U}\times_{\widetilde{X}} X$ de $X$ définit une équivalence entre les sites étales de $\widetilde{X}$ et de $X$, ce qui permet d'identifier aussi les faisceaux étales. Il faut faire nettement plus attention avec la topologie lisse-étale. On peut en effet vérifier facilement que le foncteur ci-dessus n'est même pas fidèle. Heureusement, on peut tout de même identifier les groupes de cohomologie des faisceaux abéliens sur $X$ et sur $\widetilde{X}$ via le foncteur $i_*$ (cf.~\ref{coh_et_ext_inf}). Cette section contient également un résultat analogue pour les images directes supérieures. Nous donnons un résultat de descente pour les champs algébriques, puis nous décrivons les torseurs du topos lisse-étale et montrons que, dans le cas particulier d'un groupe lisse sur la base $S$, le $H^1$ au sens des foncteurs dérivés coïncide avec le groupe des classes de torseurs.

La section~\ref{Cohomologie_plate} est dévolue à la cohomologie plate sur les champs algébriques. Nous donnons principalement deux résultats. D'une part la suite spectrale qui relie la cohomologie plate à la cohomologie lisse-étale, et d'autre part la généralisation aux champs algébriques du théorème de Grothendieck suivant lequel dans le cas d'un groupe lisse, la cohomologie plate coïncide avec la cohomologie étale (cf. \cite{Dix}, exposé~VI,
paragraphe~11).

\section{\hskip-2.3pt Cohomologie lisse-\'etale sur les champs alg\'ebriques}

\subsection{Topos lisse-étale et cohomologie des faisceaux}

Rappelons bri\`evement, pour la
commodit\'e du lecteur, la d\'efinition du site lisse-\'etale d'un champ
alg\'ebrique donn\'ee au chapitre 12 de \cite{LMB}.

\begin{sousdefi}
\label{def_lisse_etale}
Soit $\X$ un $S$-champ alg\'ebrique. On appelle site \emph{lisse-\'etale} de
$\X$ et on note \emph{Lis-\'et}$(\X)$ le site d\'efini comme suit.

Les ouverts lisses-\'etales de $\X$ sont les couples $(U,u)$ o\`u $U$ est un
$S$-espace alg\'ebrique et $u:U\fleche \X$ un morphisme repr\'esentable et
lisse. Une fl\`eche entre deux tels ouverts $(U,u)$ et $(V,v)$ est un couple
$(\varphi,\alpha)$ faisant 2-commuter le diagramme suivant :
$$\shorthandoff{!;:?}
\xymatrix@R=0.9pc@C=0.9pc{U \ar[rr]^{\varphi} \ar[rd]_u  & \raisebox{-3ex}{$^{\alpha} \FlecheNE$} & V \ar[ld]^v\\
& \X&}$$
Une famille couvrante de $(U,u)$ est une collection de morphismes
$$((\varphi_i,\alpha_i):(U_i,u_i)\flechelongue (U,u))_{i\in I}$$ telle que le
1-morphisme d'espaces alg\'ebriques
$$\coprod_{i\in I} \varphi_i : \coprod_{i\in I} U_i \flechelongue U$$
soit \'etale et surjectif.

Le site \'etale de $\X$, not\'e $\et(\X)$, est la sous-cat\'egorie pleine de
$\liset(\X)$ dont les
objets sont les couples $(U,u)$ o\`u $u$ est un morphisme \'etale, munie de la
topologie induite par celle de $\liset(\X)$.
\end{sousdefi}

On notera $\Oc_{\X}$ le faisceau structural de $\X$ défini de manière évidente. Si $\Ac$ est un faisceau d'anneaux sur $\liset(\X)$, on notera
$\Mod_{\Ac}(\X)$ la cat\'egorie des faisceaux de $\Ac$-modules sur le site
lisse-\'etale de $\X$, ou plus simplement $\Mod(\X)$ lorsque $\Ac=\Oc_{\X}$.
La cat\'egorie des faisceaux ab\'eliens sera not\'ee $\Ab(\X)$.
Les références usuelles (\cite{SGA4_1} II (6.7) et \cite{Tohoku} th\'eor\`eme (1.10.1)) permettent de vérifier que la cat\'egorie $\Mod_{\Ac}(\X)$ est une cat\'egorie ab\'elienne avec suffisamment
d'objets injectifs. En particulier il en est ainsi des cat\'egories $\Mod(\X)$
et $\Ab(\X)$.

On rappelle que le foncteur \og sections globales \fg\  est d\'efini de la
mani\`ere suivante. Si $\Fc$ est un faisceau lisse-\'etale sur $\X$, l'ensemble
$\Gamma(\X, \Fc)$ est l'ensemble des familles $(s_{(U,u)})$ de sections de $\Fc$
sur les $(U,u)\in \ob\liset(\X)$ qui sont compatibles aux fl\`eches de
restriction en un sens \'evident (cf. \cite{LMB} (12.5.3)). On v\'erifie
imm\'ediatement que le foncteur $\Gamma(\X, .) : \Ab(\X)\fleche \Ab$ est exact
\`a gauche. On d\'efinit alors $H^i(\X, .)$ comme \'etant le
$i^{\textrm{\`eme}}$ foncteur d\'eriv\'e \`a droite de $\Gamma(\X, .)$.

Il résulte maintenant de \cite{SGA4_2} V (3.5) que sur la cat\'egorie $\Mod(\X)$, le foncteur $H^i(\X,.)$ co\"{i}ncide avec le
$i^{\textrm{\`eme}}$ foncteur d\'eriv\'e \`a droite de $\Gamma(\X, .) : \Mod(\X)
\fleche \Ab$.

\begin{souslem} Soit $\X$ un $S$-champ de Deligne-Mumford et
soient $\Fc$, $\Gc$ des faisceaux lisses-étales d'ensembles (resp. de groupes
abéliens) sur $\X$. On suppose que $\Fc$ est cartésien (voir la définition dans \cite{LMB} (12....)). Alors l'application
$$\Hom(\Fc, \Gc) \flechelongue \Hom(\Fc_{\text{\rm ét}},\Gc_{\text{\rm ét}})$$
induite par le foncteur d'inclusion $\xymatrix@C=1pc{\et(\X) \ar@{^{(}->}[r] &
\liset(\X)}$  est bijective.
\label{lemme_de_prolongement_des_morphismes_etales}
\end{souslem}

\begin{sousremarque}\rm
Le résultat est faux si on ne suppose pas $\Fc$ cartésien. Il suffit de considérer un faisceau lisse-étale
abélien $\Fc$ non nul tel que $\Fc_{\text{ét}}$ soit nul. Alors, dans $\Hom(\Fc, \Fc)$ on a au moins deux
éléments distincts, à savoir l'identité et le morphisme nul, tandis que $\Hom(\Fc_{\text{ét}}, \Fc_{\text{ét}})$
est réduit à zéro. Pour exhiber un tel faisceau, on peut par exemple prendre $\X=\Spec k$, le spectre d'un corps,
et poser $\Fc(U,u)=\Omega_{U/k}$.
\end{sousremarque}

\begin{souscor}
Soient $\X$ un $S$-champ algébrique et $\Fc$ un faisceau lisse-étale sur $\X$. Alors le morphisme canonique
$$\Gamma_{\text{\rm lis-ét}}(\X,\Fc) \flechelongue \Gamma_{\text{\rm ét}}(\X,\Fc)$$
est un isomorphisme.
\end{souscor}

\begin{souscor}
Soit $\X$ un $S$-champ de Deligne-Mumford, et soit $\Fc$ un objet injectif de la
catégorie des faisceaux lisses-étales abéliens sur $\X$. Alors la restriction
$\Fc_{\text{\rm ét}}$ de $\Fc$ au site étale de $\X$ est un objet injectif de la
catégorie des faisceaux étales abéliens.
\end{souscor}

\begin{sousprop}
\label{coh_et_egale_coh_liset}
Soient $\X$ un $S$-champ de Deligne-Mumford et
$\Fc\in\Ab(\X)$ un faisceau lisse-\'etale ab\'elien sur $\X$.
On note $\Fc_{\text{\rm \'et}}$ la restriction de $\Fc$
au site \'etale de $\X$. Alors pour tout $q\geq 0$ on a un isomorphisme
canonique :
$$\flechen{H^q_{\text{\rm lis-\'et}}(\X,\Fc)}{\sim}
{H^q_{\text{\rm \'et}}(\X,\Fc_{\text{\rm \'et}}).}$$
\end{sousprop}

\noindent {\sc Fonctorialité du topos lisse-étale}
\medskip

\label{fonctorialite_lisse_etale}
Si $f: \X \fleche \Y$ est un 1-morphisme de $S$-champs alg\'ebriques, on lui associe un couple de foncteurs
adjoints $(f^{-1},f_*)$ comme dans \cite{LMB}, (12.5).
Comme on pourra le lire bient\^ot dans la prochaine \'edition de \cite{LMB}, ou
d\`es aujourd'hui dans \cite{Olsson_Sheaves_on_Artin_stacks}, il serait erron\'e
de penser que le couple $(f^{-1},f_*)$ est toujours un
morphisme de topos. Il peut en effet arriver que le foncteur $f^{-1}$ ne soit
pas exact, même lorsque $\X$ et $\Y$ sont des schémas. Pour s'en convaincre on
consultera les r\'ef\'erences cit\'ees. C'est toutefois le cas d\`es que $f$
est lisse.

Le foncteur $f_* : \Ab(\X)\fleche \Ab(\Y)$ est exact \`a
gauche puisqu'il a un adjoint
\`a gauche. On peut donc d\'efinir les foncteurs d\'eriv\'es $R^qf_* : \Ab(\X)
\fleche \Ab(\Y)$.

\bigskip
\noindent {\sc Le site lisse-lisse champêtre d'un champ algébrique}
\medskip

\label{site_llc}
Pour un certain nombre de considérations techniques, le site lisse-étale
défini dans \cite{LMB} ne contient pas suffisamment d'ouverts pour être vraiment
commode. En effet, lorsque $f$ n'est pas repr\'esentable, le foncteur $f_*$ n'est pas
induit par une application continue $\liset(\X) \fleche \liset(\Y)$, ce qui pose problème
par exemple lorsque l'on essaye de calculer les foncteurs images directes supérieures $R^qf_*$
(voir le paragraphe (\ref{images_directes_supérieures})). C'est la
raison pour laquelle nous introduisons un site un peu plus gros, qui ne
pr\'esentera plus
les m\^emes inconv\'enients. Nous démontrons ensuite (\ref{topos_llc_equiv_topos_liset}) que le topos qu'il définit est équivalent
au topos lisse-étale. Le choix de la topologie lisse plut\^ot qu'\'etale
pour ce site est essentiellement d\^u au fait que pour la topologie
\'etale, les \og ouverts lisses champêtres\fg\ 
ne sont pas toujours recouverts par un espace alg\'ebrique.

\medskip

Soit $\X$ un $S$-champ alg\'ebrique. On d\'efinit la 2-cat\'egorie des
ouverts lisses champêtres de la mani\`ere suivante. Les objets sont les couples $(\Uc,u)$
o\`u $\Uc$ est un $S$-champ alg\'ebrique et $u : \Uc \fleche \X$ est
un morphisme repr\'esentable et lisse. Un 1-morphisme entre deux tels ouverts
$(\Uc,u)$ et $(\Vc,v)$ est un couple $(\varphi,\alpha)$ o\`u $\varphi : \Uc
\fleche \Vc$ est un 1-morphisme (automatiquement représentable !) de $S$-champs alg\'ebriques et $\alpha : u
\Rightarrow v\circ \varphi$ est un 2-isomorphisme.
$$\shorthandoff{!;:?}
\xymatrix@R=0.9pc@C=0.9pc{\Uc \ar[rr]^{\varphi} \ar[rd]_u  & \raisebox{-3ex}{$^{\alpha} \FlecheNE$} & \Vc \ar[ld]^v\\
& \X&}$$
Si $(\varphi,\alpha)$ et $(\psi,\beta)$ sont deux 1-morphismes de $(\U,u)$ dans
$(\Vc,v)$, un 2-morphisme entre $(\varphi,\alpha)$ et $(\psi,\beta)$ est un
2-isomorphisme $\gamma : \varphi \Rightarrow \psi$ tel que
$\beta=(v_*\gamma)\circ \alpha$.

\begin{souslem}
\label{lemme_site_champetre}
Soient $(\U,u)$ et $(\Vc,v)$ deux ouverts lisses champêtres. Alors la cat\'egorie
des morphismes de $(\U,u)$ dans $(\Vc,v)$ est \'equivalente \`a une cat\'egorie discr\`ete.
\end{souslem}

\`A l'avenir on identifiera $\Hom((\U,u),(\Vc,v))$ \`a une cat\'egorie
discr\`ete \'equivalente et on parlera de l'\emph{ensemble} des morphismes de
$(\Uc, u)$ dans $(\Vc,v)$, et de la \emph{cat\'egorie} des ouverts lisses champêtres.
Il est clair que cette cat\'egorie admet des produits fibr\'es.

\begin{sousdefi}
On appelle site lisse-lisse champêtre, et on note $\gll(\X)$, la cat\'egorie des
ouverts lisses champêtres de $\X$ munie de la topologie engendr\'ee par la pr\'etopologie
pour laquelle les familles couvrantes sont les familles de morphismes
$$((\Uc_i,u_i)\flechelongue (\Uc,u))_{i\in I}$$ telles que le morphisme (automatiquement représentable)
$$\xymatrix{\disp \coprod_{i\in I} u_i : \coprod_{i\in I} \Uc_i \ar[r] &\Uc}$$
soit lisse et surjectif.
\end{sousdefi}

\begin{sousprop}
\label{topos_llc_equiv_topos_liset}
\begin{itemize}
\item[1)] La topologie de $\liset(\X)$ est \'egale \`a la topologie engendr\'ee par la
pr\'etopologie dite lisse, pour laquelle les familles couvrantes sont les
familles de morphismes $((U_i,u_i)\fleche (U,u))_{i\in I}$ telles que le
morphisme $\coprod_{i\in I} u_i : \coprod_{i\in I} U_i \fleche U$
soit lisse et surjectif.
\item[2)] Le foncteur d'inclusion $\xymatrix@C=1pc{\liset(\X) \ar@{^(->}[r]& \gll(\X)}$
induit une \'equivalence
de topos de la cat\'egorie des faisceaux sur le site lisse-lisse champêtre vers la
cat\'egorie des faisceaux sur $\liset(\X)$.
\end{itemize}
\end{sousprop}

\begin{sousremarque}
\label{rem_foncteurs_compatibles}\rm
Notons $\X_{\gll}$ le topos des faisceaux sur $\gll(\X)$. Via l'équivalence de topos ci-dessus, si $f : \X \fleche \Y$
est un 1-morphisme de $S$-champs alg\'ebriques, le foncteur $f_*$ est simplement donn\'e
par $(f_{*}\Fc)(\Uc,u)=\Fc(\X\times_{\Y} \Uc, \pr_{\X})$ pour tout
ouvert lisse champêtre $u : \Uc \fleche \Y$ de $\Y$. Le foncteur \og sections globales \fg\ est quant à lui donné par 
$\Gamma(\X,\Fc)=\Fc(\X, \Id_{\X})$.
\end{sousremarque}

\begin{sousremarque}\rm
\label{restriction_faisceaux_injectifs}
Si $u : \Uc \fleche \X$ est un morphisme repr\'esentable et lisse de $S$-champs
alg\'ebriques, et si $\Ac$ est un anneau du topos $\X_{\gll}$, le foncteur
$$\fonction{u^*}{\Mod_{\Ac}(\X)}{\Mod_{\Ac_{(\Uc,u)}}(\Uc)}
{\Fc}{\Fc_{(\Uc,u)}}$$
o\`u $\Fc_{(\Uc,u)}$ (resp. $\Ac_{(\Uc,u)}$) d\'esigne la restriction de $\Fc$
(resp. $\Ac$) au site lisse-lisse champêtre de $\Uc$, commute
aux limites projectives et inductives arbitraires. Il a donc un adjoint \`a
gauche $u_!$, et de plus cet adjoint est exact (tout ceci r\'esulte de
\cite{SGA4_1}~IV~(11.3.1)). Par suite $u^*$ transforme les objets injectifs en
objets injectifs (voir aussi \cite{SGA4_2} V (2.2)).
\end{sousremarque}

\subsection{Suites spectrales}

\noindent {\sc La suite spectrale relative à un recouvrement}
\medskip

\label{annexe_desc_coh}
Soit $P : \Uc^0 \fleche \X$ un morphisme représentable, lisse et surjectif de
$S$-champs algébriques.
On note
\begin{eqnarray*}
\Uc^1 &=& \Uc^0\times_{\X}\Uc^0\\
\Uc^2 &=& \Uc^0\times_{\X}\Uc^0\times_{\X}\Uc^0\\
&\vdots& \\
\Uc^n &=& \Uc^0\times_{\X} \dots \times_{\X}\Uc^0 \quad \text{($n+1$ fois)}\\
&\vdots&
\end{eqnarray*}
Les $\Uc^n$ forment avec les diagonales partielles et les projections un
$S$-champ
alg\'ebrique simplicial muni d'une augmentation vers $\X$ :
$$\xymatrix{\dots & \Uc^n & \dots & \Uc^2 \ar[r] \ar@<2ex>[r] \ar@<-2ex>[r] &
  \Uc^1 \ar@<1ex>[l] \ar@<-1ex>[l] \ar@<1ex>[r] \ar@<-1ex>[r] &
  \Uc^0 \ar[l] \ar[r] & \X}$$

Soit $\Fc$ un faisceau ab\'elien sur le site lisse-\'etale de $\X$. On
note $\Fc^i$ la restriction de $\Fc$ au site lisse-étale de $\Uc^i$.
On cherche alors \`a calculer
la cohomologie de $\Fc$ en fonction de la cohomologie des $\Fc^i$ sur les
$\Uc^i$.
On peut associer au champ alg\'ebrique simplicial ci-dessus un complexe de
\v{C}ech de la mani\`ere suivante. Pour $n\geq 2$, on note $p_{1\dots\hat{l}\dots n}$, o\`u la
notation $\hat{l}$ signifie que l'indice $l$ est omis, la projection
$\Uc^{n-1}\fleche \Uc^{n-2}$ qui \og oublie \fg\ le facteur d'indice $l$ de
$\Uc^{n-1}$. Par exemple $p_1$ et $p_2$ d\'esignent respectivement les premi\`ere
et seconde projections de $\Uc^0\times_{\X} \Uc^0$ sur $\Uc^0$. On d\'efinit alors le
complexe de \v{C}ech $S(H^q)$ comme \'etant le complexe :
$$\xymatrix{H^q(\Uc^0,\Fc^0) \ar[r]^{d^0}& H^q(\Uc^1,\Fc^1) \ar[r]^-{d^1}&
  \dots \ar[r] & H^q(\Uc^p,\Fc^p) \ar[r]^-{d^p} & \dots}$$
avec
\begin{eqnarray*}
d^0 &=& p_2^*-p_1^*\\
d^1 &=& p_{23}^*-p_{13}^*+p_{12}^*\\
&\vdots& \\
d^p &=& \sum_{l=1}^{p+2} (-1)^{l+1} p_{1\dots \hat{l}\dots (p+2)}^*.
\end{eqnarray*}

On d\'esigne par $\check{H}^p(H^q(\Uc^{\bullet},\Fc^{\bullet}))$ l'homologie
en degr\'e $p$ de ce complexe :
$$\check{H}^p(H^q(\Uc^{\bullet},\Fc^{\bullet}))=\frac{\Ker d^p}{\Im d^{p-1}}.$$

Le r\'esultat suivant, qui donne la suite
spectrale reliant la \og cohomologie de \v{C}ech \fg\ relativement \`a la
famille couvrante $P : \Uc^0 \fleche \X$, \`a la cohomologie lisse-\'etale de $\Fc$
sur $\X$, est essentiellement trivial (\emph{cf.} \cite{Brochard_these}). Il ne fait pas appel aux techniques de descente cohomologique de Deligne
présentées par Saint-Donat dans l'exposé V~bis de \cite{SGA4_2} : grâce à l'introduction du site lisse-lisse champêtre, il s'agit juste de la \og suite spectrale
de Cartan-Leray relative à un recouvrement\fg\ (\cite{SGA4_2} V (3.3)).

\begin{sousprop}
\label{suite_spectrale_de_descente}
Reprenons les hypoth\`eses et notations pr\'ec\'edentes. Il existe une suite
spectrale (fonctorielle) :
$$E_2^{p,q}=\check{H}^p(H^q(\Uc^{\bullet},\Fc^{\bullet})) \Rightarrow H^{p+q}(\X,\Fc).$$
\end{sousprop}

\bigskip
\noindent {\sc Faisceaux acycliques}
\medskip
\label{Faisceaux_acycliques}

Soit $f : \X \fleche \Y$ un morphisme de $S$-champs algébriques. Dans la mesure où le couple de foncteurs
adjoints $(f^{-1},f_*)$ n'est pas un morphisme de topos, il n'est pas évident \emph{a priori} que le foncteur
$f_*$ transforme les faisceaux injectifs en faisceaux injectifs. Il nous faudra pourtant montrer, pour obtenir la suite spectrale de Leray relative à un morphisme
de champs algébriques (cf. paragraphe (\ref{paragraphe_ss_Leray})), que le foncteur $f_*$
transforme les faisceaux injectifs en faisceaux acycliques pour le foncteur \og sections globales\fg.
L'usage de faisceaux \og flasques\fg\ en un certain sens va nous permettre de résoudre ce problème. La principale difficulté réside dans le choix de la
classe de faisceaux que l'on considère (\emph{cf.} \cite{Brochard_these}).

\begin{sousdefi}[\cite{SGA4_2} V 4.2]
Soit $\X$ un $S$-champ algébrique et soit $\Fc$ un faisceau lisse-étale abélien. On dit que $\Fc$ est
$\gll$-acyclique si pour tout morphisme représentable et lisse $u : \Uc \fleche \X$ et pour tout $q>0$, le
groupe $H^q(\Uc, \Fc_{(\Uc,u)})$ est nul.
\end{sousdefi}

\begin{sousremarque} \rm S'il est évident que les faisceaux flasques au sens de SGA~4 sont $\gll$-acycliques, il n'y a aucune
raison \emph{a priori} pour que la réciproque soit vraie. Pour quelques commentaires sur ce genre de questions, on pourra
consulter \cite{SGA4_2}~V~4.6 et~4.13. Par ailleurs il est évident que les faisceaux injectifs sont flasques,
donc aussi $\gll$-acycliques.
\end{sousremarque}

La proposition suivante, quoique fortement inspirée de la proposition~V~4.3 de \cite{SGA4_2}, en diffère légèrement.

\begin{sousprop}
\label{prop_acyclicite}
Soient $\X$ un $S$-champ algébrique, et $\Fc$ un faisceau sur $\gll(\X)$. Les propositions suivantes sont équivalentes :
\begin{enumerate}
\item[(1)] $\Fc$ est $\gll$-acyclique ;
\item[(2)] pour toute famille couvrante $((\Uc',u') \fleche (\Uc,u))$ dans $\gll(\X)$, et pour tout $q>0$,
le groupe $H^q(\Uc'/\Uc,\Fc)$ est nul (où $H^q(\Uc'/\Uc,\Fc)$ désigne le $q^{\text{ième}}$ groupe de cohomologie de
\v{C}ech de $\Fc$, aussi noté $\check{H}^q(H^0(\Uc'^{\bullet},\Fc_{(\Uc,u)}^{\bullet}))$ dans la section précédente).
\end{enumerate}
\end{sousprop}

\bigskip
\noindent {\sc Images directes suprieures et suite spectrale de Leray relative à un morphisme de champs algébriques}
\medskip
\label{paragraphe_ss_Leray}
\label{images_directes_supérieures}

La proposition suivante montre que les images directes supérieures se calculent comme l'on pense. Leurs restrictions au site étale de $\Y$ aussi (voir l'énoncé, plus précis, donné dans \cite{Brochard_these}).

\begin{sousprop}
\label{prop_image_directe_sup}
\begin{enumerate}
\item[1)] Soit $\Fc$ un faisceau ab\'elien sur $\liset(\X)$. Alors le faisceau
$R^qf_*\Fc$ est le faisceau associ\'e au pr\'efaisceau qui \`a tout ouvert
lisse-\'etale $(U,u)$ de
$\Y$ associe $H^q(\X\times_{\Y} U, \Fc_{(\X\times_{\Y} U, \pr_{\X})})$.
$$\xymatrix{\X\times_{\Y}U \ar[r] \ar[d]_{\pr_{\X}} \cartesien& U\ar[d]^u \\
\X \ar[r]^f & \Y}$$
\item[2)] Si $\Y=\Spec A$ est un schéma affine et si $\X$ est quasi-compact, alors pour
tout faisceau quasi-cohérent $\Fc$ sur $\X$ :
$$H^0(\Spec A, R^qf_*\Fc)\simeq H^q(\X,\Fc).$$
\end{enumerate}
\end{sousprop}

\begin{sousprop}
\label{images_directes_sup_compatibles}
La restriction du foncteur $R^qf_*$ \`a la cat\'egorie
$\Mod(\X)$ co\"{i}ncide avec le $q^{\textrm{i\`eme}}$ foncteur d\'eriv\'e \`a
droite de $f_* : \Mod(\X) \fleche \Mod(\Y)$.
\end{sousprop}

\begin{souslem}
\label{faisceaux_gll_acycliques_preserves}
Le foncteur $f_* : \Ab(\X)\fleche \Ab(\Y)$ préserve les faisceaux
$\gll$-acy\-cli\-ques. En particulier il transforme les faisceaux injectifs en faisceaux acycliques pour
le foncteur $\Gamma(\Y,.)$. 
\end{souslem}

\begin{sousthm}[suite spectrale de Leray relative à $f : \X \fleche \Y$]
\label{ss_Leray}
Soit $$f : \X\flechelongue \Y$$ un morphisme de $S$-champs algébriques, et soit $\Fc$ un faisceau lisse-étale abélien sur
$\X$. Il existe une suite spectrale :
$$H^p(\Y,R^qf_*\Fc) \Rightarrow H^{p+q}(\X,\Fc).$$
\end{sousthm}

\begin{sousprop}
Soient $f : \X \fleche \Y$ et $g : \Y \fleche \Zc$ des morphismes de $S$-champs algébriques et soit
$\Fc$ un faisceau lisse-étale abélien sur $\X$. Alors on a une suite spectrale :
$$R^pg_*R^qf_* \Fc \Rightarrow R^{p+q}(g\circ f)_* \Fc.$$
\end{sousprop}

\subsection{Comparaison de cohomologies}

\noindent {\sc Cohomologie et changement de base}
\medskip
\label{coh_et_chgt_de_base}

La démonstration de la proposition~\ref{lemme_coh_et_chgt_base} ci-dessous (donc celle de son corollaire~\ref{images_directes_commutent_chgt_base_plat}) utilise des hyperrecouvrements. La définition que, pour la commodité du lecteur,
nous rappelons ci-dessous, est plus simple en apparence que celle présentée dans \cite{SGA4_2}. Elle lui est
équivalente dès que la catégorie sous-jacente à $C$ admet des sommes directes arbitraires
(en effet les objets semi-représentables sont alors représentables), ce que nous supposons
dès à présent.

\begin{sousdefi}[cf. \cite{SGA4_2}~V (7.3.1.2)]
Soient $C$ un site (admettant des sommes directes arbitraires), $X$ un objet de $C$ et $U^{\bullet}\fleche X$
un objet simplicial de $C$ muni d'une augmentation vers $X$, ou de manière équivalente un objet simplicial
de $C/X$. On dit que $U^{\bullet}$ est un hyperrecouvrement de $X$ s'il possède les propriétés suivantes :
\begin{enumerate}
\item[(1)] Pour tout entier $n\geq 0$ le morphisme canonique
$$U^{n+1} \flechelongue (\cosq_n\,(\sq_n\,U^{\bullet}))_{n+1}$$
est un morphisme couvrant.
\item[(2)] Le morphisme $U^0 \fleche X$ est couvrant.
\end{enumerate}
\end{sousdefi}

\begin{sousremarque}\rm
Dans la définition précédente, $\sq_n$ désigne le foncteur qui à un objet simplicial associe son tronqué à l'ordre
$n$. Le foncteur $\cosq_n$ est l'adjoint à droite de $\sq_n$. Enfin $(\cosq_n\,(\sq_n\, U^{\bullet}))_{n+1}$ désigne
le terme d'indice $n+1$ de l'objet simplicial $\cosq_n\,(\sq_n\,U^{\bullet})$.
\end{sousremarque}

Si $U^{\bullet}\fleche X$ est un hyperrecouvrement de $X$ et si $\Fc$ est un faisceau abélien sur $X$, on lui
associe un complexe
$$0 \flechelongue H^q(U^0,\Fc^0) \flechelongue H^q(U^1,\Fc^1) \flechelongue  \dots$$
de manière tout à fait analogue à ce qui a été fait dans la section (\ref{annexe_desc_coh}) (où $\Fc^i$ désigne
la restriction de $\Fc$ à $U^i$). On note encore $\check{H}^p(H^q(U^{\bullet},\Fc^{\bullet}))$ l'homologie
en degr\'e $p$ de ce complexe. On a alors une suite spectrale (\cite{SGA4_2}~V~(7.4.0.3)) :
$$\check{H}^p(H^q(U^{\bullet},\Fc^{\bullet}))\Rightarrow H^{p+q}(X,\Fc).$$

\begin{sousprop}Soit $\X$ un $S$-champ alg\'ebrique quasi-compact, avec $S=\Spec A$ affine, et
soit $\Fc$ un faisceau quasi-coh\'erent sur $\X$.
Soient $A'$ une $A$-alg\`ebre plate, $S'=\Spec A'$, et $\X'=\X\times_S S'$.
On note $\Fc'$ l'image inverse de $\Fc$ sur $\X'$.
Alors pour tout $q\geq 0$ le morphisme naturel
$$\flechen{H^q(\X,\Fc)\otimes_A A'}{}{ H^q(\X', \Fc')}$$
est un isomorphisme.
\label{lemme_coh_et_chgt_base}
\end{sousprop}

\begin{sousprop}
\label{images_directes_commutent_chgt_base_plat}
Soit $f : \X \fleche \Y$ un morphisme quasi-compact de $S$-champs alg\'ebriques, et
soit $\Fc$ un faisceau quasi-coh\'erent sur $\X$. Soit $u : \Y' \fleche \Y$ un
morphisme plat de changement de base.
$$\xymatrix{\X' \ar[r]^v \ar[d]_g \cartesien& \X\ar[d]^f \\
\Y' \ar[r]^u & \Y}$$
Alors pour tout $q\geq 0$ le morphisme naturel
$$u^*R^qf_*\Fc \flechelongue (R^qg_*)(v^*\Fc)$$
est un isomorphisme.
\end{sousprop}

\bigskip
\noindent {\sc Cohomologie et extensions infinit\'esimales}
\medskip
\label{par_coh_et_ext_inf}
\begin{souslem}
\label{coh_et_ext_inf}
Soit $i : \X \fleche \Xt$ une immersion ferm\'ee de champs alg\'ebriques
d\'efinie par un id\'eal quasi-coh\'erent de $\Xt$ de carr\'e nul, et soit $\Fc$
un faisceau ab\'elien sur $\X$. Alors pour tout $q>0$ le faisceau $R^qi_*\Fc$ est nul. En particulier
pour tout $n$ le morphisme naturel :
$$ H^n(\Xt, i_*\Fc)\flechelongue H^n(\X,\Fc) $$
est un isomorphisme.
\end{souslem}

\begin{souscor}
Soit
$$\xymatrix{\X \ar[r]^i \ar[d]_f \cartesien &\Xt\ar[d]^{\widetilde{f}}\\ \Y \ar[r]^j &\Yt}$$
un carré 2-cartésien de $S$-champs algébriques, dans lequel on suppose que $i$ et $j$ sont des immersions fermées
définies par des idéaux quasi-cohérents de carrés nuls. Soit $\Fc$ un faisceau abélien sur $\X$. Alors on 
a pour tout $p$ un isomorphisme canonique :
$$\xymatrix{R^p\widetilde{f}_*(i_*\Fc) \ar[r]^{\sim} &j_*R^pf_*\Fc.}$$
\end{souscor}

\bigskip
\noindent {\sc Un résultat de descente}
\medskip
\label{Un_résultat_de_descente}

Le théorème suivant n'est autre qu'un résultat de descente pour les champs algébriques. Il est utilisé au paragraphe suivant.
\begin{sousthm}
\label{descente_pour_les_champs}
Soit $\Fc$ un faisceau sur $\gll(\X)$ qui, localement pour la topologie lisse sur $\X$, est représentable.
Alors $\Fc$ est représentable par un unique objet $(\Pc, p)$ de $\gll(\X)$. Autrement dit il existe un morphisme représentable et lisse
de $S$-champs algébriques $p : \Pc \fleche \X$ qui représente $\Fc$.
\end{sousthm}
\begin{sousremarque}\rm
L'hypothèse sur $\Fc$ signifie qu'il existe une présentation $x_0 : X_0\fleche \X$ de $\X$ 
telle que la restriction de $\Fc$ au site lisse de $X_0$ soit représentable par
un $X_0$-espace algébrique lisse $F_0$.
\end{sousremarque}

\bigskip
\noindent {\sc Cohomologie et torseurs}
\medskip
\label{Cohomologie_et_torseurs}

\begin{sousdefi}
\label{def_torseur}
Soient $S$ un schéma, $G$ un $S$-espace algébrique en groupes lisse sur $S$ et $\X$ un $S$-champ algébrique.
Soit $p : \Pc \fleche \X$ un 1-morphisme représentable et lisse de $S$-champs algébriques. Une action de $G$
sur $(\Pc,p)$ est un quadruplet $(\mu,\varphi_{\mu}, \varphi_e, \varphi_{\text{ass}})$, où
$\mu$ est un 1-morphisme de $G\times_S \Pc$ vers $\Pc$, et où $\varphi_{\mu}, \varphi_e$ et $\varphi_{\text{ass}}$
sont des 2-isomorphismes faisant 2-commuter les diagrammes suivants :
$$\xymatrix{G\times_S \Pc \ar[r]^-{\mu} \ar[d]_{\pr_2} \ar@{}[dr]|{\varphi_{\mu} \FlecheNE}& \Pc \ar[d]^p \\
\Pc\ar[r]_p & \X}$$
$$\shorthandoff{!;:?}
\xymatrix@R=0.9pc@C=0.9pc{\Pc \ar[rr]^-{e\times\Id_{\Pc}} \ar[rd]_{\Id_{\Pc}}
& \raisebox{-3ex}{$^{ \varphi_e} \FlecheNE$}
&G\times_S\Pc \ar[ld]^{\mu}\\ & \Pc&}$$
$$\xymatrix{G\times_S G \times_S \Pc \ar[r]^-{\Id_G\times \mu} \ar[d]_{m_G\times \Id_{\Pc}}
 \ar@{}[dr]|{\varphi_{\text{ass}} \FlecheNE}&
G\times_S \Pc  \ar[d]^{\mu} \\ G\times_S \Pc \ar[r]_{\mu} & \Pc}$$
les 2-isomorphismes $\varphi_e$ et $\varphi_{\text{ass}}$ étant assujettis aux conditions de compatibilité, que nous nous dispenserons d'écrire, assurant
que ce sont des 2-morphismes dans la 2-catégorie des ouverts lisses champêtres de $\X$ (voir juste avant~\ref{lemme_site_champetre}).

On dit qu'un tel couple $(\Pc, p)$ muni d'une action de $G$ est un $G$-torseur sur $\X$ si les deux conditions supplémentaires
suivantes sont vérifiées :
\begin{itemize}
\item[(a)] $p$ est surjectif ;
\item[(b)] le morphisme naturel
$$G\times_S \Pc \flechelongue \Pc \times_{\X}\Pc$$
induit par le triplet $(\mu, \pr_2, \varphi_{\mu})$ est un isomorphisme.
\end{itemize}
\end{sousdefi}

La proposition suivante est alors une conséquence des travaux de Giraud dans~\cite{Giraud}.

\begin{sousprop}
\label{thm_coh_torseur}
Soient $S$ un schéma, $G$ un $S$-espace algébrique en groupes commutatifs lisse sur $S$ et $\X$ un $S$-champ algébrique.
Alors le groupe $H^1(\X,G)$ est canoniquement isomorphe au groupe des classes d'isomorphie de $G$-torseurs sur $\X$
(muni de la loi de groupe induite par le produit contracté de torseurs, cf.~\cite{Giraud}~III~2.4.5).
\end{sousprop}

\begin{sousremarque}\rm
Comme d'habitude, un $G$-torseur $p : \Pc \fleche \X$ est trivial si et seulement si le morphisme structural $p$ a une section.
\end{sousremarque}

\section{Cohomologie fppf}
\label{Cohomologie_plate}
\subsection{Sorites sur la cohomologie plate}

Soit $\X$ un $S$-champ algébrique. On définit le gros site \emph{fppf} de $\X$,
noté $\fppf(\X)$, de la manière suivante. Les ouverts sont les couples $(U,u)$,
où $U$ est un espace algébrique et $u : U \fleche \X$ est un morphisme
localement de présentation finie. Un morphisme entre deux tels ouverts $(U,u)$
et $(V,v)$ est un couple $(\varphi, \alpha)$ où $\varphi : U \fleche V$ est un
morphisme d'espaces algébriques et $\alpha$ est un 2-isomorphisme de $u$ vers
$v\circ \varphi$. Une famille couvrante est une collection de morphismes
$((\varphi_i, \alpha_i) : (U_i,u_i)\fleche (U,u))_{i\in I}$ telle que le
morphisme
$$\coprod_{i\in I} \varphi_i : \coprod_{i\in I} U_i \flechelongue U$$
soit fidèlement plat et localement de présentation finie.

\begin{sousremarque}\rm
Il est évident que l'on aurait obtenu un topos équivalent en ne prenant pour
ouverts que les couples $(U,u)$ où $U$ est un schéma (appliquer le lemme de
comparaison de SGA~4, \cite{SGA4_1}~III~4.1). En particulier si $\X$ est un
schéma on retrouve le topos des faisceaux sur le gros site \emph{fppf}
usuel (considéré par exemple dans \cite{Dix}, exposé~VI, paragraphe~5, p.124).
\end{sousremarque}

On définit de manière évidente les faisceaux d'anneaux $\Oc_{\X}$ et $\Z$ sur
$\fppf(\X)$. On note $\X_{\pl}$ le topos des faisceaux sur $\fppf(\X)$. Si
$\Ac$ est un anneau du topos $\X_{\pl}$, on note $\Mod_{\Ac}^{\pl}(\X)$ la
catégorie des faisceaux de modules sur le site annelé $(\fppf(\X),\Ac)$. Elle
sera notée $\Mod^{\pl}(\X)$ lorsque $\Ac=\Oc_{\X}$, et $\Ab^{\pl}(\X)$
lorsque $\Ac=\Z$.

On définit le foncteur \og sections globales\fg\ en posant pour tout faisceau
$\Fc$ sur $\fppf(\X)$ :
$$\Gamma_{\pl}(\X,\Fc)=\lpro \Fc(U,u)$$
où la limite projective est prise sur l'ensemble des couples $(U,u)$ de
$\fppf(\X)$. Il est clair que ce foncteur commute aux limites projectives
quelconques. En particulier il est exact à gauche. De même que dans le cas des
faisceaux lisses-étales, il résulte de SGA~4 (\cite{SGA4_1}~II~6.7) que la
catégorie $\Mod_{\Ac}^{\pl}(\X)$ est une catégorie abélienne avec suffisament
d'objets injectifs. On définit alors $H^i_{\pl}(\X,.)$ comme étant le
$i^{\textrm{ième}}$ foncteur dérivé à droite de $\Gamma_{\pl}(\X,.) : 
\Ab^{\pl}(\X)\fleche \Ab$. Il coïncide sur $\Mod_{\Ac}^{\pl}(\X)$ avec le
$i^{\textrm{ième}}$ foncteur dérivé à droite de $\Gamma_{\pl}(\X,.) : 
\Mod_{\Ac}^{\pl}(\X)\fleche (\Gamma_{\pl}(\X,\Ac)$-$\Mod)$
(\cite{SGA4_2} V 3.5).

\medskip
\noindent
{\sc Fonctorialité.}

Soit $f: \X\fleche \Y$ un 1-morphisme de $S$-champs algébriques. On lui associe
un couple de foncteurs adjoints $(f_{\pl}^{-1},f_*^{\pl})$ d'une manière tout à
fait analogue à ce qui a été fait pour les faisceaux lisses-étales. Le foncteur
$f_*^{\pl}$ est alors exact à gauche, et on note $R^if_*^{\pl}$ ses foncteurs
dérivés à droite. 

\begin{sousremarque}\rm
Contrairement à ce qu'il se passe dans le cas des faisceaux lisses-étales, la limite inductive qui définit le
foncteur image inverse $f_{\pl}^{-1}$ \emph{est} filtrante. Ceci est essentiellement dû au fait qu'un morphisme
entre deux objets localement de présentation finie est lui-même localement de présentation finie (ce qui se
produisait aussi pour le site étale, mais n'était plus vrai en remplaçant étale par lisse). On en déduit
(cf. par exemple \cite{Milne_Etale_coh} annexe~A) que le foncteur $f_{\pl}^{-1}$ est exact, et donc que le
couple $(f_{\pl}^{-1},f_*^{\pl})$ est un morphisme de topos.
\end{sousremarque}

Notons ici aussi deux cas particuliers dans lesquels les foncteurs image directe
ou image inverse ont une expression plus simple. Lorsque $f$ est localement de
présentation finie, le foncteur $f_{\pl}^{-1}$ est simplement le foncteur de
restriction au site \emph{fppf} de $\X$ via le foncteur $(U,u)\mapsto (U,f\circ
u)$. Si $f$ est représentable, le foncteur $f_*^{\pl}$ provient d'une
application continue $\fppf(\X) \fleche \fppf(\Y)$, qui à un ouvert \emph{fppf}
$u : U\fleche \Y$ associe l'ouvert formé de l'espace algébrique $U\times_{\Y}\X$
muni de la projection sur $\X$. Naturellement ceci n'est pas vrai si $f$ n'est
pas représentable, puisqu'alors $U\times_{\Y}\X$ n'est pas nécessairement un
espace algébrique. Ce défaut s'avère gênant dans le calcul des images directes
supérieures, et motive l'introduction du site \emph{fppf} champêtre ci-dessous.

\medskip
\noindent
{\sc Le gros site \emph{fppf} champêtre d'un $S$-champ algébrique.}

Soit $\X$ un $S$-champ algébrique. On définit le gros site \emph{fppf} champêtre
de $\X$ de la manière suivante. Un ouvert de $\X$ est un couple  $(\Uc,u)$
o\`u $\Uc$ est un $S$-champ alg\'ebrique et $u : \Uc \fleche \X$ est
un morphisme repr\'esentable et localement de présentation finie. Un 1-morphisme
de $(\Uc,u)$ vers $(\Vc,v)$ est un couple $(\varphi,\alpha)$ o\`u
$\varphi : \Uc \fleche \Vc$ est un 1-morphisme de $S$-champs alg\'ebriques et
$\alpha$ est un 2-isomorphisme de $u$ dans $v\circ \varphi$.
Si $(\varphi,\alpha)$ et $(\psi,\beta)$ sont deux 1-morphismes de $(\U,u)$ dans
$(\Vc,v)$, un 2-morphisme entre $(\varphi,\alpha)$ et $(\psi,\beta)$ est un
2-isomorphisme $\gamma : \varphi \Rightarrow \psi$ tel que
$\beta=(v_*\gamma)\circ \alpha$.

L'analogue \emph{fppf} du lemme (\ref{lemme_site_champetre}) est alors valable,
de sorte que la 2-catégorie que nous venons de décrire est en fait équivalente à
une catégorie, que nous appellerons la catégorie des ouverts \emph{fppf}
champêtres de $\X$. On définit maintenant le gros site \emph{fppf}
champêtre de $\X$ de manière évidente, et on le note $\fppfc(\X)$.

\begin{sousprop}
Le foncteur d'inclusion $\fppf(\X) \fleche \fppfc(\X)$ induit une équivalence de
topos de la catégorie des faisceaux sur $\fppfc(\X)$ vers la catégorie des
faisceaux sur $\fppf(\X)$.
\end{sousprop}

On voit alors facilement que via cette équivalence de catégories, le foncteur
image directe $f_*^{\pl}$ provient d'une application continue
$$\fppfc(\X) \flechelongue \fppfc(\Y)$$
qui envoie un ouvert $u : \Uc\fleche \Y$ de $\Y$ sur l'ouvert $\X\times_{\Y}
\Uc\fleche \X$ de $\X$. 

\begin{sousprop} 
\label{images_directes_superieures_fppf}
Soit $f : \X \fleche \Y$ un morphisme de $S$-champs
algébriques.
\begin{itemize}
\item[1)] Soit $\Fc$ un faisceau abélien sur $\fppf(\X)$. Alors le faisceau
$R^qf_*^{\pl} \Fc$ est le faisceau associé au préfaisceau qui à tout ouvert
\emph{fppf} $u : U\fleche \Y$ de $\Y$ associe $H^q_{\pl}(\X\times_{\Y} U,
\Fc_{(\X\times_{\Y} U,\pr_{\X})})$.
\item[2)] La restriction du foncteur $R^qf_*^{\pl}$ à la catégorie
$\Mod^{\pl}(\X)$ coïncide avec le $q^{\textrm{ième}}$ foncteur dérivé à droite
du foncteur $f_*^{\pl} : \Mod(\X) \fleche \Mod(\Y)$.\ $\square$
\end{itemize}
\end{sousprop}

\subsection{Comparaison avec la cohomologie lisse-étale}
Soit $\X$ un $S$-champ algébrique. On a une application continue évidente :
$$p : \fppf(\X) \flechelongue \liset(\X)$$
induite par le foncteur d'inclusion de $\liset(\X)$ dans $\fppf(\X)$.
En particulier $p$ induit un couple de foncteurs adjoints :
$$\begin{array}{c} p_* : \X_{\pl} \flechelongue \X_{\text{lis-ét}}\\
p^{-1} : \X_{\text{lis-ét}}\flechelongue \X_{\pl}.
\end{array}$$
Le foncteur $p^{-1}$ peut être décrit de la manière suivante. On définit d'abord un adjoint à gauche au foncteur
$p_*$ pour les préfaisceaux, que l'on note $\widehat{p^{-1}}$, en associant à tout préfaisceau $\Fc$ sur $\liset(\X)$
le préfaisceau sur $\fppf(\X)$ qui à tout ouvert \emph{fppf} $u : U \fleche \X$ associe
$$(\widehat{p^{-1}}\Fc)(U,u) = \lind \Fc(V,v)$$
où la limite inductive est prise sur l'ensemble des diagrammes 2-commutatifs :
$$\shorthandoff{!;:?}
\xymatrix@R=0.9pc@C=0.9pc{U \ar[rr]^{\varphi} \ar[rd]_u  & \raisebox{-3ex}{$^{\alpha} \FlecheNE$} & V \ar[ld]^v\\
& \X&}$$
avec $(V,v)\in \ob\liset(\X)$. On définit alors $p^{-1}$ en posant $p^{-1}=\underline{a}\widehat{p^{-1}}$ où
$\underline{a}$ est le foncteur \og faisceau associé\fg. Encore une fois, on voit que $\widehat{p^{-1}}$ est défini par
une limite inductive sur $\liset(\X)$ qui n'est pas filtrante. Le couple $(p^{-1},p_*)$ n'est donc pas un morphisme de
topos a priori, et il n'y a aucune raison formelle pour que $p_*$ transforme les objets injectifs de $\Ab^{\pl}(\X)$ en
injectifs de $\Ab(\X)$, ce qui nous oblige à travailler un peu plus pour obtenir la suite spectrale (\ref{ss_comparaison_fppf_vs_liset})
ci-dessous.

\begin{souslem}
Si $\Fc$ est un objet injectif de $\Ab^{\pl}(\X)$, alors $p_*\Fc$ est un faisceau acyclique pour le foncteur
$\Gamma_{\text{\rm lis-ét}}(\X,.)$.
\end{souslem}

\begin{souslem}
Le morphisme naturel
$$\Gamma_{\pl}(\X,\Fc)\flechelongue\Gamma_{\text{\rm lis-ét}}(\X,p_*\Fc)$$
est un isomorphisme.
\end{souslem}

On déduit des deux résultats précédents une suite spectrale de Leray, qui résume les relations générales
entre cohomologie lisse-étale et cohomologie \emph{fppf}.

\begin{sousthm}
\label{ss_comparaison_fppf_vs_liset}
Soit $\X$ un $S$-champ algébrique, et soit $\Fc$ un faisceau abélien sur $\fppf(\X)$. On a alors une suite
spectrale (fonctorielle en $\Fc$) :
$$H^p_{\text{\rm lis-ét}}(\X,R^qp_*\Fc) \Rightarrow H^{p+q}_{\pl}(\X,\Fc). $$
\vskip-8mm\hfill $\square$
\end{sousthm}

\begin{sousremarque}\rm
\label{rmq_comparaison_coh_fppf_vs_liset}
On en déduit en particulier des morphismes canoniques :
$$H^p_{\text{lis-ét}}(\X,p_*\Fc) \flechelongue H^{p}_{\pl}(\X,\Fc).$$
De plus ces morphismes sont des isomorphismes si pour tout $q>0$, le faisceau lisse-étale $R^qp_*\Fc$ est nul.
\end{sousremarque}

Nous allons à présent généraliser aux champs algébriques le résultat de Grothendieck (\cite{Dix}, exposé~VI,
paragraphe~11) selon lequel si $G$ est un groupe lisse sur un schéma $X$, le morphisme canonique
$$H^q_{\text{ét}}(X,p_*G) \flechelongue H^q_{\pl}(X,G)$$
est un isomorphisme pour tout $q$. Dans la mesure où nous n'utiliserons ce résultat que pour le groupe $\gm$, nous
n'avons pas cherché à le démontrer pour les \og champs en groupes\fg\ lisses sur $\X$ (notion qui resterait d'ailleurs
à définir), mais nous nous sommes contenté de considérer un groupe lisse sur la base $S$. Dans ce cadre élémentaire, le résultat se déduit assez facilement du cas des schémas (\emph{cf.} \cite{Brochard_these}).

\begin{sousthm}
\label{coh_fppf_groupe_lisse}
Soient $S$ un schéma, $G$ un schéma en groupes lisse sur $S$, et $\X$ un $S$-champ algébrique.
On note encore $G$ le faisceau défini sur le site $\fppf(\X)$ par :
$$\forall (U,u)\in\ob\fppf(\X)\quad G(U,u):=G(U,f\circ u)=\Hom_S(U,G).$$
Alors pour tout $q$ le morphisme canonique
$$H^q_{\text{\rm lis-ét}}(\X, p_* G) \flechelongue H^q_{\pl}(\X,G)$$
est un isomorphisme.
\end{sousthm}

\bibliographystyle{../../../../latex/plain-fr}
\addcontentsline{toc}{section}{Bibliographie}
\bibliography{../../../../latex/mabiblio}
\end{document}